\documentclass[12pt]{article}
\newcommand{\N}{\mbox{I$\!$N}}
\newcommand{\R}{\mbox{I$\!$R}}

\newcommand{\qed}{{\hfill {$\rlap{$\sqcap$}\sqcup$}}\\[0.2in]\hspace*{0.5in}}
\newcommand{\qedwh}{{\hfill {$\rlap{$\sqcap$}\sqcup$}}\\[0.2in]}
\newcommand{\bk}{\\[0.1in] \hspace*{0.5in} }
\newcommand{\btd}{\bigtriangledown}
\newcommand{\mfor}{\ \ \ \ {\mbox{for}} \ \ }

\begin{document}

\vspace*{0.4in}

\centerline{\LARGE   {\bf  Refined Estimates  for  Simple Blow\,-\,ups  }}

\vspace*{0.2in}

\centerline{\LARGE {\bf   of the Scalar Curvature Equation on $S^n$}}

\vspace{0.5in}

\centerline{\Large {Man Chun  {\LARGE L}EUNG}}

\vspace*{0.15in}

\centerline{\large {National University of Singapore${\,}^*$}}

\vspace*{0.1in}

\centerline{\tt m\,atlmc\,@\,nus.edu.sg}

\vspace{0.42in}

\begin{abstract}
\vspace{0.25in}
\noindent In their work \cite{Compact} on a sharp compactness theorem for the Yamabe problem, Khuri, Marques and Schoen apply a refined blow\,-\,up analysis (what we call `\,second order blow\,-\,up argument\,' in this article) to obtain highly accurate approximate solutions for the Yamabe equation.
As for the conformal scalar curvature equation on \,$S^n$\, with \,$n \ge 4$\,,\, we examine the second order blow\,-\,up argument and obtain refined estimate for a blow\,-\,up sequence near a simple blow\,-\,up point\,.\, The estimate involves local effect from the Taylor expansion of the scalar curvature function,\,   global effect from   other blow\,-\,up points, and the balance formula as expressed in the  Pohozaev identity in an essential way.

\end{abstract}

\vspace*{1.1in}

{\bf 2010 AMS MS Classification}\,:\ \ Primary 35J60\,;\ \ Secondary 53C21. \\[0.05in]
$\overline{\ \ \ \ \ \ \ \ \ \ \ \ \ \ \ \ \ \ \ \ \  \ \ \ \ \ \ \ \ \ \ \ \ \ \ \ \ \ \ \ \ \ \ \ }$\\[0.01in]
\noindent{$\!\!{\!}^* {\small{\mbox{\, Department of Mathematics, National University of Singapore, 10, Lower Kent Ridge Rd.,}}}$}\\[0.03in]
 {\small{\mbox{\ \  Singapore 119076, Republic of Singapore\,.}}}\\[0.04in]
 \noindent{$\!\!{\!}^{**} {\small{\mbox{\, e\,-\,Appendix starts at page  44 onward. }}}$}\\[0.05in]

\pagenumbering{arabic}

{\bf \large {\bf  \S\,1. \ \ Introduction.}}

\vspace*{0.2in}

In this article, we expound local and global contributions to a refined `second order' estimate for simple  blow\,-\,ups (\,or {\it simple  isolated blow\,-\,ups\,} as known in some literature) of the prescribed scalar curvature equation

\vspace*{-0.25in}

$$  \Delta_1\, u \, - \,{\tilde c}_n\,n(n \,-\,1)\,u \, + \, ({\tilde c}_n \,{\cal K} )\, u^{{n + 2}\over {n \,-\,2}} \ \,= \ 0 \ \ \ \ \ \ \ \ {\mbox{on}} \ \ \ \ S^n. \leqno (1.1)
$$

\vspace*{-0.1in}

Here $\,{\cal K}\,,\,$ fixed once it is given, is assumed to be smooth enough [\,say, in $\,C^{\,n \,+\, 4}\,(S^n)$\,]\,,\,  $\Delta_1\,$  is the
Laplacian on $\,S^n\,$ with the standard metric $\,g_1\,,\,$ and $\,
{\tilde c}_n = {{n
- 2}\over {4\, (n \,-\,1)}}\,$ (\,$n \ge 3)\,.\,$
Via the stereographic {\it projection}
$
\,\dot{\cal P} \,:\ S^n \setminus \{\, {\bf N} \} \ \,  \longrightarrow \    {\R}^n\,,\,
$
which sends the north pole $\,{\bf N} \,\in\, S^n$ to infinity, equation (1.1) can be expressed in the simple form
$$
\Delta_o  \,v + \,({\tilde c}_n \,K)\,v^{{n + 2}\over {n \,-\,2}}  = 0\,,   \leqno (1.2)
$$

\vspace*{-0.15in}

$$
\,v \, (y) :=   u  \, ({\dot{\cal P}}^{-1} (y))\, \cdot \left( {2\over {1+ |\,y|^{\,2}}} \right)^{\!\!{{n \,-\,2}\over 2}}  \   {\mbox{and}} \ \
  K\,
(y) := {\cal K}  \,(\,{\dot{\cal P}}^{-1} (y)) \mfor y \in  \R^n\,.\ \ \  \leqno (1.3)
$$
In (1.2),  $\,\Delta_o\,$ is the
Laplacian on $\,\R^n\,$ with the standard Euclidean metric $\,g_o\,.\,$ Considered as a `\,dual\,' to  the Yamabe equation, the study on the non\,-\,compact behavior (or blow\,-\,up phenomenon) of equation (1.1)  is found to be rich and vibrant. See for examples   \cite{Chen-Lin-2} \cite{Chen-Lin-3} \cite{Leung-Supported} \cite{Li-1} \cite{Li-2},    and the references therein.

\vspace*{0.15in}

{\bf{\S\,1\,a. \  $\!\!\!$
Simple blow\,-\,up.}}\ \ When $\,\,({\tilde c}_n \,K)\,$ is equal to a constant, say, $\,n\,(n \,-\,2)\,,\,$ equation (1.2) has a family of solutions\,:

\vspace*{-0.3in}

$$
{\bf A}_{\,\epsilon_i\,,\,\, \zeta_i} \, (y) \ = \, \left(   {{\epsilon_i}\over {\epsilon_i^2 \,+ |\, y \,-\, \zeta_i|^2}} \right)^{\!{{n \,-\,2}\over 2}}\ . \leqno (1.4)
$$

For non\,-\,constant $\,K$,\, a sequence of positive solutions $\,\{ v_i\}\,$ of (1.2) which blows up at $0$\, is shown to be `close' to a sequence found in (1.4)\,.\, Precisely,
$$
\bigg\vert \ v_i \, (y) \ - \ {\bf A}_{\,\epsilon_i\,,\,\, \zeta_i} \, (y) \, \bigg\vert \ \,\le\, \  \varepsilon_i \cdot  \epsilon_i^{-\,{{n \,-\,2}\over 2} }    \ \ \ \  \mfor \ \ |\,y \,-\, \zeta_i|\,\le\,  \epsilon_i \, R_i  \ \ \ {\mbox{and}} \ \ \ i \,\gg\, 1\,, \leqno (1.5)
$$
with parameters $\,\epsilon_i \,\to\, 0$\,,\, $|\,\zeta_i| \,\to\, 0\,$  and $\,R_i \,\to\, \infty\,$  specific to $\,\{ v_i\}\,$ [\,cf. (2.21)\, in \,{\bf \S\,2\,d}\,]\,.\, Here (via a rescaling), we assume throughout this article that
$$
({\tilde c}_n\,K)\, (0) \ = \ n \, (n \,-\,2)\,. \leqno (1.6)
$$
Estimate (1.5) is rather weak -  its accuracy in general deteriorates when $\,i \to \infty\,.\,$  Moreover, (1.5) is valid (generally) in a sequence of  shrinking balls $\,B_{\zeta_i}\, (\epsilon_i\,R_i)\,.\,$  [\,The order of shrinkage  $O\, (\epsilon_i)$ makes space for the bubbles described in (1.4) to be stacked up (developed vertically, cf. the Delaunay solution \cite{K-M-P-S}\,), or  be put in juxtaposition (developed horizontally)\,.\, See \cite{Leung-Supported}  for a classification of blow\,-\,ups for equation (1.1)\,]\,.\, \bk
One can characterize  simple  blow\,-\,up  in a geometric manner - when the bound in (1.5) can be stabilized in terms of scale and accuracy,  namely,\\[0.1in]
(1.7)
$$
 {1\over C}\cdot  {\bf A}_{\,\epsilon_i\,,\,\, \zeta_i} \, (y)  \ \,\le\, \ v_i\, (y) \ \,\le\, \ C\cdot  {\bf A}_{\,\epsilon_i\,,\,\, \zeta_i} \, (y)   \!\!\!\!\!    \mfor \ {\mbox{all}}  \ \ |\,y \,-\, \zeta_i|  \ \,\le\, \ \rho_o \ \ \,{\mbox{and}}\  \ i \,\gg\, 1\,.
$$
Here $\rho_o$ and $C$ are fixed positive numbers\,.\,
See Proposition 2.24 for the precise statement. Cf. also the notion of quasi\,-\,isometry\,.\, Simple blow\,-\,ups is by far the most common non\,-\,compact behavior we encounter in equation (1.1)\,.\, In \cite{I} \cite{II} \cite{III},  blow\,-\,up sequence  with a fixed nonconstant $\,{\cal K}$\, (may not be symmetric) are constructed using the Lyapunov\,-\,Schmidt reduction method. See also \cite{Wei-Yan}\,.

\vspace*{0.15in}

{\bf{\S\,1\,b. \  Description of the main result.}} \ \ In this article, we identify three factors affecting the fixed scale behavior of simple blow\,-\,ups. \\[0.1in]
{\bf I)} \ \ {\it The local behavior of $\,K$ in terms of the Taylor expansion}
$$
({\tilde c}_n\,K) \, (y) \ = \ n\, (n \,-\,2)\ + \ [-\,{\bf P}_\ell \, (y)]\  + \ R_{\ell \, +\, 1}\, (y) \ \ \mfor \ y \in B_o\, (\rho_o)\,.   \leqno (1.8)
$$
Here $\,{\bf P}_\ell\,$ is a {\it homogeneous polynomial}\, of degree $\ell \, \in\,  \N\,,$\, and $\, R_{\ell \, +\, 1} \,$ the remainder in the Taylor expansion\,.\, [\,See (2.29) and (3.7) for the sign convention we use on ${\bf P}_\ell\,.$\,]\, We know that if $\,0\,$ is a blow\,-\,up point for equation (1.2), then $\ell \,\ge\, 2\,$ (that is, $\btd\, K \, (0) = \vec{\,0}\,;\,$ see  Theorem 5.1 in \cite{Leung-Supported} for the precise statement\,; cf. also \cite{Chen-Lin-2})\,.\,  Hence
\begin{eqnarray*}
(1.9)    & \ &  {\mbox{\   number \ \   of \  \ critical \ \  points \ \   of \ \  }} {\cal K}  \ \ {\mbox{is \ \ finite}}\\      & \ & \!\!\!\!\!\!\!\!\!\!\!\!\!\!\!\!\!\!\!\!\Longrightarrow \ \ {\mbox{equation \     (1.1) \ \ has \ \ at \ \ most \ \ finite  \   number  \   of  \  blow\,-\,up \   points}}\,. \ \ \  \ \ \ \ \ \ \ \  \ \ \ \ \ \ \ \
\end{eqnarray*}
The leading polynomial term $\,{\bf P}_\ell$\, comes into the picture when we find the difference between $\,v_i$\, and the standard solutions given in (1.4)\,.\, See (3.11)\,.\, The second order blow\,-\,up argument allows us to discern the  central information  enveloped in $\,{\bf P}_\ell\,$.\, We discuss this point more in {\bf \S\,1\,c}\, and  {\bf \S\,1\,d}\,.

\vspace*{0.15in}

{\bf II)} \ \ {\it `Flexibility'\, of the simple blow\,-\,up as measured by }\,$|\, \xi_i| \,=\, O \, (\lambda^{\,\alpha}_i)$\,,\,
\begin{eqnarray*}
(1.10) & \ &  {\mbox{where}} \ \ \ v_i\,(\xi_i) \ = \ \max \ \left\{ \ v_i\, (y) \ \,| \ \ y \,\in\  \overline{ B_o\, (\,\rho_o)}   \ \right\} \!\!\!\!\!\mfor i \gg 1\,, \ \ \ \ \ \ \ \ \ \ \ \ \ \ \ \  \\[0.15in]
(1.11) &\ &     {\mbox{and}} \ \ \   \lambda_i \ := \ [ \,v_i\,(\xi_i)\,]^{-\, {2\over {n \,-\,2}} }\,, \ \ \ \ \   \xi_i \,\to\, 0 \ \   \ \ ({\mbox{the  \ \ blow\,-\,up \ \ point}})\,. \ \ \ \ \  \ \ \ \ \ \ \ \  \ \ \ \ \ \ \ \
\end{eqnarray*}
Here $\,\rho_o$\, is a small fixed number (\,its size is related to other blow\,-\,up points)\,.\, $v_i\,$ could have other maximal points near $\,\xi_i$\,,\, but their distances to $\,\xi_i\,$ are at most $\,o\, (\lambda_i)\,$ for $\,i \gg 1\,.\,$ Refer to \,{\bf \S\,2\,g}\,.\, The position parameter \,$\xi_i$\, appears in the expression for the difference $\,[\,v_i \,-\,A_{\lambda_i\,,\,\,\xi_i}]\,,\,$ see (3.11)\,.\, Thanks to  the work of Chen and Lin \cite{Chen-Lin-2} \cite{Chen-Lin-3}, one can impose      conditions, including the following main ones (see \,{\bf \S\,2\,g\,} for the full detail)\,:
\begin{eqnarray*}
\Vert \btd\, {\bf P}_\ell\, (y) \Vert & \ge & C\, |\,y|^{\,\ell \,-\,1} \ \ \ \ \ \ \ \ \ \ \ \ \mfor y \in B_o\,(\,\rho_o)\\[0.15in]
{\mbox{and}} \ \ \ \ \   \int_{\R^n} \btd\, {\bf P}_\ell\, (y + {\cal X}) \!\!\!\!\!\!  &\cdot & \!\!\!\!\!\!  [\,A_1\, (y)]^{{2n}\over {n \,-\,2}}\  dy \ \not= \ {\vec{\,0}} \mfor {\mbox{all}} \ \ {\cal X} \in \R^n \setminus \, \{ 0 \}\,, \ \ \ \ \ \ \ \ \ \ \  \ \ \ \ \ \ \
\end{eqnarray*}
resulting
$$ \ \
|\, \xi_i| \ = \ o \, (\lambda_i) \ \  \ \ \,{\mbox{modulo \ \ a \ \ subsequecne}} \ \ \ \ \ \ \ \    \left( \!{\mbox{ that \ \,is}}\,, \ \ \  \lambda_i^{-1} \cdot \xi_i  \,\to\, 0 \right). \ \ \leqno (1.12)
$$

\vspace*{0.1in}

{\bf III)}  \ \ {\it Interaction with other blow\,-\,ups\,.\,} This  is expressed by a global harmonic function (or Green's function)

\vspace*{-0.25in}

$$
\sum_{j = 0}^k \ {{A_j}\over {|\, y \,-\, {\hat Y}_j|^{\, n- 2}}}  \mfor\    y \,\in\, \R^n \,\setminus \ \{ {\hat Y}_o \,=\,0\,, \ {\hat Y}_1\,, \ \cdot \cdot \cdot\,, \ \ {\hat Y}_k \ \}\,, \leqno (1.13)
$$
`effective' outside a neighborhood containing all the
  blow\,-\,up points $\{ {\hat Y}_j \}_{j = 0}^{\,k}\,.\,$ \,See \,{\bf \S\,2\,e}\,.\, In (1.13), \,$A_j\,$ are positive numbers.\, A major challenge here is to match the information expressed in (1.13) (\,the `collapsed region'\,)   with the one in (1.5) (\,the `blow\,-\,up' region)\,.\, See \,{\bf \S\,2\,d.1}\, and  {\bf \S\,6\,a} \,for a fuller discussion.

\vspace*{0.15in}

{\bf Main Theorem 1.14.} \ \ {\it For $\,n \,\ge\, 4\,,\,$ let $\,u_i \in C^{\, n + 4}\, (S^n)\,$ be a sequence of positive  solutions of equation\,} (1.1)\,,\, {\it with $\,{\cal K} \in C^{\, n + 4}\, (S^n)\,,\,$ and $\,v_i\,$ and $\,K\,$ be associated to $\,u_i\,$ and $\,{\cal K}\,$ via\,} (1.3)\,,\, {\it respectively}\,.\, {\it Assume that $\,\{ u_i \}\,$ has a finite number of blow\,-\,up points\, --  one of them is at the south pole, but  none at the north pole.\,  Take the following conditions\ } (1.15)\,--\,\,(1.19)\, {\it into account}\,.\\[0.125in]
(1.15) \ \  $0$ {\it is a simple blow\,-\,up point for $\,\{ v_i\}\,.$}\\[0.1in]
(1.16) \ \   ${\cal K} \,>\, 0$\, {\it in\,} $S^n$,\  {\it and}\, $\,K $ {\it is given by the Taylor expansion in\,} (1.8) \,{\it in\,} $B_o \, (\rho_o)$\,.\,  \\[0.1in]
(1.17) \ \   $(\,2 \ \le\, ) \  \ell \,\le \, n \,-\, 2$\,. \\[0.1in]
(1.18) \ \   {\it The parameters $\,\lambda_i\,$ and $\,\xi_i\,$ corresponding to the simple blow\,-\,up point at $\,0\,$ } \\[0.075in]
 \hspace*{0.51in}  \,[\,{\it via\,} (1.10) {\it and\,} (1.11)\,,\, {\it respectively}\,]\, {\it satisfy\,} (1.12)\,,\, {\it that is\,,\,} $|\,\xi_i| = o\, (\lambda_i)\,.$ \\[0.1in]
(1.19)\ \  \, {\it $\!$When \,$\ell \,=\, n \,-\,2$\,  is even and there are  more than one blow\,-\,up point}\\[0.075in]
\hspace*{0.5in}{\it or      when $\ell$ is odd,\, we   require that $\,\Delta_o^{\!(h_\ell)} \, {\bf P}_\ell\, (y) \,\equiv\, 0\,.\,
$ Here $\,h_\ell\,$ is biggest}\\[0.075in]
\hspace*{0.5in}{\it integer    less than or equal to\,} $\ell/\,2\,.\,$

\vspace*{0.1in}

{\it Then we can determine a polynomial  $\,\Gamma\,$}  ({\it constructible from $\,{\bf P}_\ell\,$ via a fixed procedure}\,)\,,\, {\it so that the following estimate holds}\, ({\it modulo a subesequence}\,)\\[0.1in]
(1.20)
\begin{eqnarray*}
  & \ &  \!\!\!\!\!\!\!\!\!\bigg\vert\   v_i \,(y) \ -\,    {\bf A}_{\lambda_i\,,\, \xi_i} \, (y) \, - \,  \left[\, \lambda^{\,\ell + 1}_i \times \Gamma   \left({\cal Y}\right)\,\right]  \!\cdot\! \large[\ {\bf A}_{\lambda_i\,,\, \xi_i} \, (y)\,\large]^{\, n \over {n \,-\,2}} \ - \ \, O_{\mbox{H}}  \! \left(\lambda_i^{{n \,-\,2}\over 2} \right)\, \bigg\vert\\[0.15in]
    & \ & \!\!\!\!\!\!\!\!\!\!\!\!=  o \!  \left(\lambda_i^{\ell  \ - \ {{n-2}\over 2}}\right) \ \     {\it for} \ \   y \,\in\, B_o\, (\rho_1) \ \ \ ( \, \rho_1 \,\le\, \rho_o \ \ {\it{is \ \ fixed}}\,)\,, \  {\it where} \ {{\cal Y} \, =\,  {{y \,-\,\xi_i}\over {\lambda_i}}}\ . \ \ \ \ \ \ \ \ \
  \end{eqnarray*}
{\it Here the term\,} $\displaystyle{O_{\mbox{H}} \! \left(\lambda_i^{{n \,-\,2}\over 2} \right)}$ {\it is defined via  the global harmonic term\,} (1.13)\,,\, {\it and its precise  expression is  found   in\,} (6.58)\,.\, (\,{\it The precise construction of   \,$\Gamma$\, is given in Proposition\,} 4.49.) \\[0.2in]
{\bf{\S\,1\,c. \ Necessity of the condition\,}} $\,\Delta_o^{\!(h_\ell)} \, {\bf P}_\ell\, (y) \equiv 0\,.\,$ \ \
At first sight the condition
\begin{eqnarray*}
(1.21) \ \ \ \ \ \ \ \ `` \  \Delta_o^{\!(h_\ell)} \, {\bf P}_\ell\, (y) & = & \Delta_o\ (\cdot \cdot \cdot\, [\,\Delta_o \,(\Delta_o  \, {\bf P}_\ell)\,]\cdot \cdot )\,(y) \  \equiv \ 0      \ \ \ \ \ \ \ \ \ \ \ \ \ \ \ \ \ \ \ \ \\
& \ & \,\leftarrow \ \ \ \ \,h_\ell\ \  \ \ \,\rightarrow
\end{eqnarray*}
for all $\,\,y \,\in\, \R^n$\," appears to be technical. Closer examination reveals that it is an integrated part of the discussion. In fact, under the conditions in Main Theorem 1.14,  when \,$\ell$\, is even, we obtain $\,\Delta_o^{\!(h_\ell)} \, {\bf P}_\ell\, (y) \equiv 0\,$ with the help of   the    Pohozaev  identity (see Proposition 6.1). The vanishing of $\,\Delta_o^{\!(h_\ell)} \, {\bf P}_\ell\,$   allows us to construct the polynomial $\,\Gamma\,$ in Main Theorem 1.14 via a reduction method,  which we begin to expound.

\newpage

{\bf{\S\,1\,d. \  Key features of the proof.}} \ \ In \cite{Compact},   Khuri,  Marques and  Schoen  introduce refined blow\,-\,up estimates for the Yamabe equation. The method is based on a second order approximation coupled with a  second order blow\,-\,up argument. We apply these methods to the scalar curvature equation (1.1), and  highlight the following  differences.\bk
The second order inhomogeneous equation is given by

\vspace*{-0.25in}

$$
\Delta_o\,  \Phi  + n \,(n + 2) \,{ A}_1^{4\over {n \,-\,2}} \cdot \Phi  = {\bf P}_\ell \cdot { A}_1^{{n + 2}\over {n \,-\,2}} \ \ \ \ \ \ \  {\mbox{in}} \ \ \R^n. \leqno (1.22)
$$

\vspace*{-0.1in}

Here $\,A_1 \,=\,A_{1\,, \, 0}$\, as given in (1.4)\,.\,
We observe that the linear operator appeared on the left hand side of (1.22)\, --\, it is used extensively in the  Lyapunov\,-\,Schmidt reduction method, see for examples \cite{Ambrosett-Scalar} \cite{Brendle-1} \cite{Brendle-2} \cite{I} \cite{II} \cite{III}\,.\, In \cite{Compact},  a solution of (1.22) is found by a linear algebra  method. The method does not disclose the precise form of the solution, which is desirable when we construct shaper estimates for simple blow\,-\,ups.\,  In this manuscript, we introduce a {\it reduction method\,,\,} which explores the recursive relations in  equation (1.22)\, (\,expounded in \,{\bf \S\,4}\,).  The condition $\,\Delta_o^{\!(h_\ell)}\, {\bf P}_\ell \,\equiv\, 0\,$ comes into the picture when we terminate the recursive process.\, As a consequence, we can determine in a step\,-\,by\,-\,step manner the exact form of the solution $\Phi\,.\,$  Although the detail is shown in \,{\bf \S\,4}\,,\, we indicate here that we know precisely what is $\,\Gamma\,$ in (1.23)\,,\, once $\,{\bf P}_\ell\,$ is given.\bk
Another unique  feature here is that the global influence from other blow\,-\,up points has to be taken into account when we  estimate in the accuracy of $\,O\,(\lambda_i^{n \,-\,2})$\,  (see \,{\bf \S\,5}\,)\,.\, To do so, we have to extend the information given by the harmonic function in (1.13) to the whole neighborhood of the blow\,-\,up point at $\,0\,,\,$ in a manner so that the second blow\,-\,up argument still works. See (6.7)\,,\, {\bf \S\,5\,b}\, and  {\bf \S\,6}\,.

 \vspace*{0.15in}


{\bf{\S\,1\,e. \ Applications\,}}: {\bf limitation on `flexibility' of simple blow\,-\,up\,,\, and locations of the blow\,-\,up points.} \ \ Consider the parameters $\,\lambda_i\,$ and $\,\xi_i$\, corresponding to the simple blow\,-\,up point at \,$0$\, via (1.10) and (1.11)\,.\, Suppose that

 \vspace*{-0.35in}

$$
   \xi_i \ = \ \lambda_i^\alpha \cdot {\vec{\,X}}  \ \  \ {\mbox{for \   a \  fixed \   vector}} \ \, {\vec{\,X}} \ \  {\mbox{and \   a \   fixed \  number }} \ \, \alpha \ > \ 0\,, \leqno (1.23)
$$

 \vspace*{-0.1in}

where $1 \, < \,\alpha \, < \, 2\,.\,$   Assume also that $\,0\,$ is the only simple blow\,-\,up point and  $\,\ell \,=\, n \,-\,2\,$.\, Then we  have

 \vspace*{-0.3in}

$$ \ \,
{\bf P}_\ell \, (\vec{\,X}) \ = \ 0 \ \ \ \ \ \    (\,{\mbox{here}} \ \ \vec{\,X} \ \ {\mbox{is \ \ considered \ \ as \ \ a \ \ point \ \ in}} \ \ \R^n)\,. \leqno (1.24)
$$

\vspace*{-0.1in}

See \,Theorem A.\,6.\,62 in the e\,-Appendix\, for the precise statement and   full layers of information available, as well as the conditions for (1.24) to hold.    In case of multiple simple blow\,-\,up points with Taylor expansions at each blow\,-\,up point  given as in (1.8), where uniformly \,$\ell\, = \,n \,-\,2\,,\,$   similar limitations exist, and they involve  the {\it locations\,} of the simple blow\,-\,up points. See {\bf \S\,A.7}\,\,\, in the e\,-\,Appendix for the exact formulas.\bk
The information  should be helpful when one seeks examples and investigates situations  with multiple simple blow\,-\,up points. Cf. \cite{Twin} on using the interaction of two close bubbles to find solutions of equation (1.1) for certain functions $\,{\cal K}$\,.\,  As a footnote,  only recently blow\,-\,up sequence with a {\it single} simple blow\,-\,up point is constructed for a fixed and none identically constant $K$ \cite{I} \cite{II} \cite{III}\,.\, Cf. also \cite{Leung-blow-up}\,,\, and \cite{Brendle-1} \cite{Brendle-2} for the Yamabe equation\,.

\vspace*{0.15in}

{\bf \S\,1\,f. \ General conditions, assumptions and conventions.}\\[0.05in]
To keep the notation clean, and without losing sight of the technical details, we assume that
$$ \,\,u_{\,i} \ \ {\mbox{and}} \ \  {\cal K} \ \ {\mbox{are \ \ in}} \ \  C^{\,n \,+\, 4} \, (S^n)\,,\, \ \ u_i \ \ {\mbox{is \ \ a  \ \ positive \ \ solution \ \ of \ \   (1.1)\,.}} \ \ \ \ \ \ \ \ \  \leqno (1.25)
$$

$$
 ``\,v_{\,i} \ \, {\mbox{and}} \ \,  K \ \, {\mbox{ descend \ \, from \   }} u \ \,  {\mbox{and}} \ \,  {\cal K} \ \, {\mbox{via\ \, (1.3)}}\,. \ \, {\mbox{Moreover,}} \ \, {\cal K} > 0 \ \, {\mbox{in}} \ \, S^n,  \ \ \ \ \ \ \ \ \  \leqno (1.26)
$$

\vspace*{-0.1in}

\hspace*{1in}{\mbox{and}} $\,({\tilde c}_n\, K)\,(0) \ = \ n\, (n \,-\,2)\,.\,"$\\[0.1in]
The degree of smoothness assumed on $\,u\,$ and $\,{\cal K}\,$ can be reduced according to the content (especially in \,{\bf \S\,2}\,)\,.

$\bullet_1$ \ \, Throughout this work,\, the dimension $\,n \,\ge \,3\,,\,$ except when otherwise is specifically mentioned, and \, ${\tilde c}_n = (n \,-\,2)/\,[\,4\,(n \,-\,1)]\,.\,$
We observe the practice on  using $\,C,\,$ possibly with sub\,-\,indices, to denote various positive constants,
which may be rendered {\it differently\,} from line to line according to the contents. {\it Whilst we use $\,{\bar c}\,$ and  $\,{\bar C},\,$ possibly with sub\,-\,index, to denote a\,}   fixed\,
{\it positive constant which always keeps the same value as it is first defined\,}.\,  \\[0.05in]
%
%
%
%
$\bullet_2$ \ \, Denote by $B_y \,(r)$ the open ball in $\,(\R^n\,, \ g_o)\,$  with center at $\,y\,$ and radius $\,r\, > \,0\,.\,$ Likewise, let $\,{\cal B}_{x}\, (\rho)\,$  be the open ball in $\,(S^n,\ g_1)\,$ with center at $\,x \,\in\, S^n\,$ and radius $\,\rho\, \in\, (0\,, \ \pi]\,.\,$ We also use the standard notation \,$\langle \ \,, \ \rangle$\, to denote the inner product in $(\R^n\,, \ g_o)\,.\,$
%
%
%
\\[0.05in]
$\bullet_3$ \ \, Given a sequence of positive numbers $\,\{\,\lambda_i\,\}\,,\,$ and a positive number $m$\,,\, we say that a sequence of numbers $\,\{\,\gamma_i\,\}\,$ satisfies

\vspace*{-0.2in}

$$
  \gamma_i \ = \ O_{\lambda_i} \,(m)\ \      \Longleftrightarrow \ \ |\,\gamma_i\,  |\ \,\le\, \   C\, \lambda_i^m \ \ \ \ \mfor\ \  i \,\gg\,  1\,. \leqno (1.27)
$$

Likewise,

\vspace*{-0.25in}

$$
\gamma_i \ = \ o_{\,\lambda_i} \,(m) \ \   \Longleftrightarrow \ \    |\,\gamma_i\,  |\ \,\le\, \ c_i\, \,\lambda_i^m \!\!\!\!\!\mfor i \,\gg\,  1\,, \ \  {\mbox{where}} \ \ c_i \,\ge\, 0 \ \ \,{\mbox{and}} \ \ c_i \,\to\, 0
$$
as $\,i \to \infty\,.\,$ The notations help to highlight the order and manage longer expressions inside the brackets.

$\bullet_4$ \ \, A statement involving a sequence is said to hold ``{\it modulo a subsequence\,}"\, if we can select a subsequence (from the original sequence in the statement) so that the statement is valid for this subsequence. As a rule, we assume  that the statement is true for the original sequence so that the notations remain clean. \\[0.3in]
{\bf \large {\bf  \S\,2. \ \ Simple blow\,-\,up.}}\\[0.1in]
{\bf \S\,2\,a. \ \  Simple blow\,-\,up and its  analytic definition\,.\,}    \ \
Intuitively, simple   blow\,-\,up develops precisely one bubble  in a neighborhood\,.\, Its  analytic definition is given by R. Schoen in  \cite{Schoen-Notes}\,.\, See also \cite{Compact} and \cite{Li-1}\,.\,
  Via a rotation, we assume without loss of generality that the blow\,-\,up point is at the south pole.   Let $\,\{ v_i \}\, $ be given as in (1.3). For a simple blow\,-\,up point,   there exists  a sequence $\,\{ \,\xi_{m_i} \, \}\, \to 0$ such that\\[0.1in]
  (2.1)
  $$
  \,{\mbox{for \   each \   }} i \,\gg\,  1, \      \xi_{m_i}\ \   {\mbox{is \   a \  {\it local}\  \,  maximum   \,  of \ }} v_i \,,\, \  {\mbox{with}} \ \,\lim_{i \to \infty}v_i \,(\xi_{m_i})\, =\, \infty\,,\,
$$

 and the rescaled average

\vspace*{-0.25in}

$$
 r \ \,\longmapsto \ \, \ r^{{n \,-\,2}\over 2} \cdot \left[ {{ \int_{ \partial B_{\,\xi_{m_i}} (r)} \ v_i\, d S }\over { \int_{ \partial B_{\,\xi_{m_i}} (r)} \ 1\, d S}}\right] \leqno (2.2)
$$

has precisely one critical point in $\,(0,\ \rho_o)\,.\,$ Here $\rho_o > 0$ is  fixed (independent on $\,v_i\,$ for $\,i\,\gg\, 1\,)\,.\,$

\vspace*{0.2in}

{\bf \S\,2\,b. \  Proportionality of simple blow\,-\,ups\,.} \ \ %
The following estimate is essentially taken from Proposition 2.3 in \cite{Li-1}\,.\, We present it in the setting of this article.

\vspace*{0.2in}

{\bf Proportionality Proposition 2.3.}\ \ {\it Under the standard conditions}\, (1.6)\,,\, (1.25) {\it and\,} (1.26)\,,\, {\it
let $\,0\,$ be a simple blow\,-\,up point for $\,\{ v_i
\}\,,\,$ and  the sequence $\, \xi_{m_i} \to 0$ carries the meaning as
in}\, (2.1) {\it and\,} (2.2).   {\it Then there exist positive constants $\,{\bar C}_1$
and $\,{\bar \rho}_o$\,}   {\it such that}

\vspace*{-0.2in}

$$
v_i \,(y) \ \,\le\, \ {{{\bar C}_1} \over{ v_i \,(\xi_{m_i})}} \,\cdot {1\over  { |\,y \,-\,\xi_{\,i}|^{\,n \,-\,2} }} \ \ \ \ \ \
{\it{for}} \ \ \ \ 0 \,< \,|\,y \,-\, \xi_{\,i}| \ \le\   {\bar \rho}_o\ \ \ {\it{and \ \ for \ \ all }} \ \ i \,\gg\,
1\,. \leqno (2.4)
$$
 {\it In addition, there is a number\, ${\bar \rho}_1 \,\in\, ( 0\,,\  {\bar \rho}_o\,)$ such that } \,({\it modulo a subsequence\,})

 \vspace*{-0.25in}

$$
[\,v_i \,(\xi_{m_i})] \cdot v_i \, (y) \ \to \  {1\over {|\, y|^{\, n \,-\,2} }} \  + \ h \, (y)\ \ \ \ \ \ \ \ \ {\it{in}} \, \ \ \ \ C^{\,2}_{\mbox{loc}} \, (\,B_o \, ({\bar \rho}_1) \setminus\, \{0\}\,)\,,\leqno (2.5)
$$
{\it where $\,h$ is a harmonic function in $\,B_o \, ({\bar \rho}_1)\,.\,$ } [{\it \,Recall that $({\tilde c}_n \,K)\,(0) \,=\, n\,(n \,-\,2)$\,.\,}]

 \vspace*{0.2in}

{\bf \S\,2\,c. \  Harmonic expression of the  collapsed part\,.} Consider a blow\,-\,up sequence of positive solutions $\,\{u_i\}\,$ of equation (1.1)\,.\,
Consider the  situations where\\[0.1in]
(2.6) \ \ ``\,the number of blow\,-\,up points is finite, say at
$
 \beta_o \ = \ {\bf S}\,, \ \cdot \cdot \cdot \,,\, \ \beta_k \,\in\, S^n \setminus \{ {\bf N} \},\,
 $\\[0.1in]
\hspace*{0.65in}and at least one of them is a simple  blow\,-\,up point (say, $\beta_{\,o}$)\,.\,"\\[0.05in]
Take a point
$$
 x_c \,\not\in\,  \{ \beta_o, \,
\cdot \cdot \cdot\,, \, \beta_k\,, \, {\bf N}\,\} \,.\, \leqno (2.7)
$$
Under the general conditions (1.25), (1.26), and also (2.6), a subsequence of

 \vspace*{-0.25in}

 $$\,\displaystyle{\left\{ \,{{u_i }\over { u_i
\,(x_c)}}\,\right\}\,}  \leqno (2.8)$$ converges to a positive $C^{\,2}$-\,function $\,{\cal H}\,$ defined on
$S^n \setminus \{ \,\beta_1, \, \cdot \cdot \cdot\,, \, \beta_k\}\,.$\, See \cite{Leung-Supported}\,.\,
With the stereographic projection $\,\dot{\cal P}\,$ onto  $\R^n$, which sends $\,{\bf N\,}$ to infinity\,,\,
  \,$\,{\cal H}\,$ can be expressed   as \,[\,cf.   the transformation in (1.3)\,]

  \vspace*{-0.35in}

\begin{eqnarray*}
(2.9) \ \ \ \ \ \ \ \ \    H \,(y) &:= &[\,{\cal H} \circ {\dot{\cal P}}^{-1} (y)] \cdot \left({2\over {1 + |\,y|^{\,2}}} \right)^{\!\!{{n \,-\,2}\over 2}}\,,\\
(2.10) \ \ \ \ \ \ \   \,H\, (y)  &= &\sum_{j = 0}^k \ {{A_{\,j}}\over { |\, y \,-\,{{\hat {\!Y}}}_{\!j}|^{\,n \,-\,2} }} \ \ \ \ \ \ {\mbox{ for}} \ \  \  y \,\in\, \R^n \setminus
 \{\ {{\hat {\!Y}}}_o\,,\,\cdot \cdot \cdot\,, \ {{\hat {\!Y}}}_k \}\,.\,
 \ \ \ \  \ \ \ \ \ \ \ \ \ \  \ \ \ \ \ \ \ \ \ \
\end{eqnarray*}

\vspace*{-0.25in}

 $$\!\!\!\!\!\!\!\!\!\!\!\!\!\!\!\!\!\!\!\!
 {\mbox{Here}}   \ \ \ \ \ \ {{\hat {\!Y}}}_{\!j}:= {\cal P} \,(\,\beta_j) \ \ \ \mfor \ \ 0 \,\le\, j \,\le \,k\,,\,  \ \ \ \ \ \ \ \ \ \ \ \ \ \ \ \ \ \      \leqno (2.11)
 $$
 and $\,A_{\,j}\,$ are positive numbers\,.\, Refer to \,{\bf \S\,4\,} in \cite{Leung-Supported}\,.\,\bk
The convergence can be quantified in the following manner. Given a sequence of positive numbers $\varepsilon_j \downarrow 0\,,\,$ and a sequence of compact sets $\{ \,{\cal C}_{\,j} \}$\, such that

 \vspace*{-0.25in}

 $$
 {\cal C}_{\,1} \,\subset\,   {\cal C}_{\,2} \,\subset\,  \cdot \cdot \cdot \ \,,\ \ \ \ \ \ \ \ \bigcup_{j\,=\, 0}^\infty\  {\cal C}_{\,j} \ = \ \R^n \setminus
 \{\ {{\hat {\!Y}}}_1\,,\,\cdot \cdot \cdot\,, \ {{\hat {\!Y}}}_k \}\,, \leqno (2.12)
 $$

  \vspace*{-0.05in}

  there exists a sequence of natural numbers $\,N_i \,\uparrow\, \infty\,$ so that
 $$
 \bigg\vert \  v_i\, (y) \,-\,[\,u_i\, (x_c)\,] \cdot H\, (y) \bigg\vert \ \,\le\, \ \varepsilon_j \cdot  [\,u_i\, (x_c)\,] \ \ \ \ \mfor i \,\ge\, N_j \ \ \ \ {\mbox{and}} \ \ y\, \in\, {\cal C}_{\,j}\ .\leqno (2.13)
 $$
 We point out when $i \to \infty\,,\,$
\begin{eqnarray*}
(2.14)     & \ &    \!\!\!\! \!\!\!\! \!\!\!\!{\mbox{ ``\,right \ hand \ side  \ of \ (2.13)\,"}} \ = \  \varepsilon_j \cdot  [\,u_i\, (x_c)]  \ \to\  0\,,  \ \ \\[0.1in]  (2.15)    & \ &    \!\!\!\! \!\!\!\!\!\!\!\!\!{\mbox{ \,\,the \  domain \ in \ which \ (2.13) \ holds \ is }} \    {\cal C}_{\,j}  \,``\to"  \, \R^n \setminus
 \{\, {{\hat {\!Y}}}_1\,,\,\cdot \cdot \cdot\,, \, {{\hat {\!Y}}}_k \}\,. \ \ \ \  \ \ \ \ \ \
\end{eqnarray*}
Cf. (2.22) and (2.23) in the \,{\bf \S\,2\,d\,.}

\vspace*{0.15in}

{\bf \S\,2\,c.1.} \  \ {\it Change of the base point.} \ \ We observe that, in (2.8), one can replace the base point $\,u_i\, (x_c)\,$ by a sequence of numbers $\,\{\,\gamma_i\,\}\,$ so that
$$
C^{-1}\cdot \gamma_i \ \,\le\, \ u_i\, (x_c)\ \,\le\, \ C\, \gamma_i \ \ \ \ \mfor i \ \,\gg\,  \ 1\ . \leqno (2.16)
$$
 A subsequence of $\ \displaystyle{\left\{ \, \gamma_i^{-1}\ \cdot u_i   \,\right\}\,} $ converges to a positive $\,C^{\,2}$\,-\,function $\,\tilde{\cal H}\,$ defined on
$S^n \,\setminus\, \{ \,\beta_o, \, \cdot \cdot \cdot\,, \, \beta_k\}\,.$\,
With the stereographic projection $\,\dot{\cal P}\,$ onto  $\R^n$,
  \,$\,\tilde{\cal H}\,$ can be expressed   as  in \,(2.10) and (2.11), with a scaling factor $\displaystyle{\,\lim_{\i \to \infty} \ \gamma_i^{-1} \cdot u_i\, (x_c)}\,$\, inserted\,.\,

\vspace*{0.15in}

{\bf \S\,2\,d. \ Renormalization and first order approximation\,}.
Let $0$ be a simple blow\,-\,up point for the sequence of positive solutions $\,\{v_i\}\,$ of equation (1.2)\,.\, With the notations in (2.1) and (2.2), define

\vspace*{-0.2in}

$$
{\cal V}_{\,i}\, ({\cal Y})\ :=\ {{v_i \,(\,\xi_{m_i} + \lambda_{m_i}\cdot{\cal Y} )}\over {v_i \,(\,\xi_{m_i})}}  \ \ \mfor \ \  {\cal Y} \,\in\, \R^n \ \ \ \ {\mbox{with}}\  \ \ \lambda_{m_i} \cdot {\cal Y}  \,\in\, B_o \,(\,\rho_o)\,. \leqno (2.17)
$$

\vspace*{-0.05in}

where $\lambda_{m_i} \ := \ [\,v_i \,(\,\xi_{m_i} ]^{-{2\over {n \,-\,2}}}\ .\,$\, Here \,${\cal V}_{\,i}\,$ satisfies the equation (extendable to $\,\R^n$\,)

 \vspace*{-0.2in}

$$
\Delta_o \,{\cal V}_{\,i} \,+\, \left[ \, ({\tilde c}_n \,K)\,(\,\xi_{m_i}  \,+\, \lambda_{m_i} \cdot {\cal Y}  )  \right] \,{\cal V}_{\,i}^{{n + 2}\over {n \,-\,2}} \ = \,0 \ \ \ \  {\mbox{in}}  \ \ \ B_o \,(\lambda_{m_i}^{-1} \cdot \,\rho_o) \leqno (2.18)
$$

Assuming (1.6), under  the conditions (1.25) and (1.26), we invoke Proposition 2.1 in \cite{Li-1}  (pp. $\!$333)   to conclude that, modulo a subsequence,  $\,\{ {\cal V}_{\,i} \}\,$ converges to $\,A_1\,=\, A_{1\,,\,\,0}\,$ as given in (1.4)\footnote{\,The job is made easier as we explain in {\bf \S\,2\,f}\,,\, for simple blow\,-\,up, we can take $\,\xi_i\,$ to be a global maximum point of $\,v_i\,$ in $\,\overline{B_o\, (\rho_o)\!}\ $\,.}.  Cf. \cite{Caffarelli-Gidas-Spruck} and \cite{Gidas-Ni-Nirenberg.1}\,.\,
The convergence happens in $\,C^1$\,-\,sense, uniformly in compact subsets in $\,\R^n$\, (for the  variable $\,{\cal Y}$\,)\,.\, This translates into a weak approximation of $\,v_i\,$,\, which can be described in the following manner. Given sequences of positive numbers $\,\{ \varepsilon_i \}\,$ and $\{ R_i \}$ with
 $\varepsilon_i \ \downarrow \ 0\,$ and  $\,R_i \ \uparrow \ \infty\,,\,$    via the Cantor diagonal argument on subsequences\,,\,  we have
$$
  |\,{\cal V}_{\,i} \, ({\cal Y}) \,-\,  A_1\, ({\cal Y})|\ \,\le\, \ \varepsilon_i   \leqno (2.19)
$$
for all $\,{\cal Y} \,\in\, B_o \, (R_i)\,$ and $\, i \,\gg\, 1\,$ (modulo a subsequence)\,.
Moreover, by choosing $R_i$ to be {\it smaller}\, if necessary, we can take it that
$$
\varepsilon_i \cdot \!R_i^{\,2(n -1)} \ \to \ 0  \ \ \ \ \,{\mbox{and}} \ \ \ \ \lambda_i \cdot \!R_i \   \to \ 0 \ \ \ \ {\mbox{as}} \ \ \,i \to \infty\,. \leqno (2.20)
$$
[\,See also \,\S\,3\,a\, in \cite{Leung-Supported}\,,\, and the proof of Proposition A.6.34 in the e\,-\,Appendix for the application of (2.20)\,] Via the change of variables
$$
y  \ = \ \xi_{m_i} + \lambda_{m_i} \, {\cal Y}\,; \ \  \ \ \ \ {\cal Y} \,\in\, B_o\, (R_i) \ \ \Longleftrightarrow \ \ y  \,\in\, B_{\, \xi_{m_i}} \, (\lambda_{m_i} \cdot R_i)\,,\,
$$
 (2.17), (2.19) and (1.4)  yield
$$
| \ v_i \, (y)   \,-\, {\bf A}_{\lambda_{m_i}\,,\  \xi_{m_i}\,} (y)\  |  \ \,\le\, \ {{\varepsilon_i}\over { \lambda_{m_i}^{{n \,-\,2}\over 2}   }} \ \  \mfor \ |\,y \,-\, \xi_{m_i}|\,\le \, \lambda_{m_i} \cdot R_i  \ \ \,{\mbox{and}} \ \ i \,\gg\, 1\,. \leqno (2.21)
$$

\vspace*{-0.1in}

Cf. (1.5)\,.\, However, we do not know, a priori, how small we can take $\varepsilon_i$ (relative to $\lambda_i$\,) and how large we can choose $\,R_i\,$ (relative to $\lambda_i^{-1}$\,)\,.\, In particular, the following scenario can occur.
\begin{eqnarray*}
(2.22) \ \ \   & \ &  \!\!\!{\mbox{``right \ hand \ side  \ of \ (2.21)\,"}} \ = \   \varepsilon_i \cdot \lambda_{m_i}^{-\,{{n \,-\,2}\over 2}   } \   \,\to\, \infty\,,\ \ \ \ \  \\[0.15in]  (2.23) \ \ \      & \ & \!\!\!\!\!  {\mbox{ \,the \  radius \   of \   \ the  \ ball \ in \ which \ (2.21) \ holds }} \, =  \,     \lambda_{m_i} \cdot R_i \to 0\,. \ \ \ \ \ \ \ \ \ \ \ \ \ \ \
\end{eqnarray*}
Our goal is to introduce bubble estimates that are accurate up  to      $\,O \,(\lambda_{m_i}^\tau)\,$ for $\,\tau\, >\, 0\,$ (as big as possible)\,,\, and to ``\,stabilize"  the domain in which the estimates hold.

 \vspace*{0.15in}

{\bf \S\,2\,d.1.} \ \ {\it Joint between shrinking bubble estimate and the expanding global harmonic  term.}\ \
As mentioned, there  are diametric contrasts between bubble estimate (2.21) and the global harmonic estimate presented in (2.13)\,.\, Adding to the list, observe that  $\Delta_o\, H = 0$\,,\, whereas $$\Delta_o \, {\bf A}_{\lambda_{m_i}\,,\,\, \xi_{m_i}} \, = \, -\,n\,(n \,-\,2)\, [\,{\bf A}_{\lambda_{m_i}\,,\, \,\xi_{m_i}\,}]^{\,{{n + 2}\over {n \,-\,2}}}\   \ (\,< 0)\,.$$   These two estimates do not immediately link to each other. We demonstrate their intricate relation when we present  estimates that are accurate up  to order $\,O_{\lambda_{m_i}} (n \,-\,2)\,.\,$

\vspace*{0.2in}

{\bf \S\,2\,e. \ An equivalent geometric expression\,.}\ \
Before we proceed to a closer relation between ${\cal V}_{\,i}\,$ and $\,A_1\,$,\, we examine a simpler  estimate  here. Not only the estimate is useful in later discussion,  it is interesting in its own right. As for the proof, we present it in \,{\bf \S\,A.1}\, in the e\,-Appendix.

{\bf Proposition 2.24.} \ \ {\it Under the standard conditions in\,} (1.6), (1.25) {\it and}\, (1.26)\,,\, {\it modulo a subsequence,\, $0$ is a simple blow\,-\,up point for \,$\{v_i\}$\, if and only if there exist a sequence\,} $\,\zeta_i\, \in\, \R^n$,\, {\it with }

\vspace*{-0.25in}

$$
\zeta_i\,  \to\,  0 \ \ \ \ \,  {\it{and}} \ \ \ \ \epsilon_i \, := \, {1\over{[\, v_i\, (\zeta_i)]^{\, {2\over {n \,-\,2}} }}} \ \to \ 0\,, \ \ \ \ \ \ \ {\it so \ \  that}\leqno (2.25)
$$

\vspace*{-0.25in}

$$
 {1\over C}\cdot  {\bf A}_{\,\epsilon_i\,,\,\, \zeta_i} \, (y)  \ \,\le\, \ v_i\, (y) \ \,\le\, \ C\cdot {\bf A}_{\,\epsilon_i\,,\,\, \zeta_i} \, (y)   \ \ \ \ \ \ \ \ {\it for \ \ all} \ \ \  |\,y \,-\,\zeta_i|  \ \,\le\, \ \rho_1\,. \leqno (2.26)
$$
{\it Here $\,C\, \ge \,1$\,  and $\,\rho_1$ are  positive constants  independent on $\,i\,.\,$}

\vspace*{0.15in}

{\bf \S\,2\,f. \ Shifting to the maximal point\,.} \ \ Let $\,\xi_i \, \,\in\, \, B_o \,(\rho_o)\,$ be given in (1.10)\,.\,
 We can
 take $\,\xi_{m_i} = \xi_i\,$ in (2.1) and (2.2). Moreover, suppose that there exists another sequence of points $\,\{ {\tilde{\xi}}_i \}\,$ which also satisfies (2.1) and (2.2) in the definition of simple blow\,-\,up points\,.\, We have  (modulo a subsequence)
$$
\ \ \ \ \ \   \ \ \ \ \   |\,{\tilde{\xi}}_i \,-\, \xi_i| \,= o\, (\lambda_i) \ \ \ \ \ \ \ \ \ \ \ \ \ \ \ \ \ \  [\,\lambda_i \ \ {\mbox{as\ \ in \ \, }} (1.11)\,]\leqno (2.27)
$$
The proofs of the above statements, which require only standard techniques, can be found in \,{\bf \S\,A.2}\, in the e\,-\,Appendix.

\vspace*{0.2in}

{\bf \S\,2\,g. \  Non\,-\,degenerate  conditions and $\,o\, (\lambda_i)\,$ restriction on flexibility.}\\[0.05in]
Non\,-\,vanishing derivatives at the blow\,-\,up point tend  to post restriction on the blow\,-\,up flexibility. One good example can be find in    \cite{Chen-Lin-3}\,,\, which  we highlight here\,,\, using  the setting of the present article. Via Taylor's expansion,
$$
(\,{\tilde c}_n\,K)\, (y) = n \,(n \,-\,2) \ +\  [\,-\,{\bf P}_\ell\,(y)] \ + \  R_{\,\ell+1}\,(y) \mfor y\in B_o \,(\rho)\,. \leqno (2.28)
$$
Here we use multi-index $\,\alpha =(\alpha_1\,, \ \cdot \cdot \cdot\,, \alpha_n)\,$,\, and
$$
{\bf P}_\ell\, (y) \ = \    \sum_{|\,{ \alpha}| = \ell}\, \left[\,   D^{\,(\ell)}_\alpha \,(-\,{\tilde c}_n \, K)\bigg\vert_{\,y\, =\,0} \cdot {{y^{\,\alpha} }\over {\alpha\,!}} \right]\,,\leqno (2.29)
$$

\vspace*{-0.15in}

$$
R_{\,\ell+1} \,(y) \ =\  O \,  \left(\, \max_{B_o (\,\rho_o + \varepsilon')} |\btd^{(\ell + 1)} \,K|  \times |\, y  |^{\,\ell + 1} \right)\,.
\leqno (2.30)
$$
[\,Negative sign is introduced for later matching. See (3.7).]
One can verify that
$$
{{|\  R_{\ell+1} \,(y)|}\over {|\,y |^\ell}} \ \to\  0 \ \ \ \ {\mbox{and}} \ \ \ \ {{\Vert \btd\, R_{\ell+1} \,(y)\Vert }\over {|\,y |^{\ell-1}}}\ \to\ 0\ \ \ \ {\mbox{as}} \ \ |\,y| \ \to\  0\,. \leqno (2.31)
$$
A more demanding condition is the lower bound \\[0.1in]
(2.32)
$$
C^{-1} \, |\,y|^{\,\ell \,-\,1}     \,\le\, \!  \sqrt{ \left( {{\partial\, {\bf P}_\ell \, (y)} \over {\partial y_{|_1}}} \right)^{\,2} + \ \cdot \cdot \cdot \ + \left( {{\partial \,{\bf P}_\ell \, (y)} \over {\partial y_{|_n}}} \right)^{2}\,} \   \,\left(   \,\le\, \   C  \, |\,y|^{\,\ell \,-\,1} \right)
$$
for  $\, y\,\in\, B_o \,(\rho)\,.$\,
Cf. the example below\,.\, From (2.32), we have
$$
0 < c\, (\varepsilon) \,\le\, \Vert\, \btd K \, (y)\Vert \,\le\, C \mfor \varepsilon \,\le\, |\, y| \,\le\, \rho\,, \ \ {\mbox{where}} \ \ \rho \ \ {\mbox{is \ \ small \ \ enough}}\,. \leqno (2.33)
$$
The following is a direct application of Lemma 3.6 in \cite{Chen-Lin-3}, after  checking  (1.2), (1.6), and  (3.2), and the conditions stated at the beginning of \,{\bf \S\,3\,} in \cite{Chen-Lin-3} (in particular, $\alpha_i \,\le\, n \,-\,2$\,, pp. 127\,,\, loc. cit.)\,,\, also verifying the conditions stated in Lemma 3.4 and Lemma 3.6 (loc. cit.),\, and  taking $\,p_i\, \equiv \,{{n + 2}\over {n \,-\, 2}}\,.\,$

\vspace*{0.15in}

{\bf Proposition 2.34.} \ \ {\it Granted the general conditions in\,} (1.6)\,,\, (1.25) {\it and\,} (1.26)\,,\, {\it suppose that $\,0$\, is a simple blow\,-\,up point for} $\{v_i\}\,.\,$ {\it Assume also\,} (2.28) {\it and\,} (2.32) {\it \,for $\,2 \,\le\, \ell \,\le\, n \,-\,2\,.\,$ If }
$$
\int_{\R^n} \btd\, {\bf P}_\ell \, (  y + {\cal X} )\, [\,A_1\,(y)]^{{2n}\over {n \,-\,2}}\ dy \  \not= \ \vec{\,0} \ \ \ \ {\it for \ \ all} \ \ {\cal X} \,\in\, \R^n \setminus \{ 0\}\,, \leqno (2.35)
 $$
 {\it then\,,\,} {\it modulo a subsequence\,,\, we have\,}  $|\,\xi_{\,i}| = o \, ( \lambda_i\,)\,$ [\,{\it Recall that $\,A_1\,=\,A_{1\,,\ 0}\,$ is given in\,} (1.4)\,,\, {\it  $\,\xi_i\,$ fulfills\,} (1.10)\,,\, {\it  $\lambda_i$ is given in}\,  (1.11)\,,\, {\it and $\,{\bf P}_\ell\,$ in\,} (2.29)\,.\,] \\[0.25in]
{\bf \S\,2\,g.1} \ {\it Examples on  $\,K\,$ with local expansions  fulfilling  the conditions $\Delta_o^{\!h_\ell} \, {\bf P}_\ell \,\equiv \, 0$\,,\,}
(2.32) {\it and\,} (2.35)\,.\, \ \ Recall that \,$h_\ell$\, is defined as the largest integer that is less than or equal to $\ell/2\,.\,$  Consider $\,n\,$ and $\,\ell\, \ge \,2\,,\,$ both even numbers\,,\, and
$$
({\tilde c}_n\, K)\, (y) \ = \ n\, (n \,-\,2) \ + \  \left[\,\left(y_{|_1}^\ell \,-\,y_{|_2}^\ell\right) \  + \ \cdot \cdot \cdot \ + \ \left(y_{|_{n \,-\,1}}^\ell \,-\,y_{|_n}^\ell\right)\,\right] \leqno (2.36)
$$
for $\, y \,\in\, B_o\, (\rho_o)\,.\,$
Using H\"older's inequality, one can verify (2.32)\,.\, Moreover,
\begin{eqnarray*}
 \!\!(y_{|_1} + {\cal X}_1)^{\ell \,-\,1} \!\!\!& = &\!\!\! y_1^{\ell -1} +  C(\ell \,-\,1\,,\, 2)  \cdot  y_1^{\ell -\,1 \,-\,2} \,{\cal X}_1^2  +  \cdot \cdot \cdot  +  C(\ell \,-\,1\,,\, \ell -2)\cdot  y_1^{\ell -1 \,-\,2} \,{\cal X}_1^{\ell \,-\,2}\ \ \\[0.15in]
 & \ & \!\!\!\!\!\!\!\!\!\!\!\!\!\!\!\!\!\!\!\!\!\!+\,{\cal X}_1 \cdot \!\left[ \,C(\ell \,-\,1\,,\, 2)\cdot  y_1^{\ell -1 \,-\,1}  \,+\, \cdot \cdot \cdot \,+\, C(\ell \,-\,1\,,\, \ell -3)\cdot y_1^{2} \cdot {\cal X}_1^{\ell \,-\,4} \,+\, {\cal X}_1^{\ell \,-\,2}\ \right].
\end{eqnarray*}

\vspace*{-0.3in}

$$
{\cal X} = ({\cal X}_1\,,\ \cdot \cdot \cdot\, \ {\cal X}_n)\,,  \ \ \ \ \ {\mbox{and}} \ \ \ \ C({j\,,\, \,k}) \ = \ {{j\,! }\over {(j \,-\,k)\,!\,\,k\,!}} \ \ \ \ \ \ \ (j \ \ge \ k) \leqno {\mbox{Here}}
$$
is the binomial coefficient.
Note that all the powers in ${\cal X}_1$ inside the  brackets are even   numbers.
As
$$
\int_{\R^n} y_{|_1}^{2j + 1} \ [\,A_1\,(y)]^{{2n}\over {n \,-\,2}}\, dy \ = \ 0\,,
$$
we have
$$
\int_{\R^n} (y_{|_1} + {\cal X}_1)^{\ell \,-\,1} \ [\,A_1\,(y)]^{{2n}\over {n \,-\,2}}\, dy \ = \ 0 \ \ \ \ \Longleftrightarrow \ \ \ \ {\cal X}_1 \ = \ 0\,.
$$
It follows that (2.35) is fulfilled with the form in (2.36)\,.\, In addition, observe   that
$$
\Delta_o^{\!h_\ell}\, \left[\,\left(y_{|_1}^\ell \,-\,y_{|_2}^\ell\right) \  + \ \cdot \cdot \cdot \ + \ \left(y_{|_{n \,-\,1}}^\ell \,-\,y_{|_n}^\ell\right)\,\right]  \ = \ 0 \ \ \ \ \ \ \ \ (\ell \ \ {\mbox{being \ \ even}})\,.
$$
One can generalized (2.36) by introducing positive multipliers onto each $\,\left(y_{|_{2j-1}}^\ell \,-\, y_{|_{2j  }}^\ell\right)$\,.\,

\newpage

{\bf \large {\bf  \S\,3. \ \  Difference between the normalization $\,{\cal V}_{\,i}\,$ and  $\,A_{\,1}\,.$}}  \\[0.15in]
 After shifting from $\,\xi_{m_i}\,$ to $\,\xi_i\,$\, as described in {\bf \S\,2\,f}\,,\, for the sake of simplicity, we continue to use the notation

 \vspace*{-0.3in}

 $$
 {\cal V}_i \, ({\cal Y}) \ := \ {{  v_i\, (\,\xi_i + \lambda_i\, {\cal Y})  }\over {  M_i }} \ \ \  \mfor {\cal Y} \,\in\, \R^n. \leqno (3.1)
 $$

  \vspace*{-0.1in}

Here

 \vspace*{-0.2in}

$$
M_i \ := \ v_i\, (\xi_i) \ \ \  \ {\mbox{and}} \ \ \ \ \lambda_i \ = \ M_i^{\,- {2\over {n \,-\,2}}   }\,, \ \ \  \ \xi_i \ \ {\mbox{is \ \ given  \ \ in \ \ (1.10)}}\,. \leqno (3.2)
$$

 Cf. (2.17) and \,{\bf \S\,2\,f\,}.\, As $\,A_1\,=\,A_{1\,, \ 0}\,$ satisfies the equation
$$
\Delta_o\, A_1 + n \,(n \,-\,2) \,A_1^{{n + 2}\over {n \,-\,2}} = 0 \ \ \ \ {\mbox{in}} \ \ \R^n, \leqno (3.3)
$$
together with equation (2.18)\,,\, which holds after the changes $\,\xi_{m_i}\, \to \,\xi_i\,$  and  $\,\lambda_{m_i}\, \to \,\lambda_i\,,\,$ it can be  seen that
\begin{eqnarray*}
(3.4)   & \ & \Delta_o \,({\cal V}_{\,i} \,-\,A_1)\, ({\cal Y} )\\[0.1in] & \  & \!\!\!\!\!\!\!\!\!\!   =\    n \,(n \,-\,2)  \!\left\{\, [\,A_1\, ({\cal Y} )]^{{n + 2}\over {n \,-\,2}} \,-\,[\,{\cal V}_{\,i}\, ({\cal Y} )]^{{n + 2}\over {n \,-\,2}}\,\right\} \\[0.1in]
& \ & \ \ \ \ +  \ \left[\,n \,(n \,-\,2) \,-\,{\tilde c}_n  \, K( \lambda_i\,{\cal Y} + \xi_{\,i})\, \right] [\,A_1\, ({\cal Y} )]^{{n + 2}\over {n \,-\,2}}\ \\[0.1in]
& \ &    \ \ \ \ \ \ \ \ \  +\, \left[\,n \,(n \,-\,2) \,-\,{\tilde c}_n \, K( \lambda_i\,{\cal Y} + \xi_{\,i})\, \right]  \left\{\,[\,{\cal V}_{\,i} \, ({\cal Y} )]^{{n + 2}\over {n \,-\,2}} \,-\,[\,A_1\,({\cal Y} )]^{{n + 2}\over {n \,-\,2}}\right\}
\end{eqnarray*}
for $\,{\cal Y} \,\in\, \R^n\,$\,.

\vspace*{0.1in}

{\bf \S\,3\,a.  \  Linear approximation to}\, $\left( A_1^{{n + 2}\over {n \,-\,2}} \,-\,\,{\cal V}_{\,i}^{{n + 2}\over {n \,-\,2}}  \right)\,$ {\bf in case of simple blow\,-\,up\,.}\ \
It follows from Proposition 2.24  that   (see \,{\bf \S\,A.10}\, in the e\,-\,Appendix for details)
\begin{eqnarray*}
(3.5)   & \ &\!\!\!\!\!\!|
[\,A_1\, ({\cal Y})]^{{n + 2}\over {n \,-\,2}} \ -\ [\,{\cal V}_1\, ({\cal Y})]^{{n + 2}\over {n \,-\,2}}   =  \left(\, {{n+2}\over {n \,-\,2}} \right)  [\,A_1\, ({\cal Y})]^{4\over {n \,-\,2}} \cdot  [\, A_1 \, ({\cal Y}) \,-\,{\cal V}_{\,i} \, \, ({\cal Y})] \\[0.1in]
 & \ &  \ \ \ \ \ \ \ \ \ \ \ \ \ \  \ \ \ \ \ \ \ \ \   \ \ \ \ \    + \ O \,(1) \,[\, A_1 \, ({\cal Y}) \,-\,{\cal V}_{\,i} \, \, ({\cal Y})]^{\,2}\cdot  [\,A_1\, ({\cal Y})]^{\,{{4\over {n \,-\,2}} -1}}
\end{eqnarray*}
for $\,|\,{\cal Y}|\ \,\le\, \ \rho_o\, \lambda_i^{-1}\,$.

\newpage

{\bf \S\,3\,b.   \ Taylor expansion of $\,\,(\,{\tilde c}_n\,K)\,$ above $\,0$\,.}\ \
Since we know that $\btd \,K (0) = 0$\,,\, and by (1.6), $({\tilde c}_n \,K) \,(0) = n\,(n \,-\,2)$\,,\, we assume that all the derivatives of $K$ vanish at $0$ up to (and equal to) order $\,\ell \,-\,1\,.\,$ Here\, $\ell \,\ge\, 2$\, is an integer. Using multi\,-\,index $\,\alpha = (\alpha_1\,, \ \cdot \cdot \cdot\,, \ \alpha_n)$\,,\, and the Taylor expansion of $({\tilde c}_n \,K)$\,,\, we obtain
\begin{eqnarray*}
 (3.6) \!\!\!\!\!& \ & n\,(n \,-\,2) \,-\,\,{\tilde c}_n \,K (\lambda_i\,{\cal Y} + \xi_{\,i})\
 =\ \sum_{|\, \alpha | = \ell} D^{ (\ell)}_\alpha \, (-\,{\tilde c}_n  K) \,\bigg\vert_{\,0} \cdot {{ ( \,\lambda_i \,{\cal Y} + \xi_{\,i})^{\, \alpha} }\over {\alpha\,!}} \ + \ \ \ \ \ \ \ \ \ \ \ \\[0.1in]
  & \ & \ \ \ \ \ \  \ \ \ \ \ \ \ \ \
  \  +\  O \,(1) \left[\,  \max_{|\,\lambda_i \, {\cal Y} \,+\, \xi_i|\,\le \,\rho_o^+}  \Vert \ \btd^{(\ell + 1)} \,K\Vert\  \right]\cdot | \,\lambda_i \,{\cal Y} + \xi_{\,i}|^{\,\ell \,+ 1}\\[0.15in]
 & \ &  \ \ \ \ \ \ \ \ \ \ \ \ \  = \ \lambda_i^\ell\cdot {\bf P}_\ell \, ({\cal Y})  \ + \ \sum_{k = 1}^\ell O \, \left( |\,\xi_i|^{\,k} \cdot (\lambda_i |\,{\cal Y}|\,)^{\, \ell \,-\, k}\,\right)  \ + \ {\cal R}_{\,3} \,({\cal Y})  \ \ \ \
 \end{eqnarray*}
 for $\, \lambda_i \cdot |\,{\cal Y}|\,\le \, \rho_o\,,\ \ i \,\gg \,1\,$.
 In the above $\rho_o^+\,$ is slightly bigger than $\rho_o\,.\,$
Moreover\,,
\begin{eqnarray*}
(3.7) \ \ \ \  \ \ \ \ \ \ \ \ \ \ \ \ {\bf P}_\ell\, ({\cal Y}) \!& = &\!   \sum_{|\,{ \alpha}\,| = \ell}\, \left[\,   D^{\,(\ell\,)}_\alpha \,(-\,{\tilde c}_n \, K)\bigg\vert_{\, 0} \cdot {{{\cal Y}^{\,\alpha} }\over {\alpha\,!}} \,\right]\,, \ \ \ \ \ \ \ \ \ \ \ \ \ \ \ \ \ \ \  \ \ \ \ \ \ \ \ \  \ \ \ \ \ \ \ \ \ \ \ \ \ \ \ \ \ \  \\[0.1in]
(3.8) \ \ \ \ \ \ \ \ \ \ \ \ \  \ \ {\cal R}_{\,3} \,({\cal Y})\!& = &\! O \,  \left(\, \max_{|\,\lambda_i \, {\cal Y} \,+\, \xi_i|\,\le\, \rho_o^+ } \Vert \ \btd^{(\ell + 1)} \,K \Vert  \cdot  \! |\,\lambda_i\,{\cal Y} + \xi_{\,i}|^{\,\ell \,+ 1} \right)\,. \ \ \ \  \ \ \ \  \ \ \ \ \ \ \ \ \ \  \ \ \ \  \ \ \ \ \ \ \ \ \ \  \ \ \ \  \ \ \ \ \ \ \ \ \ \
 \end{eqnarray*}

\vspace*{0.2in}

{\bf \S\,3\,c.  \  The mixed term\,.}\ \
Consider the last term in (3.4). Using Taylor expansion as in (3.6), and the inequality
 $$
 a > b > 0 \ \ \ {\mbox{and}} \, \ \  p \ge 1  \ \ \Longrightarrow \ \ a^p -\, b^p \ \,\le\, \ {1\over p} \cdot (a \,-\,b) \cdot  a^{\,p \,-\,1}\ , \leqno (3.9)
 $$

 \vspace*{-0.1in}

we obtain

 \vspace*{-0.25in}

\begin{eqnarray*}
(3.10) \ \   & \ & \left[\,n \,(n \,-\,2) \,-\,{\tilde c}_n \, K( \lambda_i\,{\cal Y} + \xi_{\,i})\, \right]  \left\{\,[\,{\cal V}_{\,i} \, ({\cal Y} )]^{{n + 2}\over {n \,-\,2}} \,-\,[\,A_1\,({\cal Y} )]^{{n + 2}\over {n \,-\,2}}\right\} \\[0.1in]
&  \ & \!\!\!\!\!\!\!\!\!\! =\
  O \,  \left(\, \max_{|\,\lambda_i \, {\cal Y} \,+\, \xi_i|\,\le \,\rho_o^+)}  \Vert \, \btd^{(\ell)}  K \, \Vert \cdot |\,\lambda_i\, {\cal Y} + \xi_{\,i}|^{\,\ell } \right) \, {\bf \times}\\[0.1in]
 & \ &  \ \ \ \  \times  \left[\,\,O \,(1) \, |\,{\cal V}_{\,i} \,-\,V| \times   \max \, \left\{ \ [\,{\cal V}_{\,i} \, ({\cal Y} )]^{{4}\over {n \,-\,2}}\ , \  [\,A_1\,({\cal Y} )]^{{4}\over {n \,-\,2}}  \right\} \ \right]  \ \ \ \ \ \ \ \ \ \
\end{eqnarray*}
for $\,|\,{\cal Y} |\ \le\  \lambda_i^{-1}\ \rho_o\,.$

 \vspace*{0.2in}

{\bf \S\,3\,d.  \ Isolating the key  terms and the remainder\,.}\ \
It follows from (3.4), (3.5), (3.6) and (3.10) that
$$ \Delta_o \,[\,{\cal V}_{\,i} \,-\,A_1\,]\, ({\cal Y}) + n\,(n + 2)  \,A_1^{4\over {n \,-\,2}} \,[\ {\cal V}_{\,i} \,-\,A_1\,]\, ({\cal Y}) = \lambda_i^\ell \cdot {\bf P}_\ell \, ({\cal Y}) \cdot \,A_1^{{n + 2}\over {n \,-\,2}} \ +    \  {\bf RM}\, ({\cal Y}) \leqno  (3.11)
$$
$$
\hspace*{1.8in}{\mbox{for}} \ \ \ |\,{\cal Y} | \ \,\le\, \ \lambda_i^{-1}\ \rho_o \ \ \ \ \ \ \ \left( {\cal Y} \ = \ {{ y \,-\,\xi_i}\over {\lambda_i}} \right)\,.
$$

\vspace*{-0.1in}

Here (refer to \,{\bf \S\,6\,b}\, \,and \, {\bf \S\,6\,c}\,)
 \begin{eqnarray*}
 (3.12)  \ \  \ \   {\bf RM} &= & {\bf RM}_1 \ + \  {\bf RM}_{\,1} \ + \ {\bf RM}_3 \ + \ {\bf RM}_{\,4}\,, \\[0.15in]
{\bf RM}_{\,1} \, ({\cal Y})& ``="&  \sum_{k = 1}^\ell \ O \, \left( |\,\xi_i|^{\,k} \cdot (\, \lambda_i \, |\,{\cal Y}|)^{\,\ell \,-\,k} \right) \,,\\[0.1in]
{\bf RM}_{\,2} \, ({\cal Y}) & ``="&  O \,  \left(\, \max_{|\,\lambda_i \, y \,+\, \xi_i|\, \,\le\, \,\rho_o^+}\, \Vert \ \btd^{(\ell + 1)} \,K \Vert   \times |\,\lambda_i\,{\cal Y} + \xi_{\,i}|^{\,\ell + 1} \right)\,,\\[0.15in]
{\bf RM}_{\,3}\, ({\cal Y}) & ``="&  O \,(1)\, \left\{ \,\,[\,A_1\, ({\cal Y})]^{ {{4}\over {n \,-\,2}}  -1} \cdot[\,{\cal V}_{\,i} \, ({\cal Y})- A_1\, ({\cal Y})]^{\,2} \ \right\}\,,\\[0.15in]
{\bf RM}_{\,4}\, ({\cal Y})  & ``=" &
  O \,  \left(\, \max_{|\,\lambda_i \, {\cal Y} \,+\, \xi_i|\,\le \,\rho_o^+}  \,\Vert \ \btd^{(\ell)}  K \ \Vert \times |\,\lambda_i\, {\cal Y} + \xi_{\,i}|^{\,\ell } \right) \, {\bf \times}\\[0.1in]
 & \ & \!\!\!\!\!\!\!\!\!\!\!\!\!\!\!\!\!\!\!\!\!\!\!\!\!\!\!\!\!\!\!\!\!\!\!\!     \times \ \left[\ O \,(1) \, |\,{\cal V}_{\,i}\, ({\cal Y}) \,-\,A_1\, ({\cal Y})| \times   \max \, \left\{ \ [\,{\cal V}_{\,i} \, ({\cal Y} )]^{{4}\over {n \,-\, 2}}\ , \ \ \  [\,A_1\,({\cal Y} )]^{{4}\over {n \,-\,2}}  \right\} \ \right]   \ \ \ \ \ \  \ \  \ \ \ \  \
  \end{eqnarray*}
  for $\,|\,{\cal Y} | \ \,\le\, \ \lambda_i^{-1}\ \rho_o\,.$\, (\,We use $``="$ to tell us that the right hand side is  the order of the term.)

\vspace*{0.3in}


{\bf \large {\bf  \S\,4. \ \  Cancelation of the $O\, (\lambda_i^\ell)\,$ term in (3.11).}}\\[0.2in]
We first ignore the order $\,\lambda_i^\ell\,$ in equation (3.11) and consider the linear inhomogeneous equation

\vspace*{-0.3in}

$$
\Delta_o\,  \Pi  \,+\, n \,(n + 2) \,A_1^{4\over {n \,-\,2}} \cdot \Pi  \,= \,{\cal P}_\ell \cdot A_1^{{n + 2}\over {n \,-\,2}} \ \ \ \ \ \ \  {\mbox{in}} \ \ \R^n  \leqno (4.1)
$$
(with unknown $\,\Pi\,$)\,.\, Here \,${\cal P}_\ell$\, is a {\it homogeneous\,} polynomial defined on $\,\R^n\,$  of degree $\ell \ge 1\,,\,$ and \,$A_1 \,=\, A_{1\,,\ 0}\, \, = \, \left({1\over {1 + |\,{\cal Y}|^2}} \right)^{\!\!{{n \,-\,2}\over 2}}\!.\ $
As in \cite{Compact}, potential solutions $\,\Pi\,$ can be expressed in the following    form
$$
\ \ \ \ \ \ \  \Pi\, ({\cal Y}) \,=\, {{\Gamma\, ({\cal Y})}\over { (1 + {\cal R}^{\,2})^{n\over 2} }} \ \ \ \ \ \ \mfor {\cal Y} \,\in\, \R^n \ \ \ \ \ \ \ \  \ (\,{\cal R} \,=\, |\,{\cal Y}|\,)\,. \leqno (4.2)
$$

Putting (4.2) into (4.1),   we obtain
$$
(1 + {\cal R}^{\,2})  \cdot [\,\Delta_o\, \Gamma \,] \ \,-\,\ 2n\, [\,{\cal Y} \cdot \btd \,\Gamma\,]  \ + \ 2n \,\Gamma\, =\, {\cal P}_\ell \ \ \ \ \ \ \ \ {\mbox{in}} \ \ \ \R^n\,. \leqno (4.3)
$$
Our goal is to find polynomial  solutions $\,\Gamma$\, of equation (4.3)\,,\,  and to keep  the degree of $\,\Gamma$\, as close to \,$\ell$\, as possible (this is for the second order blow\,-\,up argument), and to make explicit the dependence on $\,{\cal P}_\ell\,.\,$ We first look at some selected examples.\\[0.1in]
\noindent{\it Example} 4.4. \ \ {\it When\,}  ${\cal P}_\ell \,\equiv\, 0\,.$\ \ We can take
$$
\Gamma_1 \,({\cal Y}) \ := \ \sum_{j = 1}^n \, c_j\, {\cal Y}_{|_j}\,, \ \ \ \ {\mbox{or}} \ \ \ \
\Gamma_2 \,({\cal Y}) \ := \ {\cal R}^2 \,-\,1\,. \leqno (4.5)
$$
Here $c_j$ are any constants\,,\, and $\,{\cal Y}\,=\, ({\cal Y}_{|_1}\,, \ \cdot \cdot \cdot\,, \ {\cal Y}_{|_n} ) \,\in \R^n\,.\,$   Interestingly, up to linear combinations of $\Gamma_1$ and $\Gamma_2$\,,\, these are the  only possible solutions when ${\cal P}_\ell \equiv 0\,$ and when  we restrict $\Gamma$ to be a polynomial  of degree less than $\,n\,.\,$ Cf. Theorems 4.16 and 4.21.\\[0.1in]
\noindent{\it Example} 4.6. \ \ {\it When  $\,\ell \,\ge\, 2\,$ and $\,\Delta_o\, {\cal P}_\ell \,=\, 0\,.$}\, \ \
In this case we simply take $\,\Gamma_3 = c\, {\cal P}_\ell\,:$
\begin{eqnarray*}
(4.7) \!\!\!\!\!\!\!\!\!  &\ &
(1 + {\cal R}^{\,2})  \,\Delta_o\, \Gamma_3 \,-\,2n\, [\, ({\cal Y} \cdot \btd \, \Gamma_3) \,-\,\Gamma_3]   \ =\ -\,  2\,n\, (\ell \,-\,1) \,(c\,{\cal P}_\ell) = {\cal P}_\ell\ \ \ \ \ \ \ \ \ \ \  \ \  \ \ \ \  \   \ \  \ \ \ \  \    \\[0.1in]
&\Longrightarrow & c \ = \ -\,  {1\over {2\,n\, (\ell \,-\,1)}} \ \ \Longrightarrow \ \ \ \ \Gamma_3 \ = \ -\, {{{\cal P}_\ell}\over {2\,n\,(\ell \,-\,1)}}\ .
\end{eqnarray*}

\noindent{\it Example} 4.8. \ \ {\it When\,} $\ell \,=\, 1\,.\,$\ \
   As the left hand side of (4.3) is linear, we may assume that ${\cal P}_\ell\, ({\cal Y}) = {\cal Y}_{|_1}\,.\,$
Consider
$$
\Gamma_4\, ({\cal Y}) \ := \ a \,{\cal R}^{\,2}\, {\cal Y}_{|_1} \ + \ b\, {\cal R}^{\,4} \, {\cal Y}_{|_1}\ \ \ \ \ \ \ (\,{\mbox{maximum \ \ degree \ }} = 5)\,.\leqno (4.9)
$$
Direct calculation shows that when
$$
n \ = \ 4\,, \ \ \ \ a \ = \  {1\over {2n + 4}}\,, \ \ \ \ {\mbox{and}} \ \ \ b \ = \ {{2n-4}\over {2n + 4}} \cdot {1\over {4n + 16}}\,,
$$
 $\Gamma_4$   is a solution of (4.3) with ${\cal P}_\ell\, ({\cal Y})\, = \,{\cal Y}_{|_1}\,.\,$ The example demonstrates that, in general, $\ell = 1$  make it harder to solve equation (4.3).  Cf. an existence result for equation (1.1) obtained by Aubin in \cite{Aubin}\,.

\vspace*{0.15in}


{\bf \S\, 4\,a. \ \  Solving\,} (4.1) {\bf via the linear method\,} \,(\,$\ell\, <\, n)\,.$  \ \
  As in \cite{Compact} (pp.152)\,\,, we introduce the collection ($h_\ell$\, is the biggest integer less than or equal to \,$\ell/2$)

  \vspace*{-0.25in}

  \begin{eqnarray*}
(4.10) \ \ \ \ \ \ {\cal F} \,(\,{\cal P}_\ell\,) & \!:= &\!\! \left\{\,  {\mbox{linear \ combinations \ of}} \ ({\cal R}^{\,2})^j \,\Delta^{\,(k)}_o \, {\cal P}_\ell\,, \right.\\[0.1in]
 & \ & \ \ \ \ \ \ \ \ \ \ \ \ \ \ \ \ \ \ \ \ \ \ \ \ \ \left. 0 \,\le\, j \,\le \,k\,, \   k = 0\,, \ 1,  \ \cdot \cdot \cdot\,,\, h_\ell\,\right\},\ \ \ \ \ \ \ \ \ \  \ \  \ \ \ \  \  \ \  \ \ \ \  \
\end{eqnarray*}

\vspace*{-0.05in}

where $\,{\cal R}\, = \,|\,{\cal Y}|\,.$\, Note that $\Delta^{\!(k)}_o \, {\cal P}_\ell \equiv 0\,$ for $\,k \ \ge \  h_\ell + 1\,.\,$ Comparing to the one introduced in \cite{Compact}, (4.10) has the index \,$j$\, limited more strictly from above\,.\, Assuming that
$\,\Delta^{\!(j)}_o \, {\cal P}_\ell \not\equiv 0\,$ for $\, 0 \ \,\le\, \ j \ \,\le\, \  h_\ell\,,\,$
   ${\cal F}\,(\,{\cal P}_\ell\,)$ is a vector space with dimension,   in general,  equal to  $\, {1\over 2}\,(h_\ell \,+ 1) \,(h_\ell \,+ 2)\,$ [\,$ = \,O \,(n^2)$ when $\ell$ is close to $n$\,;\, an exceptional case is when ${\cal P}_\ell\, ({\cal Y}) = {\cal R}^\ell$\,,\, where $\ell$ is an even number\,]\,.\, We list down some simple properties concerning the operator in (4.3) and ${\cal F}\,(\,{\cal P}_\ell\,)$\,. All these can be checked readily (see \,{\bf \S\,A.3}\, in the \,e\,-\,Appendix for a proof).\\[0.1in]
%
%
%
{\bf Lemma 4.11.} \ \  ${\cal P}_\ell \,\in\, {\cal F}\,(\,{\cal P}_\ell\,)$\,.\, {\it Moreover,}

\vspace*{-0.25in}

 $$
{\it if} \ \ \ell \ \ {\it {is \ odd\  and}} \ \ \Delta_o^{\!(h_\ell)}\, {\cal P}_\ell \not\equiv 0 \ \ \Longrightarrow \ \ \sum_j\, c_j\, {\cal Y}_{|_j} \,\in\,  {\cal F}\,(\,{\cal P}_\ell\,)\,;\, \leqno (4.12)
 $$

\vspace*{-0.35in}

 $$
{\it if} \ \ \ell \ \ {\it {is \ odd\  and}} \ \   \Delta_o^{\!(h_\ell)}\, {\cal P}_\ell \not= 0 \ \ \Longrightarrow \ \ ({\cal R}^2 \,-\,1) \,\in\,  {\cal F}\,(\,{\cal P}_\ell\,)\,.\,\leqno (4.13)
 $$

 \vspace*{-0.05in}

 {\it Here at least one of the coefficient $c_j \not= 0\,.\,$} \\[0.075in]
(4.14) \ \  {\it The degree of each term in} \,${\cal F}\,(\,{\cal P}_\ell\,)$ {\it is at most\,}  $\ell\,.\,$\\[0.075in]
(4.15) \ \
$
 (1 + {\cal R}^{\,2})  \,\Delta_o\, \bullet \ \,-\,\ 2n (\,{\cal Y} \cdot \btd \,\bullet)  + 2n \, \bullet \ :    \ {\cal F}\,(\,{\cal P}_\ell\,) \to {\cal F}\,(\,{\cal P}_\ell\,)$   \ {\it{is \ \ linear}}\,.

\vspace*{0.15in}

\hspace*{0.5in}In principle, one can express the linear operator in (4.15) of  Lemma 4.11 into a matrix by using  the  basis of $\,{\cal F}\, ({\cal P}_\ell)\,$ as shown in (4.10)\,,\,  and determine whether  there is a solution or not. However, for genuine cases, and when $\,\ell\,$ is close to $\,n$\,,\,  the matrix is of the size (number of entries)   in the order $\,O\,( n^4 )\,$.\, In \cite{Compact}, the authors observes that  one can make use of the following Liouville\,-\,type theorem (shown in \cite{Chen-Lin-3}\,;\, see also \cite{Progress-Book}) to demonstrate that a solution exists\,.

\vspace*{0.15in}

 {\bf Theorem 4.16.} \ \ {\it Suppose $\psi$ is a smooth solution of\, the  equation\,}
 \vspace*{-0.2in}

$$
(1 + {\cal R}^{\,2})  \cdot [\,\Delta_o\,\psi \,] \ \,-\,\ 2n\, [\,{\cal Y} \cdot \btd \,\psi\,]  \ + \ 2n \,\psi \ = \ 0\ \ \ \ \ \ \ \ {\it{in}} \ \  \ \ \R^n\,. \leqno (4.17)
$$
$$
\lim_{r \to \infty} {{ \psi \, ({\cal Y}) }\over {{\cal R}^n}} \ = \ 0 \ \ \ \ \ \ \ \ \ \ \ \ \ (\,{\cal R}\, = \,|\,{\cal Y}|\,)\,.   \leqno (4.18) \ \ \   \ \ \ \ \ \ \  {\it{with}}
$$

$$
\psi \, ({\cal Y}) \ = \   \  c_o\, ({\cal R}^2 \,-\,1) \ + \ \sum_{j = 1}^n  c_j\, {\cal Y}_{|_j}    \ \ \ \  \ \ {\it{for}} \ \ \  {\cal Y} \,\in\, \R^n \ . \leqno (4.19) \ \ \ \ \ \ {\it Then}
$$
$$
{\it{Moreover\,,}} \ \ \ \ \  \psi \, (0)\ =\ 0 \ \ \ {\it and} \ \ \btd\, \psi\, (0)\ =\ 0 \ \ \Longrightarrow \ \ \phi \ \equiv \ 0 \ \ \ \ \ {\it in} \ \ \R^n. \leqno (4.20)
$$

 \vspace*{0.1in}

 \hspace*{0.5in}We now describe the linear method use in \cite{Compact} to find a polynomial  solution to equation (4.3)\,.\, The following result begins to reveal that the condition $\,
\Delta_o^{\!(h_\ell)} \ {\cal P}_\ell \,\equiv\, 0\,$ is tightly knitted together with the refined estimate we seek.

\vspace*{0.2in}

 {\bf Proposition 4.21.} \ \ {\it Assume that $\,{\cal P}_\ell\,$ is a homogeneous polynomial  of degree $\,\ell$\,,\, with }\, $\,2 \,\le \,\ell\, <\, n\,.\,$ {\it The  linear operator which appears in\,}  (4.15) {\it of Lemma\,} 4.11\,  {\it is a bijection if\, and only\, if}\, $\,
\Delta_o^{\!(h_\ell)} \ {\cal P}_\ell \,\equiv\, 0\,.$

\vspace*{0.15in}

 {\bf Proof.} \ \  For the sufficient part, the proof goes  in an essential manner  as in the  proof of Proposition 4.1 in \cite{Compact}, using Lemma 4.11 and Theorem 4.16 to show that the linear operator is an injection, and hence a bijection.
As for the necessary part\,,\, it follows from (4.12) and  (4.13) of Lemma 4.11\,,\, and Example 4.4 [\,the kernel contains a non\,-\,identically zero element in ${\cal F}\,({\cal P}_\ell)$\,]\,.\qedwh
\vspace*{0.2in}

{\bf \S\,4\,b. \ \   The reduction method.} \ \
Concerning the solution we find via Proposition 4.21, besides  being in $\,{\cal F}\, ({\cal P}_\ell)\,$ [\,in particular, we have property (4.14) in Lemma 4.11\,]\,,\,  there is  little we know about the solution itself. When we come to the  bubble estimates, it is natural to ask for the precise form of the solution $\Gamma$.\,   In this section we introduce a constructive method which allows us to determine each coefficient in \,$\Gamma$\,.\,    We present the precise result.

 \vspace*{0.15in}


{\bf Lemma 4.22.} \ \ {\it Let $\,{\cal P}_\ell\,$ be a homogeneous polynomial  of degree\,} $\ell  \,\ge\, 2\,$  ({\it defined on\,} $\R^n$\,)\,.
  {\it When $n$ is even.} {\it assume also that $\,\ell < n + 2\,$} ({\it no such condition when $n$ is odd}\,)\,.\, {\it Define a polynomial  $\ {\cal G}\,$ via }

  \vspace*{-0.25in}

  $$
  {\cal G}\,({\cal Y}) \ =   \ \sum^{k \,\le \,h_\ell \,-\,1}_{0 \,\le \,j \,\le \,k}  C^j_k \cdot ({\cal R}^2)^j \,[\, \Delta_o^{\!(k)} \, {\cal P}_\ell\,({\cal Y})\,] \ \ \ \ \ \  \ \ \ \ \ \ \  \  (\,{\cal R} \,=\,|\,{\cal Y}|\,)\,, \leqno (4.23)
  $$

{\it where the coefficients $\,C^{\,j}_k$ can be determined by using\,} (4.48)\,.\, {\it Then ${\cal G}$ satisfies}
\begin{eqnarray*}
(4.24)  \ \ & \ &\!\!\!\!\!\!\!\!(1 + {\cal R}^{\,2})  \,\Delta_o\, {\cal G} \ \,-\,\ 2n\, ({\cal Y} \cdot \btd \,{\cal G} )  \ + \ 2n \, {\cal G } \  =   \ {\cal P}_\ell \ +\\[0.15in]
    &  & \!\!\!\!\! \!\!\!\!\! \!\!\!\!\! \!\!\!\!\! \!\!\!\!\!    +\   [\, \Delta^{\!(h_\ell)}_o \, {\cal P}_\ell\,] \cdot  \left\{ \,a_{h_\ell} \cdot ({\cal R}^2)^{\,h_\ell} \ + \
a_{h_\ell -1} \cdot({\cal R}^2)^{\,h_\ell -1}\ + \cdot \cdot \cdot \ + \ a_1\cdot ({\cal R}^2)^1 \,+\, a_o \right\}. \ \ \ \ \ \ \ \ \ \ \ \ \ \
\end{eqnarray*}
{\it Here the numbers $\,a_k\,$ can be found by using\,} (4.45)\,.\,\bk
%
%
%
%
The precise definitions of \,$C_k^{\,j}$\, and \,$a_k$\, are obtained on the way toward the proof of Lemma 4.22. The key property is that they depend only on \,$n$\,,\, $\ell$\,,\, $j$ \,and $k$  only, and are formed by an algebraic  iteration process, which we start to describe.

\vspace*{0.1in}

{\bf \S\,4\,b\,.1\,.} \ \ {\it First step in the proof of Lemma} 4.22\,: {\it  the recurrent and reduction to powers of $\,({\cal R}^2)\,.$}\ \
 Consider first the situation where
$$
\Delta_o^{\!(k)} \, {\cal P}_\ell \not\equiv 0 \ \ \mfor \ \ 1\, \le\, k \,\le\, h_\ell \,-\,1\,. \leqno (4.25)
$$
As in Example 4.4, we take

\vspace*{-0.25in}

$$
  [\,C^{\,o}_o\,\cdot{\cal P}_\ell\,]\,, \ \ \ \ {\mbox{where}} \ \ C^{\,o}_o \ = \ {{1} \over {2n\, (1 \,-\, \ell )}}\,,\leqno (4.26)
$$

\vspace*{-0.15in}

and obtain  \\[0.075in]
(4.27)
$$
(1 + {\cal R}^{\,2})  \,\Delta_o\,  [\,C^{\,o}_o\,\cdot{\cal P}_\ell\,]  \ + \ 2n \,[1-  ({\cal Y} \cdot \btd \,)]\, [\,C^{\,o}_o\,\cdot{\cal P}_\ell\,] = {\cal P}_\ell \ + \   C^{\,o}_o \cdot \left[\,({\cal R}^{\,2}) \,\Delta_o\, {\cal P}_\ell + \,\Delta_o\, {\cal P}_\ell  \,\right]\,.
 $$
That is, we obtain $\,{\cal P}_\ell$\, in the right hand side, but  ``pay the price" by introducing $\,[\,(\,{\cal R}^{\,2}) \,\Delta_o\, {\cal P}_\ell\,]\,$ and $\,[\, \Delta_o\, {\cal P}_\ell\,]\,$.\,
 Observe that the degree of $\,\Delta_o\, {\cal P}_\ell\,$ is lowered to $\,\ell \,-\,2\,,\,$ while the degree of $\,[\,({\cal R}^2)\,\Delta_o\, {\cal P}_\ell\,]\,$ is still equal to $\,\ell\,$,\, but its structure appears simpler in the sense that 2 of the degree is taken over by $(\,{\cal R}^2)\,.$\,\bk
If $\,\ell \,=\, 2$\, or \,$3$\,,\, we are done. Assume that $\,\ell\, \ge\, 4\,.\,$ To proceed, we
 simplify the notations and highlight the change in order, and  introduce
$$
{\bf R} = ({\cal R}^{\,2}) \ \ \Longrightarrow \ \ {\bf R}^j = ({\cal R}^{\,2})^j\,, \ \ \ \ {\bf D} = [\,\Delta_o\, {\cal P}_\ell\,] \ \ \longrightarrow  \ \ {\bf D}_k := [\,\Delta_o^{\!(k)}\, {\cal P}_\ell\,]\,. \leqno (4.28)
$$
Using these notations, we have (cf. \,{\bf \S\,A.3} in the e\,-\,Appendix)\\[0.1in]
(4.29)
\begin{eqnarray*}
& \ &\!\!( 1 + {\bf R}) \cdot \Delta_o\, [\,{\bf R}^j \ {\bf D}_k] \ \ \ \ \left(\ = \ (1 + {\cal R}^2)\,\Delta_o\,[ \, ({\cal R}^2)^j \cdot \Delta_o^{\!(k)}\, {\cal P}_\ell\,]\, \right)\\[0.1in] & \ &\!\!\!\!\!\!\!\!\!\!\!\!\!\!= \ A_{\ell\,,\, j\,,\, k} \cdot (\,{\bf R}^{\,j }  \ {\bf D}_{\,k }) \  + \ {\bf R}^{\,j+1}  \ {\bf D}_{\,k+1} \  \   [\,{\mbox{degree}} = \ell + 2\,(j \,-\,k) \ \ {\mbox{on \ \ both \ \ terms}}\,]\\[0.1in]
& \ & \ \ \ \ \ \ \   +\  A_{\ell\,,\, j\,,\, k} \cdot (\,{\bf R}^{\,j-1}  \ {\bf D}_{\,k }) \ + \ {\bf R}^j  \ {\bf D}_{\,k+1} \  \ \ \ \ \ \ \ \ \\
  & \ & \hspace*{1.6in}[\ \uparrow \ \ {\mbox{degree}} = \ell + 2\,(j \,-\,k) \,-\,2 \ \ {\mbox{on \ \ both \ \ terms}}\,]\,. \ \ \ \ \ \  \ \  \ \ \ \  \ \
 \end{eqnarray*}

 \vspace*{-0.35in}

 $$
 A_{\ell\,,\,j\,,\, k\,} \ = \ (2j\,) \cdot (2j +n-2 + 2\,\ell \,-\,4 k) \leqno{\mbox{Here}}
 $$
 We realize that the process produces one same term $\,[\,{\bf R}^j \ {\bf D}_k]\,$ (times a constant), one same order term, plus two lower order terms. We illustrate the procedure when the linear operator [\,the right hand side of (4.3)] acts on $\,[\,{\bf R}^j \ {\bf D}_k\,]\,$ via the following diagram (\,$j \ge 1$\,)\,.

\begin{eqnarray*}
& \! &  \!\!\!\!\!\!\!\!\!\!\!\!\!\! \!\!\!\!\!\!\! \!\!\!\!\!\!\!\!\!\!\!\!\!\!\!\!\!\!\!\!\!\!\!\ \ \ \ (\times\,A_{\ell\,,\, j\,,\, k}   \ \rightarrow \,) \ \ \ \ \   {\bf R}^{\,j-1} \ {\bf D}_k \ \ \ \ \ \ \  \  \\
  \!\!\!\!\!& \! &\!\!\!\!\! \nearrow \\
(\times\,  \{\,A_{\ell\,,\, j\,,\, k} \,-\,2n\,[ \,\ell + 2\,(j \,-\,k) -1 )\,]\,\}\, )  \ \  \hookrightarrow \  {\bf R}^j \ {\bf D}_k & \!&  \!\!\! \!\!\!\! \  \longrightarrow \  \,{\bf R}^j \ {\bf D}_{\,k + 1} \ \ \ \ \ \ (\leftarrow \  \times\,1)\\[0.1in]
\downarrow \ (\leftarrow \   \times\,1)\!\!\!\!\!\!\!\!\!\!\!\!\!\!\!\!& \! & \\[0.1in]
{\bf R}^{\,j+1 } \ {\bf D}_{\,k+ 1}\!\!\!\!\!\! & \ &\\[0.15in]
{\mbox{Degree \ \ of \ \ the \ \ term}}\ (\,\le\,)\,:\ \   \ \ \ \ \ \  \ \ell + 2\,(j \,-\,k)  \!\!\!\!& \ & \ \ \ \   \ell + 2\,(j \,-\,k) \,-\,2
 \end{eqnarray*}

 \vspace*{-0.05in}

\centerline{\underline{Diagram 4.30}\,. \ \ The four terms, their degrees, and the multipliers.}

\newpage

We represent schematically part of the reduction procedure in the following diagram, showing the terms produced (indicated by the arrows, including itself) when the term is acted upon by  the operator $\,(1 + {\bf R})\cdot \Delta_o\,.\,$

\begin{eqnarray*}
  {\cal P}_\ell & \ &\ \ \  \  \longrightarrow \ \ \   \  \ {\bf D}\ \   \ \ \  \ \  \longrightarrow \ \ \   {\bf D}_2 \ \ \cdot \cdot \cdot    \ \ \ \ \ \ \ \  \ \ \ \  \ \  \ \     \longrightarrow \ \   [\,\Delta_o^{\!(h_\ell)} \  {\cal P}_\ell\,]  \\[0.1in]
     \downarrow \ & \ & \   \ \ \    \nearrow \ \ \ \ \ \   \ \downarrow  \ \ \ \ \ \ \ \ \ \     \nearrow   \\[0.1in]
   (1^{\mbox{st}}) \ \ \ \  \ \ \  \,{\bf R} & \! &\!\!\!\!\!\!\!\!\!\!\!{\bf D} \ \  \ \  \longrightarrow \ \  \  \ \ {\bf R}\,{\bf D}_2\ \   \ \ \   \longrightarrow \ \ \  {\bf R}\,{\bf D}_3 \ \ \cdot \cdot \cdot \ \ \       \longrightarrow \ \  ({\cal R}^{\,2})\,[\,\Delta_o^{\!(h_\ell)} \  {\cal P}_\ell\,]  \ \ \ \ \ \ \\[0.1in]
 \downarrow & \ &\!\!\!\! \ \ \ \ \ \  \nearrow \ \ \ \ \ \ \ \downarrow \ \ \ \ \ \ \ \ \nearrow \\[0.1in]
   (2^{\mbox{nd}}) \ \ \ \ \  \   {\bf R}^{\,2} & \! &\!\!\!\!\!\!\!\!\!\!\!{\bf D}_2 \ \   \longrightarrow \ \  \  \ \ \, {\bf R}^{\,2}\,{\bf D}_3  \ \   \ \    \longrightarrow \ \    \cdot \cdot \cdot \ \   \longrightarrow \ \   ({\cal R}^{\,2})^{\,2}\,[\,\Delta_o^{\!(h_\ell)} \  {\cal P}_\ell\,]  \\
& : & \ \ \  \hspace*{0.75in}\downarrow\ \ \ \ \ \ \ \ \ \ \ \ \cdot \cdot \cdot  \ \ \ \ \\[-0.1in]
   \downarrow   & \  & \!\!\!\!\!\!\! :   \ \ \ \ \  \ \ \ \ \ \   \\[0.1in]
  {\bf R}^{\,h_\ell \,-\,1}\,{\bf D}_{h_\ell \,-\,1}\!\!\!\! \!\!\!\!   \!\!\!\!\!\!\!\!\!& \  &   \ \ \ \ \ \rightarrow\  \ \ {\bf R}^{\,h_\ell-1} \,{\bf D}_{h_\ell}  \ = \ ({\cal R}^2)^{\,h_\ell -1} \,[\,\Delta_o^{\!(h_\ell)} \  {\cal P}_\ell\,]   \\[0.1in]
  \downarrow  &\  &  \ \ \ \ \ \ \ \ \ \ \ \ \  \\[0.1in]
 (h_\ell^{\mbox{th}}) \ \ \ \ \   {\bf R}^{h_\ell}& \! &\!\!\!\!\!\!\!\!\!\!\! {\bf D}_{h_\ell}   \ \ \ \ \ \ \ \     \\
   ||   & \  &    \\[0.1in]
 ({\cal R}^{\,2})^{h_\ell} & \! &\!\!\!\!\!\!\!\!\!\!\![\,\Delta_o^{\!(h_\ell)} \  {\cal P}_\ell\,]  \\[0.15in]
{\mbox{  Degree }}\  (\,\le\,)   \ \  \  \ \ell & \ &  \ \ \ \ \ \ \  \ \ \ \  \ \ \ \ell \,-\,2 \ \ \ \ \cdot \cdot \cdot \ \ \ \ \ \ \ \  \ \ \ \ \  \ \  \ \ \ \ \   \ \ \ \ \  \ \  \ \ \ \ \ \ \  \ \ \ \    \ 1\ \,/ \ \, 0\ \ \ \ \ \ \ \  \ \ \ \ \  \ \\
& \ & \hspace*{3.2in} \ell \ \,{\mbox{is \ \,odd/even}}\\
\end{eqnarray*}

\centerline{\underline{Diagram 4.31}\,. Showing   the cancelation order (top $\to$ down) on the first column.}

\hspace*{0.5in}Back to the case when $\,\ell \,=\,4\,.\,$ In the second step, we seek to eliminate the term
$$
     C^{\,o}_o \cdot [\,({\cal R}^{\,2}) \,\Delta_o\, {\cal P}_\ell\,]\ = \   C^{\,o}_o\cdot{\bf R}\,{\bf D}\,,
$$
which appears in (4.27)\,.
From (4.29) and (4.27)\,:\,  we take
$$
C^1_1 \ = \  {{-\,C^{\,o}_o} \over {A_{\ell\,,\,1\,,1} \,-\,2n\,( \ell -1 ) }} \ = \  {{-\,C^{\,o}_o} \over { 2\,(n \,-\,2) \,(2 \,-\,\ell)}}\   \ \ \ \ \  \ \ \ \ (\,{\mbox{here}} \ \ \ell \,\ge \, 4)\,.\,\leqno (4.32)
$$

It follows that

\vspace*{-0.25in}

\begin{eqnarray*}
(4.33) \!\!\!\!\!\!\!\!\!& \ & (1 + {\cal R}^{\,2})  \,\Delta_o\,  \{\,[\,C^{\,o}_o\,\cdot{\cal P}_\ell\,] + [\,C^{\,1}_1\cdot({\cal R}^2)\,\Delta_o\, {\cal P}_\ell\,] \,\} \ \\[0.1in]
& \ & \ \ \ \ \ \ \ \ \ \  + \ 2\,n \,[1-  ({\cal Y} \cdot \btd \,)]\, \{\,[\,C^{\,o}_o\,\cdot{\cal P}_\ell\ ] + [\,C^{\,1}_1\cdot({\cal R}^2)\,\Delta_o\, {\cal P}_\ell\ ] \,\}\\[0.1in]
& \ &\!\!\!\!\!\!\!\!\!\!\!\!\!\!\!\!\!\!=  \  {\cal P}_\ell \ + \   [\,C^{\,o}_o+A_{\,\ell\,,\, 1\,,\, 1} \cdot C^{\,1}_1\ ] \cdot \left[\, \,\Delta_o\, {\cal P}_\ell \ \right] + C^{\,1}_1 \cdot [\, ({\cal R}^2)^2 \, \Delta_o^2 \, {\cal P}_\ell + ({\cal R}^2)  \, \Delta_o^2 \, {\cal P}_\ell \ ]\,. \ \ \ \ \ \ \ \ \ \ \ \ \ \
\end{eqnarray*}
Inductively, we find [\,refer to (4.23)]
$$
C^{\,j}_j \ = \ {{-\,C^{j \,-\,1}_{j \,-\,1}} \over {A_{\ell\,,\,j\,,\,\, j} \,-\,2n\,( \ell -1 ) }} \ = \ {{-\,C^{j \,-\,1}_{j \,-\,1}} \over  {  (2j)\,[\,n \,-\,2 + 2 \,(\ell \,-\,j)\,] \,-\,2n\,(\ell -1)]}}\ , \leqno (4.34)
$$
where $\,1 \,\le \,j \,\le\, h_\ell\,$.\, This enables us to cancel the terms in the first column in Diagram 4.31, ending with  the term $$
\,C^{\,h_\ell \,-\,1}_{h_\ell \,-\,1}\,\cdot [\,({\cal R}^2)^{\,h_\ell}\, \Delta_o^{\!(h_\ell)}\,{\cal P}_\ell\,]\,,\, \leqno (4.35)
$$
which is present on the right hand side of  equation (4.24). \bk
Next, we proceed to cancel the terms in the second column (Diagram 4.31), starting from top toward the bottom. Gradually we move right to the next column, always proceeding from top to bottom. We summarize the cancelation in the following two cases.\\[0.2in]
{\bf *} \,{\it Cancelation of terms in the  top row in  Diagram\,} 4.31\,. \ \
From (4.33), and also from  Diagram 4.31\,,\, we have the term $$
[\, C^{\,o}_o + C_1^1 \cdot A_{\ell,\, 1\,,\, 1}\,]\cdot {\bf D}
$$
to be canceled. This is done by adding the term
$$
-\, {{C^{\,o}_o + C_1^1 \cdot A_{\ell,\, 1\,,\, 1}}\over {-\,2n\,[\, (\ell \,-\,2) -1])}}\cdot  {\bf D} \ \ \Longrightarrow \ \ C_1^{\,o} \ =\  -\, {{C^{\,o}_o + C_1^1 \cdot A_{\ell,\, 1\,,\, 1}}\over {-\,2n\,[\, (\ell \,-\,2) \,-\,1]\,)}}\ \,.
$$
With the help of the information depicted in  Diagram 4.30, and via induction, we have
$$
C_k^{\,o} = -\, {{C^{\,o}_{k \,-\,1} + C^1_k \cdot A_{\ell,\, 1\,,\  k}}\over {-\,2n\,[\, (\ell \,-\,2k) -1]\,)}} \ \ \ \ \ \ \mfor \ \ 1 \,\le\, k \,\le\, h_\ell\, -\,1\,. \leqno (4.36)
$$
The numbers $\,C^1_2$\,,\, $C^1_3\,,\, \cdot \cdot \cdot\,,\,\ C^1_{h_\ell \,-\,1}\,$ are obtained below -- see Remark 4.40.
Note that
$$
\ell \,-\,2k \ \ge \ \ell \,-\,2\, (h_\ell \,-\,1) \ \ge \ 2 \ \ \Longrightarrow \ \ (\ell \,-\,2k) -1 \ \not= \ 0 \ \ \ \ \ \ [\,{\mbox{recall \ \ (4.25)\,]}}\,.
$$
This enables us to cancel the terms in the first row, ending with
$$
\,C^{\,o}_{h_\ell \,-\,1}\cdot \,\Delta_o^{\!(h_\ell)}\,{\cal P}_\ell\,,\,
$$
which appears in the right hand side of (4.24).

 \vspace*{0.15in}

{\bf *} \,{\it Cancelation of the }\,`{\it inside}' {\it terms\,.}\ \
Finally, consider any `inside' term $\,{\bf R}^j \, {\bf D}_k$\,.\, We observe that $\,k \,> \,j\,$\, ($\,k\, =\, j\,$ appears in the first column only)\,.

 \vspace*{-0.15in}

\begin{eqnarray*}
& \! &\!\!\!\!\!\!\!\!\!\!\!\!\!\!\!\!\!\!\!\!\!\!\!\!\!\!\!\!\!\!\!\!\!\!\!\!\!\!\!\!\!\!\!\!C_{k-1}^{j -1} \  {\bf R}^{\,j-1} \ {\bf D}_{k-1} \ \ \ \ \ \ \ \ \ \ \ \ \ \  \ \ (\,j \,\ge\, 1)\\
  (\,\times\,1\  \ \rightarrow ) \ \ \downarrow  \ \ \ & \ &  \\[0.1in]
   C^j_{k \,-\,1} \,{\bf R}^j\ {\bf D}_{k \,-\,1}  \  (\,\times\,1) \ \rightarrow \ \  \ \  {\bf R}^j \ {\bf D}_k & \!&   \\
\nearrow \ (\leftarrow \,\times\,A_{\ell\,,\,j + 1\,,\  k} )\!\!\!\!\!\!\!\!\!\!\!\! \!\!\!\!  & \! & \\[0.05in]
C_k^{j + 1} \ {\bf R}^{\,j+1 } \ {\bf D}_{\,k}\ \ \ \ \ \ \ \ \ \ \ \ \ \ \ \ \ \  \  & \ &
 \end{eqnarray*}

 \vspace*{0.1in}

 \centerline{{\underline{Diagram 4.37}}. \,The three terms which give rise to an inside term (with multipliers).}

 \vspace*{0.2in}

Via induction and the discussion in (4.34) and (4.36), we may assume that the coefficients $C^{\, j \,-\, 1}_{\,k \,-\,1}$\,,\, \ $C^{\, j }_{\,k \,-\,1}$\, and \, $C^{\, j + 1}_{\,k}$\, are determined.
The term $\,{\bf R}^j \, {\bf D}_k\,$  makes its presence on the right hand   given by $$
[\ C^{j \,-\,1}_{k \,-\,1} \ + \ C^j_{k \,-\,1}
\ + \ C^{j+1}_k \cdot A_{\ell\,,\,j+ 1\,,\, k} \ ] \cdot {\bf R}^j \ {\bf D}_k \leqno (4.38)
$$
To cancel it, we introduce the term
$$
-\, {{\ C^{j \,-\,1}_{k \,-\,1} \ + \ C^j_{k \,-\,1}
\ + \ C^{j+1}_k \cdot A_{\ell\,,\,j+ 1\,,\  k}  }\over { A_{\ell\,,\,j\,,\, k } \,-\,2n\,[\,(\ell \,-\,2k + 2j) \,-\,1]}}\cdot {\bf R}^j \, {\bf D}_k
$$
to the {\it left}\, hand side\,,
$$
\!\!\!\!\!\!\!\!\!\!\cdot \cdot \cdot \cdot \cdot \cdot \ \ \Longrightarrow \ \ \ \  C^j_k \ =\ -\, {{\ C^{j \,-\,1}_{k \,-\,1} \ + \ C^j_{k \,-\,1}
\ + \ C^{j+1}_k \cdot A_{\ell\,,\,j+ 1\,, \  k}  }\over { A_{\ell\,,\,j\,,\ k } \,-\,2n\,[\,\ell +2\,(j \,-\, k)   \,-\,1]}} \ . \leqno (4.39)
$$

\vspace*{0.1in}

{\it Remark}\, 4.40. \ \ Concerning the usage of $\,C^1_2$\,,\, \  $C^1_3\,,\,\ \cdot \cdot \cdot\,,\,\ C^1_{h_\ell \,-\,1}\ $ in (4.36), we remark that, based on (4.39), in order to determine $\,C^1_2$\,,\, we need only $\,C^o_1\,,\ \,C^1_1\,$ and $\,C^2_2\,$,\, which are known via (4.34) and  (4.36)\,.\, Afterward, we can determine $C^{j \,-\,1}_j$\, for \,$3 \,\le\, j \,\le \,h_\ell \,$ (the other coefficient in the second column in Diagram 4.31)\,.\, $C^2_3$\,,\, together with \,$C_2^o$\, and \,$C^1_2$\, help to determine $\,C^1_3$\,,\, and so on.

 \vspace*{0.15in}

{\bf \S\,4\,b\,.2.} \ \ {\it Non\,-\,zero characters for\,} $({\cal R}^{\,2})^{\,j} \,\Delta_o^{\!(k)} \,=\,({\cal R}^{\,2})^{\,j} \,\Delta_o^{\!(j \,+\,  \sqcup)}$\ \ \   (\,$j \,\ge\, 1$ and $\sqcup \,\ge\, 0\,$)\,. \ \
 In order to finish the proof for Lemma 4.22, we are required to show that the denominators  in (4.34) and (4.39) are non\,-\,zero \,--\, under the conditions on $\,\ell\,$ as stated in Lemma 4.22\,.\,
Note that
$$ {\mbox{Degree}} \ \{ \,({\cal R}^{\,2})^{j} \cdot [\,\Delta^{\!(j \,+\, \sqcup)}_o \, {\cal P}_\ell\,]\,\} \ =\ \ell \,+\, 2\,[\,j \,-\,  (j + \,\sqcup\,)\,]=\, \ell \,-\,2 \,\sqcup \ \,\le\, \  \ell \ \ (\,\sqcup \ \ge \ 0\,)\,.
$$
Moreover, as   the process stops when $\,j + \sqcup = h_\ell\,,\,$ we need only to consider the situation where  $j + \sqcup \ \le\  h_\ell \,-\,1\,.\,$ It follows that
\begin{eqnarray*}
(4.41) \ \  k \,-\, j & = & \sqcup \ \ge\  0 \ \ {\mbox{and}} \ \ k \ = \ j \,+ \,\sqcup \  \,\le\, \  h_\ell\,-\, 1  \ \ \ \Longrightarrow \ \
\  j  \  \,\le\,  \ h_\ell \, -\, 1\,; \ \ \ \ \  \ \ \ \ \ \ \ \ \ \ \ \ \  \\[0.1in]
j \ \ge\  1 \!\!\!\!& \ &    \Longrightarrow \ \ \,\sqcup\  \,\le\, \ h_\ell \,-\, 2\,.
\end{eqnarray*}
We investigate the {\it characteristic equation\,,\,} which is given by the denominator in (4.39) [\,note that the denominator in (4.34) corresponds to $\,k \,=\, j\,$]\,: \begin{eqnarray*}
  & \ &  A_{\,\ell\,,\,j\,,\ k} \ \,-\,\ 2n \,[\,\ell + 2\, (j \,-\,k) \,-\,1]\ = \ 0\\[0.1in]
\Longleftrightarrow  & \ &  (2j)\, [\,2j + (n \,-\,2) + 2\, (\ell \,-\,2k)\,] \,-\,2n\, [\,\ell \,-\,2\,(k \,-\,j) \,-\,1]   \ = \  0\ \ \ \ \ \ \ \ \ \ \ \\[0.1in]
\Longleftrightarrow   & \ & (2j)^{\,2} \,-\,(2j) [\, (n \,-\,2) + 2\,  \ell  \,-\,4 \,\sqcup\,] +  n\, (2\ell \,-\,4\,\sqcup \,-\,\,2\,)   \ = \ 0 \ \ \ \ \ \\[0.1in]
\Longleftrightarrow  & \ & [\,(2j) \,-\,n] \cdot [\,(2j) \,-\,(\,2 \,\ell \,-\,4\, \sqcup \,-\,2)\,]\ = \ 0\\[0.1in]
(4.42) \ \cdot \cdot \cdot   \ \Longleftrightarrow \ \ & \ & j \ = \ {n\over 2} \ \ \  \ {\mbox{or}} \ \ \ \  j \ = \ (\ell \,-\,1) \,-\,2\,\sqcup\,.
\end{eqnarray*}

 \vspace*{0.1in}

{\bf *}{\it When $n$ is even.} \ \  Here $\,n/2\,$ is an integer, (4.41) requires us to post the restriction
$$
j \ \,\le\, \ h_\ell \,-\,1 \ < \ {n\over 2} \ \ \   \Longleftrightarrow   \ \  h_\ell  \,-\,1 \ < \ {n\over 2} \ \ \Longleftrightarrow \ \ 2\cdot  h_\ell \ < \ n + 2\,.  \leqno (4.43)
$$
That is, when $\,\ell\,$ is even, we require
$$
 \ell \,<\, n + 2\,. \leqno (4.44)
$$
Similarly, when $\,\ell\,$ is odd, we need
$$
  2 \cdot {{ \ell \,-\,1}\over 2}  \ < \ n + 2  \ \ \Longleftrightarrow \ \ \ell  \,< \, n + 3 \ \ \Longleftrightarrow \ \ \ell \,< \,n + 2 \,,
$$
as $n$ is even $\,\Longrightarrow n + 3\,$ is odd, and $\ell$ is also odd in this case.
For the second root in (4.42), since $k \,\le\, h_\ell \,-\, 1$, we have
$$
  j + \sqcup \,<\, {\ell\over 2} \ \ \Longrightarrow \ \  j + \sqcup \,\le \,{\ell\over 2} \,-\,1 \ \ \Longrightarrow \ \   \sqcup \,\le\, {\ell\over 2} \,-\,2 \ \ \ (j \,\ge\, 1) \ \ \Longrightarrow \ \   j + 2\,\sqcup\, \le\, \ell \,-\,3\,.
$$
Thus the solution $\,j \,=\, (\ell \,-\,1) \,-\,2\,\sqcup \ \ \Longleftrightarrow \ \ j + 2\,\sqcup \,=\, \ell \,-\,1\,$ is too big to happen.\,

\vspace*{0.15in}

{\bf *}{\it When $n$ is odd.} \ \  In this case,
$
\,n/2\,
$ is not an integer.\,
We need only to consider the second root in (4.42).
  As we want the term $({\cal R}^{\,2})$ to be present\,,\, and  $(\Delta_o^{\!(j + \,\sqcup)}\ {\cal P}_\ell\,)$ is not yet reduced  to  first order, therefore
\begin{eqnarray*}
 & \ &  j + \sqcup \ < \ {{\ell-1}\over 2} \ \ \Longrightarrow \ \  j + \sqcup \ \,\le\, \ {{\ell-1}\over 2} \,-\,1 \ \ \Longrightarrow \ \   \sqcup \ \,\le\, \ {{\ell-1}\over 2} \,-\,2 \ \ \ (j \,\ge\, 1) \\[0.1in] \Longrightarrow  \!\!  & \ & j + 2\,\sqcup \ \,\le\, \ \ell \,-\,3  \ \ \Longleftrightarrow \ j \ < \ (\ell \,-\,3) \,-\, 2 \sqcup\,.
\end{eqnarray*}
Once again, the solution $\,j \ = \ (\ell \,-\,1) \,-\,2\,\sqcup\,$ is too big.\, This completes the showing that the denominators in (4.34) and (4.39) are non\,-\,zero under the conditions of Lemma 4.22.\\[0.2in]
{\bf *}{\it The residue}\,.\ \ As the `pure' $\,({\cal R}^2)\,,\  ({\cal R}^2)^2\,,\, \ \cdot \cdot \cdot\,, \ ({\cal R}^2)^{\,h_\ell}\,$ terms are obtained as the by\,-\,products of the last cancelations in each column (see Diagram 4.31), except the last coefficient $\,a_o\,$,\, all the others are combination the horizontal arrow and the downward arrow (refer to Diagram 4.31)\,.\, Hence [\,together with Diagram 4.30\,;\, cf. also (4.35)  and (4.23)] we have \\[0.1in]
(4.45)
\begin{eqnarray*}
& \ & \!\!\!\!\!\!\!\!\!\!a_{h_\ell} = C^{\,h_\ell \,-\,1}_{h_\ell \,-\,1}\,, \   a_{h_\ell \,-\,1} =
\left[\ C^{\,h_\ell \,-\,1}_{h_\ell \,-\,1} \,+\, C^{\,h_\ell \,-\,2}_{h_\ell \,-\,1} \,\right]\,, \cdot \cdot \cdot \,,\,   a_1 =  \left[ \ C^{\,1}_{h_\ell \,-\,1}  \, + \, C^{\,o}_{h_\ell \,-\,1} \right]\,,\\[0.15in]
& \ & {\mbox{and}} \ \ \ a_o = C^{\,o}_{h_\ell \,-\,1}\ .
\end{eqnarray*}
The argument is completed under condition (4.25)\,.\,
Finally, suppose
$$
\Delta_o^{\!(k_o)} \, {\cal P}_\ell  \equiv 0 \mfor {\mbox{an \ \ integer \ \ }}   k_o \,\in\, [\,1\,, \ h_\ell -1\,]\,. \leqno (4.46)
$$
The process described in Diagram 4.31 ends earlier. In this case
$$
  {\cal G}\  =   \sum_{0 \,\le\, j \,\le \,k \,\le\, \,k_o \,-\,1} \!\!\!\!C^j_k \cdot ({\cal R}^2)^j \,[\, \Delta_o^{\!(k)} \, {\cal P}_\ell\,]\,, \leqno (4.47)
  $$
 where the coefficients $\,C^{\,j}_k\,$ are the same as the above.  Moreover, in this case [\,that is, with (4.46)] \,${\cal G}$\, satisfies
$$(1 + {\cal R}^{\,2})  \,\Delta_o\, {\cal G} \ \,-\,\ 2n \,({\cal Y} \cdot \btd \,{\cal G} )  \ + \ 2n \, {\cal G } \  =   \ {\cal P}_\ell\,.
$$
This completes the proof of Lemma 4.22\,.\qed
We summarize the coefficient in the following.
\begin{eqnarray*}
 \uparrow \ \    & \ &  C^{\,o}_o = {{-\,1}\over {2n\,(\ell \,-\,1)}}\,, \  \cdot \cdot \cdot\,, \ \ C^{\,j}_j = {{-\, C^{j \,-\,1}_{j \,-\,1}}\over {A_{\,\ell\,,\, j\,,\, j} \,-\,2n\,(\ell \,-\,1)}} \ ,\, \cdot \cdot \cdot\ ,
 \\[0.1in]
 \!\!(4.48)    & \ & C^{\,o}_1 = -\,{{C^{\,o}_o + C^{\,1}_1\cdot A_{\,\ell\,,\, 1\,,\, 1}  }\over {-\,2n\,[\,(\ell \,-\,2) \,-\,1]}}\,, \  \cdot \cdot \cdot\,, \ \  C^{\,o}_k = -\,{{C^{\,o}_{k \,-\,1} + C^1_k\cdot A_{\,\ell\,,\, 1\,,\, k}  }\over {-\,2n\,[\,(\ell \,-\,2k) \,-\,1]}} \ , \ \cdot \cdot \cdot\ , 
  \\[0.15in]
 \downarrow \ \  \ &\ &  C^{\,j}_k = -\, {{\ C^{j \,-\,1}_{k \,-\,1} \ + \ C^j_{k \,-\,1}
\ + \ C^{j+1}_k \cdot A_{\ell\,,\,j+ 1\,,\, k}  }\over { A_{\ell\,,\,j\,,\, k } \,-\,2n\,\{\,[\,\ell +2\,(j \,-\,k)]   \,-\,1\}}} \ \ \ \  \  {\mbox{for}} \ \ {1 \,\le\, j \ < \ k \,\le \,h_\ell \,-\, 1}\,.\ \ \ \  \ \ \ \ \  \ \ \ \ \  \ \ \ \ \  \ \ \ \ \  \ \ \ \ \  \ \ \ \ \  \ \ \ \ \  \
  \end{eqnarray*}

\newpage

{\bf Proposition 4.49.} \ \ {\it Let $\,{\cal P}_\ell\,$ be a homogeneous polynomial  of degree $\,\ell\,$ defined on} $\,\R^n.\,$  {\it Assume that}\\[0.1in]
{\bf (i)} \ \ \,{\it when $n$ is even}\,: $\,2 \,\le \,\ell  < n  + 2\,$ {\it \, and\ } $
\Delta_o^{\!(h_\ell)} \ {\cal P}_\ell = 0\,$ (\,$\Delta_o^{\!(h_\ell)} \ {\cal P}_\ell$\, {\it is degree zero\,})\,;
\\[0.1in]
{\bf (ii)} \  \,{\it when $n$ is odd}\,:  $\,2 \,\le \,\ell\,   $ {\it  and\,} $
 \Delta_o^{\!(h_\ell)} \ {\cal P}_\ell \equiv 0$ \, \,(\,{\it here\,}  $\Delta_o^{\!(h_\ell)} \ {\cal P}_\ell$\, {\it is degree one\,})\,.\,\\[0.1in]
 {\it Then   equation}\, (4.3) {\it has a polynomial  solution  $\,{\cal G}_o\,$ given by }
 $$
  {\cal G}_o \ =   \! \sum_{0 \,\le\, j \,\le \,k \,\le\, \,h_\ell\, \,-\,1} \!\!\!C^j_k \cdot ({\cal R}^2)^j \,[\, \Delta_o^{\!(k)} \, {\cal P}_\ell\,]\,. \leqno (4.50)
 $$
 {\it The coefficients $\,C^j_k\,$ are  presented in\,} (4.48)\,.\, {\it In particular, the constant and the linear terms are not present in the solution, and the   degree of each term in $\,{\cal G}\,$ is at most $\,\ell\,.\,$}

\vspace*{0.1in}

\hspace*{0.5in}Refer to \,{\bf \S\,A.4}\, in the e\,-\,Appendix for the case when $\,\Delta_o^{\!(h_\ell)} \ {\cal P}_\ell \not\equiv 0$\,.\,

\vspace*{0.3in}

{\bf \large {\bf  \S\,5. \ \  Mezzo\,-\,scale effect of the global harmonic term.}}\\[0.15in]
In this section we show that for estimates of $v_i$ with accuracy of order $O_{\lambda_i}(n \,-\,2)$ or better, the contribution from other blow\,-\,up points has to be taken into account. We continue to assume \\[0.1in]
(5.1) \ \ ``{\it the standard conditions\,} (1.6)\,,\, (1.25) {\it and} (1.26), {\it plus}\, (2.16)\,,\, {\it with the notations in\,} (1.10) {\it and} (1.11)\,"\,.

 \vspace*{0.2in}

{\bf \S\,5\,a. }   {\bf Rescaled harmonic part.} \ \
From Proposition 2.24, we can find small positive numbers $c_o$ and $c_1$ such that for $\,i\, \,\gg\, \, 1\,,\,$
$$
C^{-1} \cdot \lambda_i^{ {{n \,-\,2}\over 2}} \ \,\le\, \ v_i\, (y) \ \,\le\, \ C\, \lambda_i^{ {{n \,-\, 2}\over 2}} \mfor i \,\gg\,  1 \ \ \ {\mbox{and}} \ \ \ \ c_o \ \,\le\, \ |\,y|\ \,\le\, \ c_1\,.
$$
Together with the Harnack inequality (cf. Theorem 8.20 and Corollary 8.21 in \cite{Gilbarg-Trudinger}\,,\, p. 199), and the discussion in \,{\bf \S\,2\,c.1}\,,\, a
subsequence of

 \vspace*{-0.2in}

$$
\ \ \ \ \ \ \ \ \ \   \left\{\,M_i^{-1} \cdot u_i \right\} \ = \ \left\{ \  \lambda_i^{-{{n \,-\,2}\over 2}} \cdot u_i  \right\}  \ \ \ \ \ \ \ \ \ \  [\,M_i \ \ {\mbox{is \ \ given \ \ in}} \ \ (3.2)]\leqno (5.2)
$$
converges to a positive $\,C^{\,2}$ -\,function $\,{\cal H}_{\lambda}\,$ in
$S^n \setminus \{ \beta_o, \, \cdot \cdot \cdot\,, \, \beta_k\}\,.$\,  (  The convergence is uniform in every compact set in $\,S^n \setminus \{ \beta_o, \, \cdot \cdot \cdot\,, \, \beta_k\}\,.$)\,
With the stereographic projection $\,\dot{\cal P}\,$ onto  $\R^n$,
  \ ${\cal H}_{\lambda}\,$ can be expressed [\,cf.    (1.3)\, and \, (2.9)\,]\,   as
\begin{eqnarray*}
(5.3) \ \ \ \ \ \ \ \ \ \ H_\lambda \,(y) &:= &[\,{\cal H}_\lambda \circ {\dot{\cal P}}^{-1} (y)] \cdot \left({2\over {1 + |\,y|^{\,2}}} \right)^{\!\!{{n \,-\,2}\over 2}}\,,\\
(5.4) \ \ \ \ \ \ \ \ \ \ H_\lambda \,(y)  &= &\sum_{j \,=\, 0}^k \ {{{\cal A}_{\,l}}\over { |\, y \,-\, {{\hat {\!Y}}}_l|^{\,n \,-\,2} }} \ \ \
 \ \ \ \  \ \ \ \ \ \ \ \ \ \  \ \ \ \ \ \ \ \ \ \  \ \ \ \  \ \ \ \ \ \ \ \ \ \  \ \ \ \ \ \ \ \ \ \
\end{eqnarray*}
for $\,y \,\in\, \R^n \setminus
 \{\  {{\hat {\!Y}}}_o\,=\,0\,,\,\cdot \cdot \cdot\,, \ {{\hat {\!Y}}}_k \}\,$.
 Here  $\,{{\hat {\!Y}}}_j\,:= \,{\cal P} \,(\beta_j)\,,\,$ $0 \, \le\, j \,\le \,k\,,\,$ and $\,{\cal A}_{\,j}\,$ are positive constants [\,a constant times $A_{\,j}$ which appears in (2.10)\,].
 From the Proportionality Proposition 2.3, (5.3) and \,{\bf \S\,2\,f}\,,\, together with (1.10) and (1.11), we obtain
 $$
[\,v_i\,(\xi_i)]\cdot   v_i\, (y) \ = \  \lambda_i^{-{{n \,-\,2}\over 2}} \cdot v_i\, (y)   \ \longrightarrow \  {1\over {|\,y|^{\,n \,-\,2} }} \ + \ h\, (y) \ \ \  \ {\mbox{in}} \ \  C^2_{\mbox{loc}} (B_o \ ({\bar \rho}_1) \setminus \{ 0 \})\,.
 $$
 Recall that we assume (without loss of generality) $({\tilde c}_n\, K )\, (0)\, =\, n\, (n \,-\,2)\,.$\, Hence we know that

 \vspace*{-0.15in}

$$
{\cal A}_{\,o} \ = \ 1\,. \leqno (5.5)
$$
$$
{\mbox{Define}} \ \ \ \ \ \ \  {\mbox{H}}_{\lambda_{\ge 1}} \,(y)\ := \ \sum_{j\, =\, 1}^k \ {{{\cal A}_{\,j}}\over { |\, y \,-\,{{\hat {\!Y}}}_{\!j}|^{\,n \,-\,2} }} \ \ \ \ \ \ {\mbox{for}} \ \  \  y \,\in\, \R^n \setminus
 \{\  {{\hat {\!Y}}}_1\,,\,\cdot \cdot \cdot\,, \ {{\hat {\!Y}}}_k \}\,.\, \leqno (5.6)
$$
Note that $\,{\mbox{H}}_{\lambda_{\ge 1}}\,(y)\,=\,H_\lambda\, (y) \,-\,|\,y|^{-\,(n \,-\,2)}\ \ $ is well\,-\,defined and smooth on a neighborhood of $\,0\,.\,$

\vspace*{0.15in}

{\bf \S\,5\,b.}   \ {\bf Estimating}  $\,|\ {\cal V}_{\,i} \, ({\cal Y}) \,-\,A_1 \,({\cal Y})\,|\,$ {\bf on the mezzo\,-\,scale}\, $\,C_o \,\le\, \lambda_i\, |\,{\cal Y}| \,\le\, C_1\,.$

\vspace*{0.1in}

We start with the convergence occurring in (5.2)\,:
\begin{eqnarray*}
& \ &   \!\!\!\!\!\!\!\!\!\!\!\!\!\!\!\!\!\!\!\!\!\!\!\!\!{{u_i \,(x)}\over {\lambda_i^{{n \,-\,2}\over 2} }} \ \to \  {\cal H}_{\lambda} \,(x) \mfor x \,\in\, S^n \Big\backslash\!\left[\  \bigcup_{j\,=\,0}^k \ {\cal B}_{\beta_l}\, (\rho) \right]\,, \ \  \ \ \ \rho > 0 \ \ {\mbox{small \ \ but \ \ fixed}}\\[0.1in]
\Longrightarrow & \ &\!\!\!\!\!\!\Bigg\vert \ {{u_i\, (x)}\over {\lambda_i^{{n \,-\,2}\over 2} }} \ \,-\,\  {\cal H}_\lambda\, (x) \Bigg\vert \ \,\le\, \ \varepsilon    \mfor \ {\mbox{all}} \ \ \ i \,\gg\,  1 \ \ {\mbox{and}} \ \ x \,\in\, S^n \Big\backslash\!\left[\  \bigcup_{j\,=\,0}^k \ {\cal B}_{\beta_l}\, (\rho) \right]\\[0.1in]
\Longrightarrow & \ &\!\!\!\!\!\!\Bigg\vert \ {{u_i\, (x)}\over {\lambda_i^{{n \,-\,2}\over 2} }} \cdot \!\left( {2\over {1 + |\,y|^{\,2}}} \right)^{\,{{n \,-\,2}\over 2}} \! \,-\,   {\cal H}_\lambda\, (x) \cdot\! \left( {2\over {1 + |\,y|^{\,2}}} \right)^{\,{{n \,-\,2}\over 2}} \!\Bigg\vert \ \,\le\, \ \varepsilon \! \cdot \!\left( {2\over {1 + |\,y|^{\,2}}} \right)^{\,{{n \,-\,2}\over 2}}\\[0.1in]
& \ &  \ \ \ \ \  \mfor \ {\mbox{all}} \ \ \  i \,\gg\,  1 \ \ {\mbox{and}} \ \ x \,\in\, S^n \Big\backslash\!\left[\ {\bf N} \ \cup \  \bigcup_{j\,=\,0}^k \ {\cal B}_{\beta_l}\, (\rho) \right]; \ \  \ y = \dot{\cal P}  \, (x)
\end{eqnarray*}
\begin{eqnarray*}
\Longrightarrow & \ &\!\!\!\!\!\!\Bigg\vert \ {{v_i\, (y)}\over {\lambda_i^{{n \,-\,2}\over 2} }}   \ \,-\,\  H_\lambda \,(y)  \, \Bigg\vert \ \,\le\, \ {{C\,\varepsilon }\over {(1 + r)^n}}  \mfor\ {\mbox{all}} \ \ \ i \,\gg\,  1 \ \ {\mbox{and}} \ \ y \,\in\, \R^n \Big\backslash   \left[\  \bigcup_{l\,=\,0}^k \ B_{{\hat Y}_l} \ (r_o) \right]\\[0.15in]
\Longrightarrow & \ & \Bigg\vert \ {{v_i\, (y)}\over {\lambda_i^{{n \,-\,2}\over 2} }} \ \,-\,\ {1\over {|\, y|^{\, n \,-\,2} }}    \ \,-\,\ {{\mbox H}}_{ {\lambda}_{\ge 1}} \  (y)   \, \Bigg\vert \ \,\le\, \ {{C\,\varepsilon }\over {(1 + r)^n}}\ \ \ \ \ \   [\,{\mbox{using}} \ \ (5.5) \ \ {\mbox{and}} \ \ (5.6)\,] \ \ \ \ \ \ \ \ \ \
\\[0.15in]
(5.7) \ \cdot \cdot \cdot \cdot &\!\!\! \cdot\, \cdot \ \ & \Longrightarrow  \ \ \Bigg\vert \ v_i\, (y)    \ \,-\,\ {{\lambda_i^{{n \,-\,2}\over 2} }\over {|\, y|^{\, n \,-\,2} }}    \ \,-\,\  [\,{{\mbox H}}_{ {\lambda}_{\ge 1}} \  (y) \,]\cdot \lambda_i^{{n \,-\,2}\over 2}   \, \Bigg\vert \ \,\le \ \,{{C\,\varepsilon \cdot \lambda_i^{{n \,-\,2}\over 2}  }\over {(1 + r)^n}}
\end{eqnarray*}
for all $\,i \,\gg\,  1\,$ and $\,c_o \ \,\le\, \ |\,y|\ \,\le\, \  c_1\,.$\, As usual\,,\, $r \ = \ |\,y|\,.\,$ Note that we may take
$$
\displaystyle{c_1\ \,\le\, \ {1\over 2}\cdot \min_{1 \,\le\, j \,\le\, k}} \,|\,{\hat Y}_j| \,. \leqno (5.8)
$$
In (5.7), we replace $\,y\, \to \,\xi_i + \lambda_i \, {\cal Y}$\,.\, For $\,i \,\gg\,  1\,$ and $\,c_o \ \,\le\, \ |\ \xi_i + \lambda_i \, {\cal Y}|\ \,\le\, \  c_1\,$\,:

(5.9)
$$\Bigg\vert \ v_i\, (\xi_i + \lambda_i \, {\cal Y})    \ \,-\,\ {{\lambda_i^{{n \,-\,2}\over 2} }\over {|\  \xi_i + \lambda_i \, {\cal Y}|^{\, n \,-\,2} }}    \ \,-\,\  [\,{{\mbox H}}_{ {\lambda}_{\ge 1}} \      (\xi_i + \lambda_i \, {\cal Y})\,]\cdot \lambda_i^{{n \,-\,2}\over 2}   \, \Bigg\vert  \ \,\le\, \  {{C\,\varepsilon \cdot \lambda_i^{{n \,-\,2}\over 2}  }\over {(1 + |\,\xi_i + \lambda_i \, {\cal Y}|\,)^{\,n}}}
 $$
 \begin{eqnarray*}
  &\Longrightarrow &  \Bigg\vert \ {{ v_i\, (\xi_i + \lambda_i \, {\cal Y})}\over {v_i\,(\xi_i)}}    \ \,-\,\ {{\lambda_i^{\,n \,-\,2} }\over {|\  \xi_i + \lambda_i \, {\cal Y}|^{\, n \,-\,2} }}    \ \,-\,\  [\,{{\mbox H}}_{ {\lambda}_{\ge 1}} \     (\xi_i + \lambda_i \, {\cal Y})\,]\cdot \lambda_i^{\,n \,-\,2}   \, \Bigg\vert \ \,\le\, \ C\,\varepsilon   \cdot \lambda_i^{\,n \,-\,2}   \\[0.15in]
 & \Longrightarrow  & \Bigg\vert \ {\cal V}_{\,i}\, ({\cal Y})    \ \,-\,\ {1\over {|\  [\, \lambda_i^{-1}\cdot \xi_i\,] + \, {\cal Y}|^{\, n \,-\,2} }}    \ \,-\,\ \lambda_i^{\,n \,-\,2}  \cdot  [\,{{\mbox H}}_{ {\lambda}_{\ge 1}} \       (\xi_i + \lambda_i \, {\cal Y})\,]   \,      \Bigg\vert \  \,\le\,  \  C\,\varepsilon \cdot \lambda_i^{\,n \,-\,2}  \ .\\
\end{eqnarray*}
Next, we seek to show that the second term in the last inequality above is ``close" to $A_1$ in the mezzo\,-\,range,  under the condition that $|\, \lambda_i^{-1}\cdot \xi_i\,| = O \, (1)\,.$\, We first note that
\begin{eqnarray*}
& \ &  \ \ y \ = \ \xi_i + \lambda_i \, {\cal Y}\,, \ \ \ \ {\mbox{where}} \ \ \ c_o \ \,\le\,  \ |\,y| \ \,\le\, \ c_1  \\[0.1in]
& \ &  \!\!\!\!\!\!\!\!\!\!\!!\!\!\!\!\!\!\!\!\!\!\!\!\!\!\!\!\!\!\!\!\!\!\!\!\!\!\!\!\!\!\!\!\!\!\!\! \Longrightarrow  \ \      c_o \, \le\,  |\ \xi_i + \lambda_i \, {\cal Y} |\ \,\le\, \   c_1 \
  \Longrightarrow \   c_o \cdot \lambda_i^{-1}    \,\le\,   |\ ( \lambda_i^{-1}\cdot \xi_i\,) \,+   \, {\cal Y} |  \,\le\,    c_1\cdot \lambda_i^{-1} \ \ \ \ \ \ \ \ \ \ \ \ \ \ \ \ \ \ \ \ \ \\[0.1in]
  & \ &  \   \ \ \ \ \ \ \ \ \ \ \   \ \ \ \ \ \ \ \ \ \ \   \   \ \ \ \ \ \ \ \ \ \ \     [\,{\mbox{assuming}} \  \lambda_i^{-1} \cdot |\,\xi_i| \,=\, O \, (1) \,] \ \ \\[0.1in]
 (5.10) \cdot \cdot \cdot  \  & \Longrightarrow  & [\,c_o \,-\,o\, (1)\,] \cdot \lambda_i^{-1}  \ \,\le\, \ |   \, {\cal Y} |\ \,\le\, \    [\,c_1 + o\, (1)\,] \cdot \lambda_i^{-1}\ ,  \end{eqnarray*}
where \,$o\, (1) \to  0 \,$ as $\,i \to \infty\,.$\,
Thus again, for $\,i \,\gg\,  1\,$ and $\,c_o \ \,\le\, \ |\ \xi_i + \lambda_i \, {\cal Y}\,|\ \,\le\, \  c_1\,$,\,  we continue with
\begin{eqnarray*}
(5.11)& \ & \Bigg\vert \ A_1\, ({\cal Y} )    \ \,-\,\ {1\over {|\  [\, \lambda_i^{-1}\cdot \xi_i\,] + \, {\cal Y} |^{\, n \,-\,2} }}      \  \Bigg\vert  \\[0.15in]
& = & \Bigg\vert \ \left({1\over {1 + |\,{\cal Y} |^{\,2}}}  \right)^{\!\!{{n \,-\,2}\over 2}}    \ \,-\,\ \left({1\over {  |\,{\cal Y} |^{\,2} + 2\, {\cal Y}  \cdot [\, \lambda_i^{-1}\cdot \xi_i\,] + O \, (1)}}  \right)^{\!\!{{n \,-\,2}\over 2}}        \, \Bigg\vert \\[0.15in]
 & \,\le\, & \ {2\over {n \,-\,2}} \cdot  \Bigg\vert\  {1\over {1 + |\,{\cal Y} |^{\,2}}} \  \,-\,\ {1\over {  |\,{\cal Y} |^{\,2} + 2\, {\cal Y}  \cdot [\, \lambda_i^{-1}\cdot \xi_i\,] + O \, (1)}}
\ \Bigg\vert \times\\[0.15in]
 & \ & \   \times \ \max \, \left\{ \ \left({1\over {1 + |\,{\cal Y} |^{\,2}}}  \right)^{\, {{n \,-\,2}\over 2} \,-\, 1}\  \  {\bf ,} \ \     \ \left({1\over {  |\,{\cal Y} |^{\,2} + 2\, \langle{\,\cal Y}\,, \    [\, \lambda_i^{-1}\cdot \xi_i]\,\rangle \, + \, O \, (1)}} \right)^{\, {{n \,-\,2}\over 2} \,-\, 1} \,\right\} \ \ \\[0.15in]
 & \,\le\, & C' \cdot  \Bigg\vert\  {{ 2\, \langle{\,\cal Y}\,, \    [\, \lambda_i^{-1}\cdot \xi_i]\,\rangle  + O \, (1) \,-\,1} \over { [\,1 + |\,{\cal Y} |^{\,2}] \cdot [\, |\,{\cal Y} |^{\,2} + 2\, \langle{\,\cal Y}\,, \    [\, \lambda_i^{-1}\cdot \xi_i]\,\rangle  + O \, (1)]}}
\ \Bigg\vert \,\cdot  \left( {1\over {|\, {\cal Y} |^{\,2}}}\right)^{\, {{n \,-\,2}\over 2} \,-\,1}  \\[0.15in]
 & \,\le\, & C' \cdot  \left[{1\over {|\,{\cal Y} |^{\,4} }}\  + \ {{O\, (1)}\over {|\,{\cal Y} |^{\,4} }} \  + \ {{ |\,{\cal Y} | \cdot |\, \lambda_i^{-1} \cdot \xi_i\,| }\over {|\,{\cal Y} |^{\,4} }}  \right] \cdot  \left( {1\over {|\, {\cal Y} |^{\,2}}}\right)^{\, {{n \,-\,2}\over 2} \,-\,1}  \\[0.15in]
 & \,\le\, & C' \cdot  \left[{1\over {|\,{\cal Y} |^{\,4} }}\  + \ {{O \, (1)}\over {|\,{\cal Y} |^{\,4} }} \  + \ {{  O\, (1) }\over {|\,{\cal Y} |^{\,3} }}  \right] \cdot  \left( {1\over {|\, {\cal Y} |^{\,2}}}\right)^{\, {{n \,-\,2}\over 2} \,-\,1} \\[0.15in]
 & \,\le\, & C' \cdot   {{  O\, (1) }\over {|\,{\cal Y} |^{\,3} }}  \ * \ \left( {1\over {|\, {\cal Y} |^{\,2}}}\right)^{\, {{n \,-\,2}\over 2} \,-\,1} \   \le\  C' \cdot   {{  O\, (1) }\over {|\,{\cal Y} | }}  \cdot \left( {1\over {|\, {\cal Y} |^{\,2}}}\right)^{\, {{n \,-\,2}\over 2}  } \  \,\le\,  \  O\, (1) \cdot \lambda_i^{\, n \,-\,1}\ .  \ \ \ \ \ \ \ \ \ \ \
\end{eqnarray*}
In the above, we apply the inequality
$$
   a \,>\, b\, >\, 0 \ \   {\mbox{and}} \, \    p \,\ge \,1  \ \ \Longrightarrow \ \ a^p -\, b^p \ \,\le\, \ p^{-1} \cdot (a \,-\,b) \cdot  a^{\,p \,-\,1}\,.\\[0.1in]
$$
Note that when $\,j\,\not= \,0\,,\,$ from (5.8) and (5.10), we have
\begin{eqnarray*}
(5.12)  & \ &   {{A}\over {|\, (\xi_i + \lambda_i\, {\cal Y}) \,-\, {\hat Y}_j\,|^{\, n \,-\,2} }}\ = \ {{A}\over {|\, (\lambda_i\, {\cal Y}  \,-\,{\hat Y}_j\,) + \xi_i \,|^{\, n \,-\,2} }} \\[0.15in]
  & = &  \!\!\!{{A}\over {|\,  \lambda_i\, {\cal Y} \,-\, {\hat Y}_j\,\,|^{\, n \,-\,2} }} \times {1\over { \left( 1 \, + \, {{ \xi_i}\over { |\, \lambda_i \,{\cal Y} \,-\,{\hat Y}_j|}}\right)^{\, n \,-\,2} }}\ = \
 {{A}\over {|\,  \lambda_i\, {\cal Y}  \,-\, {\hat Y}_j\, \,|^{\, n \,-\,2} }} \cdot \left[ \ 1 + O \, (|\,\xi_i\,|\,) \,\right] \ \ \ \ \ \ \ \ \ \ \ \ \ \ \ \ \ \ \
  \end{eqnarray*}
  for $\,[\,c_o \,-\,o\,(1)]  \cdot \lambda_i^{-1} \,\le \,|\, {\cal Y}|\,\le \,[\,c_1 + o\,(1)]  \cdot \lambda_i^{-1}\,$.\, Thus if we install the terms
$$
{\mbox H}_{\ge 1}\, ({\cal Y}) \ := \ \sum_{j \ge 1} \ {{{\cal A}_{\,j}}\over {|\,  \lambda_i\, {\cal Y}  \,-\, {\hat Y}_j \,|^{\, n \,-\,2} }}\ , \leqno (5.13)
$$
  [$\ {\mbox H}_{\ge 1}\,({\cal Y})\,$ is smooth and harmonic in $\,B_o\,(c_1 \cdot \lambda_i^{-1}\,)\,]\,,\,$ and apply the triangle inequality,   (5.8), (5.11) and (5.12) furnish us with the following mezzo\,-\,scale estimate.\\[0.1in]
{\bf Lemma 5.14.} \ \ {\it For $\,n \ge 3\,,\,$ under the conditions and notations in\,} (5.1)\, {\it for $\,\{ u_i \}\,\,$ $\,\{v_i\}\,$,\, $\,\lambda_i\,$ and $\,\xi_i\,,\,$ assume also that $\,\lambda_i^{-1}\cdot |\,\xi_i| = O\, (1)$\,.\, For any $\,\varepsilon \,>\, 0$\,,\, we have}

\vspace*{-0.25in}

$$
 |\ {\cal V}_i\, ({\cal Y}) \,-\,A_1\, ({\cal Y}) \,-\,\lambda_i^{\,n \,-\,2}\cdot {\mbox H}_{\ge 1}\, ({\cal Y})\,| \ \,\le\,  \ C\,\varepsilon \cdot  \lambda_i^{\,n \,-\,2}   \ + \     O\, (1) \cdot   \lambda_i^{\, n \,-\,1} \leqno (5.15)
 $$
{\it  for all \,$i \,\gg\, 1\,$ \,and\ } ${\bar c}_1 \cdot \lambda_i^{-1} \ \,\le\, \ |\, {\cal Y}|\ \,\le\, \ {\bar c}_2 \cdot \lambda_i^{-1}\ .$  {\it \,Here $\  {\bar c}_1 \, > \, 0$ can be taken to be any small  \,} ({\it but fixed}\,) {\it constant as long as $\ {\bar c}_1 \,<\, {\bar c}_2$\,.\,} [\,{\it In\,} (5.15)\,,\, {\it  the order in the right hand side  is} \,$\,O_{\lambda_i} (n \,-\,2)\,.$\,]

\newpage

{\bf \large {\bf  \S\,6. \ \  Second order blow\,-\,up argument and the proof of  }}\\[0.075in]
\hspace*{0.56in}{\bf \large {\bf  Main Theorem 1.14.}}\\[0.15in]
In this section, {\it we take up all the assumptions stated in Main Theorem}\, 1.14.
%
%
To begin with, we observe the following.

\vspace*{0.15in}

{\bf Proposition 6.1.} \ \ {\it For $\,n \,\ge\, 4\,,\,$  under the general conditions listed in} (5.1)\,,\, {\it we also
take the following conditions\ } (i)\,--\,(iii)\ {\it into account}\,.\\[0.075in]
(i) \ \  $\ 0$ {\it is a simple blow\,-\,up point for $\,\{ v_i\}\,.$}\\[0.075in]
(ii) \ \,    $K $ {\it is given by\,} (1.8) \,{\it in\,} $B_o \, (\rho_o)$\,,\,  {\it where\ } $2 \,\le \, \ell < \, n- 2\,.\,$ \\[0.075in]
(iii) \     {\it The parameters $\,\lambda_i\,$ and $\,\xi_i\,$ corresponding to the simple blow\,-\,up point at $\,0\,$ } \\[0.075in]
 \hspace*{0.31in}  \,[\,{\it via\,} (1.10) {\it and\,} (1.11)\,]\, {\it satisfy\,} (1.12)\,,\, {\it that is\,,\,} $|\,\xi_i| = o\, (\lambda_i)\,.$ \\[0.075in]
(iv) \  \,$\ell\,$ {\it is even}\,.\,\\[0.075in]
{\it Then $\,\Delta_o^{\!(h_\ell)} \, {\bf P}_\ell\, (y) \ = \ 0\,$}\, [$\,{\bf P}_\ell\,$ {\it is found in\,} (1.8)\,.\, {\it When $\ell$ is even}, $\,h_\ell \,=\, \ell/2\,,\,$ {\it and\,}  $\,\Delta_o^{\!(h_\ell)} \, {\bf P}_\ell$\, {\it is a number}\,]\,.\,  {\it The same conclusion also holds  when $\,\ell = n \,-\,2\,$ with an additional assumption that $\,0\,$ is\, the only blow\,-\,up point} (\,$\ell$ {\it is still required to be even\,})\,.

\vspace*{0.15in}

\hspace*{0.5in}The key of the proof is to combine  the change of center formula (see A.6.33 in the e\,-\,Appendix)  with the condition  \,$|\,\xi_i| = o\, (\lambda_i)\,$.\,   Other  arguments actually proceed in   similar fashion as those found  in \cite{Leung-Supported} and \cite{Li-1}\,.\,  For the benefit of readers, we present the estimates in \,{\bf \S\,A.6.\,d}\, in the e\,-\,Appendix. \bk
Together with condition (1.19) in the Main Theorem, Proposition 4.49 and Proposition 6.1\,,\, we  can secure a solution $\,\Pi_{\bf p}\,$ of the equation

$$
\Delta_o\, \Pi_{\bf p} \ + \ n \,(n + 2) \,A_1^{4\over {n \,-\,2}} \cdot \Pi_{\bf p} \ = \  {\bf P}_\ell \,\cdot A_1^{{n + 2}\over {n \,-\,2}} \ \ \ \ \ \ \ \ \  {\mbox{in}} \ \ \ \R^n\ .
$$

 \vspace*{-0.2in}

$$
 {\mbox{Moreover\,,}} \ \   \,\Pi_{\bf p}\, ({\cal Y}) \,= \, {{\Gamma_{\bf p} ({\cal Y})  }\over { (1 + {\cal R}^{\,2})^{\,{n\over 2}} }} \   \  \mfor {\cal Y} \,\in\, \R^n\,, \ \ \ \ {\mbox{where}} \ \  {\cal R}\, = \,|\, {\cal Y}|\,. \leqno (6.2)
$$

Thanks to Proposition 4.49, the precise form of $\,\Gamma_{\!\bf p}\,$ is known once $\,{\bf P}_{\!\ell}\,$ is given.
 It

 follows from (3.11) that

 \vspace*{-0.25in}

$$
 \Delta_o \,(\,{\cal V}_{\,i} \,-\,A_1 \,-\,\lambda_i^{\,\ell}\cdot \Pi_{\bf p} \,)+  n\,(n + 2)  \,A_1^{4\over {n \,-\,2}} \ (\,{\cal V}_{\,i} \,-\,A_1 \,-\,\lambda_i^{\,\ell}\cdot \Pi_{\bf p} \,) \ = \  {\bf  RM}   \leqno (6.3)
$$
in $\,\R^n\,,\,$ where the `remainder' $\,{\bf RM}\,$ is given in (3.12)\,,\, Cf. (3.1)\,--\,(3.3)\,.

\vspace{0.2in}

{\bf \S\,6\,a\,.} \  {\bf    Inclusion of the harmonic term via interpolation\,.} \ \
 Consider

 \vspace*{-0.25in}

\begin{eqnarray*}
 (6.4) \ \ \   {\cal D}_i^{\Pi} \, ({\cal Y}) & := & \left[ \ {\cal V}_{\,i}\, ({\cal Y}) \,-\,A_1\, ({\cal Y}) -\,\lambda_i^{\,\ell}\cdot \Pi_{\bf p} ({\cal Y}) \,-\,\ \lambda_i^{\,n \,-\,2}\cdot {\mbox H}_{{\ge 1}}\, ({\cal Y})\,\right] \ +\\[0.05in]
  & \ & \ \  \ \ \ \ \ \ \ \ \ \ \ \ \ \ \ \  \    \ \ \ \  +\  [\,\lambda_i^{\,n \,-\,2} \,\cdot\, h_o \,] \cdot  \left[\, 1 \,-\,{{\tilde {\cal R}\,({\cal Y})}\over {c\,\lambda_i^{-1} }}\,\right] \ \ \  \ \ \  \ \ \  \ \ \  \ \ \  \ \ \
\end{eqnarray*}

\vspace*{-0.1in}

in the region $\,  |\,{\cal Y}|\, \,\le\, \, c\,\lambda_i^{-1}\,.$\,\,  Here $\,c \,\le\, c_1\,$ is a positive constant to be fixed [\,(cf. (5.8)]\,,\, $\,{\mbox H}_{{\ge 1}} \,$ is given in (5.13)\,,\, and
$$ \ \ \ \ \ \ \ \
 h_o \, := \, {\mbox H}_{{\,\ge 1}}  \, (0)\ \ \ \ \ \ \ \ \ [\,{\mbox{that  \ \ is\,,\, \ \, setting \ \ }} {\cal Y} \ = \ 0 \ \ {\mbox{in \ \ (5.13)}}\, ]  \ .\leqno (6.5)
$$
In addition,  $\,\tilde {\cal R} \,\in\, C^\infty  \ ( \R^n)\,$ is a non\,-\,negative function which satisfies\\[0.05in]
(6.6)

 \vspace*{-0.25in}

$$
\tilde {\cal R} \,({\cal Y}) =  |\,{\cal Y}| \!\!\!\!\!\!\mfor |\,{\cal Y}|\ge 1\,, \ \  \tilde {\cal R}\,(0)\, = \,0\,, \ \ \btd_{\cal Y} \, \tilde {\cal R}\,(0)\, =\, \vec{\ 0}\,, \ \ {\mbox{and}} \ \ |\,\Delta_o \ \tilde {\cal R}|\le C \ \ {\mbox{in}} \ \    \R^n.
$$

 \vspace*{0.1in}

{\bf \S\,6\,a.1.  } {\it Joint between bubble estimate and the global harmonic term\,.} \ \
From (1.4)\,,\, (1.6)\, and \,(2.17), we know that $\,{\cal V}_i \, (0) = A_1\, (0) = 1\,$ for all $\,i\,.\,$ Here, we group together the terms in (6.4) which link to the harmonic function $\,{\mbox{H}}_{\ge 1}\,$,\, and note that\\[0.1in]
(6.7)

\vspace*{-0.3in}

\begin{eqnarray*}
\lambda_i^{\,n \,-\,2}\cdot {\mbox H}_{\ge\, 1}\, ({\cal Y}) \, + \, [\,\lambda_i^{\,n \,-\, 2}  \cdot  h_o \,] \!\cdot  \! \left(\, 1 \,-\,{{\tilde {\cal R}\,({\cal Y})}\over {c\,\lambda_i^{-1} }}\,\right) & = & 0 \ \ \ \ \ \   {\mbox{if}} \ \ {\cal Y} \ = \, 0\,,\\[0.1in]
\lambda_i^{\,n \,-\,2}\cdot {\mbox H}_{\ge\, 1}\, ({\cal Y}) \, + \, [\,\lambda_i^{\,n \,-\, 2}  \cdot  h_o \,] \!\cdot  \! \left(\, 1 \,-\,{{\tilde {\cal R}\,({\cal Y})}\over {c\,\lambda_i^{-1} }}\,\right) & = & \lambda_i^{\,n \,-\,2}\cdot {\mbox H}_{\ge\, 1}\, ({\cal Y})
\end{eqnarray*}
if $\,|\,{\cal Y}| = c\, \lambda_i^{-1}\,.$ That is, via $\,\tilde {\cal R}\,$,\, the bubble estimate and the global harmonic term are joint. See the comments in \,{\bf \S\,2  d.\,1}\,,\, (2.14), (2.15), (2.22) and (2.23).

\vspace*{0.15in}

{\bf \S\, 6\,b.} \   {\bf   Ingredients for the method to work.} \ \ \cite{Compact} provides  the framework for  what we call the ``second level blow\,-\,up argument\," (see also \cite{Chen-Lin-3})\,.\,  It is an exquisite method which goes to the root of the blow\,-\,up phenomenon. Here we summarize the key steps and apply it to our situation.

\vspace*{0.15in}

{\bf \S\, 6\,b.1. }  {\it First order vanishing property.}\ \
From (1.6),\, (3.1),\, (3.2), Proposition 4.49 (from there we know that $\,\Gamma_{\bf p}$ contains no constant term)\,,\, and (6.7),  we obtain
$$
 {\cal D}_i^{\Pi} \, (0)\ = \ 0 \mfor  \ i \ \,\gg\,  \  1\,. \leqno (6.8)
$$
Referring to (5.13),  (observe that presence of $\lambda_i$ in the right hand side below)
$$
 {\partial \over {\partial {\cal Y}_{|_1}}}   \left[\  {{{\cal A}_{\,l}}\over {|\,  \lambda_i\, {\cal Y}  \,-\,   {\hat Y}_j\, |^{\, n \,-\,2} }} \right]\,\Bigg\vert_{\,{\cal Y} = \,0}  \!\! = \ (n \,-\,2)\,{\cal A}_l \cdot {{ \lambda_i\,{\hat Y}_{j_{\,|_1}} }\over {|\,   {\hat Y}_j \,|^{\, n} }} \ \ \  \    [\, {\hat Y}_j = ({\hat Y}_{j_{\,|_1}}\,,\, \cdot \cdot \cdot\, \ {\hat Y}_{j_{\,|_n}}\,)\,]\,. \leqno (6.9)
$$
From the definition of $\,A_1 \,=\,A_{1\,, \ 0}$ [\,cf. (1.4)]\,,\,  (1.6),\, (3.1),\, (3.2), Proposition 4.49 (from there we recognize that $\,\Gamma_{\bf p}$ contains no first order  terms)\,,\, (6.5)\,,\, (6.6) and (6.9), we obtain
$$
 \Vert \, \btd_{\cal Y}\,{\cal D}_i^{\Pi} \, (0)\, \Vert \  = \ 0\, (\lambda_i^{\, n \,-\,1}\,) \ \  \mfor   \  i \ \,\gg\,  \  1\,. \leqno (6.10)
$$

\vspace*{0.15in}

{\bf \S\, 6\,b.2. }   {\it The maximum in $B_o\,(c \, \lambda_i^{-1})$\,.\,} \ \
Because $\,\Delta_o\, [\, {\mbox{H}}_{\ge \,1} \,-\,h_o\,] = 0\,,\,$  we have
$$
   \Delta_o\, {\cal D}_i^{\Pi} \, \ + \ n\,(n + 2)\,A_1^{4\over {n \,-\,2}} \cdot {\cal D}_i^{\Pi}  \ = \ {\bf R.H.S\,}_i\ \ \ \ \  {\mbox{in}} \ \ \ B_o\, (c\,\lambda_i^{-1}\,)\ ,\leqno (6.11)
$$
where [\,by using (6.3) and (6.6)]\,,\,
\begin{eqnarray*}
(6.12) \ \    {\bf R.H.S\,}_i \!\! & := &\!\!   {\bf RM} \ -\  \lambda_i^{\, n \,-\,\,2} \left\{  \,{{h_o  }\over {c\, \lambda_i^{-1} }} \cdot\! \Delta_o \,  \tilde {\cal R}    \ + \ n \, (n + 2) \cdot {{h_o   }\over {c\, \lambda_i^{-1} }} \cdot A_1^{4\over {n \,-\,2}} \cdot \tilde {\cal R}\  \right.\\[0.1in]
&   \ &\ \ \ \ \ \ \ \ \  \ \ \ \ \ \ \ \ \  \ \ \ \ \ \ \ \ \   \ \left. + \ n\,(n + 2)\,A_1^{4\over {n \,-\,2}} \cdot \left[ \, {\mbox H}_{\ge 1} \,-\, h_o\ \right]   \, \right\}. \ \ \ \ \ \ \ \   \ \ \   \ \ \  \ \ \
\end{eqnarray*}
From (2.26) in Proposition 2.24\,,\, Proposition 4.49 and the expression for $\,\Pi_{\bf p}\,$ [\,that is, (6.2)]\,,\,  we have
$$
\ell \ \,\le\, \ n \,-\,2 \ \ \Longrightarrow \ \ \Lambda_{\,i} := \max_{|\,{\cal Y}|\le\, c \,\lambda_i^{-1}} |\, {\cal D}_i^{\Pi}\, ({\cal Y})| \ < \  \infty  \mfor i \,\gg\, 1\,.\leqno (6.13)
$$

\vspace*{-0.15in}

We assert that

\vspace*{-0.25in}

$$
\Lambda_i  = o_{\lambda_i} \, (\ell  )\ . \leqno (6.14)
$$

The assertion  is equivalent to
$$
\Lambda_i \,=\, o\, (1) \, \lambda_i^{\ell } \ \ \Longleftrightarrow \ \
{{  \lambda_i^{\ell}}\over {\Lambda_i   }} \,=\, {1\over {o\, (1)}}\ \ \Longleftrightarrow \ \ {{  \lambda_i^{\ell}}\over {\Lambda_i   }} \,\to\, \infty\,.\leqno (6.15)
$$
Suppose that this is not the case. Then (modulo a subsequence)
$$
{{  \lambda_i^{\ell} }\over {\Lambda_i   }} \ \,\le\, \ C   \mfor {\mbox{all}} \ \ i \,\ge\, 1\ \ \ \ \Longleftrightarrow \ \ \  {{1 }\over {\Lambda_i   }} \ \,\le\, \ {C\over {  \lambda_i^{\ell} }}   \mfor {\mbox{all}} \ \ i \,\ge\, 1\,. \leqno (6.16)
$$
In what follows, we seek to derive a contradiction from (6.16)\,.\,

\vspace*{0.15in}

{\bf \S\, 6\,b.3. }    {\it Renormalization.}\ \
Consider the function

\vspace*{-0.25in}

$$
{\bf W}_i :=   {{{\cal D}_i^{\Pi}}\over {  \Lambda_i   }} \ \ \ \ \ \ \ \ \  {\mbox{defined \ \ in}} \  \ \ \ \ \  \overline{B_o \left(c\,\lambda^{-1}_i\right)\!}\ \,. \leqno (6.17)
$$

\vspace*{-0.05in}

From (6.11), we have

\vspace*{-0.25in}

$$\Delta_o \, {\bf W}_i \ + \ n\,(n + 2) \, A_1^{4\over {n \,-\,2}} \cdot {\bf W}_i   \ = \ \Lambda_i^{-1}\cdot {\bf R.H.S\,}_i   \ \ \ \ \ {\mbox{in}} \ \  \  B_o \,( \,c \, \lambda_i^{-1}\,)\,.\, \leqno (6.18)
$$

\vspace*{0.15in}

{\bf \S\, 6\,b.4. }    {\it Order of magnitude of the  remainder.}\ \ The key property we want to show about $\,{\bf R.H.S\,}_i\,$ are the following (\,under the condition in the Main Theorem).
\begin{eqnarray*}
(6.19)\  & \ & \!\!\!\!\!\!\!\!\!\!{\mbox{Given  \  any  \ }} {\cal R}_o \ > \ 0\,,\ \ \Lambda_i^{-1}  \,\cdot \!|\, {\bf R.H.S\,}_i   \,| \  \to \   0 \ \   {\mbox{uniformly}} \ \ {\mbox{in}} \ \ B_o\, ({\cal R}_o)\,.\\[0.15in]
(6.20)\  & \ &   \!\!\!\!\!\!\!\!\!\!{{|\, {\bf R.H.S\,}_i\,({\cal Y})|}\over{\Lambda_i}}   \   \,\le\, \  {{C\,  }\over {(1 + {\cal R})^4}} \ \,+ \ O\,(\lambda_i)\cdot \chi_{B_o\,(1)}\ \   \ \ {\mbox{for}} \ \  \ {\cal R}\, \,\le\, \,c\,\lambda_i^{-1}. \ \ \ \ \  \ \ \ \ \  \ \ \ \ \ \ \ \ \ \ \ \ \ \ \ \ \ \ \ \
\end{eqnarray*}
Here $\,{\cal R} \ = \ |\,{\cal Y}|\,$ and $\,\chi_{B_o\,(1)}\,$ is the characteristic function of the unit ball.
These are demonstrated in {\bf \S\,6\,c}\,.

\vspace*{0.15in}

{\bf \S\, 6\,b.5. }   {\it Vanishing on the whole.}\ \
From (6.8), (6.9), (6.10), (6.16) and (6.17)\,,\, we know that

\vspace*{-0.3in}

$$
|\, {\bf W}_i| \ \,\le\, \  1 \ \ \ \ \ {\mbox{in}} \ \ \ B_o\, (c\,\lambda_i^{-1})\,,\leqno (6.21)
$$
$$
 {\bf W}_i \,(0) \ = \ 0 \mfor i \,\gg\, 1\,, \ \ \ \mbox{and} \ \ \ \ \btd \, {\bf W}_i\, (0) \,=\,O \, (\lambda_i) \ \to \  {\vec{\,0}} \ \ \   {\mbox{as}} \ \   i \,\to \,\infty\,. \leqno (6.22)
$$
It follows from (6.11), (6.18), (6.19), (6.21)  and standard elliptic theory that (modulo a subsequence)   $$\,{\bf W}_i \longrightarrow {\bf W} \ \ {\mbox{     \  uniformly  \  in  \ every  \ compact  \ subset  \ of}} \ \ \R^n. \leqno (6.23)
 $$
Here $\,{\bf W}\,$ is a $\,C^2$\,-\,function  satisfying
$$
\Delta_o \,{\bf W} \ + \ n \,(n + 2)\,A_1^{4\over {n \,-\,2}} \cdot  {\bf W} \ = \  0  \ \ \ \ \ \ {\mbox{in}} \ \ \ \R^n. \leqno (6.24)
$$
In addition, from (6.22),
$$
{\bf W}\,(0) \ = \ 0 \ \ \ \ \& \ \ \ \btd  {\bf W}\,(0) \ = \ {\vec{\,0}}\ .\leqno (6.25)
$$
Moreover, in \,{\bf \S\,6\,b.7}\,,\, we show that
$$
|\, {\bf W} \, ({\cal Y})| \,\to\, 0 \ \ \ \ {\mbox{when}} \ \ \ |\, {\cal Y}|\,\to\, \infty\,.\leqno (6.26)
$$
A standard boot\,-\,strap argument shows that $\,{\bf W}$ is smooth in $\R^n$\,.\,
It follows from the Liouville\,-\,type theorem  for  (6.24) (see Lemma 2.4 in \cite{Chen-Lin-3}, cf. also \cite{Progress-Book}\,) that
$$
{\bf W} \ \equiv \ 0 \ \ \ \ \ {\mbox{in}} \ \ \ \R^n \ \ \Longrightarrow \ \ {\bf W}_i \ \to \  0 \ \ {\mbox{uniformly \ \ in}} \ \ B_o \, ({\cal R}_o) \,\subset\,  \R^n.\leqno (6.27)
$$
Here $\,{\cal R}_o\,$ can be any given positive number.
 On the other hand, by the definition of $\Lambda_i$\,,\, there is a point
$$
 {\cal Y}_{\mu_i} \ \,\in\, \ \overline{B_o \,(c\,\lambda_i^{-1})\!\!} \ \  \ \ \ \ {\mbox{such\ \  that}} \ \ \  {\bf W}_i \, ({\cal Y}_{\mu_i})\  = \ 1\,.\leqno (6.28)
$$
We produce a contradiction with (6.27)    by showing that we can find a positive number $\,{\cal R}_o\,$ such that

\vspace*{-0.25in}

$$
|\, {\cal Y}_{\mu_i}| \ \,\le\, \  {\cal R}_o  \ \ \mfor {\mbox{all}} \ \ \ \ i \,\gg\,  1 \ \ \ \ \ \ \left( \Longleftrightarrow \ \ \max_{\overline{B_o\,({\cal R}_o) } } \ {\bf W}_i \ = \ 1 \right)\,. \leqno  (6.29)
$$

 \vspace*{0.1in}

{\bf \S\,6\,b.6.}    {\it Smallness of $\,|\,{\bf W}_i\,|\,$ near the boundary\,} $\,\partial B_o\,(\,c\,\lambda_i^{-1}\,)$\,. \ \ Given (5.15) in Lemma 5.14, we turn our attention to $\,\lambda_i^\ell\cdot \Pi_{\bf p}\,$   in \,(6.4)\,.\, Based on Proposition 4.49 and condition (6.2), we have
   $$
   {\mbox{the \ \  degree \ \ of \ \ the \ \ polynomial

     \ \ }} \Gamma_{\bf p} \ \ \,{\mbox{in \ \ (6.2)\ \ \,is \ \ at \ \ most  }} \ \ n \,-\,2\,. \leqno (6.30)
   $$
   It follows that
\begin{eqnarray*}
(6.31)  \ \ \ \ \ \ \ \ \  \,     & \ & \!\!\!\!\!\!\!\!\!\!\!\!\!\!| \  \lambda_i^{\ell} \cdot  \Pi_{\bf p} ({\cal Y})| \ = \ \lambda_i^{\ell} \cdot {{ | \ \Gamma_{\bf p}  ({\cal Y})| }\over { (1 + {\cal R}^{\,2})^{n\over 2} }} \\[0.1in]
& \ &  \ \ \ \ \ \  \ \ \  \  \ \le\, \    C\!\cdot\! {{ \lambda_i^{\ell}\,(1 + {\cal R})^{\,n \,-\,2} }\over { (1 + {\cal R})^n}} \ \,\le\,  \ {{ C_1\,\lambda_i^{\ell  } }\over {  (1 + {\cal R})^2 }}   \ \,\le\, \   C_2\, \lambda_i^{\,\ell + 2}  \\[0.1in]
& \ &
\hspace*{1.2in} \mfor {\cal R}\, =\,|\, {\cal Y}| \,\ge\,\,(1 \,-\,\delta) \cdot c\,\lambda_i^{-1}\\[0.1in]
(6.32) \ \cdot \cdot \cdot   \cdot \cdot \cdot  & \ & \!\!\!\!\!\!\!\!\!\Longrightarrow    \ \Lambda_i^{-1} \cdot | \  \lambda_i^{\ell} \cdot  \Pi\, ({\cal Y})| \ = \ O \, (\lambda_i^2) \ \ \  [\,{\mbox{via}} \ \ (6.16)]  \ \ \ \ \ \ \ \ \ \ \ \ \ \ \ \ \ \ \ \ \ \ \ \ \ \
\end{eqnarray*}
for $\,{\cal R}\, \,\ge\,\,(1 \,-\,\delta) \cdot c\,\lambda_i^{-1}\,$.\, Together with (6.4), (6.6), Lemma 5.15   and (6.31), we have
 $$
 |\,{\bf W}_i\, ({\cal Y})| \,=\, O\, (\varepsilon) \ + O\,(\lambda_i^2) \mfor |\,{\cal Y}| \, = \, c\,\lambda_i^{-1}  \ \ \ {\mbox{and}} \ \ \ \ i \,\gg\, 1 . \leqno (6.33)
 $$

Moreover, for $\,(1 \,-\,\delta) \cdot [\, c\,\lambda_i^{-1}\,] \,\le\, |\,{\cal Y}|\, \,\le\, \,[\, c\,\lambda_i^{-1}\,]\,$,\, we have
$$
[\,\lambda_i^{\,n \,-\,2} \,\cdot\, h_o \,] \cdot  \bigg\vert \  1 \,-\,{{\tilde {\cal R}\,({\cal Y})}\over {c\,\lambda_i^{-1} }}\,\bigg\vert \, = \,[\,\lambda_i^{\,n \,-\,2} \,\cdot\, h_o \,] \cdot  \bigg\vert \  1 \,-\, {{|\,{\cal Y}|}\over {c\,\lambda_i^{-1} }}\,\bigg\vert \ \,\le\, \ [\,\delta \cdot h_o] \cdot \lambda_i^{\,n \,-\,2} \ . \leqno (6.34)
$$
Thus if we choose $\,\delta\, >\, 0\,$ to be small [\,relative to the constant $C$ in (6.16) and $\,h_o\,$ only]\,,\, and combine (6.32) with (6.33), (6.34) and  Lemma 5.15, we obtain
$$
 |\,{\bf W}_i  \, ({\cal Y})| \ \,\le\, \ {1\over 2} \ \  \mfor  (1 \,-\,\delta) \cdot [\, c\,\lambda_i^{-1}\,] \ \,\le\, \ |\,{\cal Y}| \ \,\le\, \  [\, c\,\lambda_i^{-1}\,] \ \ \ \ {\mbox{and}} \ \  \ i \,\gg\, 1\,. \leqno (6.35)
$$

\vspace*{0.1in}

{\bf \S\,6\,b.7.}    {\it The decay to \,$0$  --  proof of\ } (6.26)\,. \ \ To demonstrate (6.26), we show that for any positive number $\,{\tilde\varepsilon}\,$ (small and given), we can find a positive number $\,{\cal R}_{\tilde\varepsilon}$\, and a natural number $\,I_{\tilde\varepsilon}\,$,\, such that
$$
|\,{\bf W}_i\, ({\cal Y})| \ \,\le\, \   {\tilde\varepsilon} \mfor {\mbox{all}} \ \ i \,\ge \,I_{{\tilde\varepsilon}} \ \ \ {\mbox{and}} \ \ \ {\cal R}_{\tilde\varepsilon}\ \,\le\, \ |\, {\cal Y}| \ \,\le\, \   (1 \,-\,\delta)\cdot [\,c\, \lambda_i^{-1}\,]\,. \leqno (6.36)
$$
Here $\delta \,\in\, (0, \ 1)$ is fixed, and the integer $\,I_{{\tilde\varepsilon}}\,$ depends on $\,{\tilde\varepsilon}\,$ only\,.\, In particular, $\, (1 \,-\,\delta)\cdot [\,c\, \lambda_i^{-1}\,] \ \to \  \infty$\, as $\,i \,\to\, \infty\,.\,$  Once we have (6.36), together with the uniform convergence of $\,{\bf W}\,$ to $\,{\bf W}_i\,$ on any given compact subset of $\,\R^n$\,,\, we have
$$
|\,{\bf W} \, ({\cal Y})| \ \,\le\, \ {\tilde\varepsilon} \ \ \mfor {\mbox{all}} \ \ \   |\, {\cal Y}| \ \ge \  {\cal R}_{\tilde\varepsilon}\ \ \ \Longrightarrow  \ \ (6.26)\,.\leqno (6.37)
$$
To demonstrate the proof for (6.36),
via (6.\,33), we already know that $\,|\,{\bf W}_i|\,$ is `small' along the boundary $\,\partial B_o\,(\,c\,\lambda_i^{-1}\,)$\,.\,  The value in $\, B_o\,(\,c\,\lambda_i^{-1}\,)\,$  is governed by the equation describing $\,\Delta_o\, {\bf W}_i\,$ [\,that is, \,(6.18)\,]\,,\,  and  the
 Green representation formula, which, in the present situation, is given by
\begin{eqnarray*}
 (6.38)   & \ &  \!\!\!\!\!\!\!\!\!\!{\bf W}_i\, ({\cal Y}_{out} ) \ = \
 \int_{B_o (\,c\,\lambda_i^{-1}\,)} \!\! G_i\, ({\cal Y}, \  {\cal Y}_{out}) \left\{-\, n \,(n + 2) [\, A_1\, ({\cal Y})]^{4\over {n \,-\,2}} \cdot {\bf W}_i\,({\cal Y}) \right.\\[0.1in]
  & \ &\hspace*{1.6in} \ \ \ \ \ \  \ \ \ \ \ \ \ \ \ \left. + \  \Lambda_i^{-1} \, \cdot  {\bf R.H.S\,}_i\, ({\cal Y}) \right\}\, d\,{\cal Y}\\[0.1in]
 & \ & \ \ \   \ \ \ \ \ \ \ \ \ \ \ \ \ +\  \int_{\partial B_o (\,c\,\lambda_i^{-1}\,)}\,   [\,{\bf n}  \,\cdot \btd_{\cal Y} \, G_i ({\cal Y}, \  {\cal Y}_{out})]\, {\bf W}_i \,({\cal Y})\ dS_{\cal Y}
\end{eqnarray*}
for $\,{\cal Y}_{out} \,\in\, B_o\, (\,c\,\lambda_i^{-1}\,)\,.$
Here $\,G_i\,$ is the Green function for \,$\Delta_o$\, in $\,B_o\, (\,c\,\lambda_i^{-1}\,)\,$ with the Dirichlet boundary condition.
See, for example, \cite{McOwen} \,\, Note that   $$\displaystyle{G_i\, ({\cal Y}\,, \ {\cal Y}'\,) \ \approx \ -\, {1\over {(n \,-\,2)\,\Vert \,S^{n \,-\,1} \Vert }} \cdot {1\over {|\,{\cal Y}\,-\, {\cal Y}' |^{\,n \,-\, 2} }}\,}$$ when ${\cal Y}$ is close to ${\cal Y}'$\,.\,\, The sign is negative of the one used in \cite{Compact}\,.\, \,Using (1.4) and (6.21),   we obtain
$$
  \bigg\vert \,[\,A_1\,({\cal Y})]^{4\over {n \,-\,2}}  \,\cdot\, {\bf W}_i \,({\cal Y}) \bigg\vert \ \,\le\, \  \left( {1\over {1 + {\cal R}^2}}\right)^{2} \ \,\le\, \ {C\over {1 + {\cal R}^4 }} \ \ \mfor {\cal Y} \in\R^n\,. \leqno (6.39)
$$

Consider points $\,{\cal Y}_{out}\,$ so that
$$
|\,{\cal Y}_{out}| \ \le\  (1 \,-\,\delta) \cdot [\,c\, \lambda_i^{-1}]\,.
$$
Via proportional property on the Green function (see  \,{\bf \S\,A.5}\,  in the e\,-\,Appendix\,;\, cf. also pp. 157 in  \cite{Compact}\,)\,,\, we have
$$
\ \ |\,G_i\, ({\cal Y}, \ \,{\cal Y}_{out})| \,\le\, \left[ \,C_1 + {{C_2}\over {\delta^{\,n \,-\,2} }} \right] \cdot {1\over {|\,y \,- {\cal Y}_{out}|^{\, n \,-\,2} }} \mfor {\cal Y} \,\in\, B_o \, (c\,\lambda_i^{-1}\,) \,\setminus \{ {\cal Y}_{out}\}\,, \leqno (6.40)
$$
$$
|\ {\bf n} \cdot \btd\, G_i \,({\cal Y}, \ \, {\cal Y}_{out})|\ \,\le\, \ {{C_3 }\over {\delta^{\, n  }}} \cdot  \lambda_i^{\,n \,-\,1} \ \ \mfor {\cal Y} \,\in\, \partial B_o \, (c\,\lambda_i^{-1}\,)\,,\,\leqno (6.41)
$$
where $\,|\,{\cal Y}_{out}| \,\le\, (1 \,-\, \delta) \cdot [\,c\, \lambda_i^{-1}\,]\,.\,$
Here $\,C_1\,,\, \ C_2$  and $C_3$ are positive constant independent on $i\,$ and $\delta\,.\,$ It follows from (6.33) and (6.41) that

\vspace*{-0.2in}

\begin{eqnarray*}
 (6.42)  & \ & \bigg\vert \, \int_{\partial B_o (\,c\,\lambda_i^{-1})}\, [\,{\bf n}  \,\cdot \btd_{\cal Y} \, G_i ({\cal Y}, \  {\cal Y}_{out})\,]\, {\bf W}_i \,({\cal Y})\ dS_{\cal Y} \,\bigg\vert \\[0.1in]
&\le  &   {{C_2 \cdot \varepsilon}\over {\delta^{\, n  }}} \,\cdot\,  \lambda_i^{\,n \,-\,1} \cdot  \int_{\partial B_o (\,c\,\lambda_i^{-1})}\, dS_{\cal Y}  \\[0.1in]
 & \le &   {{C_2 \cdot \varepsilon}\over {\delta^{\, n  }}} \,\cdot\,  \lambda_i^{\,n \,-\,1} \cdot \left[\Vert\, S^{n \,-\,1} \Vert \cdot (c\, \lambda_i)^{n \,-\,1} \right]\ \,\le\, \ {{C \cdot \varepsilon}\over {\delta^n}}\,.\ \ \ \ \ \ \ \ \ \ \ \ \ \ \ \ \ \ \ \ \ \ \ \ \ \ \ \ \ \ \ \
\end{eqnarray*}
Here the consider $C$ is independent on $i$\,.\, Putting the information into (6.38)\,,\, together with (6.20) and (6.40), we obtain\\[0.1in]
(6.43)
\begin{eqnarray*}
 |\, {\bf W}_i \,({\cal Y}_{out} )\,|
&  \le  &\!\!\!   \left[ \,C_1 + {{C_2}\over {\delta^{\,n \,-\,2} }} \right] \cdot \int_{B_o (c\,\lambda_i^{-1})}\left( {1\over {|\,{\cal Y} \,-\,{\cal Y}_{out}|^{\,n \,-\,2}}}  \cdot {1\over {(1 + |\,{\cal Y}\,|)^{\,4} }} \right) dy \\[0.1in]
  & \ & \  +   \int_{B_o \,(1)} O \,(\lambda_i) \, + \, {{C \cdot \varepsilon}\over {\delta^n}}  \ \ \ \  \ \ \ \ \    ( \,{\mbox{from\  the\ harmonic \  term}}\,)\\[0.1in]
 & \le & \ C_\delta \cdot \left[ \  {1\over {(1 + |\,{\cal Y}_{out}|\,) }}   \ + \ O\,(\lambda_i) \ + \ O\, (\varepsilon)\,  \right]   \ \ \ \ \ \ \ \ \ \  \ \end{eqnarray*}
 for $\,|\,{\cal Y}_{out}| \,\le \,(1 \,-\, \delta)\,[\,c\,\lambda_i^{-1}\,]\,.$
 Refer to \cite{Chen-Lin-3} for the estimation of the first integral in (6.43)\,.\, Thus we can find $\,R_{\tilde\varepsilon} \,>\, 0\,$ and $\,I_{\tilde\varepsilon}\,$ such that for all $\,i \ \ge \ I_{\tilde\varepsilon}\,,\,$ we have
$$
|\, {\bf W}_i \,({\cal Y}  )| \ \,\le\, \ \tilde\varepsilon \ \ \ \mfor \ \ \ R_{\tilde\varepsilon} \ \,\le\, \ |\,{\cal Y}| \ \,\le\, \ (1 \,-\,\delta)\,[\,c\,\lambda_i^{-1}\,]\,. \leqno (6.44)
$$

\vspace*{0.1in}

{\bf \S\,6\,b.8.}   \ \  {\it Further restriction on the location of the maximum\,--\, proof of\,} (6.29)\,.\ \ In view of (6.35), we actually have

\vspace*{-0.1in}

$$
 \ \ \ \ \ \ \ \ \ \ \ \ \ \ \ \ \ \ {\cal Y}_{\mu_i} \ \,\in\, \ \overline{B_o \,\left(  (1 \,-\, \delta)\,[\,c\,\lambda_i^{-1} ]   \right)}\ \ \ \ \ \ \ \ \ \ \ \ \ \ \ \ \ \  [\,{\mbox{cf}}. \ (6.28)]\,.
$$
Here $\,\delta \,>\, 0\,$ is chosen close to $\,0\,$ (as explained in \,{\bf \S\,6\,b.6}\,)\,.\, Argue as in (6.43), we arrive at a similar conclusion
 $$
1 \ = \  |\  {\bf W}_i \,({\cal Y}_{\mu_i} )| \  \,\le\, \
  C \, \left[ \  {1\over { 1 + |\,{\cal Y}_{\mu_i}|  }} \ + \ O\,(\lambda_i) \ + \ O\, (\varepsilon)\,  \right] \,. \leqno (6.45)
$$
It follows that there is a fixed positive number  \,$R_o$\, such that
$$|\,{\cal Y}_{\mu_i}|\ \,\le\, \ R_o\, \mfor i \,\gg\, 1\,.\,$$
Hence we establish (6.29), and obtaining a contradiction to (6.27). Thus (6.16) must be wrong. That is, (6.14) holds.

\vspace*{0.2in}

{\bf \S\,6\,c.  Terms in the remainder ${\bf R.H.S.}_i$\, in (6.12)\,.}\ \ Here we verify (6.19) and (6.20). Recall that ${\bf RM}$ is decomposed into four components as expressed in (3.12)\,.\,

\vspace*{0.15in}

{\bf \S\,6\,c.1.}  \  {\it First term.} \ \ Under the condition $\,|\,\xi_i| \,= \,o_{\lambda_i} (1)\,,\,$ we have\\[0.1in]
(6.46)
\begin{eqnarray*}
 \!\!\!\!\!\!\!\!\!\! \!\!\!\!\!\!\!\!\!\! \!\!\!\!\!\!\!\!\!\!& \ &   \Lambda_i^{-\,1} \cdot \left[\,  |\,\xi_i|^{\,k} \cdot \left( \lambda_i \, |\,{\cal Y}| \,\right)^{\ell \,-\,k} \,\right] \times [\,A_1\, ({\cal Y})]^{\,{{n + 2}\over {n \,-\,2}}} \ \ \    (1 \,\le\, k \,\le\, \ell\,; \ \ |\,\xi_i|\, = \,o \, (1) \cdot \lambda_i\,) \\[0.1in]
& \,\le\, & \lambda_i^{-\,\ell} \cdot o\, (1) \cdot \lambda_i^k \cdot \lambda_i^{\ell \,-\,k} \times {\cal R}^{\, \ell \,-\,k}   \cdot \! \left( {1\over {1 + {\cal R}^{\,2}}}  \right)^{{n + 2}\over 2} \ \ \  \  \ [\ {\mbox{using}} \ \ (6.16)\,; \ \ {\cal R} = |\,{\cal Y}|\,]\ \ \ \ \ \ \ \ \ \ \ \ \\[0.1in]
   & \,\le\, & o \,(1)\,  \cdot (1 + {\cal R})^{\, \ell \,-\,k}   \cdot\! \left( {1\over {1 + {\cal R} }}  \right)^{ \!n + 2 } \ \le \  {{o \, (1)}\over {     (1 + {\cal R})^{\,5}  }} \ \ \ \ \  \mfor {\cal R} \ \,\le\, \  c\,\lambda_i^{-1} \\[0.15in]
   & \to & 0 \ \ \ \ \ \ \ \ \       {\mbox{uniformly \ \ in  \ \ }} B_o\, ({\cal R}_o) \ \ \  \ \ \ \ \ \ (\,i \,\gg\, 1\,; \ \ \ell \,\le\, n -2\,, \ \ \ \ k \,\ge\, 1)\ .
 \end{eqnarray*}

\vspace*{0.1in}

{\bf \S\,6\,c.2.}    \ {\it Second term.}\\[0.1in]
(6.47)
 \vspace*{-0.1in}
\begin{eqnarray*}
\!\!\!\!\!\!& \ & \Lambda_i^{-\,1} \cdot \left( \max_{B_o (\lambda_i\,{\cal Y} + \xi_{\,i})}\!\Vert\,\btd^{\,(\ell + 1)} K \,\Vert \right) \cdot  |\,\lambda_i\,{\cal Y} + \xi_{\,i}|^{\,\ell + 1} \times [\,A_1\, ({\cal Y})]^{\,{{n + 2}\over {n \,-\,2}}} \\[0.1in]
 &  \le  & \!\!\!\!C\,  \lambda_i^{-\,\ell}\ \bigg\vert \,\lambda_i \left( {\cal Y} + {{\xi_i}\over {\lambda_i}} \right) \bigg\vert^{\,\ell + 1}  \times {1\over {(1 + {\cal R})^{\,n + 2 }}}   \ \le \ C\, \,\lambda_i \cdot (1 + {\cal R} )^{\, \ell \,+\, 1    }  \,\cdot {1\over {(1 + {\cal R})^{\,n + 2 }}} \\[0.1in]
 & \ &  \ \ \ \  \ \ \  \ \ \ \ \ \ \ \ \ \ \ \ \  \ \ \ \  \ \ \  \ \ \ \ \ \ \ \ \ \ \ \ \  [\, {\mbox{via}} \ \ (6.16)\,; \ \ |\,\xi_i| \,=\, o \,(\lambda_i)\,,\ \ \ell \,\le\, n\,] \ \ \ \ \ \ \ \ \\[0.1in]
 &  \le  & \!\!\!\! {{C\,\lambda_i}\over {(1 + {\cal R})^{\,3} }} \ \le \ {{C_1}\over {(1 + {\cal R})^4}}
  \!\!\!\!\mfor \  {\cal R} \ \,\le\, \  c\,\lambda_i^{-1} \ \ \left[\ \Longrightarrow \ \ \lambda_i\, (1 + {\cal R})\ \le\ C_2\,\right]\,.
\end{eqnarray*}

 \vspace*{0.1in}

{\bf \S\,6\,c.3.} \ \  {\it Third term.} \ \ Using
$ \
a^{\,2} \ = \  (a \,-\,b)^2 \,+\, 2\,b\,(a \,-\,b) \,+\, b^{\,2}\,,\
$
we obtain\\[0.1in]
(6.48)
  \begin{eqnarray*}
& \ & \Lambda_i^{-1} \cdot  A_1^{ {{4}\over {n \,-\,2}}  -1} \, ({\cal V}_{\,i} \,-\,A_1)^{\,2} \ \ \ \ \ \ \ \ \ \ \ \ \ \   \left[\ {\mbox{take}} \ \  a \,=\,({\cal V}_{\,i} \,-\,A_1)\,, \ \ b \,=\,\lambda_i^{\,\ell} \cdot \Pi_{\bf p} \,  \right]\\[0.15in]
& \ & \!\!\!\!\!\!\!\!\!\!\!\!\!\!\!\le  A_1^{ {{4}\over {n \,-\,2}}  -1} \, \Lambda_i^{-1}   \,\{\, [\ {\cal V}_{\,i} \,-\,A_1 \,-\,\lambda_i^{\,\ell} \cdot \Pi_{\bf p} \,]^{\,2}  + 2 \,(\lambda_i^{\,\ell} \  \Pi_{\bf p}\,) \!\cdot\!\! |\,{\cal V}_{\,i} \,-\,A_1 \,-\,\lambda_i^{\,\ell} \, \Pi_{\bf p}\,| + (\lambda_i^{\,\ell} \cdot \Pi_{\bf p} )^{2}  \}\\[0.1in]
& \ & \!\!\!\!\!\!\!\!\!\!\!\!\!\!\!\le\left( {1\over {1 + {\cal R}^{\,2}}} \right)^{{6 \,-\,n}\over 2}  \,\left\{  \,|\,   ({\cal V}_{\,i} \,-\,A_1 \,-\,\lambda_i^{\,\ell} \cdot\Pi_{\bf p}| \cdot |\,{\bf W}_i  \ + \  O_{\lambda_i}(\,(n \,-\,2) \,-\,\ell)| \right. \\[0.1in]
 & \ & \ \ \ \ \ \ \ \ \    \left. +\  2 \,|\,\lambda_i^{\,\ell}\cdot \Pi_{\bf p}\,| \cdot |\,{\bf W}_i + O_{\lambda_i}(\,(n \,-\,2) \,-\,\ell)|\, + \, (\Lambda_i^{-1}\,\lambda_i^{\,\ell} ) \cdot  [\, \lambda_i^{\,\ell} \cdot (\Pi_{\bf p}\,)^{\,2} \right\}\\
& \ & \ \ \ \ \ \ \ \ \ \ \ \ \ \ \ \ \ \ \ \ \ \ \ \ \ \ \ \
\ \ \  \ \ \ \ \ \ \ \ \ \ \ \ \ \ \  \  \ \  [\, {\mbox{using}} \ \ (6.4)\,, \ \ (6.16) \ \ {\mbox{and}} \ \ (6.17)\,]
\end{eqnarray*}
 \begin{eqnarray*}
& \ & \!\!\!\!\!\!\!\!\!\!\!\!\!\!\!\le{C\over {(1 + {\cal R})^{\,6-n} }} \cdot \left\{  \ |\, ({\cal V}_{\,i}\, ({\cal Y}) \,-\,A_1\,({\cal Y}) \,-\,\lambda_i^{\ell}\cdot \Pi_{\bf p}\ | \ + \ \lambda_i^{\,\ell} \, |\, \Pi_{\bf p} ({\cal Y})| \ + \ \lambda_i^\ell \cdot |\,\Pi_{\bf p}({\cal Y})|^{\,2}\,\right\}
\\[0.1in]
 & \ & \ \ \ \ \ \ \   \ \ \ \ \ \ \ \  [\  {\mbox{as}} \ \ |\,{\bf W}_i + O_{\lambda_i}(\,(n \,-\,2) \,-\, \ell)| \\[0.1in]
& \ & \!\!\!\!\!\!\!\!\!\!\!\!\!\!\!\le|\,{\bf W}_i| + |\,  O_{\lambda_i}(\,(n \,-\,2) \,-\, \ell\,)| \ \,\le\, \ 1 + C \ \ {\mbox{when}} \ \ \ell \,\le\, (n \,-\,2)\,]\\[0.2in]
 & \ &   \longrightarrow \  0 \ \ \ {\mbox{in}} \ \ B_o\, (R_o) \ \   {\mbox{uniformly}} \ \ \ \       [\,{\mbox{from \ \ (2.59)}}\,:  \ |\, {\cal V}_i \,-\,A_1\,| \ \to \ 0 \ \ {\mbox{in}} \ \ B_o\, (R_o)\,] \end{eqnarray*}
 \begin{eqnarray*}
\ \ \ \ \ \ \ \ \ \  \ \ & \  & \!\!\!\!\!\!\!\!\!\!\!\!\!\!\!\!\!\!\!\!  \,\le\, \ {{C_1}\over {(1 + {\cal R})^{\,6-n} }} \cdot \left\{  \ |\, ({\cal V}_{\,i}\, ({\cal Y}) \,-\,A_1\,({\cal Y})| \ + \ \lambda_i^{\,\ell} \, |\, \Pi_{\bf p}({\cal Y})| \ + \ \lambda_i^\ell \cdot |\,\Pi_{\bf p}({\cal Y})|^{\,2}\,\right\}\\[0.1in]
& \  & \!\!\!\!\!\!\!\!\!\!\!\!\!\!\!\!\!\!\!\!  \,\le\, \  {{C_2}\over {(1 + {\cal R})^{\,6-n} }} \cdot \left\{ \ {1\over {(1 + {\cal R})^{\, n \,-\,2} }}     \ + \ {{  \lambda_i^\ell \cdot (1 + {\cal R})^\ell }\over {(1 + {\cal R})^n}}   \ + \ {{  \lambda_i^\ell \cdot (1 + {\cal R})^{2\,\ell} }\over {(1 + {\cal R})^{2\,n}}} \right\} \\[0.1in]
& \ &  \hspace*{1.75in}[\,{\mbox{via}} \ \, (3.7) \ \, {\mbox{and \ \    Proposition 4.49}}\,] \\[0.1in]
& \  & \!\!\!\!\!\!\!\!\!\!\!\!\!\!\!\!\!\!\!\!  \,\le\,\   {C\over {(1 + {\cal R})^{\,4} }}    \ \ \  \ \  \mfor {\cal R} \,\le\, c\,\lambda_i^{-1} \ \ \left[\ \Longrightarrow \ \ \lambda_i\, (1 + {\cal R}) \,\le\, C_2\,\right]\,.
 \end{eqnarray*}
Here $\,\ell \,\le\, n \,-\, 2\,.$
\vspace*{0.15in}

{\bf \S\,6\,c.4.} \ \  {\it Fourth term.} \ \ From (3.12), we obtain\\[0.1in]
(6.49)
\begin{eqnarray*}
& \ & \!\!\!\!\!\!\!\!\!\!\!\!\!\!\!\!\!\!\!\! \!\!\!\!\!\!\!\!\!\!\!\!\!\!\!\!\!\!{\bf RM}_4 \, ({\cal Y})  \\[0.1in]
& \ & \!\!\!\!\!\!\!\!\!\!\!\!\!\!\!\!\!\! \!\!\!\!\!\!\!\!\!\!\!\!\!\!\!\!\!\!\!\!=\  \Lambda_i^{-\,1} \cdot \left\{ \ O \,  \left(\, \max_{|\,\lambda_i \, {\cal Y} \,+\, \xi_i|\,\le \,\rho_o)}  \Vert \, \btd^{(\ell)} \,K \, \Vert \times |\,\lambda_i\, {\cal Y} + \xi_{\,i}|^{\,\ell } \right)\right\}\, \times\\[0.1in]
& \ &  \ \ \ \ \ \ \ \ \ \   \times \left[\ O \,(1) \, |\,{\cal V}_{\,i} \,-\,A_1\,| \times   \max \, \left\{ \ {\cal V}_{\,i}^{{4}\over {n \,-\,2}}\,, \ \ A_1^{{4}\over {n \,-\,2}}   \right\}\, \right]   \\[0.15in]
& \ & \!\!\!\!\!\!\!\!\!\!\!\!\!\!\!\!\!\! \!\!\!\!\!\!\!\!\!\!\!\!\!\!\!\!\!\!\!\!\le \ \, C\,\Lambda_i^{-1} \cdot  |\,\lambda_i\,{\cal Y} + \xi_{\,i}|^{\,\ell}  \cdot \left[\ A_1^{ {{4}\over {n \,-\,2}}   } \, |\,{\cal V}_{\,i} \,-\,A_1\,|  \,\right] \ \ \ \ \  \ \ \ \ \  \ \ [\,{\mbox{cf.}} \ \ (A.1.3)\,]\\[0.15in]
   & \ &   \!\!\!\!\!\!\!\!\!\!\!\!\!\!\!\!\!\! \!\!\!\!\!\!\!\!\!\!\!\!\!\!\!\!\!\! \mfor |\,{\cal Y}| \ \,\le\, \ c\cdot \lambda_i^{-1} \ \ \ \ \ \ \ [\,{\mbox{using}} \   \ {\mbox{Propositions \ 2.3 \   and \   2.24}}  \,] \ \ \ \ \ \ \ \ \ \ \ \ \ \ \\[0.1in]
& \ & \!\!\!\!\!\!\!\!\!\!\!\!\!\!\!\!\!\! \!\!\!\!\!\!\!\!\!\!\!\!\!\!\!\!\!\!\!\!\le \ \,C_1\, \Lambda_i^{-1} \,\lambda_i^{\,\ell} \cdot (\,1 + {\cal R}\,)^{\ell}  \left[ A_1^{ {{4}\over {n \,-\,2}}   } \cdot \! \left( |\, {\cal V}_{\,i} \,-\,A_1 \,-\,\lambda_i^{\,\ell}\ \Pi_{\bf p}\,| + \lambda_i^{\,\ell} \, |\,\Pi_{\bf p}\,|\,   \right)   \right]\\
& \ & \!\!\!\!\!\!\!\!\!\!\!\!\!\!\!\!\!\! \!\!\!\!\!\!\!\!\!\!\!\!\!\!\!\!\!\!\!\!\le \ \,C_2\,  \lambda_i^{\,\ell} \cdot (\,1 + {\cal R}\,)^{\,\ell}\, \times  {1\over {(1 + {\cal R})^{\,4} }} \times  |\, {\bf W}_i \, + \, O_{\lambda_i} ((n \,-\,2) \,-\,\ell)|\ + \ \\[0.1in]
 & \ &   \ \ \ \ \ \   +\  C_3\,  (\Lambda_i^{-1}\, \lambda_i^{\,\ell}\,) \cdot (\,1 + {\cal R}\,)^{\ell}\, \cdot  {1\over {(1 + {\cal R})^{\,4} }} \cdot     {{\lambda_i^{\,\ell} \, (1 + {\cal R})^{ \ell }}\over {(1 + {\cal R})^n}} \\[0.1in]
& \ & \!\!\!\!\!\!\!\!\!\!\!\!\!\!\!\!\!\! \!\!\!\!\!\!\!\!\!\!\!\!\!\!\!\!\!\!\!\!\le \ \,{{ C_4\, \lambda_i^\ell \cdot (1 + {\cal R})^\ell }\over { (1 + {\cal R})^4}} \cdot \left[ \ 1 + {1\over {(1 + {\cal R})^{n \,-\,\ell} }}\right] \ \   \longrightarrow \  0 \ \ \   {\mbox{uniformly \ \ in}} \ \ B_o \, (R_o)\\[0.15in] \ \ \ \ \ \ \ \ \ \ \ \
& \ & \!\!\!\!\!\!\!\!\!\!\!\!\!\!\!\!\!\! \!\!\!\!\!\!\!\!\!\!\!\!\!\!\!\!\!\!\!\!\le \ \, {{C_5}\over {  (1 + {\cal R})^{\,4}  }} \ \ \  \mfor \ \  {\cal R}  \,= \, |\,{\cal Y}| \,\le \,c\,\lambda_i^{-1} \   \ \ \ \ ({\mbox{when}} \ \ \ell \,\le\, n \,-\,2)\,. \ \ \ \ \ \ \ \  \ \ \  \ \ \
\end{eqnarray*}

\vspace*{0.2in}

 {\bf \S\,6\,c.\,5.}   {\it  The inserted harmonic term.}\ \ The last couple of terms to be considered in $\,\Lambda_i^{-1} \cdot {\bf R.H.S.}_i\,$ [\,cf. (6.12), (6.19) and (6.20)\,]\, are
$$
\   {1\over {\Lambda_i}} \,\cdot \bigg\vert \  \lambda_i^{\,n \,-\,2} \cdot \,{{h_o}\over {c\,\lambda_i^{-1} }} \cdot  \, \Delta_{\cal Y}\, {\tilde{\cal R}}\,({\cal Y}) \,  \bigg\vert  \ \,\le\, \ C\, \lambda_i  \,\cdot\, \chi_{B_o\, (1)} \ \ \ \  \ \  [\,{\mbox{via}} \ \ (6.16) \ \ {\mbox{and}} \ \ \ell \,\le \, n \,-\,2\,]\ ,  \leqno (6.50)
$$
(6.51)
\begin{eqnarray*}
  & \ &  {{ \lambda_i^{\,n \,-\,2} }\over {\Lambda_i}}  \cdot n\,(n + 2)  \, {{h_o}\over {c\,\lambda_i^{-1} }} \cdot   [\,A_1\, ({\cal Y})]^{\, {4\over {n \,-\,2}} } \cdot {\tilde{\cal R}}\,({\cal Y}) \\ & \,\le\, & C\, \lambda_i  \cdot \chi_{B_o\, (1)} \ + \ C_1 \cdot \lambda_i  \cdot  {{\cal R}\over {(1 + {\cal R})^4}} \ \ \longrightarrow \ 0 \ \ {\mbox{uniformly \ \ in}} \ \ B_o\, (R_o)\\
 & \,\le\, & C\, \lambda_i  \cdot \chi_{B_o\, (1)} \ + \ C_2  \cdot  {1\over {(1 + {\cal R})^4}} \ \ \ \ \ \ \ \mfor {\cal R} = |\,{\cal Y}| \,\le\, c\,\lambda_i^{-1}\ , \ \ \ \ \ \ \ \ \ \  \ \ \ \ \ \ \ \ \ \ \
\end{eqnarray*}
(6.52)
\begin{eqnarray*}
& \ & {{\lambda_i^{\,n \,-\,2} }\over {\Lambda_i}} \cdot  n\,(n + 2) \,[\,A_1\,({\cal Y})]^{\, {4\over {n \,-\,2}} }\cdot |\, {\mbox{H}}_{\ge 1}\,({\cal Y}) \ \,-\,\ h_o \, |\\[0.15in]
&    \,\le\, & C \cdot \left( {1\over {1 + {\cal R}^2}} \right)^4 \ \times
\sum_{j \ge 1} \   \bigg\vert  \, \left[ \ {{{1}}\over {|\,  \lambda_i\, {\cal Y}  \,-\, {\hat Y}_j \,|^{\, n \,-\,2} }} \ \,-\,\ {{1}\over {|\,   {\hat Y}_j \,|^{\, n \,-\,2} }} \right] \, \bigg\vert \times {\cal A}_{\,j}   \\[0.15in]
 & \ &  \ \ \    \longrightarrow \  0 \ \ {\mbox{uniformly \ \ for }} \ \  {\cal Y} \,\in\, B_o\,(R_o) \ \ \ \\[0.1in]
    & \ &  \hspace*{-0.3in}[\,{\mbox{recall \   (5.13) \  and  \ (6.5)\,,\, \  observe \   also \   that }} \  \lambda_i\  {\cal Y} \,\to\, 0 \!\!\!\!\!\!\!\!\!\!\!\!\mfor \!\!\!\!|\,{\cal Y}| \ \,\le\, \ R_o\,] \ \ \ \ \ \ \ \ \ \ \\
 & \,\le\, & {C\over {(1 + {\cal R})^4}}  \ \ \ \ \ \ \ \  \ \ \ \ \ \ \ \   \ \ \ \ \  \ \ \ \ \ \ \ \ \ \ \ \ \  \ \ \ \ \  \mfor {\cal R} \,\le\, c\, \lambda_i^{-1}\ .
\end{eqnarray*}
Thus we   estimate  each term in the remainder and confirm the right orders in (6.19) and (6.20).
Combining the  above  discussion, we obtain the following.\\[0.2in]
{\bf Theorem 6.53.} \ \ {\it Under the conditions in Main Theorem}\, 1.14\,,\, {\it we have}\,
$$
|\, {\cal D}^\Pi_i\, ({\cal Y})| \ = \ o_{\lambda_i} (\ell) \  \ \ {\it{for}} \ \ |\,{\cal Y}| \ \,\le\, \ c\,\lambda_i^{-1} \     \ {\it{and}} \ \ i \,\gg\,  1 \ \  ({\it modulo \ a \ subsequence}\,) \,. \leqno (6.54)
$$

\vspace*{0.2in}

 {\bf \S\,6\,d.} \     {\bf Proof of Main Theorem 1.14}\,  --  \,{\bf zooming out to the original  scale}\,.\\[0.1in]
   As in the transformation  (2.17) (see also \,{\bf \S\,2\,f}\,), we note that,
from the definitions of $\,{\cal V}_i\,,\, \ A_1$\,,\,  and (6.2)\,,\, estimate (6.54) in Theorem 6.53 leads  to

\vspace*{-0.2in}

\begin{eqnarray*}
(6.55)  & \ & \!\!\!\bigg\vert\  {{v_i \,(\,\xi_{\,i} + \lambda_i \,{\cal Y})}\over {v_i\,(\xi_i)}}\  -\  \left( {1\over {1 + |\,{\cal Y}|^{\,2}}}\right)^{\!\!{{n \,-\,2}\over 2}} \ \,-\,\ \lambda_i^{\,\ell} \cdot {{ \Gamma_{\bf p} \,({\cal Y})}\over {(1 + |\,{\cal Y}|^{\,2})^{n\over 2} }} \  \,-\, \\[0.15in]
 & \ &  \ \ \ \ \  \,-\,\ \lambda_i^{\,n \,-\,2}\cdot {\mbox H}_{\,\ge\, 1}\, ({\cal Y})\   + \  \lambda_i^{\,n \,-\,2} \,\cdot\, h_o \,\cdot\, \left( 1 \,-\,{{\tilde {\cal R}}\over {c\,\lambda_i^{-1} }}\right)  \,\bigg\vert \
  =   \ o\, (\lambda_i^{\ell}\,) \ \ \  \ \ \ \   \ \ \ \
  \end{eqnarray*}

  for $\,|\,{\cal Y}|\,\le \,c\, \lambda_i^{-1}\,.\,$ Via the transformation
  $
  \,y = \lambda_i\, {\cal Y} \, + \, \xi_i\,
  $ and the definition $\,M_i \,:=\, v_i\,(\xi_i)\,$,\,
 (6.55) is rewritten a\\[0.1in]
 (6.56)
  \begin{eqnarray*}
 & \ & \!\!\!\!\!\!\!\!\!\!
  \bigg\vert\  v_i \,(y) \   -\  M_i \cdot \left( {1\over {1 + \lambda_i^{-2}\,|\,y \,-\, \xi_{\,i}|^{\,2}}}\right)^{\!\!{{n \,-\,2}\over 2}} \ \,-\,\ M_i \cdot \lambda_i^{\,\ell} \cdot {{ \Gamma_{\bf p} \ (\lambda_i^{-1}\,(y \,-\,\xi_{\,i}))}\over {[\ 1 + \lambda_i^{-2}\,|\,y \,-\, \xi_{\,i}|^{\,2}\,]^{n\over 2} }} \ \,-\, \ \ \  \\[0.15in]
 & \ &     \!\!\!\!\!\!\!\!\!\!\!\!\!\!-\,\ M_i \cdot \lambda_i^{\,n \,-\,2}\cdot {\mbox H}_{\ge 1}\, (\lambda_i^{-1}\,(y \,-\,\xi_{\,i}))\  + \ M_i \cdot \lambda_i^{\,n \,-\,2} \,\cdot\, h_o \,\cdot  \left( 1 \,-\,{{
 {\tilde {\cal R}} \, ( {\lambda_i^{-1}\,(y \,-\,\xi_{\,i})) } }\over {c\,\lambda_i^{-1} }}\right)  \,\Bigg\vert \\[0.15in]
  & \ & \!\!\!\!\!\!\!\!\!\!\!\!\!\!\!= \ M_i \cdot o\, (\lambda_i^{\ell}\,)
  \end{eqnarray*}
  for $\, |\,y| \ = \ |\,\lambda_i \,  {\cal Y} + \xi_{\,i} | \ \,\le\, \ c \,-\,o\, (1) \ \  (i \,\gg\,  1)\,.$
  Recall that $M_i \ = \ \lambda_i^{-\,{{n \,-\,2}\over 2}}\ \,$ [\,see (1.10), (1.11), (2.17) \,and\, \,{\bf \S 2\,f}\,]\,,\, and also the form of $\,{\mbox H}_{\ge 1}\,$ as expressed in (5.13)\,.\, Hence we come to the conclusion that\\[0.1in]
  (6.57)
   \begin{eqnarray*}
& \ &  \!\!\!\!\!\!\!\!\!\!\!\!\!\Bigg\vert\  v_i \,(y)\,-   \left( {{\lambda_i}\over { \lambda_i^{\,2} \ + \ |\,y \,-\,\xi_{\,i}|^{\,2}}}\right)^{\!\!{{n \,-\,2}\over 2}} \!\!\!  \,-\,  \left[\, \lambda^{\,\ell \,+ 1}_i \cdot \Gamma_{\bf p} \!\left({{y \,-\,\xi_i}\over { \lambda_i }}\right) \right]  \!\cdot\! \left( {{\lambda_i}\over {\lambda_i^2 + |\,y \,-\,\xi_{\,i}|^{\,2} }}\right)^{\!{n \over 2}} \!- \\[0.15in]
 & \ &  \!\!\!\!\!\!    \,-\, \ \, \lambda_i^{\,{{n \,-\,2}\over 2}} \left[ \ \sum_{j\,\ge\, 1} \left( {{{\cal A}_{\,j}}\over { |\ (y \,-\,\xi_i) \,-\,{\hat Y}_{\,j}|^{\, n \,-\,2} }} \   \,-\,{{{\cal A}_{\,j}}\over { |\   {\hat Y}_{\,j}|^{\, n \,-\,2} }}   \right)\ + \ \lambda_i \cdot {{ h_o \cdot {\tilde {\cal R}} \, ({\cal Y}) } \over c}    \right] \Bigg\vert \ \ \ \\[0.2in]
   & \ & \!\!\!\!\!\!\!\!\!\!\!\!\!\!\!= \  o_{\lambda_i} \! \left(\ell  \ \,-\,\ {{n-2}\over 2}\ \right) \ \ \ \ \ \ \ \ \ \   \ \ \    \mfor  |\,y\,| \, \le\, \rho_2 \ \ \ \ \ \  \ \ \ \ \ \ \ ( i \,\gg\, 1)\,. \ \ \ \  \ \ \ \  \ \ \ \
\end{eqnarray*}
Here $\,\rho_2 > 0\,$ is a number slightly less than $c$\,.\, With  \,(1.4), and  \\[0.05in]
(6.58)

\vspace*{-0.15in}

$$
O_{\mbox{H}}\left( \lambda_i^{{n \,-\,2}\over 2} \right) \, := \  \lambda_i^{\,{{n \,-\,2}\over 2}} \left[ \ \sum_{j\,\ge\, 1} \left( {{{\cal A}_{\,j}}\over { |\ (y \,-\,\xi_i) \,-\,{\hat Y}_{\,j}|^{\, n \,-\,2} }} \   \,-\,{{{\cal A}_{\,j}}\over { |\   {\hat Y}_{\,j}|^{\, n \,-\,2} }}   \!\right)\, + \, \lambda_i  \cdot {{h_o \cdot {\tilde {\cal R}} \, ({\cal Y}) } \over c}    \right] \,,
$$
we arrive at (1.20)\,.\, [\,As usual, \,${\cal Y} \,=\,\lambda_i^{-1} \,(\,y \,-\,\xi_i)$\,.\,]

\vspace*{0.5in}

 {\bf e\,-Appendix}\,  is available at page 44 onward.

\vspace*{0.15in}

{\large D}EPARTMENT OF {\large M}ATHEMATICS\,, \,\\[0.1in]
{\large N}ATIONAL {\large U}NIVERSITY OF {\large S}INGAPORE\,,\,\\[0.1in]
{\large 10,  \,L}OWER {\large K}ENT {\large R}IDGE  {\large R}OAD\,,
\\[0.1in]
{\large S}INGAPRE {\large 119076}\,,
\\[0.1in]
{\large R}EPUBLIC OF {\large S}INGAPORE \ \ \ \
\\[0.2in]
{\tt matlmc@math.nus.edu.sg}\\[0.1in]

\newpage

\centerline{\bf \LARGE {\bf e\,-\,Appendix to the Article  }}

\medskip \smallskip

\begin{center} {\Large   {\bf ``\,Refined Estimates for Simple Blow\,-\,ups  }}
\medskip \medskip \smallskip \smallskip \\ {\Large {\bf  of the  Prescribed Scalar Curvature on $S^n$\,''}}.

\vspace{0.43in}

{\Large {Man Chun  {\LARGE L}}}{EUNG}\\[0.165in]

{\large {National University of Singapore}} \\[0.1in]
{\tt matlmc@nus.edu.sg}\\
\end{center}
\vspace{0.33in}


\vspace*{0.2in}

{\it In this appendix we follow the notations, conventions, equation numbers, section numbers, lemma, proposition and theorem numbers as  used in the main article\,}  \cite{Refine}\,,\,  {\it unless otherwise is specifically mentioned} (\,{\it  for instances, those equation numbers starting with}\, `A'\,).  \\[0.4in]
{\large{\bf \S\,A.1} \ \ {\bf  A proof of Proposition}\, 2.24.}\\[0.2in]
{\it Proof of the necessary part\,}  ($\Longrightarrow$)\,.\, Here we can take $\,\zeta_i \,=\, \xi_{m_i}\,$ and
$
\epsilon_i \,  \,=\, \,  \lambda_{m_i} \,,
$
according to the analytic definition of simple blow\,-\,up point as in (2.1) and (2.2)\,.\,
With the notations in  (2.20), (2.26) is equivalent to
$$
C^{-1} \cdot A_1 \,({\cal Y})\ \le \ {\cal V}_{\,i}\,({\cal Y}) \ \le \ C\, A_1\,({\cal Y}) \ \ \ \ \ \ \mfor \ \ \  \,|\,{\cal Y}| \ \le \ \rho_o \,\lambda_{m_i}^{-1}\,. \leqno (A.1.1)
$$
[\,C\,f. (2.17) and (2.19)\,.\,]\, It follows from  (1.4) and $\,A_1 \,=\, A_{1\,, \ 0}\,$ that
$$
{1\over {2^{{n \,-\, 2}\over 2} }} \cdot {1\over {|\,{\cal Y}|^{\, n \,-\, 2} }} \ \le \  A_1\, ({\cal Y}) \ \le \ {1\over {|\,{\cal Y}|^{\, n \,-\, 2} }} \ \ \mfor |\,{\cal Y}|\ge 1 \ \ \ [\,{\mbox{i.e.}} \ \ {\cal Y} \not\in B_o \,(1)]\,.\leqno (A.1.2)
$$
 When $i \gg 1\,,\,$ $\,|\,\lambda_{m_i}\,{\cal Y} | \ \le \ \rho_o \ \Longrightarrow \ |\,\lambda_{m_i}\,{\cal Y} \,+\, \xi_{\,i}| \ \le \ {\bar \rho}_o\ $  for $\,0\, < \,\rho_o\, <\, {\bar \rho}_o\,,\,$
From  (2.19) and   Proportionality Proposition 2.3 we obtain\\[0.1in]
(A.1.3)

\vspace*{-0.25in}

\begin{eqnarray*}
{\cal V}_{\,i} \,({\cal Y}) & \,=\, & {{v_i \,(\lambda_{m_i}\,{\cal Y} \,+\, \xi_{\,i})}\over{M_i}} \
 \le  \ {{1}\over{M_i^{\,2}}}\cdot {{C_1}\over{|\,(\lambda_{m_i}\,{\cal Y} \,+\, \xi_{\,i}) \,-\, \xi_{\,i} |^{\,n \,-\, 2}  }} \\[0.1in]
& \le& {{1}\over{M_i^{\,2}}}\cdot \!{{C_1}\over{\lambda_{m_i}^{\,n \,-\, 2} \, |\,{\cal Y}  |^{\,n \,-\, 2}  }} \  \le \  {{C_1}\over{|\,{\cal Y}|^{\,n \,-\, 2} }} \ \ \ \ \ \mfor 0 < |\,\lambda_{m_i}\,{\cal Y} |\, \le \,\rho_o \ \ \ {\mbox{and}} \ \ i \gg 1\,.
\end{eqnarray*}
Here we use (2.4)\,.\,
As we already know that ${\cal V}_{\,i} \,({\cal Y}) \to A_1\,({\cal Y})$ \, uniformly for\, ${\cal Y}\, \in \,B_o\, (1)\,,\,$ together with the above estimate and (A.1.2), we have
$$
{\cal V}_{\,i} \,({\cal Y})\, \le\,  C  \,A_1\,({\cal Y}) \ \ \mfor \ i \, \gg\,  1 \ \ \ {\mbox{and}} \ \ \ |\,{\cal Y}| \,\le\, \rho_o \,\lambda_{m_i}^{-1}\ . \leqno (A.1.4)
$$
As for the lower bound in (A.1.1), from  Proportionality Proposition 2.3, we have
$$
M_i \cdot v_i \, (y) \, \ge\,  {a\over {|\,y|^{n \,-\, 2}}} \ +\ h \, (y) \ \,-\, \ o\, (1) \mfor 0 \, <\,  \rho_1 \, \le\,  |\, y| \, \le\,  \rho_o\,.\leqno (A.1.5)
$$
Here $\,o\,(1) \, \to\,  0\,$ when $\,i \, \to\,  \infty\,.\,$
Choosing $\, \rho_o\, $ to be small enough (correspondingly adjusting \,$\rho_1 \, <\,  \rho_o$\,)\,,\, we  obtain
$$
M_i \cdot v_i \, (y) \ \ge \ {{2^{-1}\cdot a}\over {|\,y- \xi_{m_i} |^{\,n \,-\, 2}}}  \ \  \mfor 0 < \rho_1 \le |\, y- \xi_{m_i}  | \le \rho_o \ \ \ {\mbox{and}} \ \ i \gg 1\,.\leqno (A.1.6)
$$
Repeat the argument in (A.1.3), we have
$$
{\cal V}_{\,i} \,({\cal Y}) \, \ge\,  \  {{C_2}\over{|\,{\cal Y}|^{\,n \,-\, 2} }} \ \ \ \ \mfor \ \   \rho_1 \,\lambda_{m_i}^{-1}  \, \le\,  |\,{\cal Y}| \, \le\,  \rho_o \,\lambda_{m_i}^{-1} \ \ \ \ {\mbox{and}} \ \ \ i \, \gg \, 1\,. \leqno (A.1.7)
$$
Again using the uniform convergence  ${\cal V}_{\,i} \,({\cal Y})\,  \to \, A_1\,({\cal Y})$ \,for\, ${\cal Y}\,\in\,B_o\, (1)\,,\,$ and (A.1.2), we obtain
$$
{\cal V}_{\,i} \,({\cal Y})\ \ge \  {{C_3}\over{|\,{\cal Y}|^{\,n \,-\, 2} }} \ \ \ \mfor |\, {\cal Y} | \,=\, 1\,.\leqno (A.1.8)
$$
Let us pay attention to the region $\,1 \,\le\,  |\,{\cal Y}| \,\le\,   \rho_o \,\lambda_{m_i}^{-1}\ $ and
consider the function
$$
\left( {\cal V}_{\,i} \,({\cal Y}) \ \,-\, \ {{C_4}\over{|\,{\cal Y}|^{\,n \,-\, 2} }} \right)\!. \leqno (A.1.9)
$$
Choose $C_4 \,=\, \min \, \{ \, C_2\,, \ C_3\}\,.\,$
It follows from equation (2.18) that
$$
\Delta_o   \!\left( {\cal V}_{\,i} \,({\cal Y}) \,-\, {{C_4}\over{|\,{\cal Y}|^{\,n \,-\, 2} }} \right) \ \,=\, \ \Delta_o\,{\cal V}_{\,i} \,({\cal Y}) \ < \ 0 \ \ \ \ \mfor  1 \, \le\,   |\,{\cal Y}| \, \le \, \rho_o \,\lambda_{m_i}^{-1}. \leqno (A.1.10)
$$
For the boundary, we apply  (A.1.7) and (A.1.8)\,,\, then we use the maximum principle to obtain
$$
{\cal V}_{\,i} \,({\cal Y}) \ \ge \ {{C_4}\over{|\,{\cal Y}|^{\,n \,-\, 2} }}   \ \ \ \ \mfor  \ \ 1\,  \le\,   |\,{\cal Y}| \, \le\,  \rho_o \,\lambda_{m_i}^{-1}.\leqno (A.1.11)
$$
As before, we already know that ${\cal V}_{\,i} \,({\cal Y}) \, \, \to \, A_1\,({\cal Y})$\, for \,${\cal Y} \,\in\,B_o\, (1)\,.\,$ Combining with the above estimate and (A.1.2), we obtain
$$
{\cal V}_{\,i} \,({\cal Y}) \ \ge \ C^{-1} \cdot  \,A_1\,({\cal Y}) \ \ \ \  \mfor i\,  \gg\,  1 \ \ \ {\mbox{and}} \ \ \ |\,{\cal Y}| \, \le\,  \rho_o \,\lambda_{m_i}^{-1}
$$
once we choose $C$ to be small enough.\,
This complete the proof of ($\,\Longrightarrow$)\,.\,\\[0.2in]
{\it Proof of the sufficient part\,}  ($\Longleftarrow$)\,.
\ \  Assuming that we have (2.26), that is,
$$
 {1\over C}\cdot \left( {{\epsilon_i}\over {\epsilon_i^2 \,+\, |\,y \,-\, \zeta_i|^2 }}\right)^{\!{{n \,-\, 2}\over 2}} \ \le \ v_i\, (y) \ \le \ C\cdot \left( {{\epsilon_i}\over {\epsilon_i^2 \,+\, |\,y \,-\, \zeta_i|^2 }}\right)^{\!{{n \,-\, 2}\over 2}} \leqno (A.1.12)
$$
for   $\,|\,y \,-\, \zeta_i|  \ \le \ \rho_o\,.\,$  It becomes clear that for $i \gg 1$\,,\, there is a point $\,\xi_i\,$ such that
$$
v_i\, (\xi_i) \,=\, \max \ \left\{ \ v_i\, (y) \ | \ \ y \,\in\  \overline{ B_o\, (\,\rho_o)}   \ \right\} \ \ \ \  {\mbox{and}} \ \ \ |\,\zeta_i \,-\, \xi_i|\ \le\  B\cdot\epsilon_i\,. \leqno (A.1.13)
$$
Here $\,B\,$ is a fixed positive number [\,one can take \,$ B^2 \,=\, C^{\,{4\over {n \,-\, 2}}}\, \,-\, \,1\,$,\, where \,$C\,$ is the constant in (A.1.12)]\,.\, Moreover,
$$
\lambda_i \ := \ {1\over {[\,v_i\, (\xi_i)]^{2\over {n \,-\, 2}} }} \ \ \Longrightarrow \ \ c^{-1}\cdot  \epsilon_i \ \le \ \lambda_i \ \le \  c \cdot \epsilon_i \,\mfor \,i \gg 1\,. \leqno (A.1.14)
$$
Here $\,c \,\ge\, 1\,$ is a constant.\,
In addition, via the triangle inequality $|\,y \,-\, \xi_i| \,\le\, |\,y \,-\, \xi_i|\, +\, |\, \xi_i \,-\, \xi_i|$\,,\, and vice versa, we have
$$
\ \ \  (A.1.13) \ \ {\mbox{and}} \ \ (A.1.14)  \ \ \Longrightarrow \ \ {1\over D} \ \le \ {{ \epsilon_i^2 \,+\, |\, y \,-\, \zeta_i|^{\,2}}\over { \lambda_i^2 \,+\, |\, y \,-\, \xi_i|^{\,2} }} \ \le\  D \mfor y \in \R^n\,.\leqno (A.1.15)
$$
Here we can take the constant $\,D \,=\, c^2\, (B^2 \,+\, 1) \,=\, c^2\cdot  C^{4\over {n \,-\, 2}} \,.\,$
Thus
$$
    {1\over {C'}}\cdot  {\bf A}_{\lambda_i\,,\  \xi_i} \, (y)  \ \le \ v_i\, (y) \ \le \ C'\cdot {\bf A}_{\lambda_i\,,\  \xi_i} \, (y)  \mfor  |\,y \,-\, \xi_i|  \ \le \ \rho_1\,, \leqno (A.1.16)
$$
where $\,\rho_1 > 0\,$ is slightly less than $\,\rho_o\,.\,$
Thus, without loss of generality, we may assume that
$$\zeta_i \ \,=\, \ \xi_i \ \ \  \ {\mbox{and}} \ \ \  \epsilon_i \ \,=\, \ \lambda_i \mfor i \gg 1\,\,.\leqno (A.1.17)$$
Directly,

\vspace*{-0.3in}

$$
 \left( {{\lambda_i}\over {\lambda_i^2 \,+\, |\,y \,-\, \xi_i|^2 }}\right)^{\!{{n \,-\, 2}\over 2}} \ \le \ {1\over {|\,y \,-\, \xi_i|^{{n \,-\, 2}\over 2} }} \ \ \ \Longrightarrow \ \  \ \ v_i\, (y)\ \le \ {C\over {|\,y \,-\, \xi_i|^{{n \,-\, 2}\over 2} }} \leqno (A.1.18)
$$
for $\, 0 \,<\, |\,y \,-\, \xi_i|  \, \le \, \rho_1\,.\,$
That is, $0$ is an isolated blow\,-\,up point. \bk
Next, consider the rescaled average

\vspace*{-0.25in}

 $$
 r \ \,\longmapsto \ \, {\bar w}_i\, (r) \ := \ r^{{n \,-\, 2}\over 2} \cdot \left[ {{ \int_{ \partial B_{\,\xi_i} (r)} \ v_i\, d S }\over { \int_{ \partial B_{\,\xi_i} (r)} \ 1\, d S}}\right] \leqno (A.1.19)
$$
Via the change of variables $\,r \,=\, e^{-t}\,,\,$ and (A.1.16), we have
\begin{eqnarray*}
(A.1.20)   & \ &
{1\over {C'}} \cdot \left[ \ {1\over {e^{(\,t \,-\, t_i)}  \,+\, e^{-\,(\,t \,-\, t_i)} }} \right]^{\,{{n \,-\, 2}\over 2}}  \ \le \ {\bar w}_i\, (t) \ \le \  C' \cdot \left[ \ {1\over {e^{(\,t \,-\, t_i)}  \,+\, e^{-\,(\,t \,-\, t_i)} }} \right]^{\,{{n \,-\, 2}\over 2}}\ \ \ \ \ \ \  \\[0.15in]
& \ & \ \ \ \ \ \   \mfor \ \ r \,=\, e^{-\,t} \  \le \ \rho_o \ \ \Longleftrightarrow \ \  t \ \ge \ T_1\,, \ \ \ \ {\mbox{where}} \ \ T_1 \ := \ -\, \ln \, \rho_1\,.
\end{eqnarray*}
In the above
$$
\,\,e^{-\, t_i} \ \,=\, \ \lambda_i \ \ \Longleftrightarrow \ \   t_i \,=\, -\, \ln\,\lambda_i \ \ \mfor i \,=\, 1\,, \ 2\,, \cdot \cdot \cdot \ \ \    (t_i \to \infty \ \ {\mbox{as}} \ \ i \to \infty). \,\,\leqno (A.1.21)
$$
By performing a blow\,-\,up analysis as in  Theorem 4.2 in \cite{Asy-Grad-Est}\,  (cf. also \S\,7\,c in \cite{Leung-Supported}), and using (A.1.20), we obtain
$$
{\tilde w}_i\, (t) \ \,=\, \ {\bar w}_i\, (t \,+\, t_i) \ \to\  \left[ \ {1\over {e^{\,t} \,+\, e^{-\,t} }} \right]^{\,{{n \,-\, 2}\over 2}}\ \,=\, \  {1\over {2^{{n \,-\, 2}\over 2} }} \cdot \left({1\over{ \cosh \,t}} \right)^{\!{{n \,-\, 2}\over 2}} \ . \leqno (A.1.22)
$$
The convergence is in $C^2$-\,sense, uniform on any given  bounded interval in $\R\,.\,$
Directly,\\[0.1in]
(A.1.23)
$$
{d\over {dt}} \left({1\over{ \cosh \,t}} \right)^{\!{{n \,-\, 2}\over 2}}   \,=\, \,0 \ \ \ {\mbox{iff}} \ \ t \,=\, 0\,, \ \ \ \ {{d^2}\over {dt^2}}\!  \left({1\over{ \cosh \,t}} \right)^{\!{{n \,-\, 2}\over 2}}  \ \le \ -\,c^2 \ < \ 0  \mfor |\,t|\le \delta\,.
$$
Here $c\,$ is a constant which  depends on the small number $\delta > 0\,.\,$
It follows that $\,{\tilde w}_i\,$ has only one critical point  in $\,[\, \,-\, \,\delta\,, \ \delta\,]\,$ for all \,$i \gg 1\,.\,$ Likewise, the first statement in (A.1.23) shows that ${\tilde w}_i$ has no critical point in $[\,-\,T_2\,, \ T_2] \,\setminus\, [\, \,-\, \delta\,, \ \delta\,]\ $       for $i \gg 1\,.\,$ Here $T_2 \,\in  \,[\,-\,(t_i \,-\, T_1)\,, \ \infty)\,$ is a  (fixed) large positive number. Any potential critical point in  $[\,-\,(t_i \,-\, T_1)\,,\, \ \infty) \, \setminus \, [\, \,-\, T_2 \,, \ T_2\,]\,$ can be ruled out by using
(A.1.20), together with Lemma 5.1 in \cite{Chen-Lin-3} (cf. also Lemma 7.16 and Lemma 7.25 in \cite{Leung-Supported}\,,\, and the proof of Theorem 4.1 in \cite{Asy-Grad-Est})\,.\, Note that (A.1.20) implies   $$ \ \,
C^{-1} \cdot \left[\ {1\over{ \cosh \,t}} \right]^{{n \,-\, 2}\over 2}   \le \ {\tilde w}_i\, (t) \ \le \ C\cdot \left[\ {1\over{ \cosh \,t}} \right]^{{n \,-\, 2}\over 2} \mfor t \in [\,-\,(t_i \,-\, T_1)\,, \ \infty)\ . \leqno (A.1.24)
$$
Any ``\,small\," critical value for ${\tilde w}_i$ comes from a  local minimum\,,\, and according to (5.3) in \cite{Chen-Lin-3}  and Lemma 5.1 (loc. cit.)\,,\, ${\tilde w}_i$ has to increase in either direction, which eventually contradicts (A.1.24).
Via a translation back to $\,{\bar w}_i$\, as defined in (A.1.19),  it follows that\,,\, for each $ i \gg 1\,,\,$ ${\bar w}_i\, (r)\,$ has only one critical point (around $\,r_i := e^{-\,t_i}\,$)\, for $\,r \,\in\, (0\,, \ \rho_1)\,$.\, This completes the checking that $\,0\,$ is a simple blow\,-\,up point for $\,\{ v_i \}\,.\,${{\hfill {$\rlap{$\sqcap$}\sqcup$}}

\newpage

{\large{\bf \S\,A.2} \ \ {\bf  Shifting to the maximal point\,.}}\\[0.2in]
From \S\,A.1\,,\, in particular, (A.1.13) and  (A.1.17), we can take $\,\xi_{m_i} \,=\, \xi_i\,$ in (2.1) and (2.2)  in the definition of simple blow\,-\,up points\,. Suppose there exists another sequence of points $\,\{ {\tilde{\xi}}_i \}\,$ which also satisfies (2.1) and (2.2)\,.\, We show that (modulo a subsequence)
$$
|\,{\tilde{\xi}}_i \,-\, \xi_i| \,=\, o\, (\lambda_i) \leqno (A.2.1)
$$
The proof, which requires only standard argument, can be readily  recognized by people working on the area. For the benefit of general readers, we present   the argument, and refer to available papers for selected technical details.
Set
$$
{\tilde \lambda}_i \ \,=\, \ {1\over{ [\,v_i\, ({\tilde{\xi}}_i ) ]^{2\over {n \,-\, 2}} }} \ \ \ \ \Longrightarrow \ \ \ \ {\tilde \lambda}_i \ \ge \ \lambda_i  \ \ \left(  \ \,=\, \ {1\over{ [\,v_i\, ({{\xi}}_i ) ]^{2\over {n \,-\, 2}} }} \right). \leqno (A.2.2)
$$
To receive a contradiction, suppose
$$
 \ \ \ \ \ \   {\tilde \lambda}_i^{\,-1}   \cdot |\, {\tilde{\xi}}_i \,-\, \xi_i | \ \to \ \infty \ \ \ \ \ \ \ \ \  \ ({\mbox{modulo \ \ a \ \ subsequence}})\,. \leqno (A.2.3)
$$
As in the proof of Lemma 3.10 in \cite{Leung-Supported}, both $\,\{ {\tilde{\xi}}_i \}\,$ and $\,\{ {{\xi}}_i \}\,$ can be used as in (2.17)   to form (distinct) bubbling sequences (see \cite{Leung-Supported}), contradicting that the blow\,-\,up is isolated (Proposition 3.32 in \cite{Leung-Supported})\,.\, Hence there is a positive constant $B$ such that
$$
  {\tilde \lambda}_i^{-1}   \cdot |\, {\tilde{\xi}}_i \,-\, \xi_i | \ \le \ B \ \ \ \mfor \ \ i \,\gg\, 1\,.\leqno (A.2.4)
$$
As in \S\,2\,d\, (see also the proof of Lemma 3.10 in \cite{Leung-Supported}), \\[0.1in]
(A.2.5)$$
{\tilde {\cal V}}_i\, (Y) \ := \ {{ v_i\, ({\tilde \lambda}_i \,Y \,+\, {\tilde{\xi}}_i) }\over { v_i\, ( {\tilde{\xi}}_i)  }} \ \to \ \ A_1\, (Y) \ \ \ \ {\mbox{in \ \ }} C^1\,-\,{\mbox{sense\,, \ \ uniformly \ \ for }} \ \ Y \in B_o\, (R)\,.
$$
Once we take $R$ to be large enough, the point $Y_i$ defined by $$
{\tilde \lambda}_i \,Y_i \ \,+\, \ {\tilde{\xi}}_i \ \,=\, \ \xi_i \ \ \ \
{\mbox{satisfies}}  \ \ \ |\, Y_i| \ \le \ B \ < \ R/2\,.\,\leqno (A.2.6)
$$
Hence $\,{\tilde {\cal V}}_i\,$ has a critical point  at $\,Y \,=\, Y_i\,$ for\, $i \gg 1\,.\,$ In particular, $\,\btd \ {\tilde {\cal V}}_i \,(Y_i)\ \,=\, \ 0\,$ for $\,i\gg 1\,.\,$ On the other hand,
$$
\min \ \{\, |\, \btd \, A_1\, (Y)| \ : \ \ \delta \le |\, Y| \le R \ \} \ \ge \ c^2_{\delta,\ R}\  > \  0\,. \leqno (A.2.7)
$$
Here $\,c^2_{\delta,\ R}\,$ is a positive number depending on $\delta$ and $R$ ($\,c^2_{\delta,\ R}\, \to 0\,$ as $\delta \to 0^+$)\,.\,
The convergence in (A.2.5) together with   (A.2.6) and (A.2.7) shows that $|\, Y_i|\,\le \,\delta\,.\,$ Moreover, we can let $\,\delta \,\to\, 0\,$ as $\,i \,\to\, \infty\,.\,$ Hence
$$
|\,  {\tilde{\xi}}_i \,-\,  {{\xi}}_i| \ \,=\, \ o\, ({\tilde \lambda}_i )\,. \leqno (A.2.8)
$$
Using again the convergence in (A.2.5) and (A.2.8), we have
$$
 \lambda_i \ \le \ {\tilde \lambda}_i \ \le \  [\,1 \,+\, o\,(1)] \cdot \lambda_i \ \ \Longrightarrow \ \ |\,  {\tilde{\xi}}_i \,-\,  {{\xi}}_i| \ \,=\, \ o\, ({\lambda}_i )\,.
$$
In the above, $o\, (1) \to 0^{\,+}$ as $i \to \infty\,.\,$

\vspace*{0.5in}

{\large{\bf \S\,A.3} \ \ {\bf  Proof of Lemma}\, 4.11.}\\[0.2in]
The first conclusion in (i)   follows directly from  definition (4.10), and (ii) from  the limitation $j \ \le \  k\,.\,$ As for the second and third conclusions in (i), observe that  when $\,\ell\,$ is odd, $\,\Delta_o^{\!(h_\ell)}\, {\cal P}_\ell\,$ is a degree one polynomial\,.\, If it is not equivalent to zero, then $\,\Delta_o^{\!(h_\ell)}\, {\cal P}_\ell \ \,=\,  \ \sum_j\, c_j\, {\cal Y}_{|_j}   \, \in \, {\cal F}\,(\,{\cal P}_\ell\,)\,$ as claimed\,.\bk
Similarly, when $\,\ell\,$ is even, $\,\Delta_o^{\!(h_\ell)}\, {\cal P}_\ell\, \,=\, c \not=0\,$ is a number, then $c \,  \in \, {\cal F}\,(\,{\cal P}_\ell\,)\,,\,$ and so does $\,c\,{\cal R}^2   \, \in \, {\cal F}\,(\,{\cal P}_\ell\,)\,.\,$ Their difference is also in $\, {\cal F}\,(\,{\cal P}_\ell\,)\,.\,$\bk
 As for (iii), we first observe that, via direct calculation, we  have
 $$
   {\cal Y} \cdot \btd \, Q_l \,=\, l \cdot Q_{\,l}  \ \ \ \ {\mbox{  for \ \ any \ \ homogeneous \ \ polynomial  \ \ with \ \ degree \ \ }} l\,.\leqno (A.3.1)
$$
 For $\,j\, \ge\, 1$\,,\, using the product formula
 $$
 \Delta_o\,(f\cdot g) \ \,=\, \ f \cdot ( \Delta_o\,g) \ \,+\, \ 2\, \langle\,\btd \,f\,,\, \btd\, g\,\rangle \ \,+\, \ g \cdot ( \Delta_o\,f)\,,
 $$
 we obtain
\begin{eqnarray*}
 (A.3.2) \ \ \ \ \ \ \ \ \ \ \ \ \ \ \     \Delta_o\, [\,({\cal R}^{\,2})^j \,\Delta^k_o \, {\cal P}_\ell\,]
& \,=\,  & ({\cal R}^{\,2})^{\,j } \, \Delta^{k+1}_o \, {\cal P}_\ell \, \ \,+\,  \ \,A_{\ell\,,\, j\,,\, k} \cdot ({\cal R}^{\,2})^{\,j-1}\,\,\Delta^k_o \, {\cal P}_\ell\,, \ \ \ \ \ \\[0.15in]
(A.3.3) \ \ \ \  \ \ \ \    ({\cal R}^{\,2})\,\Delta_o\, [\,({\cal R}^{\,2})^{\,j} \,\Delta^k_o \, {\cal P}_\ell\,]
& \,=\,  & ({\cal R}^{\,2})^{\,j+1} \, \Delta^{k+1}_o \, {\cal P}_\ell \ \,+\, \  A_{\ell\,,\, j\,,\, k} \cdot ({\cal R}^{\,2})^{\,j}\,\,\Delta^k_o \, {\cal P}_\ell\,,\\[0.15in]
{\mbox{(A.3.4)}} \ \ \ \ \ \ \ \ \  \ \ \ \ \  \ \ \ \ \  \ \ \ \ \ \  \ \   \ \ \,  A_{\ell\,,\,j\,,\, k\,} & \,=\, & (2j\,) \cdot (2j +n-2 \,+\, 2\,\ell \,-\, 4 k)\,.\ \ \ \ \ \ \ \ \  \ \ \ \ \ \ \ \  \ \ \ \ \ \\
 \end{eqnarray*}

 \vspace*{-0.2in}

As $\,j \,\le\, k \ \Longrightarrow \ (j \,+\, 1) \le (k \,+\, 1)$, the terms which appear on the right hand side above belong to $\, {\cal F}\,(\,{\cal P}_\ell\,)\,.\,$ \qed

\newpage

{\large{\bf \S\,A.4} \ \ {\bf    The case when $\,\Delta_o^{\!(h_\ell)}\, {\cal P}_\ell \not\equiv 0$\,.\, }}\\[0.2in]
For  $\,\ell \,\le\, n \,-\,2$\,,\, where $\,n\,$ is even, we discuss how to eliminate the condition\,  $\,\Delta_o^{\!(h_\ell)}\, {\cal P}_\ell \,= \,0\,$\,  by adding higher order (up to $n$\,-\,th order) terms.
We start with
\begin{eqnarray*}
(A.4.1)
 \!\!\!  \!\!\! & \ & (1 \,+\, {\cal R}^{\,2})  \,\Delta_o\, ({\cal R}^\ell\,) \ \,-\, \ 2n \,[\,{\cal Y} \cdot \btd \, ({\cal R}^\ell\,) \,]  \ \,+\, \ 2n \, ({\cal R}^\ell\,)   \\[0.15in]
& \ & \!\!\!\!\!\!\!\!\!\!\!\!\!\!\!\!\!\!     \,=\, \!\left( [\, 1 \,+\, {\cal R}^{\,2}\,] \, \Delta_o \,-\,2n \,{\cal R} \!\cdot\! {{\partial}\over {\partial {\cal R}}} \,+\,  2\,n\right) {\cal R}^\ell\     \,=\,   (\ell \,-\, 2) \,(\ell \,-\, n)\,{\cal R}^\ell \,+\, \ell\, (\ell \,+\, n \,-\, 2) \, {\cal R}^{\, \ell \,-\, 2}\,.
\end{eqnarray*}
In particular, when $\,\ell \,=\, 2\,$ or $\,n\,$\,,\, we have
\begin{eqnarray*}
(A.4.2)\!\!\!\!\!\!\!\!\!& \ & (1 \,+\, {\cal R}^{\,2})  \,\Delta_o\, ({\cal R}^n\,) \,-\, 2n \,[\,{\cal Y} \cdot \btd \, ({\cal R}^n\,) \,]  \, \,+\,  2n \, ({\cal R}^n\,)   \,=\, 0 \!\cdot\! {\cal R}^n \,+\, [\,2n\,(n \,-\, 1)]\,{\cal R}^{n \,-\, 2},\\[0.1in]
(A.4.3)\!\!\!\!\!\!\!\!\!& \ & (1 \,+\, {\cal R}^{\,2})  \,\Delta_o\, ({\cal R}^2\,) \,-\, 2n \,[\,{\cal Y} \cdot \btd \, ({\cal R}^2\,) \,]  \, \,+\,  2n \, ({\cal R}^2\,) \,=    0 \cdot {\cal R}^{\,2} \ \,+\, \   [\,2\,n]\,.
\end{eqnarray*}
Consider finding a radial function $F\, (r)$ so that
 \begin{eqnarray*}
(A.4.4) \!\!\!\!\!\!\!\!\!& \ & \ \ (1 \,+\, {\cal R}^{\,2})  \,\Delta_o\, F\, (r) \ \,-\, \ 2\,n\, [\,{\cal Y} \cdot \btd \, F\, (r) \,]  \ \,+\, \ 2n \, F\, (r) \\[0.1in]
& \,=\, & \!\!- \   [\, \Delta^{h_\ell}_o \, {\cal P}_\ell\,] \cdot \left\{ a_o \ \,+\, \ a_1\cdot ({\cal R}^2)  \ \,+\, \ \,+\, \cdot \cdot \cdot \ \,+\, \ a_{h_\ell -1} \cdot ({\cal R}^2)^{\,h_\ell -1} \ \,+\, \  \,a_{h_\ell} \cdot ({\cal R}^2)^{\,h_\ell}   \right\}. \ \ \ \ \
\end{eqnarray*}

\vspace*{0.1in}

({\bf i}) \ \ \,We start with using a $({\cal R}^{\,2})$ term to cancel the constant term. By (A.4.3) above, we won't introduce any new $({\cal R}^{\,2})$ term.\\[0.2in]
({\bf ii}) \ \
The  $({\cal R}^{\,2})$\,-\,term in the right hand side can be canceled by introducing an $({\cal R}^2)^2$\ -\,term [\,using (A.4.1), and \,$\ell\,(\ell \,+\, n \,-\, 2) \not= 0$\,]\,.\, By doing so, a new\, $({\cal R}^2)^2$\,-\,term is introduced to the right hand side.\\[0.2in]
({\bf iii}) \  The  combined $({\cal R}^2)^2$\ -\,term in the right hand side can be canceled by introducing an $({\cal R}^2)^{\,3} $\ -\,term\,.\, By doing so, an $({\cal R}^2)^{\,3} $\,-\,term  is introduced to the right hand side. The process goes on until we reach  the ${\cal R}^{\,n-2}$ term ($n\,\ge\, 4$ is even). Introducing a ${\cal R}^n$\ -\,term cancels the ${\cal R}^{\,n-2}$\,-\,term, and via (A.4.2), it {\it does not}\,  re\,-\,introduce itself to the right hand side (\,that is,\, ${\cal R}^n$ is not present)\,.

\newpage

\begin{eqnarray*}
({\cal R}^2)^{\,{n\over 2}} \ \ \rightarrow ||||  & \ &   \ \ \ \ \ \ \ \ \ \   \\
\searrow \ \ & \ & \hspace*{1.3in}  \\
({\cal R}^2)^{\,{{n-2}\over 2}}  \ \ \rightarrow \ \ & \ & \!\!({\cal R}^2)^{\,{{n-2}\over 2}}    \ \ \ \ \ \  \ \ \ \   \ \ \ \ \ \  \\
\!\!\searrow \ \ & \ & \    \ \ \ \ \ \ \ \ \ \    \\
& \ &   \cdot    \\
& \ &   \cdot    \\
& \ &   \cdot    \hspace*{1.2in}   \\
& \ &  \cdot    \hspace*{1.2in}  \\
({\cal R}^2)^2 \ \ \rightarrow \ \  & \ & \!\!({\cal R}^2)^2     \\
\!\!\searrow \ \ & \ & \     \hspace*{1.2in} \\[0.1in]
{\cal R}^{\,2} \rightarrow |||| \ \ & \ & \!\!{\cal R}^{\, 2}     \ \ \ \ \ \  \ \ \ \    \ \ \ \  \ \   \\[0.1in]
\searrow \ \   & \ &       \\[-0.1in]
  & \ & \!\!c_o     \
\end{eqnarray*}

\vspace*{0.2in}

\centerline{\underline{Diagram A.4.5}. \ \ The cancelation order from bottom upward (when \,$n$\, is even)\,.}

\vspace*{0.2in}

In summary, when $\,n \,\ge\, 4\,$\,,\, $\ell$\, even with $\ell \le n \,-\, 2\,,\,$ we can find a polynomial \\[0.05in]
(A.4.6)

\vspace*{-0.25in}

$$
F\, (r) := [\, \Delta^{h_\ell}_o \, {\cal P}_\ell\,] \cdot
 \left\{ B_1\cdot({\cal R}^2)^1 \,+\, B_2\cdot({\cal R}^2)^2 \,+\, \cdot \cdot \cdot \,+ B_k\cdot({\cal R}^2)^k \,+\, \cdot \cdot \,+\, \,B_{n\over 2}\cdot({\cal R}^2)^{\,{{n  }\over 2  }} \right\}\, ,
$$
which satisfies (A.4.4)\,,\, where

\vspace*{-0.3in}

\begin{eqnarray*}
 (A.4.7) \!\!\!\!\!\!\!\!& \ &
 B_1 \,=\, -{{a_o}\over {2n}}\,, \ \   B_2 \, \,=\, \, -\,{{a_1}\over {4\, (4 \,+\, n \,-\, 2)}}\,, \ \ \ B_3 \,=\, -\, {{a_2 \,+\, (4 \,-\, 2) \, (4 \,-\, n)\cdot B_2}\over  {6\, (6 \,+\, n \,-\, 2)}}\,, \\[0.15in]
  & \ & \!\!\!\!\!\!\!\!\!\!\!\!\!\!\!\!\!\!\cdot \cdot \cdot\,, \ \   B_k \,=\,-\, {{ a_{ k  \,-\, 1} \,+\, \left[ \ 2 \left(k \,-\,1 \right) \,-\, 2\, \right]\cdot \left[\ 2 \left( k \,-\,1 \right) \,-\, n \,\right] \cdot B_{k \,-\, 1}}\over {(2k)\, [\,(2k) \,+\, n \,-\, 2\,]}}  \ \   \  {\mbox{ for} } \ \  3 \,\le\, k \,\le\, {n\over 2}\ .
\end{eqnarray*}
In particular,

\vspace*{-0.3in}

$$
B_{n\over 2} \ \,=\, \ -\, {{ \ a_{ {n\over 2}  \,-\, 1} \ \,+\, \ (n \,-\, 4)\,(-2) \cdot B_{{n\over 2} \,-\, 1}\ }\over {n\, [\,2\, (n \,-\, 1)\,]}}\ . \leqno (A.4.8)
$$

\newpage

{\bf Proposition A.4.9.} \ \ {\it For $n \ge 4\,$ and $\,\ell \le n \,-\, 2\,,\,$ both being  even\,,\,  let $\,{\cal P}_{\ell}\,$ be a homogeneous polynomial  of degree $\,\ell\,.\,$ Then   equation\,} (4.3) {\it has a solution given by}
$$
 \sum_{0 \le j \,\le \,k \le \,h_\ell \,-\, 1} \!\!\!\!C^j_k \cdot ({\cal R}^2)^j \,[\, \Delta_o^{\!(k)} \, {\cal P}_\ell\,] \ \,+\, F\, (r) $$
 {\it where \,$F$ is given in } (A.4.6)\,,\, (A.4.7) {\it and}\, (A.4.8)\,.\\[0.2in]
 \hspace*{0.5in}Although Proposition A.4.9 allows us to find a solution of equation (4.3) without the condition $\,\Delta_o^{\!(h_\ell)}{\cal P}_\ell \,\equiv\, 0\,$,\, the presence of an order $n$ term (that is, ${\cal R}^n$) hinders the application of the second order blow\,-\,up argument, cf. \S\,6\,b.6\,.\\[0.2in]
{\it Remark on uniqueness}\,.\ \
Let $\Gamma_a$ and $\Gamma_b$ be two polynomial solutions to equation (4.3), with maximum degrees $\, < n\,.\,$ Via Theorem 4.16,  [\,mindful of condition (4.18) which requires the maximum degrees $\, < n\,$]\,,\, we have
$$
[\,\Gamma_a \,-\, \Gamma_b\,]\,({\cal Y}) \,=\, c_o\,(1 \,-\, {\cal R}^2) \ \,+\, \ \sum\, c_j\, {\cal Y}_{|_j} \ \ \ \ \ \ \ \ \ \ (\,{\cal R} \,=\, |\,{\cal Y}|\,)\,. \leqno (A.4.10)
$$
Thus when $\,\ell\, <\, n\,$, the solution found in Proposition 4.49 are ``unique" in the case of (A.4.10)\,.\, This cannot be extended to the solution found in Proposition A.4.9, even though $\,\ell \,<\, n\,$\,,\, as the presence of $\,{\cal R}^n$\,-\,term voids the application of Theorem 4.16.\,\\[0.075in] $\displaystyle{\left(\, {\mbox{For \ a \ generic}} \ \,  n\,,\, \ B_{n\over {\,2}} \not= 0\,.\,\right)}$

\newpage

{\large{\bf \S\,A.5} \ \ {\bf Bounds on the Green  function on $B_o\, (a)\,$ with a point}}}\\[0.075in]
\hspace*{0.62in}  {\large{\bf not too close to the boundary.}}\\[0.2in]
The Green's function for \,$\Delta_o$\, on $\,B_o\, (a)\,$ is given by
$$
G\, (y, \  \xi) \ =\ -\, {1\over {(n \,-\, 2)\, \Vert\, S^{\, n \,-\, 1} \Vert }} \cdot \!\left[ \ {1\over {|\, y \,-\, \xi|^{\,n \,-\, 2}}} \ \,-\, \  \left( {a\over {|\, \xi|}}\right)^{\, n \,-\, 2} \!\!\cdot {1\over {|\, y \,-\, \xi^*|^{\,n \,-\, 2}}}\right]\,, \leqno (A.5.1)
$$
where $\,\xi^*\,$ is the reflection of the point $\,\xi\,$ upon the sphere $\,\partial B_o\,(a)\,$,\, given by
$$
\xi^* \ \,=\, \ {{a^2}\over {|\,\xi|^2}}\cdot \xi \ \ \ \Longrightarrow \ \  \ |\,\xi^*| \ \,=\, \ {{a^2}\over {|\,\xi|}} \ \ \ \mfor\ \,  \xi \,\in\, B_o\, (a) \setminus \{ 0 \}\,.\leqno (A.5.2)
$$
See for example \cite{McOwen}\,.\, Here $\,\Vert\, S^{\, n \,-\, 1} \Vert\,$ denotes the volume/meassure of $\,S^{n \,-\, 1}\,$ with the standard metric\,.\, In order to obtain the bound in (6.40), we need only to consider the second term in (A.5.1). Let
$$
|\, \xi| \ \,=\, \ \gamma \, a \ \ \Longrightarrow \ \ |\,\xi^*| \ \,=\, \  {a\over \gamma} \ \left( \,> \ a \right)\,, \ \ \ \ {\mbox{where}} \ \ \gamma \ \le \  ( 1 \,-\, \delta)\,.
$$
It follows that, for $|\, \xi| \ \le \ (1 \,-\, \delta)\, a\,,\,$
\begin{eqnarray*}
(A.5.3) \ \  \ \left( {a\over {|\, \xi|}}\right)^{\, n \,-\, 2}\!\!\!\! \cdot {1\over {|\, y \,-\, \xi^*|^{\,n \,-\, 2}}} &  \le & \left( {1\over {\gamma}}\right)^{\, n \,-\, 2}\!\! \cdot {1\over {\big\vert \ a \,-\, {a\over \gamma} \,\big\vert^{\,n \,-\, 2}}}\  \,=\, \  \left( {1\over { 1- \gamma}}\right)^{\, n \,-\, 2}\!\!\!\! \cdot {1\over {a^{\,n \,-\, 2}}}\ \ \ \ \ \ \\[0.15in]
& \le & {1\over {\delta^{\, n \,-\, 2} }} \cdot {{2^{\,n \,-\, 2}}\over {|\, y \,-\, \xi |^{\, n \,-\, 2} }}\ \  \ \ \ \ ({\mbox{since}} \ \ |\  \xi \,-\, y|\le 2\,a)\,.
\end{eqnarray*}
We obtain
$$
|\,G\, (y, \ \xi)| \,\le\, \left[ \,C_1 \,+\, {{C_2}\over {\delta^{\,n \,-\, 2} }} \right] \cdot {1\over {|\,y \,-\, \xi|^{\, n \,-\, 2} }} \mfor |\,\xi| \le (1 \,-\, \delta) \cdot \!a \ \  {\mbox{and}} \ \ y \in B_o \, (a) \setminus \{ \xi \}.
$$
As for the  bound in (6.41), we note that the term
$$
{\bf n} \cdot \btd_y\, G_i \,(y, \ \xi)
$$
is indeed the Poisson kernel for $\,\Delta_o$\, on $B_o\, (a)\,,\,$ which is given by (e.g. \cite{McOwen}  pp. 116)
\begin{eqnarray*}
(A.5.4) \ \ \ & \ & {1\over {a \, \Vert\,S^{\, n \,-\, 1} \Vert}} \cdot {{ \,a^2 \ \,-\, \ |\,\xi|^{\,2}}\over { |\, y \,-\, \xi|^{\,n} }} \ \ \ \ \ \ \ \mfor y \in \partial B_o\,(a)\\[0.15in]
& \le & {{ a^2  }\over {a \, \Vert\, S^{\, n \,-\, 1} \Vert\cdot (\, \delta\, a \, )^{\,n} }} \ \ \ \ \ \ \ \ \ \  \ \ \ \mfor |\,\xi|\le (1 \,-\, \delta)\,a \ \ \ {\mbox{and}} \ \ |\,y| \ \,=\, \ a \ \ \ \ \ \ \ \ \ \ \ \ \  \ \ \ \\[0.15in]
& \le & {{C }\over {\delta^n}} \cdot {1\over {a^{\,n \,-\, 1}}}\ .
\end{eqnarray*}

\newpage

{\large{\bf \S\,A.6} \ \ {\bf  Balance and cancelation\,.}}\\[0.2in]
%
%
%
%
We first describe the classic balance formula for equation (1.2)\,,\, due to Stanislav I. Pohozaev. See for examples \cite{Li-1} and \cite{Leung-Supported}\,.\\[0.15in]
{\bf Theorem A.6.1.} \ \ {\it   Let $\,v\,$  and $\,K\,$ satisfy the general conditions in\,} (1.27) {\it and\,} (1.28)\,.\, {\it
We have the following.}\\[0.1in]
{\bf (I)\ }\ \ \,{\bf Global Pohozaev's identity.} \ \
$$
\int_{\R^n}\  \langle\,y\,, \    \btd \, K\,  (y)\,\rangle \, [\,v\,(y)\,]^{{2n}\over {n \,-\, 2}}  \  dy \   \,=\, \ 0\,. \leqno (A.6.2)
$$

\vspace*{0.1in}

({\bf II}) \ \ {\bf Mezzo\,-\,scale Pohozaev's identity.} \ \ {\it For a fixed number $\rho_o > 0$\,,\,  \,we have }
\begin{eqnarray*}
(A.6.3)&\,& \ \ \ \int_{B_{\,o} (\,\rho_o) } \ \langle\,y\,, \    \btd \, K\,  (y)\,\rangle \,[\,v\,(y)\,]^{{2n}\over {n \,-\, 2}}  \  dy   \ \,=\, \ {{1}\over {{\tilde c}_n}}\cdot {{2n}\over {n \,-\, 2}}  \,\int_{\partial B_{\,o} (\,\rho_o) }    \langle\!\vec {\ {\bf V}}\,, \ {\bf n} \,\rangle\,  dS\,,   \ \ \ \ \ \ \ \ \ \ \ \ \ \ \ \ \ \ \ \ \ \ \ \ \ \ \ \\[0.15in]
(A.6.4)&\,&\vec {\ {\bf V}}  (y) \ \,=\, \  {{n \,-\, 2}\over 2}\, v (y) \,\btd v (y)\,- \, {{|\btd v (y)|^{\,2}}\over
2} \ y  \ \,+\, \ \left[\,\langle\,y\,, \    \btd \, v\,  (y)\,\rangle \,
\right] \,\btd v\, (y)\\[0.1in]
 & \ & \ \ \ \ \ \ \ \ \ \ \ \ \ \ \  \ \ \ \ \ \ \ \ \ \ \ \ \ \ \ \ \ \ \ \ \ \ \ \ \ \ \ \ \ \ \ \ \ \,+\,  \ {{n \,-\, 2}\over {2n}} \cdot  {\tilde c}_n \cdot \left\{  \, [\,v\, (y)]^{{2n}\over {n \,-\, 2} }   \  K\, (y)\,\right\}\,y\,.  \ \ \ \
 \ \ \ \  \end{eqnarray*}
{\it In}\, (A.6.3) {\it and}\, (A.6.4)\,,\, $y$ {\it is treated as a vector\,,\, and\,  ${\bf n}$\, is the unit outward normal on $\,\partial B_{\,o} (\,\rho_o)\,.\,$ }

\vspace*{0.1in}

{\bf \S\,A.6.a.}\ \ {\it   Order of vanishing outside the blow\,-\,up points.} \ \ Observe that, via (A.1.2) and gradient estimate \cite{Asy-Grad-Est}, we obtain

\vspace*{-0.2in}

\begin{eqnarray*}
& \ & \!\!\!\max_{\partial B_{\,o}\, (\,\rho_o)} |\,\langle\,y\,, \    \btd \, v_i\,  (y)\,\rangle \,|\cdot v_i \ \,=\, \  O_{\lambda_i}(n \,-\, 2)\,,\ \ \ \ \ \ \ \ \
\max_{\partial B_{\,o}\, (\,\rho_o)} | \btd   v_i|^{\, 2} \   \,=\, \ O_{\lambda_i}(n \,-\, 2)\,,\\[0.1in]
& \ &  \!\!\!\max_{\partial B_{\,o}\, (\,\rho_o)} |  \btd \,v_i|\cdot v_i \ \ \,=\, \  O_{\lambda_i}(n \,-\, 2)\,,\ \ \ \ \     \
\max_{\partial B_{\,o}\, (\,\rho_o)}   [\,v\, (y)]^{{2n}\over {n \,-\, 2} }   \ |\, K\, (y)|\cdot |\, y|  \ \,=\, \ O_{\lambda_i}(2n)\ \ \ \ \ \\[0.15in]
& \ & \!\!\!\!\!\!\!\!\!\!\!\!\!  (A.6.5)  \ \cdot \cdot \cdot \cdot \cdot \cdot  \cdot \,\cdot \cdot \ \   \Longrightarrow \ \ \ \ \int_{B_{\,o} (\,\rho_o) }    \langle\,y\,, \    \btd \, K\,  (y)\,\rangle \,\, [\,v_i\, (y)]^{{2n}\over {n \,-\, 2}}   \,dy \ \,=\, \ O_{\lambda_i} (n \,-\, 2)\ .
\end{eqnarray*}
See also \cite{Leung-Supported}\,.\, Whereas from (2.12) and (2.19), we have
\begin{eqnarray*}
    \bigg\vert \ \int_{\R^n \, \setminus \,\Omega} \langle\,y\,, \    \btd \, K\,  (y)\,\rangle \,\, [\,v_i\, (y)]^{{2n}\over {n \,-\, 2}} \,\, dy \bigg\vert & \,=\, &  O_{\lambda_i}\, (n ) \ \ \ \ \ \ \ \ \ \ \ \ \ \ \ \ \ \ \ \ \ \ \ \  \\[0.15in]
(A.6.6) \ \cdot \cdot \cdot \cdot \cdot \cdot \cdot \cdot \ \  \Longrightarrow \ \   \int_{\Omega}   \langle\,y\,, \    \btd \, K\,  (y)\,\rangle \,\, [\,v_i\, (y)]^{{2n}\over {n \,-\, 2}} \,\, dy  & \,=\, &   O_{\lambda_i} \,(n )\,.
\end{eqnarray*}

In the discussion,
 $$
 \Omega \ \,=\, \ \bigcup_{j \,=\, 0}^k \ B_{ {\hat Y}_j}\, (\rho) \ \ \ \ \ \ \ \left[ \,  B_{ {\hat Y}_j}\, (\rho)\, \cap\, B_{ {\hat Y}_{\,l}} \, (\rho) \ \,=\, \ \emptyset \mfor j \,\not=\, l \,\right]\,, \leqno (A.6.7)
 $$
 where (as usual) $\,\{\  {{\hat {\!Y}}}_o\,=\,0\,,\,\cdot \cdot \cdot\,, \ {{\hat {\!Y}}}_k \}\,$ is the collection of all blow\,-\,up points.

\vspace*{0.2in}

{\bf \S\,A.6.b.}\ \ {\it    Linking the Pohozaev integral to the condition\,} $\Delta_o^{\!(h_\ell)}\, {\bf P}_\ell\,\equiv\,0\,.\,$\\[0.1in]
In the consideration of the  integrals in the Pohozaev identities\,,\,
we often encounter integral in the expression (A.6.9) below. We first record down the following observation.\\[0.2in]
 {\bf Lemma A.6.8.} \ \ {\it For a homogeneous polynomial $\,{\cal Q}_{\ell}\,$} ({\it defined on\,} $\R^n$) {\it of degree\,} $\ell \, \le \, n \,-\, 1$\,.\,  {\it If\,} \,$\ell$\, {\it is even\,,\, then   the following equivalence holds.}
 $$
\int_{\R^n} \   {\cal Q}_{\ell}\, (y) \cdot  \left( {1\over {1+ |\,y|^{\,2}}} \right)^n\,d y \ \,=\, \ 0 \ \ \ \ \Longleftrightarrow \ \ \ \Delta_{\,o}^{\!(h_\ell)}\ Q_\ell \ \,=\, \ 0\ . \leqno (A.6.9)
$$
({\it Recall that $\,h_\ell \,=\, \ell /2\,$ when $\,\ell\,$ is even\,.})
\vspace*{0.15in}

{\bf Proof.} \ \ We observe that, as $\,\ell \,\le\, n \,-\, 1\,,\,$   the integral in (A.6.9) is absolutely convergent. Keeping the notation $\,y \,=\, (y_{|_1}\,,\, \cdot \cdot \cdot\, \ y_{|_n} ) \in \R^n\,$,\, consider a typical term in $\,{\cal Q}_{\ell}\ $:

\vspace*{-0.3in}

$$ \ \
y_{|_1}^{\alpha_1} \cdot  y_{|_2}^{\alpha_2}   \cdot \cdot \cdot \   y_{|_{n  }}^{\alpha_{n  }}\,, \ \ \ \ \ \ \ \ {\mbox{where}} \ \ \alpha_j \,\ge\, 0 \ \ \  {\mbox{and}}  \ \ \ \sum_{j \,=\, 1}^n \alpha_j  \ \,=\, \ \ell \ \le \ n \,-\, 1\,.   \leqno (A.6.10)
$$
If one of the indices (say, $\alpha_j$) is an odd natural number, via symmetry, we have
$$
\, \int_{\R^n} \  [\  y_{|_1}^{\alpha_1} \cdot  y_{|_2}^{\alpha_2}   \cdot \cdot\, y_{|_j}^{\alpha_j} \cdot \cdot \   y_{|_{n  }}^{\alpha_{n  }}] \cdot  \!\left( {1\over {1+ |\,y|^{\,2}}} \right)^n\,d y \ \,=\, \ 0\ \ \ \ \    (\alpha_j \ \ {\mbox{is\ \  odd}}\,)\,.  \leqno (A.6.11)
$$
Direct calculation also shows that in this situation
$$
 \ \ \ \ \ \ \ \ \ \ \ \ \ \ \ \ \ \ \ \ \Delta_{\,o}^{\!({h_\ell})}\, \left[\ y_{|_1}^{\alpha_1} \cdot  y_{|_2}^{\alpha_2}   \,\cdot \cdot\  y_{|_j}^{\alpha_j} \cdot \cdot \     y_{|_{n  }}^{\alpha_{n  }} \right] \ \,=\, \ 0 \ \ \ \ \ \ \  \ \ \ \ \ \ \ \  \ \ \ \ \ \ \ \ \ \ \ (\alpha_j \ \ {\mbox{is\ \  odd}}\,)\,.
$$
(\,Recall that $\,\Delta_o^{\!(h_\ell)}\, {\cal Q}_\ell\,$ is a number when $\ell$ is even.\,)
Thus we are left with the case where any one index in (A.6.10) is an  even natural number  or zero. Let us introduce the following notion\,: a multi\,-\,index

\vspace*{-0.15in}

$$\ \ \  \ \ \ \ \ \
\alpha \,=\, (\alpha_1\,,\, \ \alpha_2\,,\, \ \cdot \cdot \cdot\,, \ \alpha_n) \ \ \ \ \ \ \ \ \  \ \ \ \ \ \ \bigg( \ |\, \alpha| \ \,=\, \ \sum_{j \,=\, 1}^n \,\alpha_j \ \,=\, \ \ell > 0 \, \bigg)
$$

\vspace*{-0.05in}

is {\it even} if each $\alpha_j\ \, (\,1 \,\le\, j \,\le\, n\,)\,$ is either  an even natural number or zero\,.\, With respect to this, the simplest case that can happen to the integral in (A.6.11) is
$$
J \, := \, \int_{\R^n}    y_{|_1}^2 \,\cdot \cdot \cdot  \  y_{|_{h_\ell}}^2 \!\left( {1\over {1+ r^{\,2}}} \right)^n\! d y\,.\leqno (A.6.12)
$$
We seek to reduce other even multi\,-\,index cases to that in (A.6.12)\,.\,
As $\,\ell \,\le\, n \,-\, 1\,$,\, we arrange in this way
$$
y_{|_1}^{\,k \,+\, 2} \cdot  y_{|_2}^{\alpha_2}   \cdot \cdot \cdot \   y_{|_{n \,-\, 1}}^{\,\alpha_{n \,-\, 1}}\,, \ \ \ \ \ \ \ {\mbox{where}} \ \ \ \ \alpha \,=\, (k \,+\, 2\,,\, \ \alpha_2\,,\, \ \cdot \cdot \cdot\,, \ \alpha_{n-1}\,, \ 0) \ \ \ {\mbox{is \ \ even}}\,.
$$
Here $\,k \,\ge \,2\,$ is an even number\,.\,
Via symmetry, the ordering is not important when we compute the integral in (A.6.11)\,.\,
One obtains the following reduction formula.
\begin{eqnarray*}
(A.6.13)& \ & \ \ \ \int_{\R^n} \ \,  y_{|_1}^{\,k \,+\, 2} \cdot  \left[ \ \ y_{|_2}^{\alpha_2}   \cdot \cdot \cdot \   y_{|_{n \,-\, 1}}^{\alpha_{n \,-\, 1}}\, \right]   \left( {1\over {1 \,+\, |\,y|^{\,2}}} \right)^{\,n} dy \\[0.15in]
 &\ & \!\!\!\!\!\!\!\!\!\!\!\!\!\!\!\!\!\!\!\!\!\!\!   \,=\, \ (k \,+\, 1)    \int_{\R^n}   \ \,   y_{|_1}^{\,k  }  \cdot y_{|_n}^{\,2  }  \cdot \left[ \ \ y_{|_2}^{\alpha_2}   \cdot \cdot \cdot \   y_{|_{n \,-\, 1}}^{\alpha_{n \,-\, 1}}\, \right]  \left( {1\over {1 \,+\, |\,y|^{\,2}}} \right)^{\,n} dy \mfor 2 \,\le\, k \,\le\, n \,-\, 3\,, \ \ \ \ \ \ \ \ \ \
\end{eqnarray*}
 by using  Fubini's theorem and integration by parts formula. See \S\,A.8 below\,.\bk
In view of (A.6.11) and (A.6.13), we introduce the following notation.  For an integer $\,m\,\ge\, 0$\,,\, define
\begin{eqnarray*}
(A.6.14) \ \ \ \ \ m\,!_{-2} & \,=\, &  1 \ \ \ \ {\mbox{if}} \ \ m \,=\, 0\ \ \ {\mbox{or \  }} 2\,; \ \ \ \ \ m\,!_{-2} \ \,=\,  \ 0 \ \ \ \ {\mbox{if}} \ \ m \ \ {\mbox{is \ \ odd}}\,;
\\[0.1in]
m\,!_{-2} & \,=\,  & ( m \,-\, 1) \,( m \,-\, 3)\,(m \,-\, 5) \cdot \cdot \cdot \ \cdot \,3 \cdot 1 \ \ \ \ \ \ \ {\mbox{if}} \ \ m \,\ge\, 4 \ \ \, {\mbox{is \ \ even}}\,. \ \ \ \ \ \  \
\end{eqnarray*}
Via the vanishing formula (A.6.11) and the reduction formula (A.6.13), we have
$$ \int_{\R^n} \  [\ y_1^{\,\alpha_1} \cdot \cdot \cdot  \,  y_n^{\,\alpha_n} \,] \left( {1\over {1+ |\,y|^{\,2}}} \right)^n\,d y \ \,=\, \  (\alpha_1)\,!_{-2} \,\times \,\cdot \cdot \cdot \,\times (\alpha_n)\,!_{-2}\cdot J\ .\leqno (A.6.15)
$$
On the other side,  calculation shows that
$$
B:= \Delta_{\,o}^{\!( h_\ell)}\, \left[\ y_1^2 \cdot y_2^2 \ \cdot \cdot \cdot \ y_{h_\ell}^2\ \right]  \ \,=\,  \ \ell \, (\ell \,-\, 2) \,(\ell \,-\, 4) \,\cdot  \cdot \cdot \,2 \cdot 1\,. \leqno (A.6.16)
$$

{\bf Claim.} \ \ Let $\,\alpha_2\,, \ \cdot \cdot \cdot\,, \ \alpha_{n \,-\, 1}\,$   be even natural numbers or zero,\, and
$$
\ell \,=\, ( k \,+\, 2) \,+\, {\alpha_2}  \ \,+\, \ \cdot \cdot \cdot \ \,+\, \   \alpha_{n \,-\, 1}\,, \ \ \ \ {\mbox{where}} \ \ k \ge 2 \ \ {\mbox{is \ \ an \ \ even \ \ integer}}\,. \leqno (A.6.17)
$$
Then
\begin{eqnarray*}
(A.6.18) \ \   \Delta_{\,o}^{\!( h_\ell)}  \left\{\,   y_{|_1}^{\,k \,+\, 2}  \cdot \left[ \ \ y_{|_2}^{\alpha_2}   \cdot \cdot \cdot \   y_{|_{n \,-\, 1}}^{\alpha_{n \,-\, 1}}\, \right]   \, \right\} \, \,=\, \, (k \,+\, 1) \cdot \!\Delta_{\,o}^{\!( h_\ell)}\, \left\{\, y_{|_1}^{\,k  }  \cdot  y_{|_n}^{\,2  }  \left[ \ \ y_{|_2}^{\alpha_2}   \cdot \cdot \cdot \   y_{|_{n \,-\, 1}}^{\alpha_{n \,-\, 1}}\, \right]   \right\}.
\end{eqnarray*}
Refer to (A.6.17), set
$$
{\cal o}:= \,{\alpha_2}  \ \,+\, \ \cdot \cdot \cdot \ \,+\, \   \alpha_{n \,-\, 1}\ . \leqno (A.6.19)
$$
We demonstrate how to use induction on $\,{\cal o}\,$ to prove the assertion in \S\,A.8 in this e\,-\,Appendix.
Thus using (A.6.18) repeatedly, we are led to
$$
\Delta_o^{\!(h_\ell)} \ \left[\,  y_{|_1}^{\alpha_1}   \cdot \cdot \cdot \   y_{|_{n }}^{\alpha_{n }}\,  \right] \ \,=\,  \  (\alpha_1)\,!_{-2} \,\times \,\cdot \cdot \cdot \,\times (\alpha_n)\,!_{-2}\cdot B\,. \leqno (A.6.20)
$$
Using the linearity of the operations and symmetry, (A.6.15) and (A.6.20) yield

\vspace*{-0.25in}

$$
\int_{\R^n} \   {\cal Q}_{\ell}\, (y) \cdot \left( {1\over {1+ |\,y|^{\,2}}} \right)^n\,d y \ \,=\,  \ {J\over B}\cdot  [\ \Delta_{\,o}^{\!(h_\ell)}\ {\cal Q}_{\ell}\,]\, . \leqno (A.6.21)
$$
In particular, we establish (A.6.9)\,.{{\hfill {$\rlap{$\sqcap$}\sqcup$}}

\vspace*{0.2in}

{\bf \S\,A.6.\,c.}\ \ {\it   Change of center.}\ \ With (6.57),\, let us  write
$$
v_i\, (y) \,=\, {\bf A}\, (y) \ \,+\, \ {\bf B}\, (y) \ \,+\, \ {\bf C}\, (y) \mfor |\,y| \ \le \ \rho_1\,, \leqno (A.6.22)
$$
where
\begin{eqnarray*}
(A.6.23) \ \ \ {\bf A}\, (y) & \,=\, &\! \left( {{\lambda_i}\over { \lambda_i^{\,2} \,+\, |\,y \,-\, \xi_{\,i}|^{\,2}}}\right)^{{n \,-\, 2}\over 2} \ \ ,\\[0.15in]
(A.6.24) \ \ \ {\bf B}\, (y) & \,=\, &\! \lambda_i\cdot \left[\, \lambda^{\,\ell  }_i \cdot \Gamma_{\bf p} \, \left({{y \,-\, \xi_i}\over { \lambda_i }}\right)\,\right]   \cdot \left( {{\lambda_i}\over {\lambda_i^2 \,+\, |\,y \,-\, \xi_{\,i}|^{\,2} }}\right)^{\, n \over 2} \ ,\\[0.175in]
(A.6.25) \ \ \ {\bf C}\, (y) & \,=\, &\!  [\, v_i \,-\, {\bf A} \,-\, {\bf B}\,]\, (y) \ =\  o_{\lambda_i} \left(\ell  \ \,-\, \ {{n-2}\over 2} \right) \ \,+\, \ O_{\lambda_i} \left(  {{n-2}\over 2} \right)\ .\  \ \ \ \  \ \  \ \ \ \  \ \  \ \ \ \  \ \
\end{eqnarray*}
 For simplicity, we suppress the subindex $\,i\,$ in \,{\bf A}\,,\, {\bf B}\, and \,{\bf C}\,.\, Because of the polynomial nature of $\,\Gamma_{\bf p}\,$,\, in general  ${\bf B}$ is not rotationally symmetric,  contrasting to ${\bf A}\,.\,$ As this part of the discussion is used repeatedly in this article, we consider in general a homogeneous polynomial $Q_\ell$ defined on $\,\R^n$\, with degree $\ell \,\in \,[\,2\,, \ n \,-\, 2\,]\,.\,$  For $\rho > 0\,,\,$ consider the integral
$$
  \int_{B_o (\,\rho)}  {Q_\ell}\, (y) \cdot   [\,{\bf A}\, (y )]^{{2n}\over {n \,-\, 2}} \, dy\,.
$$

Here $\,{\bf A}\,$ is given in (A.6.23)\,.\, We  first observe that
\begin{eqnarray*}
(A.6.26) & \ & \bigg\vert \ \int_{  B_o \, (\,\rho_2)} f\,(y) \cdot \left( {{{\lambda_i} }\over {{\lambda_i} ^{\,2} \,+\, r^{\,2}}}  \right)^n \ dy \ \,-\, \ \int_{  B_{\xi_i} \, (\,\rho_2)} f\,(y) \cdot \left( {{{\lambda_i} }\over {{\lambda_i} ^{\,2} \,+\, r^{\,2}}}  \right)^n \ dy \ \bigg\vert \\[0.15in]
& \le & \int_{  B_o \, (\,\rho_2\, \,+\, \,|\,\xi_i|)\,\setminus B_o \, (\,\rho_2 \,- \,|\,\xi_i|)} |\,f\,(y)| \cdot \left( {{{\lambda_i} }\over {{\lambda_i} ^{\,2} \,+\, r^{\,2}}}  \right)^n \ dy \ \,=\, \  \   O \, (|\,\xi_i|  \cdot \lambda_i^{\,n})  \ \ \ \ \ \ \ \ \ \ \ \
\end{eqnarray*}
for $\,i \,\gg\, 1\,.\,$ Here $\,f\,$ is a bounded continuous function defined on a slightly bigger ball\,,\, and $\,\xi_i\,$ continues to find its meaning in (1.13)\,.\, In particular, $\,|\,\xi_i| \,\to\, 0\,.\,$\ It follows that
$$
\int_{B_o \,(\,\rho)}  {Q_\ell}\, (y)\cdot  [\,{\bf A}\, (y )]^{{2n}\over {n \,-\, 2}} \, dy \  \,=\, \  \int_{B_{\xi_i} (\,\rho)}  {Q_\ell}\, (y) \cdot   [\,{\bf A}\, (y )]^{{2n}\over {n \,-\, 2}} \, dy \ \,+\, \ o_{\lambda_i} \,(n)\,. \leqno (A.6.27)
$$

Let us arrange
$$
  {Q_\ell}\, (y) \ =\   {Q_\ell}\, (\xi_i \,+\, [\,y \,-\, \xi_i\,]\,)
  \ \,=\, \ {Q_\ell}\, ( \xi_i ) \ +\ M^Q \, (y \,-\, \xi_i)\ \,+\, \  {Q_\ell}\, ( y \,-\, \xi_i )\ .\leqno (A.6.28)
$$
Here the ``intermediate" term $\,M^{Q} \, (\,y \,-\, \xi_i)$ can be further broken down into $\,\ell \,-\, 1\,$ terms based on the degree on $\,\xi_i\,:$
\begin{eqnarray*}
(A.6.29)  \!\!\!\!& \ & M^Q\, (\xi_i\,;\ y \,-\, \xi_i) \\[0.15in]
  [  & =
&   O\, (|\,\xi_i| \,|\,y \,-\, \xi_i|^{\ell \,-\, 1})\ \,+\, \   O\, (|\,\xi_i|^2\, |\,y \,-\, \xi_i|^{\ell \,-\, 2})   \ \,+\, \ \cdot \cdot \cdot \  \,+\, O\,(|\,\xi_i|^{\ell \,-\, 1}\cdot |\,y \,-\, \xi_i|\, )\ ]  \ \ \ \ \  \\[0.15in]
  & \ & \,=\, \ {\Xi^{\,Q}}_{\!\!\!\!\!\!1}\  \,(\xi_i\,; \ y- \xi_i) \ \,+\, \ {\Xi^{\,Q}}_{\!\!\!\!\!\!2}\  \,(\xi_i\,, \ y- \xi_i)\,+ \ \cdot \cdot \cdot \ \,+\, \ {\Xi^{\,Q}}_{\!\!\!\!\!\!\ell \,-\, 1}\  \,(\xi_i\,, \ y- \xi_i)\ ,
\end{eqnarray*}
$$
{\mbox{where \ \ (formally)}} \ \ \ \ \   {\Xi^{\,Q}}_{\!\!\!\!\!\!h}\ \, (\,\xi_i\,, \ z) \,=\,     \sum_{|\, \alpha| \,=\, h} {1\over {\alpha\,!}} \,\cdot\,   \xi_i^\alpha \ D_\alpha^{(h)}\, {Q_\ell} \, (z)\bigg\vert_{\,z \,=\, (\,y \,-\, \xi_i)}\     \leqno (A.6.30)
$$
for $\, 1 \,\le\, h \,\le\, \ell \,-\, 1\,.\,$
We continue with
\begin{eqnarray*}
(A.6.31)& \ & \!\!\!\!\!\int_{B_{\xi_i}  (\,\rho)}  {Q_\ell}\, (y) \cdot   [\,{\bf A}\,(y )]^{{2n}\over {n \,-\, 2}} \  dy  \  \,=\, \
 \int_{B_{\xi_i}  (\,\rho)}   {Q_\ell}\, (y) \,  \left( {{\lambda_i}\over {\lambda_i^2 \,+\, |\,y \,-\, \xi_i|^2}} \right)^n\ dy  \\[0.2in]
   &\  &\!\!\!\!\!\!\!\!\!\!\!\!\!\!\!\!\!\!\!\!\!\!\!\!\!\!\!\!\! =\
 \int_{B_{\xi_i} \,(\,\rho)}   {Q_\ell}\, (y \,-\, \xi_i) \,  \left( {{\lambda_i}\over {\lambda_i^2 \,+\, |\,y \,-\, \xi_i|^2}} \right)^n\ dy \\[0.15in]
  & \ & \!\!\!\!\!\!\!\!\!\!\!\!\!\!\!\!\!\!\!\!\!\!+ \  \int_{B_{\xi_i}  (\,\rho)}   \!\!\!M^Q\, (\xi_i\,;\ y \,-\, \xi_i)   \left( {{\lambda_i}\over {\lambda_i^2 \,+\, |\,y \,-\, \xi_i|^2}} \right)^n\! dy     \ \,+\, \ \int_{B_{\xi_i} \,(\,\rho)}   \!\!{Q_\ell}\, ( \xi_i) \,  \left( {{\lambda_i}\over {\lambda_i^2 \,+\, |\,y \,-\, \xi_i|^2}} \right)^n\! dy  \\[0.2in]& \,=\, &
 \int_{B_o \,(\,\rho)}  {Q_\ell}\, (z) \,  \left( {{\lambda_i}\over {\lambda_i^2 \,+\, |\,z|^2}} \right)^n\! dz \ \,+\, \ \int_{B_o \,(\,\rho)}  M^Q\, (\xi_i\,;\ z) \,  \left( {{\lambda_i}\over {\lambda_i^2 \,+\, |\,z|^2}} \right)^n\ dz  \ \,+\, \\[0.15in]
 & \ & \ \ \ \ \ \ \ \  \ \,+\, \ {Q_\ell}\, ( \xi_i) \, \int_{B_{o} (\,\rho)}    \left( {{\lambda_i}\over {\lambda_i^2 \,+\, |\,z|^2}} \right)^n\ dz  \ \ \ \ \   \ \ \ \ \ \ \ \   \ \ \ \ \     [\,z \,=\, (y \,-\, \xi_i)\,]\, . \\[0.2in]
 & \ & \!\!\!\!\!\!\!\!\!\!\!\!\!\!\!\!\!\!\!\!\!\! \,=\, \
 \lambda_i^\ell \, \int_{B_o \,(\,\lambda_i^{-1}\cdot\,\rho)}  Q_\ell \, (y) \,  \left( {{1}\over {1 \,+\, |\,y|^2}} \right)^n\! dy \ \,+\, \\[0.15in]
 & \ & \!\!\!\!\!\!\!\!  \,+\, \ \lambda_i^\ell \left[ \ \sum_{h \,=\,1}^{\ell \,-\, 1} \int_{B_o \,(\,\lambda_i^{-1}\cdot\,\rho)} \Xi^{\,Q}_{h}\ \left( {{\xi_i}\over {\lambda_i}}\,, \ y \right) \,  \left( {{1}\over {1 \,+\, |\,y|^2}} \right)^n\ dy \right] \ \,+\, \ \\[0.15in]
 & \ &   \ \ \ \ \,+\, \  \lambda_i^\ell  \cdot\left[ \, Q_\ell \left( {{\xi_i}\over {\lambda_i}}\right) \right]\cdot \int_{B_o \,(\,\lambda_i^{-1}\cdot\,\rho)} \left( {{1}\over {1 \,+\, |\,y|^2}} \right)^n\ dy \ \ \ \ \ \ \ \ \ \ \ (\,z \  \to \ \lambda_i \cdot y)\,.\\
  \end{eqnarray*}

\vspace*{0.1in}

{\bf \S\,A.6.\,c.1}\ \ {\it  Expressing the integrals on\,} $\R^n$\,. \ \ Observe also that the above argument remains valid if we replace $\,\rho\,$ by $\,\infty\,.\,$  Indeed, consider (in general) a homogeneous polynomial $\,{\cal Q}_k$\, defined on $\,\R^n\,$ with degree $k \in [\,0\,,\, \ n \,-\, 2]\,.\,$ For  a sequence of positive numbers $\,r_i \to \infty\,:\,$
    $$
 \bigg\vert \,\int_{\R^n \,\setminus B_o \,(\,r_i)}  {\cal Q}_k \, (y) \,  \left( {{1}\over {1 \,+\, |\,y|^2}} \right)^n\! dy \bigg\vert  \ \le \ C \int_{r_i}^\infty {{r^k}\over {r^{2n} }}\cdot  r^{n \,-\, 1} \, dr \ \le \ C \int_{r_i}^\infty {{1}\over {r^{3} }} \, dr \ \le {{C_2}\over {r^2_i}}
   $$
   $$
   \!\!\!\!\cdot \cdot \cdot \, \Longrightarrow  \int_{B_o \,(\,r_i)}  {\cal Q}_k\, (y) \,  \left( {{1}\over {1 \,+\, |\,y|^2}} \right)^n\! dy  \,=   \int_{\R^n} {\cal Q}_k\, (y)  \left( {{1}\over {1 \,+\, |\,y|^2}} \right)^n\! dy \,+\ O \,(r_i^{-2})\,, \leqno (A.6.32)
   $$
  for $\,i \,\gg\, 1\,$ and $\,0 \,\le\, k \,\le\, n \,-\, \,2\,.\,$ Putting $\,r_i \, \,=\, \, \lambda_i^{-1} \cdot \rho$\,,\, and combining (A.6.27), (A.6.31) and (A.6.32) with $\,|\,\xi_i| \to 0\,,\,$  we obtain

\vspace*{-0.25in}

  \begin{eqnarray*}
 (A.6.33)   & \ & \int_{B_o \,(\,\rho)}  {Q_\ell}\, (y)\cdot  [\,{\bf A}\, (y)]^{{2n}\over {n \,-\, 2}} \, dy\  \,=\, \
 \lambda_i^\ell \,\int_{\R^n}  Q_\ell \, (y) \,  \left( {{1}\over {1 \,+\, |\,y|^2}} \right)^n\! dy \ \,+\, \ \ \ \ \ \ \ \ \ \ \\[0.2in]
  & \ & \ \ \ \ \ \  \ \ \ \ \ \ \ \ \  \ \ \  \ \ \ \ \    \,+\, \ \lambda_i^\ell \left\{ \ \sum_{h \,=\,1}^{\ell \,-\, 1} \int_{\R^n}  \left[ \, \Xi^{\,Q}_{h}\! \left( {{\xi_i}\over {\lambda_i}}\,, \ y \right) \,\right]\cdot  \left( {{1}\over {1 \,+\, |\,y|^2}} \right)^n\ dy \right\} \ \,+\, \ \ \ \ \ \ \ \ \ \ \  \\[0.1in]
& \ &  \ \ \ \ \ \   \,+\, \  \lambda_i^\ell  \cdot\left[ \, Q_\ell \left( {{\xi_i}\over {\lambda_i}}\right) \right]\cdot \int_{\R^n}  \left( {{1}\over {1 \,+\, |\,y|^2}} \right)^n\ dy \ \,+\, \ O_{\lambda_i}\, ( \ell \,+\, 2 ) \ \,+\, \ o_{\lambda_i} ( n ) \,. \ \ \ \ \  \ \ \ \  \end{eqnarray*}

\vspace*{0.25in}

{\bf \S\,A.6.d.}\ \ {\it
  $\Delta_o^{\!h_\ell}\,{\bf P}_\ell \, \equiv \, 0$ when $\,\ell$\, is even and } $\ell \,<\, n \,-\, 2$\,,\, {\it or  when $\,\ell \,= \,n \,-\, 2\,$}\\[0.075in]
   \hspace*{0.6in}(\,$n$ {\it being even}) {\it  and with only one simple blow\,-\,up point\,.}    \\[0.3in]
  {\bf Proposition A.6.34.} \ \ {\it For $\,n \ge 4\,,\,$  under the general conditions\,} (1.6)\,,\, (1.25),\, (1.26)\,,\,  {\it assume that $\,\{ u_i \}\,$ has finite number of blow\,-\,up points, one at the south pole, but  none at the north pole.\,  Take the following conditions\ } (i)\,--\,(iii)\,{\it into account}\,.\\[0.15in]
(i) \ \  $\ 0$ {\it is a simple blow\,-\,up point for $\,\{ v_i\}\,.$}\\[0.1in]
(ii) \ \,    $K $ {\it is given by\,} (1.8) \,{\it in\,} $B_o \, (\rho_o)$\,,\,  {\it where\ } $2 \,\le \, \ell < \, n- 2\,.\,$ \\[0.1in]
(iii) \     {\it The parameters $\,\lambda_i\,$ and $\,\xi_i\,$ corresponding to the simple blow\,-\,up point at $\,0\,$ } \\[0.075in]
 \hspace*{0.31in}  \,[\,{\it via\,} (1.10) {\it and\,} (1.11)\,]\, {\it satisfy\,} (1.12)\,,\, {\it that is\,,\,} $|\,\xi_i|\, =\, o\, (\lambda_i)\,.$ \\[0.1in]
(iv) \  \,$\ell\,$ {\it is even}\,.\,\\[0.15in]
{\it Then $\,\Delta_o^{\!(h_\ell)} \, {\bf P}_\ell\, (y) \ \,=\, \ 0\,$}\, ($\,\Delta_o^{\!(h_\ell)} \, {\bf P}_\ell$\, {\it is a number when $\,\ell\,$ is even}\,)\,.\,  {\it The same conclusion also holds  when $\,\ell\, \,=\, \,n \,-\, 2\,$ with an additional assumption that $\,0\,$ is\, the only blow\,-\,up point} (\,$\ell$ {\it is still required to be even\,})\,.

\vspace*{0.15in}

{\bf Proof.} \ \ The key is to combine  the change of center formula (A.6.33) with the condition  \,$|\,\xi_i|\, \,=\, \,o\, (\lambda_i)\,$,\, and observe the lower order terms.\, Other  arguments actually proceed in   similar fashion as those found  in \cite{Leung-Supported} and \cite{Li-1}\,.\,  For the benefit of readers, we present the estimates in detail. From (A.6.5) and (A.6.6), we have
$$
(A.6.35) \ \ \ \int_{B_{\,o} (\,\rho_o) }    \langle\,y\,, \    \btd \, K\,  (y)\,\rangle \,\, [\,v_i\, (y)]^{{2n}\over {n \,-\, 2}}   \,dy   \  \,=\, \  \cases{\ O_{\lambda_i} (n \,-\, 2)\,, \!& in \ general\,; \cr \cr
              \,O_{\lambda_i} (n )\,, \!&    one blow\,-\,up point. \cr}
$$
Throughout this proof we assume that the positive constant \,$\rho_o > 0\,$ is chosen to be small enough\,.\,
%
Let us pay attention to (2.21) and the number $\,R_i\,$ satisfying (2.20), together with the remark in \,\S\,2\,f\, (on shifting the center, see also \S\,A.2\,)\,.\, Note that
\begin{eqnarray*}
(A.6.36)\!\! \!\!\! & \ & \int_{ B_o (\rho_o) }  \!   \langle\,y\,, \    \btd \, K\,  (y)\,\rangle \,   [\,v_i\, (y)]^{{2n}\over {n \,-\, 2}}\,dy \ \,=\, \ \int_{B_o (\lambda_i \, R_i) } r\, {{\partial K}\over {\partial r}} \,   [\,{\bf A} \, (y)]^{{2n}\over {n \,-\, 2}}\ dy \ \cdot \cdot \cdot \ ({\bf I}) \\[0.15in]
& \ & \ \  \ \ \  \ \ \  \ \ \  \   \  \ \  \ \ \ \  \ \  +\ \int_{ B_o (\lambda_i \, R_i) }   r\, {{\partial K}\over {\partial r}} \, \left\{ \,[\,v_i\, (y)]^{{2n}\over {n \,-\, 2}} \,-\,[\,{\bf A} \, (y)]^{{2n}\over {n \,-\, 2}}\right\} \ dy\ \cdot \cdot \cdot \ ({\bf II})  \\[0.15in]
& \ & \ \  \ \ \  \ \ \  \ \   \ \  \  \ \ \  \ \  \ \ \ \ \ \  \ \ \ \ \ \  \ \ \ \     +\
\int_{B_o (\rho_o)\,\setminus \, B_o (\lambda_i \, R_i) }  r\, {{\partial K}\over {\partial r}} \,  [\,v_i\, (y)]^{{2n}\over {n \,-\, 2}}\ dy\,.\ \cdot \cdot \cdot \ ({\bf III})\ \ \ \   \\
\end{eqnarray*}

\vspace*{-0.1in}

(i) \ \ We begin  with the core term ({\bf I})\,.\, From (1.8), we have
$$
r \cdot {{\partial \,[\,{\tilde c}_n\,K]}\over {\partial r}} \ \,=\, \ \langle \,y\,, \ \btd \, [\,{\tilde c}_n\,K] \ \rangle \ \,=\, \ \ell \times [\,-{\bf P}_\ell\, (y)] \ \,+\, \ O\, \left( |\,y|^{\,\ell \,+\, 1}\,\right)  \leqno (A.6.37)
$$
for $\,y \in B_o\, (\rho_o)\,.\,$  It follows that [\,recall that \,${\bf A}$\, is given in (A.6.23)\,]

\vspace*{-0.3in}

\begin{eqnarray*}
& \ & \\
(A.6.38)\!\!\!\!\!\!\!\!\!& \ & \int_{ B_o (\lambda_i \, R_i) } r\, {{\partial K}\over {\partial r}} \, [\,{\bf A}\,(y)]^{{2n}\over {n \,-\, 2}} \ dy\ \,=\, \ -\,{{ \ell }\over {c_n}} \, \cdot \int_{ B_o (\lambda_i \, R_i) } {\bf P}_\ell\,\cdot [\,{\bf A}\,(y)]^{{2n}\over {n \,-\, 2}} \  dy \ \,+\, \\[0.15in]
 & \ &\ \ \ \ \ \ \ \  \ \ \ \ \ \ \ \ \ \  \ \ \ \ \ \ \ \ \ \  \ \ \ \ \ \  \,+\, \ \int_{ B_o (\lambda_i \, R_i) } O\, (\,|\,y|^{\,\ell \,+\, 1}\,) \cdot [\,{\bf A}\,(y)]^{{2n}\over {n \,-\, 2}} \  dy \,. \ \ \ \ \ \ \ \ \\[0.15in]
 & \  & \ \ \   =\   -\,{{ \ell }\over {c_n}} \, \cdot \int_{ B_o (\lambda_i \, R_i) } {\bf P}_\ell\,\cdot [\,{\bf A}\,(y)]^{{2n}\over {n \,-\, 2}} \  dy \ \ \,+\, \ O_{\lambda_i} \, (\ell \,+\, 1) \ \ \ \mfor \ell \,\le\, n \,- 2 \ \ \ \ \ \ \ \\[0.1in]
 & \ &    \!\!\!\![\, {\mbox{see\ \  (A.6.41)\,; \ \ this \ \ only \ \ requires}} \ \   |\,\xi_i| \,=\, O\, (\lambda_i)  \ \uparrow \,]\\[0.15in]
 & \  &  \ \ \   =\   -\,\lambda_i^\ell\cdot {{ \ell }\over {c_n}} \, \cdot \int_{ \R^n } {\bf P}_\ell\,(y)\cdot \left( {1\over {1 \,+\, |\,y  |^2}} \right)^n \  dy \ \,+\, \ o_{\lambda_i} \, (\ell) \\[0.11in]
 & \ & \!\!\!\!\!\!\! [\, {\mbox{estimating \   as \    in \   (A.6.31)\,--\ (A.6.33)\,, \ }} \uparrow  \ \ {\mbox{using \ \ }} \bigg\vert\,{{\xi_i}\over {\lambda_i}}\bigg\vert  \ =\  o\,(1) \ \ {\mbox{and}}  \\[0.1in]
    & \ & \ \ \ \   {\mbox {  replacing}} \ \, \rho_o \ \  {\mbox{by}}  \ \,\lambda_i\cdot \!R_i  \ \ {\mbox{in \ \ (A.6.31)}}\,, \ \ {\mbox{and \ \ take \ \ }} r_i \,=\, R_i  \ \ {\mbox{in \ \ (A.6.32)}} \, ]\,.
\end{eqnarray*}
Here $\,2 \,\le \, \ell\,\le \, n \,-\, 2\,.\,$

\vspace*{0.15in}

(ii) \ \ Now we turn to the term marked ({\bf II})\,.\,
Via inequality  (3.14) and (2.26) in Proposition 2.24, we have

\vspace*{-0.1in}

\begin{eqnarray*}
(A.6.39) & \ &   \bigg\vert\ \int_{ B_{\xi_i} (\lambda_i \, R_i) }  \!   \langle\,y\,, \    \btd \, K\,  (y)\,\rangle \,  \left\{ \,[\,v_i\, (y)]^{{2n}\over {n \,-\, 2}} \,-\,[\,{\bf A} \, (y)]^{{2n}\over {n \,-\, 2}}\right\} \,dy\,\bigg\vert\\[0.1in]
& \le & C \cdot {{\varepsilon }\over {\lambda_i^{{n \,-\, 2}\over 2} }}  \cdot {1\over {\lambda_i^{{n \,+\, 2}\over 2} }} \cdot \int_{ B_{\xi_i} (\lambda_i \, R_i) }  \,r^{\ell}\, \,dy \ \ \ \ \ \ \ \ \ \ \ \ \ \ \ \ \ \ \ \ \ \ \ \  (\,r\, =\, |\,y|)\\[0.1in]
& \le & C   \cdot {\varepsilon \over {\lambda_i^{n  }}} \cdot \int_{ B_{o} (\lambda_i \, (R_i\,+\,1)) }  \,r^{\ell}\, \,dy \ \ \ \ \ \ \ \ \ \ \ \ \ \ \ \ \ \  \ \ \ \  \ \ \ \ \ \ \ \ \ (\, \lambda_i\cdot |\,\xi_i| \to 0\,)\\[0.1in]
& \le & C_1 \, \cdot {\varepsilon \over {\lambda_i^{n }}} \cdot [\,\lambda_i \, (R_i\,+\,1)]^{\ell \,+\, n} \ \,=\, \ C_1 \cdot \lambda^\ell \,[\,\varepsilon\, (R_i\,+\,1)^{\ell \,+\, n} ] \ \,=\, \ o_{\lambda_i} (\ell)\,. \ \ \ \ \ \  \ \ \ \ \ \  \ \ \ \
\end{eqnarray*}
In the above we apply
$$
\ \ \ \ \ \ \  \varepsilon_i \cdot R_i^{\,2(n \,-\, 1)}\ \,=\, \ o\, (1) \ \ \ \ \ \ \ \ \ \ \ \ \ \ \  \ \ \ \ \ \  \ \ \ \ \ \  \  \ \ \  \ \ \ [\,{\mbox{see}} \ \ (2.20)\,]\,.
$$
(iii) \ \ As for the term marked ({\bf III}), using the binomial expansion on $\,(|\,z| \,+\, |\,\xi_i|\,)^\ell$\,,\, together with (2.5) or (2.26) (similar to Lemma 2.4 in \cite{Li-1}\,), we have
\begin{eqnarray*}
& \ &  \!\!\!\!\!\!\!\!\!\!\!\!\!\!\!(A.6.40) \ \ \ \ \bigg\vert \ \int_{B_o (\,\rho_o) \,\setminus B_o (\lambda_i \, R_i) }  \!\!  \langle\,y\,, \    \btd \, K\,  (y)\,\rangle \,  [\,v_i\, (y)]^{{2n}\over {n \,-\, 2}}\,dy \, \bigg\vert \\[0.15in]
  & \le &  C_1 \cdot  \int_{B_o (\,\rho_o) \,\setminus B_o (\lambda_i \, R_i) }  |\,y|^\ell \cdot \left( {{\lambda_i}\over {\lambda_i^2 \,+\, |\, y \,-\, \xi_i|^2}} \right)^{n}\,dy\ \\[0.15in]
  & \le &  C_2 \cdot  \int_{B_{\xi_i} (\,\rho_o') \,\setminus B_{\xi_i} (\lambda_i \, (R_i -\,c)) }  |\,y|^\ell \cdot \left( {{\lambda_i}\over {\lambda_i^2 \,+\, |\, y \,-\, \xi_i|^2}} \right)^{n}\,dy\\[0.1in]
   & \ &  \ \ \ \ \ \ \  \ \ \ \ \   [\,{\mbox{here}} \ \ \rho_o' \ \ {\mbox{is \ \ slightly \ \ bigger \ \ than \ \ }} \rho_o\,, \\[0.1in]
   & \ &  \ \ \ \ \ \ \  \ \ \ \ \ \ \  \ \ \ \ \ \ \    \ \   \      c \ \ {\mbox{is \  \ a \ \ \ big \ \ enough  \ \ constant \ \ (require \ \ }} |\,\xi_i| \,=\, O \,(\lambda_i))\,]\\[0.15in]
  & \le &  C_3 \cdot  \int_{B_o (\,\rho_o') \,\setminus B_o (\lambda_i \, (R_i -\,c)) }  \left( \ \sum_{j \,=\, 0}^\ell \,|\,z|^{\ell -\, j} \cdot |\,\xi_i|^j \right)\cdot \left( {{\lambda_i}\over {\lambda_i^2 \,+\, |\, z |^2}} \right)^{n}\,dz \ \ \ \ \ \    (\,y \,=\, z \,+\, \xi\,)\\[0.15in]
  &  \  &\!\!\!\!\!\!\!\!\!\!\!\!\!\!\!\le \ C_4\,\lambda_i^\ell  \cdot \!\left[ \, \sum_{j \,=\, 0}^\ell  \,  \int^{{{\rho_o'}\over {\lambda_i}}}_{R_i -\,c }
    \left( {{1}\over {1 \,+\, r^2}} \right)^{n}r^{\ell \,-\, j \,+\, (n \,-\, 1)}\,dr \right] \ \     (\,{\mbox{polar \   coordinates \   and}} \   r  \to  \lambda_i \cdot r)\\[0.15in]
  &  \  &\!\!\!\!\!\!\!\!\!\!\!\!\!\!\!=\  o\, (\lambda_i^\ell\,) \ \ \ \ \ \ \ \ \ \ \ \ \ \ \ \ \ \ \ \ \ \ \ \ \ \ \ \ \ \ \ \ \ \ \ \ \ \ \ \ \ \ \ \ \ \ \ \ \ \ \ \ \ \ \ \ \ \ \ \ \ \ \ \ \ \  \mfor 2\, \le\, \ell \,\le \, n \,-\, 2 \ \, \ \\[0.1in]
  & \ &  \!\!\!\!\!\!\!  (\,{\mbox{needs \ \ only\ \ }} |\,\xi|\,=\,O\,(\lambda_i)\,;\, \ \ \ R_i \to \infty \ \ {\mbox{and}} \ \ \lambda_i  \cdot\! R_i \to 0 \ \Longrightarrow \ \ \lambda_i^{-1} \cdot \rho_o > R_i\,)\,.
\end{eqnarray*}

Similarly,

\vspace*{-0.17in}

$$
 \int_{  B_o (\lambda_i \, R_i) } \!\! |\,y|^{\ell \,+\, 1} \left( \!{{\lambda_i}\over {\lambda_i^2 \,+\, |\, y \,-\, \xi_i|^{\,2}}} \!\right)^n\!\!dy  \ \le \  C_1 \int_{  B_{\xi_i} (\lambda_i \, (R_i \,+\, \, c)) }  \!\!\!\!\!\!\!\!|\,y|^{\ell \,+\, 1}   \left( \!{{\lambda_i}\over {\lambda_i^2 \,+\, |\, y \,-\, \xi_i|^{\,2}}} \! \right)^n \!\! dy \leqno (A.6.41)
$$
$$
\ \  \le C_2    \int_{  B_o (\lambda_i \, (R_i +\, c)) }  \left( \ \sum_{j \,=\, 0}^{\ell \,+\, 1} \,|\,z|^{\ell \,+\, 1 -\, j} \cdot |\,\xi_i|^j \right)\cdot \left( \!{{\lambda_i}\over {\lambda_i^2 \,+\, |\, z |^{\,2}}} \!\right)^{n}\,dz  \ \,=\, \ O_{\lambda_i} ( \ell \,+\, 1)\,,
$$
which requires only $\,|\,\xi_i| \,=\, O\,(\lambda_i)\,.\,$ Here (as before) $\,y \,=\, z \,+\, \xi\,,\,$ and we apply the change of variables $\,z \to \lambda_i \cdot z\,.\,$  Using (A.6.35), (A.6.36), (A.6.38), (A.6.39)\, and \,(A.6.40)
we obtain
\begin{eqnarray*}
& \ & \!\!\!\!\!\!\!\!\!\!\!\!\!\!\!\int_{B_{\,o} (\,\rho_o) }    \langle\,y\,, \    \btd \, K\,  (y)\,\rangle \cdot [\,v_i\, (y)]^{{2n}\over {n \,-\, 2}}   \,dy \ \,=\, \ \lambda_i^\ell \, \int_{\R^n}  {\bf P}_\ell \, (y) \,  \left( {{1}\over {1 \,+\, |\,y|^2}} \right)^n\! dy \ \,+\, \ o_{\lambda_i}\, (\ell)\\
(A.6.42)\  &\cdot & \!\!\!\cdot \,\cdot \,\cdot \ \Longrightarrow \ \ \int_{\R^n}  {\bf P}_\ell \, (y) \,  \left( {{1}\over {1 \,+\, |\,y|^2}} \right)^n\! dy \ \,=\, \ 0\,.
\end{eqnarray*}
Here $\,2\,\le\, \ell \,\le\, n \,-\, 2\,.$\, (A.6.42) together with Lemma A.6.8 \,($\ell$ is even)\, imply that $\,\Delta_o^{\!(h_\ell)} \, {\bf P}_\ell\, (y) \equiv 0\,.$\, Let us end the proof with the remark that the condition $\,|\,\xi_i| \,=\, o\, (\lambda_i)\,$ is only `fully' used  in the last step in (A.6.38)\,.\qed
We highlight that the smaller order term in (A.6.41) depends on convergence parameters $\,\varepsilon_i\,$ and $\,R_i$\,,\, as well as the condition $\,|\,\xi_i| \,=\, o\,(\lambda_i)\,.\,$  In the next section, we apply the refined estimate (A.6.22) to discern out the layers of information hidden in $\,o_{\lambda_i}\, (\ell)\,.\,$

\newpage


{\bf \S\,A.6.\,e. \ \ } {\it
  Isolating the key term with lowest order in $\lambda_i$\,.}   \ \ Suppose that estimate (6.57) holds for \,$\{ v_i \}$\, inside $B_o\,(\rho_2)\,,\,$ where $\,\rho_2 \,>\, 0\,$ is a constant. From (A.6.22)\,--\ (A.6.25), we obtain \begin{eqnarray*}
(A.6.43) \   \ \  \int_{B_o \,(\rho_2)} r\, {{\partial K}\over {\partial r}} \cdot [\,v_i\, (y)]^{{2n}\over {n \,-\, 2}} \  dy
 & =
 & \int_{B_o \,(\rho_2)} r\, {{\partial K}\over {\partial r}} \, [\,{\bf A}\, (y)]^{{2n}\over {n \,-\, 2}} \, dy \ \,+\, \\[0.1in]
   \,+\, \
\int_{B_o \,(\rho_2)} & \ &  \!\!\!\!\!\!\!\!\!\!\!\!\!\!\! r\, {{\partial K}\over {\partial r}} \, \left( [\,({\bf A} \,+\, {\bf B} \,+\, {\bf C})\,(y)\,]^{{2n}\over {n \,-\, 2}} \,-\, [\,{\bf A}\, (y)]^{{2n}\over {n \,-\, 2}}\right) dy\,. \ \ \ \ \ \ \ \ \ \ \ \
\end{eqnarray*}
In order to estimate the last term in the above, we make use of the inequality
 \begin{eqnarray*}
(A.6.44) \!\!\!\!\!\!  & \ & \int_{B_o \,(\,\rho_2)} \bigg\vert\,r\, {{\partial K}\over {\partial r}} \bigg\vert\cdot \bigg\vert \ [\,({\bf A} \,+\, {\bf B} \,+\, {\bf C})\,(y)]^{{2n}\over {n \,-\, 2}} \,-\, [\,{\bf A}\, (y)]^{{2n}\over {n \,-\, 2}}\ \bigg\vert \ dy \,\bigg\vert \\[0.1in]
 & \le & \varepsilon  \int_{B_o \,(\,\rho_2)} \bigg\vert \, r\, {{\partial K}\over {\partial r}} \bigg\vert  \   {\bf A}^{{2n}\over {n \,-\, 2}}\, dy \ \,+\, \ {{C_n}\over {\varepsilon^{{2n}\over {n \,-\, 2}} }} \cdot  \int_{B_o \,(\,\rho_2)}  \bigg\vert \, r\, {{\partial K}\over {\partial r}} \bigg\vert  \, \left( |\,{\bf B}|^{{2n}\over {n \,-\, 2}} \,+\, |\,{\bf C}|^{{2n}\over {n \,-\, 2}}\right) dy\ . \ \ \ \ \ \ \ \
 \end{eqnarray*}
 Here $\,\varepsilon > 0\,$ is a given (small) number, and the dimensional constant $C_n$ is independent on $\,\varepsilon\,.\,$
  We demonstrate  the argument toward (A.6.44) in   \S\,A.3 in the  Appendix\,.\,

\vspace*{0.15in}

 {\it Remark}\, A.6.45.  \ \ Suppose that we seek to find \,$\varepsilon > 0\,$ so that
$$
\varepsilon \cdot \lambda_i^\ell \ \,+\, \ {{C_n}\over {\varepsilon^{{2n}\over {n \,-\, 2}} }} \cdot  \lambda_i^{\ell\, \,+\, a} \ \,=\, \  \varepsilon \cdot \lambda_i^\ell \ \,+\, \ {{C_n}\over {\varepsilon^{{2n}\over {n \,-\, 2}} }} \cdot  \lambda_i^{\,a \,-\, t} \cdot \lambda_i^{\,\ell \,+\, t} \ .
$$
That is, we want to re\,-\,distribute some order of $\,\lambda_i\,$ to the first term so that in the end the two terms have the same order $\,O_{\lambda_i}\, (\ell \,+\, t)$\,:
$$
  \varepsilon \, \,=\, \, (\lambda_i)^t   \Longrightarrow \,     {{C_n}\over {( \lambda_i  )^{\,{{2n}\over {n \,-\, 2}} \cdot t} }} \, \,=\, \, \lambda_i^{\,a \,-\, t} \,
\Longrightarrow \, {{2n}\over {n \,-\, 2}}\cdot\, t \, \,=\, \,  (a \,-\, t)\  \Longrightarrow  \, t\,=\, {{n \,-\, 2}\over {3\,n \,-\, 2}}\,\cdot  \,a\,.  \ \ \
 $$

 \vspace*{0.1in}

{\bf \S\,A.6.f}  \ \ {\it
  Estimate on the leading order term.}   \ \
Recall (1.8) and (4.12)\,.\,
\begin{eqnarray*}
r \cdot {{\partial \,[\,{\tilde c}_n\,K]}\over {\partial r}} & \,=\, & \langle \,y\,, \ \btd \, [\,{\tilde c}_n\,K] \ \rangle \ \,=\, \ \ell \times [\,-{\bf P}_\ell\, (y)\,] \ \,+\, \ O\, (\,|\,y|^{\,\ell \,+\, 1}\,) \\[0.1in]
(A.6.46) \ \cdot \cdot \cdot \cdot  \Longrightarrow & \ &  \!\!\!\!\!\!\!\!\!\!\!\!\! \int_{B_o \,(\,\rho_2\,)} r\, {{\partial K}\over {\partial r}} \, [\,{\bf A}\,(y)\,]^{{2n}\over {n \,-\, 2}} \ dy\ \,=\, \ -\,{{ \ell }\over {c_n}} \, \cdot \int_{B_o \,(\,\rho_2\,)} {\bf P}_\ell\,\cdot [\,{\bf A}\,(y)\,]^{{2n}\over {n \,-\, 2}} \  dy \ \,+\, \\[0.15in]
 & \ &\ \ \ \ \ \ \ \  \ \ \ \ \ \ \ \ \ \  \ \ \ \ \ \ \ \ \ \  \ \ \ \  \,+\, \ \int_{B_o \,(\,\rho_2\,)} O\, (\,|\,y|^{\,\ell \,+\, 1}\,) \cdot [\,{\bf A}\,(y)\,]^{{2n}\over {n \,-\, 2}} \  dy \,. \ \ \ \ \ \ \ \ \
\end{eqnarray*}
The first term in the right hand side of the last equation in (A.6.46) can be expanded by using (A.6.33). As for the second term, it can be estimated as in (A.6.41)\, (replacing $\,\lambda_i\, R_i\,$ by $\,\rho_2$)\, showing that the term is of order $O_{\lambda_i}(\ell \,+\, 1)\,.$\, Hence we obtain the following.

 \vspace*{0.2in}

{\bf Lemma A.6.47.} \ \ {\it Under the conditions in\,} (2.63)\,,\,  (A.6.23),\, {\it $\ell \in [\,2\,,\ n \,-\, 2\,]\,,\,$   suppose that\,,\, for} $\,i \,\gg \,1\,,\,$ $\,\xi_i  \, =\,   \lambda_i^{1 \,+\, {\eta_{\,o}}}\cdot \!\vec{\,X}\,,\,$ {\it where\,} $\vec{\,X} \,\in \,\R^n\,$ {\it is a fixed vector\,.\, Then we have}
\begin{eqnarray*}
(A.6.48) & \ &  \,-\, \int_{B_o \,(\,\rho_2)} r\, {{\partial K}\over {\partial r}} \, {\bf A}^{{2n}\over {n \,-\, 2}} \ dy  \,=\, \ \lambda_i^\ell  \cdot   {{ \ell }\over {c_n}}  \, \int_{\R^n} {\bf P}_\ell \,(y) \left( {1\over {1 \,+\, |\,y|^2}}\!\right)^n \! dy \, \,+\, \\[0.25in]
& \ &   \!\!\!\!\!\!\!\!\!\!\!\!\!\!\!\!\!\!\!\!\!\!\!\!\!\!\!\!         +\  \int_{\R^n} \left[\ \lambda_i^{\ell \,+\, {\eta_{\,o}}}\cdot \Xi_1^{\bf P} \, (\vec{\,X}\,, \ y) \ \,+\, \cdot \cdot \cdot \ \,+\, \ \lambda_i^{\ell \,+\, (\ell \,-\, 1)\cdot \, {\eta_{\,o}}}\cdot \Xi_{\ell \,-\, 1}^{\bf P} \, (\vec{\,X}\,, \ y) \right]\cdot \left( {1\over {1 \,+\, |\,y|^2}}\right)^n \!dy \ \,+\, \\[0.1in]
& \ & \ \ \ \ \ \ \ \ \ \ \ \ \ \ \  \,+\, \ \lambda_i^{\ell \,+\, \ell \cdot \,{\eta_{\,o}}}\cdot {\bf P}_\ell\, (\vec{\,X})  \int_{\R^n}  \left( {1\over {1 \,+\, |\,y|^2}}\right)^n \! dy \ \,+\, \  O_{\lambda_i}\, (\ell \,+\, 1)\,.
\end{eqnarray*}
{\it Here $\ \Xi_h^{\bf P}\,$ is defined as in\,} (A.6.29) {\it and}\, (A.6.30) {\it by replacing $\,Q_\ell\,$ by\,} ${\bf P}_\ell\,.$

\vspace*{0.2in}

{\bf \S\,A.6.\,g.} \ \ {\it  Estimate on the term involving \,${\bf B}$\,.\,} \ \
We first shift the center in the term
\begin{eqnarray*}
(A.6.49)\!\!\!\!\!\!& \ & \int_{B_o\, (\,\rho_2)} \bigg\vert\, r {{\partial K}\over {\partial r}} \bigg\vert \cdot   \bigg\vert\,  {\bf B}\, (y) \bigg\vert^{{2n}\over {n \,-\, 2}} \ d y\ \le \ C\, \int_{B_o\, (\,\rho_2)} |\,y|^{\,\ell} \cdot   \bigg\vert\,  {\bf B}\, (y) \bigg\vert^{{2n}\over {n \,-\, 2}} \ d y\\[0.15in]
& \le & C\,\int_{B_{\xi_i}\, (\,\rho_2)} |\,y|^{\,\ell}  \cdot \ \bigg\vert\,  {\bf B}\, (y) \bigg\vert^{{2n}\over {n \,-\, 2}} \ d y \ \,+\, \ C_1\, |\,\xi_i| \cdot \lambda^{{n\over 2} \cdot {{2n}\over {n \,-\, 2}}  }  \\[0.15in]
& \ & \ \ \ \  \ \ \ [\ {\mbox{cf. \ \ (A.6.22)}}\,,\, \ \ {\mbox{observe \ \ that \ \ }} \lambda_i^{\,\ell} \cdot |\ \Gamma_{\bf p} \, ({\cal Y})| \le c \ \ {\mbox{in}} \ \ B_o\, (\,\rho_2)\,] \\[0.15in]
& \le & C\,\int_{B_{\xi_i}\, (\,\rho_2)} \left[ \,|\,y \,-\, \xi_i|^{\,\ell} \,+\, C_1\, |\,\xi_i| \cdot  |\,y \,-\, \xi_i|^{\,\ell \,-\, 1 } \,+\, \, \cdot \cdot \cdot \ \,+\, \  |\,\xi_i|^{\,\ell} \,\right] \cdot   \bigg\vert\,  {\bf B}\, (y) \bigg\vert^{{2n}\over {n \,-\, 2}} \ d y \\[0.1in]
& \ & \ \ \ \ \ \ \ \ \ \ \ \ \ \ \ \ \ \ \ \ \ \
  \,+\, \ \,o_{\lambda_i} \left(  {{n^2}\over {n \,-\, 2}}  \,+\,1 \right)\ \ \ \ \ \ \ \       [\,{\mbox{recall \ \ that}} \ \ |\,\xi_i| \,=\, o\, (\lambda_i)\,] \ . \ \ \ \ \ \ \ \ \ \ \ \ \
 \end{eqnarray*}
 In the above, we apply the triangle inequality  and the binomial expansion as in
 $$\ \ \ \ \  \ \ \ \ \
 |\,y|^\ell \ \le \ (  |\,y \,-\, \xi_i|  \,+\,  |\,\xi_i|\,)^\ell \ \ \ \ \ \ \  \ \ \ \ \ \ \ \  \ \ \ (\,\ell \ \ {\mbox{a  \ \ positive \ \ integer}})\,.
 $$

Introduce the change of variables
$$
|\, y \,-\, \xi_i|  \, \,=\, \,  \rho  \, \,=\, \,  \lambda_i\, \tan \, \theta \ \ \ \ \ \Longrightarrow \ \ \  \ \ \ \tan \, \theta  \, \,=\, \,  {{|\, y \,-\, \xi_i|}\over { \lambda_i  }}\ \ \ \Longrightarrow \ \ \ |\,{\cal Y}| \ \,=\, \ \tan \,\theta\,. \leqno (A.6.50)
$$
Recall that in the expression (A.6.24) for $\,{\bf B}\,$,
$$|\,  {\bf B}\, (y) |^{\,{{2n}\over {n \,-\, 2}}} \ \,=\, \ \lambda_i^{\,( \ell \,+\, 1) \cdot \,{{2n}\over {n \,-\, 2}}  } * |\,\Gamma_{\, \le \,\ell} \ ({\cal Y}) |^{\,{{2n}\over {n \,-\, 2}}  } \cdot \left[\  {1\over {\lambda_i\,(1 \,+\, \tan^2 \,\theta\,)}} \right]^{{n^2}\over {n \,-\, 2}}\  \ \ \ \left({\cal Y} \ \,=\, \ {{y \,-\, \xi_i}\over {\lambda_i}}\right)\,.
$$
Moreover,\\[0.1in]
(A.6.51)

\vspace*{-0.25in}

$$
|\,\Gamma_{\, \le \,\ell} \ ({\cal Y}) | \ \le \ C\, \left[ \, {\cal R}^2 \,+\, \cdot \cdot \cdot \,+\, {\cal R}^{\,\ell}\ \right]   \  \Longrightarrow   \cases{\ |\,\Gamma_{\, \le \,\ell} \ ({\cal Y}) | \ \le\ C_1 & for $ \    {\cal R} \,=\, |\, {\cal Y}| \,\le\, 1$\,; \cr \cr
              \ |\,\Gamma_{\, \le \,\ell} \ ({\cal Y}) | \ \le \ C_2\, {\cal R}^{\,\ell} & for   $ \  |\, {\cal Y}| \,\ge \,1$. \cr}
$$

For $\,0 \,\le\, l \,\le\, \ell\,,\,$ we have

 \begin{eqnarray*}
\!\!& \ &\!\!\!\!\!\!\!\!\!\!\!\!\!\!\!(A.6.52) \ \  \ \ \ \ \int_{B_{\xi_i}  (\,\rho_2)}  |\,y \,-\, \xi_i|^{\,l}\cdot |\, {\bf B}\, (y) |^{{2n}\over {n \,-\, 2}} \ \,d y \\[0.15in]
&  \,=\, & \left( \ \int_{B_{\xi_i}  (\,\lambda_i)} \,+\, \ \int_{B_{\xi_i}  (\,\rho_2)\,\setminus B_{\xi_i}\, (\,\lambda_i)} \right)\ |\,y \,-\, \xi_i|^{\,l}\cdot |\, {\bf B}\, (y) |^{{2n}\over {n \,-\, 2}} \ d y \\[0.2in]
& \le & C_3 \int_0^{\lambda_i}  \, r^{\,l    } \cdot \left[ \lambda_i^{\, \ell \,+\, 1 \,-\, {n\over 2} } \right]^{{2n}\over {n \,-\, 2}}\ [\, r^{ \,n \,-\, 1   }\,d\,r\,] \\
& \ & \ \ \ \ \ \ \ \ \ \ \ \ \ \ \ \ \ \ \ \ \ \ \ \  \ \ \ \ \ \ \ \ \ \ \ \   \ \ \  [\,{\mbox{where}} \ \ r \,=\, |\,y \,-\, \xi_i|\,; \ \ {\mbox{using \ first \ half \ in (A.6.40)}}\,] \\[0.15in]
& \ & \!\!\!\!+ \  C_2\, \int_{\,\arctan\,1}^{\,\arctan \, \left( {{\rho_2}\over {\lambda_i}} \right) } \ \  [\,  \lambda_i \tan\,\theta\,]^{\,l}\cdot \lambda_i^{\, (\,\ell \,+\, 1)\cdot \, {{2\,n\, }\over {n \,-\, 2}}  } \cdot  (\tan\,\theta)^{  {{2\,n\,\ell}\over {n \,-\, 2}}  }  \cdot \left[\ {1\over {\lambda_i\,(1 \,+\, \tan^2 \theta)}} \right]^{{n^2}\over {n \,-\, 2}} \times \\[0.15in]
& \ & \hspace*{3.4in} \ \ \ \ \ \times \,\{\lambda_i^n\,[\,\tan \theta]^{n \,-\, 1}\,\,\sec^2 \, \theta\} \ d \theta  \\[0.2in]
& \le &  O_{\lambda_i}  \left( l\, \,+\, \,\ell \cdot {{2n}\over {n \,-\, 2}}    \right) \ \,+\, \\[0.1in]
 & \ &   \!\!\!\!\!\!\!\!\!\!+ \ O_{\lambda_i}  \left(   \,l\ \,+\,  \ [\,\ell \,+\, 1] \cdot {{2\,n  }\over {n \,-\, 2}}  \, \,+\, \, n \, \,-\, \,  {{n^2}\over {n \,-\, 2}}   \right) \!\cdot  \!\int_0^{\arctan \, \left( {{\rho_2}\over {\lambda_i}} \right) } \   {{[\, \cos \, \theta \,]^{\, {{2n^2}\over {n \,-\, 2}} }  }\over { [\, \cos \, \theta \,]^{ \,l \,+\,  {{2\,n\,\ell}\over {n \,-\, 2}}\, +\, (n \,-\, 1) \,+\, 2  }   }} \ \ d\,\theta    \\[0.2in]
& \le & O_{\lambda_i}  \!\left( l\, \,+\,  {{2n\,\ell}\over {n \,-\, 2}}    \right) \, \,+\,  \,  O_{\lambda_i}  \left( l\, \,+\,  {{2n\,\ell}\over {n \,-\, 2}}    \right)  \!\cdot\! \int_0^{\arctan \, \left( {{\rho_2}\over {\lambda_i}} \right) } \    {{[\, \cos \, \theta \,]^{\, {{2n^2}\over {n \,-\, 2}} }  }\over { [\, \cos \, \theta \,]^{ \,l \,+\, \left( {{2\,n\,\ell}\over {n \,-\, 2}} \,-\, 4\right)\,+\, 4 \,+\, (n \,+\, 1) }   }} \ \ d\,\theta    \\[0.15in]
\\[0.2in]
& \le & O_{\lambda_i} \! \left( l\, \,+\,  {{2n\,\ell}\over {n \,-\, 2}}    \right)  \, \,+\,  \,  O_{\lambda_i} \! \left( l\, \,+\,  4  \right)  \! \cdot \!\int_0^{\arctan \, \left( {{\rho_2}\over {\lambda_i}} \right) }   {{[\, \cos \, \theta \,]^{\, {{2n^2}\over {n \,-\, 2}} }  }\over { [\, \cos \, \theta \,]^{ \,l \,+\,  n   \,+ \,5  }   }} \ \ d\,\theta   \  \ \  \       [\, {\mbox{see  \  \S\,A.6\,g.2\,}}] \\[0.15in]
& \,=\, & \ o_{\lambda_i}  \left( l\, \,+\,  4 \right) \ \,+\, \ o_{\lambda_i}  \left( l\, \,+\,  4 \right)  \hspace*{2.9in}({\mbox{as \ \ }} \ell \,\ge\, 2)\\[0.1in]
  & \ & \ \ \ \ \  \ \ \ \ \ \ \ \ \ \ \ \left[\, {\mbox{note \ \ that}} \ \ l \,+\, n \,+\, 5\ \le\  (n \,-\, 2) \,+\, n \,+\, 6 \,=\, 2n \,+\, 3 \,\le\,{{2\,n^2  }\over {n \,-\, 2}}   \ \right]\,. \\
\end{eqnarray*}

Using (A.6.49), (A.6.52)   and $\,|\,\xi_i| \,=\, O \, (\lambda_i)$\,,\,  we obtain the following.\\[0.2in]
%
%
%
{\bf Lemma A.6.53.} \ \ {\it Let $\,{\bf B}\,$ be given as in\,} (A.6.24)\,,\, {\it $\,|\,\xi_i| \,=\, O\, (\lambda_i)$\,,\, and $K$ satifies\,} (1.8)\,.\, {\it Then for $\,2 \,\le\, \ell \,\le\, n \,-\, 2\,$,\, we have}
$$
\int_{B_o\, (\,\rho_2)} \bigg\vert\  r {{\partial K}\over {\partial r}} \bigg\vert \cdot \bigg\vert\,  {\bf B}\, (y) \bigg\vert^{{2n}\over {n \,-\, 2}} \ d y \ \,=\, \ o_{\lambda_i}  \left( \ell \, \,+\,  4 \right)\ . \leqno (A.6.54)
$$

\vspace*{0.2in}

{\bf \S\,A.6.\,h. \  } {\it Estimate on the term involving \,${\bf C}$\,.\,} \ \ As
$$
\ell \,\le\, n \,-\, 2\  \ \ \Longrightarrow \ \ \ell \,+\, 1 \,-\, {n\over 2} \ \le \ {{n \,-\, 2}\over 2}\,. \leqno (A.6.55)
$$
From (A.6.25) \,and $\,|\,\langle y\,, \ \btd \,K\rangle | \,\,\le\, c\,$ in $\,B_o\, (\rho_o)\,,\,$ we obtain
$$ \int_{B_o\, (\,\rho_2)} \bigg\vert\, r {{\partial K}\over {\partial r}} \bigg\vert \cdot \bigg\vert\,  {\bf C}\, (y) \bigg\vert^{{2n}\over {n \,-\, 2}} \ d y   \ \,=\, \   O_{\lambda_i}\! \left(\,\left[\ \ell  \,-\, {{n-2}\over 2} \right]\,\cdot {{2n}\over {n \,-\, 2}}\,\right)\,. \leqno (A.6.56)
$$
Observe that
\begin{eqnarray*}
 & \ &  \left( \ell \,+\, 1 \,-\, {n\over 2} \right)\,\cdot {{2n}\over {n \,-\, 2}}\  > \ \ell\\[0.1in]
  \Longleftrightarrow & \ &  \!\!\!\!\!\!\!\!\ell \ > \ {{n\,(n \,-\, 2)}\over {n \,+\, 2}} \ \ \Longrightarrow \ \ \ell \,=\, (n \,-\, 2\,) \ \ \& \ \ n \ge 4\,; \ \ {\mbox{or}} \ \ \ \ell \,=\, (n \,-\, 3\,) \ \ \& \ \ n > 6\ .
\end{eqnarray*}
Moreover,
$$
 \ell \,=\, n \,-\, 2 \ \  \Longrightarrow \ \ \left( \ell \,+\, 1 \,-\, {n\over 2} \right)\,\cdot {{2n}\over {n \,-\, 2}}\  \,=\, \ n\,,\leqno (A.6.57)
$$
and when $\,\ell \,= \,(n \,-\, 3\,)$ and $\,n \,>\, 6\,,\,$
$$
\left[\, (n \,-\, 3) \,+\, 1 \,-\, {n\over 2} \right]\,\cdot {{2n}\over {n \,-\, 2}} \ \,-\, \ (n \,-\, 3) \ \,=\, \ {{n \,-\, 6}\over {n \,-\, 2}} \ > \ 0\ . \leqno (A.6.58)
$$

\newpage

{\bf \S\,A.6.\,i. \  } {\it Remaining estimates.} \ \ Using $\,|\,\langle\, y\,,\ \btd \, K \rangle\,| \ \le \ C\, |\,y|^\ell\,$ and $\,|\,\xi_i| \, \,=\, \, o\, (\lambda_i)$\, [\,actually we only need $\,|\,\xi_i| \, \,=\, \, O\, (\lambda_i)$]\,,\, as in (A.6.33), we have
$$
\ \ \varepsilon\,\int_{B_o\, (\,\rho_2)} \bigg\vert\, r {{\partial K}\over {\partial r}} \bigg\vert \cdot \,[\,  {\bf A}\, (y) ]^{\,{{2n}\over {n \,-\, 2}} }\ d y \ \,=\, \ \varepsilon\cdot O_{\lambda_i}   (\ell)  \ \ \   \mfor 2\ \le   \ell \ \le\  n \,-\, 2\ . \leqno (A.6.59)
$$

\vspace*{0.15in}

{\bf \S\,A.6.\,i.1. \  }   {\it Estimate on the outside.} \ \ Similar to (A.6.32),\, for $\,l\,\le\, n \,-\, 1$,\, we have
\begin{eqnarray*}
 (A.6.60) \ \ & \ & \int_{\R^n \setminus B_o \, (\,\rho_2)} r^{\,l} \cdot \left( {{\lambda_i}\over {\lambda_i^{\,2} \,+\, r^{\,2}}}  \right)^n \ dy \ \le \ \int_{\R^n \setminus B_o \, (\,\rho_2)} r^{\,l} \cdot \left( {{\lambda_i}\over {  r^{\,2}}}  \right)^n \ dy \\[0.15in]
& \,=\, & \lambda_i^n \cdot  \int_{\R^n \setminus B_o \, (\,\rho_2)} {{1}\over {  r^{\,2n -\, l}}}    \ dy \ \le \ C\, \lambda_i^n \cdot  \int_{ \rho_2}^\infty {{1}\over {  r^{\,2n -\, l \,-\, (n \,-\, 1)}}}    \ dr
\  \,=\, \  O\, (\lambda_i^n) \ . \ \ \ \ \ \ \ \ \ \ \ \ \ \ \ \ \ \
\end{eqnarray*}

\vspace*{0.2in}

{\bf \S\,A.6.\,i.2. \  }   {\it The angle.} \ \ Let
$$
\theta_{\lambda_i} := \arctan \   {{\rho_2}\over {\lambda_i}} \leqno (A.6.61)
$$
$$
\Longrightarrow \ \  \cos\ \theta_{\lambda_i}  \  \,=\, \ {1\over {\sec\, \theta_{\lambda_i} }} \  \,=\, \ {1\over {\sqrt{\sec^{\,2}\, \theta_{\lambda_i}  \,} }}\  \,=\, \ {1\over {\sqrt{1 \,+\, \tan^{\,2}\, \theta_{\lambda_i}  \,} }} \  \,=\, \ {1\over {\sqrt{ 1 \,+\, {{{\rho_2}^{\,2} }\over {{\lambda_i} ^{\,2}}} } }} \  \,=\, \ O \, ({\lambda_i} )\,.
$$

\vspace*{0.3in}

{\bf \S\,A.6.\,j\,.\,} \ \ {\it Limitation on flexibility\,.\,}\\[0.15in]
  {\bf Theorem A.6.62.} \ \ {\it Assume the conditions in Main Theorem\,} 1.14  {\it and restrict $\ell$ to be either $\,n \,-\, 2\,$ or $\,n \,- \,3\,$\,.\, Moreover, assume the following.}  \\[0.1in]
  {\bf (i)} \ \ \, {\it When  $\,\ell \,=\, n \,-\, 2\,,\,$ we take it that ${\bf S}$ is the only blow\,-\,up point\,.}\\[0.1in]
   {\bf (ii)} \  \, {\it When   $\,\ell \,=\, n \,-\, 3\,,\,$ we take it that $\,n > 6$\,,\, and it is possible to have other  }\\[0.075in]
   \hspace*{0.47in}{\it blow\,-\,up point}({\it  s }$\!$)\,.\\[0.1in]
   {\it Let $\ \Xi_h^{\bf P}\,$ be given in\,} (A.6.29) {\it and\,} (A.6.30) {\it with $\,Q_\ell\,$ replaced by $\,{\bf P}_\ell$\,.\,    Assume also that}
$$
  \xi_i \ \,=\, \   \lambda_i^{1 \,+\, {\eta_{\,o}}} \cdot \vec{\,X} \ \ \ \ {\it{for}} \ \  i \gg 1 \ \ \,{\it{and \ \ a \ \ fixed \ \ }} \  \vec{\,X} \,\in\, \R^n\,. \leqno (A.6.63)
$$

\vspace*{-0.25in}

$$
{\it Here} \ \ \ \ \ \ \ \ \ \ \ \ \ \ \ \ \eta_o \ \le \, \cases{\   {2\over {3n \,-\, 2}} & \ \ \ \ {\it when}  $ \ \   \ell \,=\, n \,-\, 2$\,; \cr \cr
              \  {{n \,-\, 6}\over {(n \,-\, 3)(3n \,-\, 2)}}   & \ \ \ \ {\it when}   $ \ \ \ell \,=\, n \,-\, \,3\ \ {\it{and}}  \ \ n > 6$\,. \ \ \ \ \ \ \ \  \ \ \ \ \ \ \  }
$$
{\it Then it is necessary that }\, ${\bf P}_\ell\, (\vec{\,X}) \,=\, 0\,,\,$ {\it and}
   \begin{eqnarray*}
 \int_{\R^n}  \Xi_1^{\bf P}\, \left( \vec{\,X}\,, \ y \right) \,  \left( {1\over {1 \,+\, |\,y|^2}} \right)^n dy \ \,=\, \, \cdot \cdot \cdot \, \,=\, \,    \int_{\R^n}  \Xi^{\bf P}_{\ell \,-\, 1}\, \left( \vec{\,X}\,, \ y \right) \,  \left( {1\over {1 \,+\, |\,y|^2}} \right)^n \!dy   \ \,=\, \ 0\,.
  \end{eqnarray*}

 \vspace*{0.1in}

 {\bf Proof.} \ \ Using estimate (A.6.43), (A.6.44), (A.6.48), (A.6.54) (A.6.56)\ --\ (A.6.58), we obtain

 \vspace*{-0.25in}

 \begin{eqnarray*}
(A.6.64)  \!\!\!\!\!\!&\ & \,-\, \!\int_{B_{\,o} (\,\rho_2) } \langle\,y\,, \    \btd \, K\,  (y)\,\rangle \,\, [\,v_i\, (y)\,]^{{2n}\over {n \,-\, 2}}   \ dy \ \,=\, \  \lambda_i^\ell \cdot {\ell\over {{\tilde c}_n}}  \, \int_{\R^n}   {\bf P}_\ell \, (y) \,  \left( {{1}\over {1 \,+\, |\,y|^2}} \right)^n\!\!\!dy \, \,+\, \\[0.1in]
   & \ &    \!\!\!\!\!\!\!\!\!\!\!\!\!\!\!\!     +\  \int_{\R^n} \left[\ \lambda_i^{\ell \,+\, {\eta_{\,o}}}\cdot \Xi_1^{\bf P} \, (\vec{\,X}\,, \ y) \ \,+\, \cdot \cdot \cdot \ \,+\, \ \lambda_i^{\ell \,+\, (\ell \,-\, 1)\cdot \, {\eta_{\,o}}}\cdot \Xi_{\ell \,-\, 1}^{\bf P} \, (\vec{\,X}\,, \ y) \right]\cdot \left( {1\over {1 \,+\, |\,y|^2}}\right)^n \!dy \ \,+\, \\[0.1in]
 & \ &         +\  \ \lambda_i^{\ell \,+\, \ell \cdot \,{\eta_{\,o}}}\cdot {\bf P}_\ell\, (\vec{\,X})  \int_{\R^n}  \left( {1\over {1 \,+\, |\,y|^2}}\right)^n \! dy \ \,+\, \  O_{\lambda_i}\, (\ell \,+\, 1) \ +\\
\end{eqnarray*}
$$
\hspace*{0.8in} \  \,+\, \ \varepsilon \cdot O_{\lambda_i}\,(\ell) \ \ \,+\, \ \  {{C_n}\over {\varepsilon^{{2n}\over {n \,-\, 2}} }} \,\cdot \,  \cases{\  O_{\lambda_i}(\,\ell \,+\, 2) & for $ \ \   \ell \,=\, n \,-\, 2$\,; \cr \cr
              \ O_{\lambda_i}\left(\ell +{{n \,-\, 6}\over {n \,-\, 2}}\right)  & for   $ \ \ \ell \,=\, n \,-\, 3\,, \ \ n > 6$. \cr}
$$

 \vspace*{-0.05in}

Referring to Remark A.6.45\,:
\begin{eqnarray*}
(A.6.65)\ \ \  {\mbox{when}} \ \  \ell & \,=\, & n \,-\, 2\,, \ \   \ell \cdot \eta_o \ < \ 2\cdot {{n \,-\, 2}\over {3n \,-\, 2}} \ \ \Longleftrightarrow \ \ \ \eta_o \ < \ {2\over {3n \,-\, 2}} \ \ \\[0.1in]
& \ & \hspace*{2.1in}\left( a \,=\, 2\,, \ \ t \,=\, {{2\,(n \,-\, 2)}\over {3n \,-\, 2}} \,<\, 1\right)\,; \\[0.1in]
(A.6.66)\ \   \  {\mbox{when}} \ \ \ell & \,=\, & n \,-\, 3\,, \ \    \ell \cdot \eta_o \ < \ {{n \,-\, 6}\over {n \,-\, 2}} \cdot {{n \,-\, 2}\over {3n \,-\, 2}} \ \ \Longleftrightarrow \  \eta_o \, < \, {{n \,-\, 6}\over {(n \,-\, 3)\,(3n \,-\, 2)}}\ \ \\[0.1in]
& \ & \hspace*{1.4in}\left( a \,=\, {{n \,-\, 6}\over {n \,-\, 2}} \,, \ \ t \,=\, {{n \,-\, 2}\over {3n \,-\, 2}}\cdot  {{n \,-\, 6}\over {n \,-\, 2}} \,<\, 1\right)\,. \ \ \ \ \ \ \
   \end{eqnarray*}
  Combining with (A.6.35), we  come to the  conclusion of the theorem.\qed

\newpage

{\large \bf \S\,A.7. \ \ $\ell \,=\,n \,-\, 2$\, and multiple simple blow\,-\,up points\ --}\\[0.1in]
 \hspace*{0.75in}{\large \bf off\,-\,center cancelation.}\\[0.15in]
We present the consideration on   global cancelation/balance with finite number of  blow\,-\,up points, say, at
 $$
 \ \ \ \ \ \  \  \ \ \ \ \ \    {{\hat{Y}}_o} \,=\, 0\,,\ \ {{\hat{Y}}_1}\,, \ \cdot \cdot \cdot\,, \ {{\hat{Y}}_k} \ \ \ \ \ \  \ \ \ \ \ \ \  \  \ \ \ \ \ \  \ \ \ \ \ \  \  \ \ \ \ \ \  \  (k \ge 1\,)\,. \leqno (A.7.1)
 $$
[\,Cf. (2.7)\,.\,]\, {\it  Throughout this section\,} (\S\,A.7) {\it we assume the general conditions\,} (1.6)\,,\, (1.25)\,,\, (1.26)\,,\, {\it  $\,n \,>\, 6\,,\,$ and } \\[0.05in]
(A.7.2)
 $$
 {{\hat{Y}}_j} \ \  {\it{is \   a \  simple \   blow\!-\!up \   point \ and}}  \sum_{2\, \le\, l\,<\, n \,-\, 2} \!\!\!|\,\btd^{(\,l)}  K\, ( {{\hat{Y}}_j})| \,=\, 0 \ \ \ {\it{for}}  \ \,0\,\le \, j \, \le \, k.
 $$
 Cf. (1.9) and the remark preceding it. For the simple blow\,-\,up point at $0$, we keep the notations on $\,\xi_i\,$ and $\,\lambda_i\,$ as introduced in (1.10) and (1.11).
 Likewise (cf. Proposition 2.24 and \S\,2\,f) we set
\begin{eqnarray*}
 (A.7.3) \ \ \    \ \ \ \ \ \ \ \xi_{_{m_i}} \ : \ \ \ v_i\,  (\xi_{_{m_i}}  )  &  \,=\, & \max \ \left\{ \, v_i\, (y) \ | \ \ y \, \in\,  \overline{B_{ {\hat Y}_m }   (\rho_3)\!} \   \right\}\,, \ \ \  \ \ \ \xi_{m_i} \ \to \ {\hat Y}_m\,, \ \ \ \ \ \\[0.1in]
(A.7.4) \    \ \ \ \ \ \ \  \    \ \ \ \ \ \ \ \    \ \ \ \ \ \ \ \  \lambda_{m_i}  & := &  {1\over{ [\, v_i\, (\xi_{_{m_i}} )\, ]^{\,{2\over {n \,-\, 2}} }}}\    \ \ \ \ \ \ \ \ \ \ \ \ \ \ \ \ \  \mfor 1 \,\le\, m \,\le \,k\,. \ \ \ \
 \end{eqnarray*}
 Here $\,\rho_3 > 0\,$ is a constant (small enough) so that Proposition 2.3 and    Proposition 2.24 hold   after a translation to each individual blow\,-\,up point,  and
 $$
 \,B_{ {\hat Y}_m }   (\rho_3)  \,\cap\, B_{ {\hat Y}_j }   (\rho_3) \ \,=\, \ \emptyset\ \ \ \ \mfor\  j \,\not=\, m\,.
 $$
Via Proposition 2.24 and the Harnack inequality \cite{Leung-Supported} , $$
 {1\over C} \cdot \lambda_{m_i}^{\,{{n \,-\, 2}\over 2}}  \  \le \   \min_{|\,y -\, {\hat Y}_m | \,=\, \rho_3}  \!\!  v_i\,(y) \  \le \  \max_{|\,y -\, {\hat Y}_m | \,=\, \rho_3}    \!\!v_i\, (y)\  \le \ C\,\lambda_{m_i}^{\,{{n \,-\, 2}\over 2}}   \leqno (A.7.5)
 $$
 for \,$1 \,\le\, m\,\le \, k $\, and $\ i \,\gg\, 1\,.\,$  As there is  no blow\,-\,up point which appears at the north pole [\,cf. (2.11)\,]\,,\,
apply  the Harnack inequality \cite{Leung-Supported} again on
 $$\ \displaystyle{
  S^n \,\setminus \ {\dot{\cal P}}^{-1} \left( \ \bigcup_{j \,=\, 0}^k \ B_{ {\hat Y}_j}\, (\rho_3) \! \right)
  }\ $$ and
obtain

$$ \ \ \ \ \
{1\over C}  \ \le  \ {{\lambda_{m_i} }\over {\lambda_i}} \ \le \ C \ \ \ \ \ \mfor 1 \,\le\, m \,\le\, k \ \ \ {\mbox{and}} \ \ i \gg 1\,. \leqno (A.7.6)
$$
 It follows that, modulo a subsequence,
 $$
 S_m \ := \ \lim_{i \to \infty} \ {{\lambda_{m_i} }\over {\lambda_i}}  \ \ \ \ {\mbox{is \ \ well\,-\,defined}}\mfor 1 \,\le \,m \,\le\, k\,. \leqno (A.7.7)
 $$
Using the global formula (A.6.6),  together with (A.7.5) and (A.7.6)\,,\, we have
$$
 \int_{\Omega}  \   \langle \,y\,,\, \ \btd_y \ K \, (y)\,\rangle \, \left[\,v_i\,(y)\right]^{{2n}\over {n \,-\, 2}}\ dy   \ \,=\, \  O_{\lambda_i}\,(\,n )\,, \leqno (A.7.8)
 $$
 where $\,\Omega \ \,=\, \ B_o\, (\,\rho_3) \  \cup \   B_{\,{\hat{Y}}_1}\, (\,\rho_3) \ \cup\, \cdot \cdot \cdot\ \cup  \ B_{\,{\hat{Y}}_k}\, (\,\rho_3)\ .$

\vspace*{0.2in}

{\bf \S\,A.7.\,a.} \  \ {\it Off\,-\,origin blow\,-\,up point.} \
Consider the simple blow\,-\,up at \,${{\hat Y}_m}\,$,\, where \,$0 \,<\, m\, \le \,k\,.\,$ Via Taylor expansion,

\vspace*{-0.3in}

\begin{eqnarray*}
(A.7.9) \ \ \ \    {\tilde c}_n \cdot K (y) & \,=\, & ({\tilde c}_n \cdot K) \,({{\hat Y}_m}) \ \,+\, \ \sum_{|\,\alpha| \,= \,{n \,-\, 2}}  {1\over {\alpha\,!}} \ [\, D^{({n \,-\, 2})}_{\,\alpha}\, (\,{\tilde c}_n K)\,({\hat Y}_m\,) \,]\cdot  (y \,-\, {{\hat Y}_m})^\alpha \  \\[0.15in]
  & \ & \ \ \ \ \ \ \ \ \ \ \ \ \ \ \ \ \ \   \ \ \ \      +\  O \,(|\,y \,-\, {{\hat Y}_m}|^{\,n \,-\, 1}\,) \ \   \mfor \ \ \,|\,y \,-\, {{\hat Y}_m}|\, \le \,\rho_3\ \ \ \ \ \ \ \ \ \ \ \ \\[0.15in] \ \ \ \ \ \ \ \ \ \ \ \
  & \,=\, & ({\tilde c}_n \cdot K )\,({{\hat Y}_m}) \ \,+\, \  [\,-\,{\bf P}_{n \,-\, 2\,,\,m}\, (y)\,]  \   +\  O \,(|\,y \,-\, {{\hat Y}_m}|^{\,n \,-\, 1}\,)\,.
\end{eqnarray*}
Here $\,{\bf P}_{n \,-\, 2\,,\,m}\, (y)\,$ is defined by the equation above [\,see also (A.7.15)\,]. For the sake of continuity, we keep the sign convention on `{\bf P}'\,,\, which we use in this article.
Assume that
$$
\Delta_o^{\!(h_{n \,-\, 2})} \  {\bf P}_{n \,-\, 2\,,\ m} \,(\,{\hat Y}_m) \ \equiv \ 0 \ \ \ \ \ \   \mfor \ \ 0 \,\le\, m \,\le\, k\,, \leqno (A.7.10)
$$
and all the corresponding conditions as in Main Theorem (1.17) hold for each individual simple blow\,-\,up point, except
  $\,({\tilde c}_n \cdot K)\,({\hat Y}_m)\,$ may not be $\, n\,(n \,-\, 2)\,.\,$ Thus the estimate contains a scaling factor
\begin{eqnarray*}
(A.7.11) \ \ \ v_i\, (y) & \,=\, & \left[ \,{{n\, (n \,-\, 2)}\over {({\tilde c}_n\, K)\, (\,{\hat Y}_m)\,}} \right]^{\,{{n \,-\, 2}\over 4}} \cdot \left( {{\lambda_{m_i} }\over {\lambda_{m_i}^2 \,+\, |\, y \,-\, \xi_{m_i}|^2   }} \right)^{\!{{n \,-\, 2}\over 2}} \ \ \ \,+\, \\[0.2in]
& \ & \!\!\!\!\!\!\! \!\!\!\!\!\!\!\!\!\!\!\!\!\!\!\!\!\!\!\!\!\!\!\!\!\!\!\!\!\!\!\!\!\!+  \ \left\{ \ {\mbox{expressions \ similar \ to \ those \ in \   (A.6.24)\ and \  (A.6.25)}}\,\right\}  \mfor \,y \in B_{ {\hat Y}_m}\, (\,\rho_3)\,.\,
\end{eqnarray*}
[\,Here $\,\rho_3\,$ is made smaller if necessary\,.\,]\,\bk
We find the first derivative by using change of variables $\,y\,=\, z \,+\, {\hat Y}_m\,,\,$
$$
{\tilde c}_n \cdot \langle \,y\,,\, \ \btd_y \ K \, (y)\,\rangle
=  {\tilde c}_n \cdot\langle \,( \,y \,-\, {{\hat Y}_m}) \,,\, \ \btd_y \ K \, (y)\,\rangle   \ \,+\, \   {\tilde c}_n \cdot\langle\, {{\hat Y}_m}\,,\, \ \btd_y \ K \, (y)\, \rangle \,,\leqno (A.7.12)
$$

\vspace*{-0.35in}

\begin{eqnarray*}
(A.7.13) \ \ \  \ \ \ \ \ \  \ \ \  \langle \,( y \,-\, {{\hat Y}_m}) \,,\, \ \btd_y \ K \,(y)\,\rangle  & \,=\, &  \big\langle \,z\,,\, \    (\btd_y \ K)\,|_{\,y \,= \,z \,+\, {\hat Y}_m}   \big\rangle\\[0.1in]
  &= &\big\langle \,z\,,\, \    \btd_z \ K_{\to} \,(z)\,\big\rangle\,|_{\,z \,=\, ( \,y \,-\, {{\hat Y}_m})}\ . \ \ \  \ \ \ \ \ \  \ \ \ \ \ \  \ \ \
\end{eqnarray*}

\vspace*{-0.4in}

$$
\ \ \ \ K_{\to} \,(z) \ \,=\, \ K \, (z\, \,+\, \,{\hat Y}_m)\,.\leqno{\mbox{Here}}
$$
Consider the second expression in the right hand side of (A.7.12)\,.\,  As in (A.7.13), we have
\begin{eqnarray*}
(A.7.14) \  \    {\tilde c}_n \cdot\langle\, {{\hat Y}_m}\,,\, \ \btd_y \ K \, (y)\, \rangle & \,=\, &    \langle\, {{\hat Y}_m}\,,\, \ \btd_y \ {\bf P}_{n \,-\, 2\,,\ m}\, (y)\, \rangle \ \,+\, \ O\, (|\,y \,-\, {{\hat Y}_m}|^{\,n \,-\, 2}\,)\\[0.15in]
& \ & \!\!\!\!\!\!\!\!\!\!=\ \langle\, {{\hat Y}_m}\,,\, \ \btd_z \ {\bf P}_{m\,\to} \, (z)\, \rangle |_{\,z \,=\, y \,-\, {{\hat Y}_m}}  \,+\, \ O\, (|\,y \,-\, {{\hat Y}_m}|^{\,n \,-\, 2}\,). \ \ \ \ \ \ \ \ \ \ \ \ \ \ \ \\[0.1in]
\!\!\!{\mbox{Let}} \ \ \ \ \ \ \ \ \  \ \ \ \ \ \ \ \ \ \ \ \ \ \ \ \  \ \  \ \  {\cal L}_m\, (z) & := &  \langle\, {{\hat Y}_m}\,,\, \ \btd_z \ {\bf P}_{m\,\to} \, (z)\, \rangle \ \ \ \ \ \ \ \  \mfor \ \ 1\, \le \,m \,\le\, k\,.
\end{eqnarray*}

${\cal L}_m\,$  is a homogeneous polynomial of degree $n \,-\, 3\,.\,$  As in (A.7.13),  from (A.7.9), we recognize
$$
\ {\bf P}_{m\,\to} \, (z)\, =\, {\bf P}_{n \,-\, 2\,,\ m} \, ( z \,+\, {\hat Y}_m ) \, \,=\, \, \,-\,  \!\!\!\sum_{|\,\alpha| \,= \,{n \,-\, 2}}    {1\over {\alpha\,!}} \ \left[\, D^{({\,n \,-\, 2})}_{\,\alpha}\, ({\tilde c}_n\,  K)\,( {\hat Y}_m) \,\right]\,\cdot  z^\alpha\,. \leqno (A.7.15)
$$
Observe that, via (A.7.10), $\,\Delta_o^{h_{n \,-\, 2}} \ {\bf P}_{m\,\to} \, (z)  \ \equiv \ 0\,.\,$ \bk
{\it In the following we assume that, for each $\,m$ with $\ 0 \le m \le k\,,\,$ there is a positive number}\, $\eta_m$ \,{\it such that}
$$
\ \ \eta_m\ <\ {{n \,-\, \,6}\over {(n \,-\, 3)\,(3n \,-\, 2)}} \ \ \ \ \ \  {\it{and \ }} \ \ \ \ \ \xi_{m_i} \ \,=\, \ \lambda_{m_i}^{1 \,+ \,\eta_m} \cdot {\vec{\,X}}_m \ \ \ \ \ {\it{for}} \ \  i \gg 1\,, \ \ \ \ \leqno (A.7.16)
$$
{\it{where}}\, ${\vec{\,X}}_m \ \in \ \R^n\,$  {\it{is \  fixed}}\,.\, Consider the integral

 \vspace*{-0.15in}

\begin{eqnarray*}
(A.7.17)\!\!\!\!\!\!& \ & {\tilde c}_n \cdot \int_{B_{\, {\hat Y}_m} \!(\,\rho_3)}  \  \langle \,y\,,\, \ \btd_y \ K \, (y)\,\rangle \, [\,v_i\,(y)\,]^{{2n}\over {n \,-\, 2}}\ dy
\\[0.1in]
 & \ & \hspace*{1.8in} [\,\downarrow \ {\mbox{via \ \ (A.7.12)\,, \ \ (A.7.13) \ \ and \ \ (A.7.14)}}\,]\\[0.1in]
& \,=\, & {\tilde c}_n \cdot \int_{B_{ \,{\hat Y}_m} \!(\,\rho_3)}   \langle \,z\,,\, \ \btd_z \ K_{\to} \, (z)\,\rangle|_{\,z \,=\, (y \,-\, \,{\hat Y}_m)} \, [\,v_i\,(y)\,]^{{2n}\over {n \,-\, 2}}\ dy \ \,+\, \\[0.15in]
 & \ &   \!\!\!\!\!\!\!\!\!\!\!\!+ \  \int_{B_{\, {\hat Y}_m} \!(\,\rho_3)}   \langle \,{\hat Y}_m\,,\, \ \btd_z \ {\bf P}_{m\,\to} \, (z)\,\rangle|_{\,z \,=\, (y -\, {\hat Y}_m)} \, [\,v_i\,(y)\,]^{{2n}\over {n \,-\, 2}}\ dy \ \,+\, \ O_{\lambda_{m_i}}\, (n \,-\, 2)\\[0.15in]
 & \ &  \!\!\!\!\!\!\!\!\!\!\!\!\!\!\!\!\!\!\!\!\!\![\,{\mbox{using \  (2.26) \ \& \ (2.27)\,, \  as \  in (A.6.41) \  for \  the \  term}} \  \, O\,(\,|\,y \,-\, {\hat Y}_m|^{\,n \,-\, 2}\,)\, \uparrow\ ]\\[0.15in]
 & \,=\, & {\tilde c}_n \cdot \int_{B_o  (\,\rho_3)}   \langle \,z\,,\, \ \btd_z \ K \, (z)\,\rangle \, [\,v_i\,(z+ {\hat Y}_m )\,]^{{2n}\over {n \,-\, 2}}\ dz \ \ \ \ \ \ \ \ \ \ ( z \,=\,  y \,-\, {\hat Y}_m)\\[0.15in]
 & \ & \ \ \ \ \ \ \ \ \ \ \ \ \ \ \ \ \,+\, \  \int_{B_o(\,\rho_3)}  {\cal L}_m\,(z) \, [\,v_i\,(z \,+\, {\hat Y}_m)\,]^{{2n}\over {n \,-\, 2}}\ dz \ \,+\, \ O_{\lambda_{m_i}}\, (n \,-\, 2) \ \ \ \ \ \ \ \ \ \ \ \ \ \ \  \\
 \end{eqnarray*}

  \vspace*{-0.15in}

\hspace*{0.45in} [\ $\downarrow$ \  using (A.7.11) as in (A.6.33), (A.6.46), and the proof of Proposition A.6.34\,]

 \begin{eqnarray*}
 & \,=\, & \lambda_{m_i}^{n \,-\, 2}\left[\,{{n\, (n \,-\, 2)}\over {{\tilde c}_n\cdot K\, ( {{\hat Y}_m} )}} \right]^{n\over 2} \int_{\R^n}    {\bf P}_{\!\to m}\,(z) \ \left( {1\over {1 \,+\, |\,z|^2}}\right)^n \ dz \, \,+\, \,   o_{\lambda_{m_i}}\, (n \,-\, 2) \, \,+\, \ \ \ \ \ \ \ \ \ \ \\[0.15in]
 & \ & \  \ \   \ \  \,+\, \ \lambda_{m_i}^{n \,-\, 3} \left[\,{{n\, (n \,-\, 2)}\over {{\tilde c}_n\cdot K\, ( {{\hat Y}_h} )}} \right]^{n\over 2}  \int_{\R^n}  {\cal L}_m\,(z)  \ \left( {1\over {1 \,+\, |\,z|^2}}\right)^n \ dz \ \,+\, \\[0.15in]
 & \ & \ \ \ \ \ \  \ \   \,+\, \    O_{\lambda_{m_i}}\, ([\,n \,-\, 3] \,+\, \eta_m) \ +\, \cdot \cdot \cdot \,+ \, O_{\lambda_{m_i}}\, ([\,n \,-\, 3] \,+\, [\,n \,-\, 4]\cdot \eta_m)\ \,+\, \\[0.15in]
 & \ & \ \ \ \ \  \ \  \     \,+\, \ \lambda_{m_i}^{(n \,-\, 3) \,+\, (n \,-\, 3)\cdot\, \eta_m }\cdot {\cal L}_m\,({\vec{\,X}}_m) \left[\,{{n\, (n \,-\, 2)}\over {{\tilde c}_n\cdot K\, ( {{\hat Y}_h} )}} \right]^{n\over 2}  \int_{\R^n}    \ \left( {1\over {1 \,+\, |\,z|^2}}\right)^n \ dz \ \,+\, \\[0.1in]
 & \ &   \  \ \ \ \ \ \ \ \  \ \    \,+\, \    O_{\lambda_{m_i}}  \left( [\, n \,-\, 3]  \ \,+\, \ {{n \,-\, 6}\over {3n \,-\, 2}} \right)  \\[0.15in]
 & \ &   \ \ \ \   [\,\uparrow {\mbox{as \ in \ (A.6.34), \ Lemma \ A.6.47, \ (A.6.64)\,--\,(A.6.66) \  with}} \ \ \ell \,= \,n  \,-\, 3\,,\ n > 6\,] \ .
 \end{eqnarray*}

In (A.7.17), we put priority on   terms with lowest order in $\,\lambda_{m_i}\,$.\, If $\,n\,$ is even, then $\,n \,-\, 3\,$ is odd, giving
 $$
 \int_{\R^n}   {\cal L}_m\,(z)  \ \left( {1\over {1 \,+\, |\,z|^2}}\right)^n \ dz \,=\, 0\ . \leqno (A.7.18)
 $$
 In case $n$ is odd, then $\,n \,-\, 3\,$ is even, and $\,h_{n \,-\, 2} \,=\, {{n \,-\, 3}\over 2}\,.\,$ Thus the condition
\begin{eqnarray*}
& \ &  \Delta_o^{ \left( h_{n \,-\, 2} \right)  }\  {\bf P}_{m\,\to}   \ \equiv \ 0 \ \ \Longleftrightarrow \ \ \Delta_o^{ \!\left({{n \,-\, 3}\over 2} \right)  }\  {\bf P}_{m\,\to}   \ \equiv \ 0 \\[0.15in]
& \Longrightarrow &
 \Delta_o^{\! \left( {{n \,-\, 3}\over 2} \right)  }\ {\cal L}_m\,(z) \ \,=\, \  \Delta_o^{\! \left( {{n \,-\, 3}\over 2} \right)  }\  \langle \,{\hat Y}_m\,, \ \btd_z\ {\bf P}_{m\,\to} \,(z)\, \rangle \\[0.1in]
  & \ & \ \ \ \ \ \  \ \ \ \ \ \ \ \ \ \ \ \ \ \ \,   \,=\, \ \bigg\langle \,{\hat Y}_m\,, \ \btd_z \ \left[  \ \Delta_o^{ \!\left( {{n \,-\, 3}\over 2} \right)  } \ {\bf P}_{n \,-\, 2\,,\ m} \,(z) \right] \bigg\rangle \ \,=\, \ 0\,.
\end{eqnarray*}
 Again we obtain (A.7.18) by using Lemma A.6.8.\bk
   When we add  up the integrals and estimates from each simple blow\,-\,up point,  for simplicity, we skip the terms with intermediate orders
 $$
 O_{\lambda_{m_i}}\, (\,[\,n \,-\, 3] \,+\, \eta_m) \ \,+\, \,\cdot \cdot \cdot \,+\,  \ O_{\lambda_{m_i}}\, (\,[\,n \,-\, 3] \,+\, [\,n \,-\, 4]\cdot \eta_m)\,,
 $$
 and draw a conclusion from (A.7.8) and (A.7.17) that

 \newpage

 (A.7.19)
$$
\ \,{\cal L}_m \, ({\vec{\,X}}_m)  \ \,=\, \  \langle\, {{\hat Y}_m}\,,\, \ \btd_z \ {\bf P}_{n \,-\, 2\,,\,m} \ ({\vec{\,X}}_m) \, \rangle \ \,=\, \ 0 \ \  \ \ ({\mbox{where}} \ \ \xi_{m_i} \ \,=\, \ \lambda_{m_i}^{1 +\, \eta_m} \cdot {\vec{\,X}}_m)\,,
$$

\vspace*{-0.25in}

$$
 (n \,-\, 3)\cdot \eta_m \ < \ {{n \,-\, 6}\over {3n \,-\, 2}} \ \ \ \ \ \  \ \left( {\mbox{observe \ \ that \ \ }} {{n \,-\, 6}\over {3n \,-\, 2}}  \,<\, 1 \right), \leqno {\mbox{(A.7.20) \ \ \  {\it provided}}}
$$
{\it and there is no interference from other blow\,-\,up points, precisely}\,:\\[0.1in]
(A.7.21)
$$
   (n \,-\, 3) \cdot \,\eta_{m} \ \not= \ h \cdot \,\eta_j \ \ \ \  {\it{for}} \ \ \   j \not= m\,, \   1 \,\le\, h \,\le\, n \,-\, \,3\ \ \ (\,h \ \ {\mbox{ is \ a \   natural \   number}})\,.
$$
In case some of the $\eta_h \ \,=\, \ \eta_m$\,,\, together with (A.7.7),  we obtain

\vspace*{-0.15in}

\begin{eqnarray*}
(A.7.22) \!\!\!\!\!\!\!& \ &  \left[\,{{n\, (n \,-\, 2)}\over {({\tilde c}_n\,K)\, ( {{\hat Y}_m} )}} \right]^{n\over 2}  \cdot \,S_m^{\,(n \,-\, 3) \,+\, (n \,-\, 3) \,\eta_m} \,\cdot  \langle\, {{\hat Y}_m}\,,\, \ \btd_z \ {\bf P}_{n \,-\, 2\,,\ m} \ ({\vec{\,X}}_m) \, \rangle \ \\[0.15in]
  & \ &  \!\!\!\!\!\!\!\!\!\!\!\!    \,+\, \! \sum_{\scriptstyle 0 \,<\, h \,   \le \,k   \,,\, h \not= m \atop \scriptstyle {\mbox{with }} \eta_h\, =\,  \eta_m } \!  \left[\,{{n\, (n \,-\, 2)}\over {({\tilde c}_n\, K)\, ( {{\hat Y}_h} )}} \right]^{n\over 2} \cdot \,S_h^{\,(n \,-\, 3) \,+\, (n \,-\, 3) \,\eta_h} \,\cdot
  \langle\, {{\hat Y}_h}\,,\, \ \btd_z \ {\bf P}_{n \,-\, 2\,,\,\,h} \ ({\vec{\,X}}_h) \, \rangle \ \,=\, \ 0\   , \ \ \ \ \  \ \ \ \ \ \ \ \
\end{eqnarray*}
where we set the condition

\vspace*{-0.25in}

$$
 \eta_j \,\not= \, \eta_m \ \ \Longrightarrow \ \  {\mbox{(A.7.21) \ \ holds}}\,.\leqno (A.7.23)
$$

\vspace*{-0.05in}

We summarize the conditions assumed in the balance and cancelation formulas (A.7.19) and (A.7.22)\,: besides the ones mentioned  next to them [\,they are (A.7.20), (A.7.21), and  (A.7.23)\,]\,,\, and the conditions found in Main Theorem 1.14 for each blow\,-\,up point, plus (A.7.1)\,--\,\,(A.7.4),\, (A.7.10) and (A.7.16)\,.

\newpage

{\large{\bf \S\,A.8.} \ \ {\bf Verification of\ } (A.6.13) {\bf and\,} (A.6.18).}\\[0.2in]
 Refer to (A.6.13) for the notation we use.
 \begin{eqnarray*}
& \ & \\[-0.25in]
& \ & \int_{\R^n} \ \, (  {\mbox{\,terms \ \ without  \ }} y_{|_1} \ \,\& \ \,y_{|_n}  ) \cdot \ y_{|_1}^{\,k \,+\, 2} \cdot \left( {1\over {1+ r^{\,2}}} \right)^{\!n} \, d y \ \ \ (\,{\mbox{absolute \ \ convergence}})\\[0.15in]
& \,=\, & \int_{-\infty}^\infty \cdot  \cdot \cdot \ \int_{-\infty}^\infty \,(  {\mbox{\,terms \ \ without  \ }} y_{|_1} \ \,\& \ \, y_{|_n}  ) \ \times \  \ \ \ \ \ \ \ \ \ \ \ \ \ \ \ \    (\,{\mbox{Fubini's Theorem}}) \\[0.1in]
 & \ &  \!\!\!\!\!\!\!\!\!\times \left[ \ \int_{-\infty}^\infty\int_{-\infty}^\infty \,y_{|_1}^{\,k \,+\, 2} \cdot  \left( {1\over {1 \,+\, [\,y_{|_2}^{\,2} \,+\, \cdot \cdot \cdot \,+\, y_{|_{n-1}}^2 \,] \,+\, y_{|_1}^{\,2} \,+\, y_{|_n}^{\,2}}} \right)^{\!n}  dy_{|_1} \,dy_{|_n} \right]dy_{|_2} \cdot \cdot \cdot dy_{|_{n \,-\, 1}}\\[0.15in]
& \,=\, & \int_{-\infty}^\infty \cdot  \cdot \cdot \ \int_{-\infty}^\infty \,(  {\mbox{\,terms \ \ without  \ }} y_{|_1} \ \& \ y_{|_n}  ) \ \times \  \\[0.1in]
 & \ & \!\!\!\!\!\!\!\times\, \left[\  \int_0^\infty\int_0^{2 \,\pi} \,y_{|_1}^{\,k \,+\, 2}  \cdot   \left( {1\over {1 \,+\, [\,y_{|_2}^{\,2} \,+\, \cdot \cdot \cdot \,+\, y_{|_{n \,-\, 1}}^{\,2}\,] \,+\, \rho^{\,2}}} \right)^{\!n} \cdot \rho \  d\theta \, d \rho \right]dy_{|_2} \cdot \cdot \cdot dy_{|_{n-1}}\\
 & \ & \ \ \ (\,{\mbox{polar \ \ coordinates \ \ on \ \ }} \R^2\,,   \ \   \rho^2 \ \,=\, \ y_1^2 \,+\, y_n^2\,,\, \ \ y_1 \, \,=\, \ \rho\cdot \sin\,\theta\,, \ \ y_n \, \,=\, \ \rho\cdot \cos\,\theta )\\[0.15in]
& \,=\, & \int_{-\infty}^\infty \cdot  \cdot \ \int_{-\infty}^\infty \,(  {\mbox{\,terms \ \ without  \ }} y_{|_1} \ \,\&\, \ y_{|_n}  ) \ \times  \   \\[0.1in]
 & \ & \!\!\!\!\!\!\!\!\!\!\!\!   \times\, \left[ \ \int_0^\infty\int_0^{2 \,\pi} \,\rho^{\,k \,+\, 3}\, (\sin^{\,k \,+\, 2} \theta)  \cdot   \left( {1\over {1 \,+\, [\,y_{|_2}^{\,2} \,+\, \cdot \cdot \cdot \,+\, y_{|_{n-1}}^{\,2}\,] \,+\, \rho^{\,2}}} \right)^{\!n}  \cdot  d\theta \, d \rho\, \right]\,dy_{|_2} \cdot \cdot \cdot dy_{|_{n-1}}\\[0.15in]
 & \,=\, & \int_{-\infty}^\infty \cdot  \cdot \ \int_{-\infty}^\infty \,(  {\mbox{\,terms \ \ without  \ }} y_{|_1} \ \,\&\, \ y_{|_n}  ) \ \times  \   \\[0.1in]
& \ & \!\!\!\!\!\!\!\!\!\!\!\!   \times\,\left[  \!\int_0^\infty \!\!\rho^{\,k \,+\, 3}  \cdot   \left( {1\over {1 \,+\, [\,y_{|_2}^{\,2} \,+\, \cdot \cdot \cdot \,+\, y_{|_{n-1}}^{\,2}\,] \,+\, \rho^{\,2}}} \right)^{\!n}  \cdot  \left\{   \int_0^{2 \,\pi} \!\! \sin^{\,k \,+\, 2} \theta  \,d\theta \right\}  d \rho \, \right]\,dy_{|_2} \cdot \cdot\,   dy_{|_{n-1}}\ .
\end{eqnarray*}
Here ``\,(terms without $\,y_{|_1} \ \,\&\, \ y_{|_n}\,$)" is a polynomial on $y_{|_2}\,,\, \cdot \cdot \cdot\,, \ y_{|_{n-1}}$\,,\, having  sufficiently low degree so that the integral is absolutely convergent.
Likewise,
\begin{eqnarray*}
& \ & \int_{-\infty}^\infty \cdot  \cdot \ \int_{-\infty}^\infty \,(  {\mbox{\,terms \ \ without  \ }} y_{|_1} \ \,\&\, \ y_{|_n}  ) \ \times \   \\
 & \ &   \!\!\!\!\!\!\!\!\!\!\!\!\!\!\!\!\times \, \left[ \ \int_{-\infty}^\infty\int_{-\infty}^\infty \,y_{|_1}^{\,k}\,y_{|_n}^{\,2} \ \left( {1\over {1 \,+\, [\,y_{|_2}^{\,2} \,+\, \cdot \cdot \cdot \,+\,  y_{|_{n-1}}\,] \,+\, y_{|_1}^{\,2} \,+\, y_{|_n}^{\,2}}} \right)^{\,n}\ dy_{|_1} \,dy_{|_n} \right]dy_{|_2} \cdot \cdot \cdot dy_{|_{n-1}}\\[0.15in]
& \,=\, & \int_{-\infty}^\infty \cdot  \cdot \ \int_{-\infty}^\infty \,(  {\mbox{\,terms \ \ without  \ }} y_{|_1}\ \,\&\,\ y_{|_n}  ) \ \times  \   \\
 & \ & \!\!\!\!\!\!\!\!\!\!\!\!\!\!\!\!   \times   \left[\,  \int_0^\infty \!\!\!\rho^{\,k \,+\, 3}    \! \left( \!{1\over {1 \,+\, [\,y_{|_2}^{\,2} \,+\, \cdot \cdot   \,+\,  y_{|_{n-1}}] \,+\, \rho^{\,2}}} \!\right)^{\!\!n} \! \left\{ \int_0^{2 \,\pi} \!(\sin^k \theta)\,(\cos^{\,2} \theta) \,d\theta \right\} \, d \rho \right]dy_{|_2} \cdot \cdot  \, dy_{|_{n-1}}\,.
\end{eqnarray*}

A direct calculation using integration by parts shows that
$$
\int_0^{2 \,\pi} (\sin^{\,k \,+\, 2} \theta)\, d\theta \ \,=\, \ \,-\, \int_0^{2 \,\pi} (\sin^{\,k \,+\, 1} \theta)\, d\,[\,\cos\,\theta] \ \,=\, \  (k \,+\, 1) \int_0^{2 \,\pi} (\sin^k \theta)\,(\cos^{\,2} \theta) \,d\theta\,.
$$
Hence we deduce  (A.6.13).\bk
To show (A.6.18), recall that
$$
{\cal o} \ \,=\, \ {\alpha_2}  \ \,+\, \ \cdot \cdot \cdot \ \,+\, \   \alpha_{n \,-\, 1}\ .
$$
We demonstrate how to use induction on $\,{\cal o}\,$ to prove the assertion. Recall that
$$
 \ell \,\in\, [\,0, \ n\,-\, 2\,) \ \ \ {\mbox{is \ \ even}}\,.
$$
{\bf (I)} \ \ When $\,{\cal o} \ \,=\, \ 0$\,,\, the term is a constant. We have
\begin{eqnarray*}
\Delta_o^{\!(h_\ell)} \  y_{|_1}^{\,k \,+\, 2}  & \,=\, & (k \,+\, 2) \, (k \,+\, 1) \cdot \cdot \cdot \, 3 \cdot 2 \cdot 1\,;
 \\[0.15in]
 \Delta_o  \ \left\{(k \,+\, 1) \, y_{|_1}^{\,k  } \,y_{|_n}^2  \right\} & \,=\, & (k \,+\, 1) \,k \, (k \,-\, 1) \, y_{|_1}^{\,k \,-\, 2  }\,y_{|_n}^2 \ \,+\, \ 2 \, (k \,+\, 1) \, y_{|_1}^{\,k  }\,,\\[0.15in]
 \Delta_o^{\!(2)}  \ \left\{(k \,+\, 1) \, y_{|_1}^{\,k  } \,y_{|_n}^2  \right\} & \,=\, & (k \,+\, 1) \,k \, (k \,-\, 1) \, (k \,-\, 2) \, (k \,-\, 3) \, y_{|_1}^{\,k \,-\, 4  }\,y_{|_n}^2 \ +\\[0.1in]
  & \ &  \ \ \ \ \ \,+\, \ 2 \times 2\, (k \,+\, 1) \,k \, (k \,-\, 1) \, y_{|_1}^{\,k \,-\, 2  }\,,\\
  & \cdot & \\
  & \cdot & \\
  & \cdot &\\
  \end{eqnarray*}

\vspace*{-0.5in}

  \begin{eqnarray*}
  \Delta_o^{ \, (h_\ell \,-\, 2)}  \ \left\{(k \,+\, 1) \, y_{|_1}^{\,k  } \,y_{|_n}^2  \right\} & \,=\, & (k \,+\, 1) \,k \, (k \,-\, 1) \, (k \,-\, 2) \, (k \,-\, 3) \, \cdot \cdot \cdot \, 3\cdot \, y_{|_1}^{\,2 }\,y_{|_n}^2 \ +\\[0.1in]
  & \ & \ \ \ \ \ \ \ \ \ \,+\, \ 2 \,\left[ \,{h_\ell} \,-\, 2 \right] (k \,+\, 1) \,k \, (k \,-\, 1) \cdot \cdot \cdot \, 5  \,\cdot y_{|_1}^{\,4  }\,,\\[0.1in]
  \Delta_o^{ \, (h_\ell \,-\, 1)}  \ \left\{(k \,+\, 1) \, y_{|_1}^{\,k  } \,y_{|_n}^2  \right\} & \,=\, & (k \,+\, 1) \,k \, (k \,-\, 1) \, (k \,-\, 2) \, (k \,-\, 3) \, \cdot \cdot \cdot \, 3\cdot 2 \cdot  \, [\,y_{|_1}^{\,2 }+ y_{|_n}^2] \ +\\[0.1in]
  & \ & \ \ \ \ \ \ \ \ \ \ \,+\, \ 2 \,\left[ \,{h_\ell} \,-\, 2\, \right] (k \,+\, 1) \,k \, (k \,-\, 1) \cdot \cdot \cdot \, 5 \cdot 4 \cdot 3\, y_{|_1}^{\,2 }\,,\\[0.1in]
  \Delta_o^{ \, h_\ell}  \ \left\{(k \,+\, 1) \, y_{|_1}^{\,k  } \,y_{|_n}^2  \right\} & \,=\, & (k \,+\, 1) \,k \, (k \,-\, 1) \, (k \,-\, 2) \, (k \,-\, 3) \, \cdot \cdot \cdot \, 3\cdot 2 \cdot  1 \cdot \, 4\ \,+\, \\[0.1in]
  & \ & \ \ \ \ \ \,+\, \  \left[\, \ell \,-\, 4 \right] (k \,+\, 1) \,k \, (k \,-\, 1) \cdot \cdot \cdot \, 5 \cdot 4 \cdot 3\,\cdot 2 \cdot 1\\[0.1in]
  & \,=\, & (k \,+\, 1) \,k \, (k \,-\, 1) \, (k \,-\, 2) \, (k \,-\, 3) \, \cdot \cdot \cdot \, 3\cdot 2 \cdot  1 \!\cdot\! [\,\ell \,-\, 4 \,+\, 4\,]\\[0.1in]
   & \,=\, & (k \,+\, 2) \,(k \,+\, 1) \,k \, (k \,-\, 1) \, (k \,-\, 2) \, (k \,-\, 3) \, \cdot \cdot \cdot \, 3\cdot 2 \cdot  1  \ \ \ \ \ \\[0.1in]
    & \ & \hspace*{1.55in}  (\,{\mbox{as}} \ \ \ell \,=\, k \,+\, 2  \ \ {\mbox{in \ this \ case}})\,.
\end{eqnarray*}
Hence the case $\,{\cal o}\, =\, 0\,$ is settled.

 \vspace*{0.2in}

{\bf (II)} \ \ As an induction hypothesis, suppose that
\begin{eqnarray*}
\Delta_o^{\! (h_\ell)} \ \left\{ y_{|_1}^{\,k \,+\, 2} \cdot [ \,\cdot \cdot \cdot \ {\mbox{degree}} \ \,=\, {\cal o} \ \cdot \cdot \cdot \,] \right\} \,=\, \Delta_o^{\!(h_\ell)} \ \left\{(k \,+\, 1) \, y_{|_1}^{\,k  } \,y_{|_n}^2 \times [ \,\cdot \cdot \cdot \ {\mbox{degree}} \ \,=\, \ {\cal o} \ \cdot \cdot \cdot \,]  \ \right\}
\end{eqnarray*}
holds for $\,\ell \,=\, ( k \,+\, 2) \ \,+\, \ {\cal o}\,,\,$ where $\,k \, \ge \,2 \,$ (variable), but $\,{\cal o} \,>\, 0\,$ (fixed). We continue to use the notations above and {\it there is no\,  $\,y_{|_1}\,$ or \,$\,y_{|_n}\,$ inside the homogeneous polynomial denoted by\,} $\,[ \,\cdot \cdot \cdot \ {\mbox{degree}} \ \,=\, \,{\cal o}   \ \cdot \cdot \cdot \,]\,  $\,.\, Let us go on to show
\begin{eqnarray*}
(A.8.1)\ \ \ \ \ \ \ \ \ \ \ \  & \ & \Delta_o^{\! (h_\ell)} \ \left\{ y_{|_1}^{\,k \,+\, 2} \cdot [ \,\cdot \cdot \cdot \ {\mbox{degree}} \ =\, {\cal o} \,+\, 2 \ \cdot \cdot \cdot \,] \right\}\\[0.1in] &  \,=\, & \Delta_o^{\! (h_\ell)} \ \left\{(k \,+\, 1) \, y_{|_1}^{\,k  } \,y_{|_n}^2 \cdot [ \,\cdot \cdot \cdot \ {\mbox{degree}} \ \,=\, \,{\cal o} \,+\, 2 \ \cdot \cdot \cdot \,]  \ \right\}\ , \ \ \ \ \ \ \ \ \ \ \ \ \ \ \ \ \ \ \ \ \ \ \ \
\end{eqnarray*}
where $k \,\ge\, 2$ is even. Let us find the first Laplacians\,:
\begin{eqnarray*}
(A.8.2) \!\!& \ & \Delta_o  \, \left\{ y_{|_1}^{\,k \,+\, 2} \cdot [ \,\cdot \cdot \cdot \ {\mbox{degree}} \ \,=\, \  {\cal o} \,+\, 2 \ \cdot \cdot \cdot \,] \right\}\\[0.1in]
& \,=\, & (k \,+\, 2) \, (k \,+\, 1) \ y_{|_1}^{\,k  } \cdot [ \,\cdot \cdot \cdot \ {\mbox{degree}} \ \,=\,  \ {\cal o} \,+\, 2 \ \cdot \cdot \cdot \,] \\[0.1in]
 & \ & \ \ \ \ \ \ \ \ \ \ \ \ \ \ \ \ \,+\, \  y_{|_1}^{\,k \,+\, 2} \cdot \left\{ \Delta_o \ [ \,\cdot \cdot \cdot \ {\mbox{degree}} \ \,=\, \ {\cal o} \,+\, 2 \ \cdot \cdot \cdot \,] \right\}\\[0.1in]
& \,=\, & k \, (k \,+\, 1) \ y_{|_1}^{\,k  } \cdot [ \,\cdot \cdot \cdot \ {\mbox{degree}} \ \,=\,  \ {\cal o} \,+\, 2 \ \cdot \cdot \cdot \,] \\[0.1in]
& \ &  \ \ \  \ \,+\, \ 2 \, (k \,+\, 1) \ y_{|_1}^{\,k  } \cdot [ \,\cdot \cdot \cdot \ {\mbox{degree}} \ \,=\,  \ {\cal o} \,+\, 2 \ \cdot \cdot \cdot \,]\\[0.1in]
 & \ & \ \ \ \ \ \ \ \ \ \ \ \ \ \ \ \ \,+\, \  y_{|_1}^{\,k \,+\, 2} \cdot \left\{ \Delta_o \ [ \,\cdot \cdot \cdot \ {\mbox{degree}} \ \,=\,  \ {\cal o} \,+\, 2 \ \cdot \cdot \cdot \,] \right\}\\[0.1in]
 & \ & \ \ \ \ \ \ \ \ \ \  \ \ \ \ \ \ \ \ \ \  \ \ \ \ \ \ \  (\ \leftarrow \ \ \ \ \ \uparrow \ \ {\mbox{degree}} \ \,=\,  \ {\cal o} \ \ \ \ \ \rightarrow \ )\,;
 \\[0.22in]
& \ &  \Delta_o  \, \left\{(k \,+\, 1) \, y_{|_1}^{\,k  } \,y_{|_n}^2 \cdot [ \,\cdot \cdot \cdot \ {\mbox{degree}} \ \,=\,  \ {\cal o} \,+\, 2 \ \cdot \cdot \cdot \,]  \ \right\}\\[0.1in]
 &= & \  (k \,+\, 1)\,k\,(k \,-\, 1) \,y_{|_1}^{\,k  \,-\, 2} \,y_{|_n}^2 \cdot [ \,\cdot \cdot \cdot \ {\mbox{degree}} \ \,=\,  \ {\cal o} \,+\, 2 \ \cdot \cdot \cdot \,]\\[0.1in]
 & \ & \ \ \ \ \ \ \ \ \ \ \ \  \,+\, \ 2\,(k \,+\, 1) \,y_{|_1}^k \cdot [ \,\cdot \cdot \cdot \ {\mbox{degree}} \ \,=\, \ {\cal o} \ +\  2 \ \cdot \cdot \cdot \,]\\[0.1in]
 & \ & \ \ \ \ \ \ \ \ \ \ \ \ \ \ \ \  \ \ \ \ \ \  \,+\, \  (k \,+\, 1) \,y_{|_1}^{\,k  } \,y_{|_n}^2 \cdot \left\{ \Delta_o \ [ \,\cdot \cdot \cdot \ {\mbox{degree}} \ \,=\, \ {\cal o} \,+\, 2 \ \cdot \cdot \cdot \,] \right\}\\[0.1in]
 & \ & \ \ \ \ \ \ \ \ \ \ \ \ \ \ \  \ \ \ \ \ \ \ \ \ \  \ \ \ \ \ \ \  \ \ \ \ \ \ \ \ \   \ \ \  (\ \leftarrow \ \ \ \ \ \uparrow \ \ {\mbox{degree}} \ \,=\, \  {\cal o} \ \ \ \ \ \rightarrow \ )\,,\\[0.1in]
& \ & \ \ \ \ \ \ \ \ \ \ \ \  \ \ \ \ \ \   [\ {\mbox{observe \ \ that \ \ }}  (k \,+\, 2)\, (k \,+\, 1) \,-\, 2\, (k \,+\, 1) \ \,=\, \ k \, (k \,+\, 1)\,]\,. \ \ \ \ \ \ \ \ \
\end{eqnarray*}
Via the induction hypothesis, the last two terms in the respective expressions are equal. After simplification, to verify (A.8.1), it suffices to  show that
\begin{eqnarray*}
& \ & \Delta_o^{\!(\,h_\ell \,-\, 1)} \ \left\{ y_{|_1}^{\,k } \cdot [ \,\cdot \cdot \cdot \ {\mbox{degree}} \ \,=\, \   {\cal o} \,+\, 2 \ \cdot \cdot \cdot \,]\, \right\}\\[0.1in]
 &= & \Delta_o^{\!(\,h_\ell \,-\, 1)} \ \left\{(k \,-\, 1) \, y_{|_1}^{\,k \,-\, 2  } \ y_{|_n}^2 \cdot [ \,\cdot \cdot \cdot \ {\mbox{degree}} \ \,=\, \ {\cal o} \,+\, 2 \ \cdot \cdot \cdot \,]  \ \right\}\ .
\end{eqnarray*}
Applying $\,\Delta_o\,$ on the  terms
$$
\left\{ y_{|_1}^{\,k } \cdot [ \,\cdot \cdot \cdot \ {\mbox{degree}} \ \,=\, \   {\cal o} \,+\, 2 \ \cdot \cdot \cdot \,]\, \right\} \ \   {\mbox{and}} \ \  \left\{(k \,-\, 1) \, y_{|_1}^{\,k \,-\, 2  } \ y_{|_n}^2 \cdot [ \,\cdot \cdot \cdot \ {\mbox{degree}} \ \,=\, \   {\cal o} \,+\, 2 \, \cdot \cdot \cdot \,]  \, \right\},
$$
using  similar calculation and cancelation as in (A.8.2), we come down gradually to verify
\begin{eqnarray*}
(A.8.3) \ \ \ \ \ \ & \ & \Delta_o^{\! \left(h_\ell \,-\, {{k \,-\, 2}\over 2} \right)} \ \left\{ y_{|_1}^{\,4 } \cdot [ \,\cdot \cdot \cdot \ {\mbox{degree}} \ \,=\, \   {\cal o} \,+\, 2 \ \cdot \cdot \cdot \,]\, \right\}\\[0.1in]
 &= &\Delta_o^{\! \left(h_\ell \,-\, {{k \,-\, 2}\over 2} \right)} \   \left\{3 \, y_{|_1}^{\,2 } \,y_{|_n}^2 \cdot [ \,\cdot \cdot \cdot \ {\mbox{degree}} \ \,=\,  \ {\cal o} \,+\, 2 \ \cdot \cdot \cdot \,]  \ \right\}\,.\ \ \ \ \ \ \ \ \ \ \ \ \ \ \ \ \ \ \ \ \ \ \ \
\end{eqnarray*}
Note that
$$
\ell \,=\, (k \,+\, 2) \,+\, ( {\cal o} \,+\, 2) \ \ \Longrightarrow \ \ {\ell\over 2} \,-\, {{k \,-\, 2}\over 2} \ \,=\, \ {{4 \,+\, ( {\cal o} \,+\, 2)  }\over 2}\ .
$$
Apply the Laplacian on the two terms inside the brackets in (A.8.3) and obtain
\begin{eqnarray*}
& \ & \Delta_o  \ \left\{ y_{|_1}^{\,4} \cdot [ \,\cdot \cdot \cdot \ {\mbox{degree}} \ \,=\,  \ {\cal o} \,+\, 2 \ \cdot \cdot \cdot \,] \right\}\
= \ 4 \cdot 3 \ y_{|_1}^{\,2 } \cdot [ \,\cdot \cdot \cdot \ {\mbox{degree}} \ \,=\,  {\cal o} \,+\, 2 \ \cdot \cdot \cdot \,]\ \,+\,  \\[0.1in]
 & \ & \ \ \ \ \ \ \ \ \ \ \ \ \ \ \ \ \,+\, \  y_{|_1}^{\,4} \cdot \left\{ \ \Delta_o \ [ \,\cdot \cdot \cdot \ {\mbox{degree}} \ \,=\,  \ {\cal o} \,+\, 2 \ \cdot \cdot \cdot \,] \,\right\}\,;\\[0.1in]
 & \ & \ \ \ \ \ \ \ \ \ \  \ \ \ \ \ \ \ \ \ \  \ \ \ \ \ \,    (\ \leftarrow \ \ \ \ \ \uparrow \ \ {\mbox{degree}} \ \,=\,  \ {\cal o} \ \ \ \ \ \rightarrow \ )
 \\[0.25in]
& \ &  \Delta_o  \ \left\{3 \, y_{|_1}^{\,2  } \ y_{|_n}^2 \cdot [ \,\cdot \cdot \cdot \ {\mbox{degree}} \ \,=\, \ {\cal o} \,+\, 2 \ \cdot \cdot \cdot \,]  \ \right\}\\[0.1in]
 &= & \  3 \cdot 2 \cdot 1 \cdot [\ y_{|_1}^2 \,+\, y_{|_n}^2\,] \cdot [ \,\cdot \cdot \cdot \ {\mbox{degree}} \ \,=\, \ {\cal o} \,+\, 2 \ \cdot \cdot \cdot \,]\ \,+\, \\[0.1in]
 & \ & \ \ \ \ \  \ \ \ \ \ \ \ \ \ \ \,+\, \  3 \,y_{|_1}^{\,2  }\, y_{|_n}^2 \cdot \left\{ \ \Delta_o \ [ \,\cdot \cdot \cdot \ {\mbox{degree}} \ \,=\,  \ {\cal o} \,+\, 2 \ \cdot \cdot \cdot \,] \, \right\}\\[0.1in]
 & \ & \ \ \ \ \ \ \ \ \ \ \ \ \ \ \  \ \ \ \ \ \ \ \ \ \  \ \, \ \ \ \  (\ \leftarrow \ \ \ \ \ \uparrow \ \ {\mbox{degree}} \ \,=\,  {\cal o} \ \ \ \ \ \rightarrow \ )\,.
 \\
\end{eqnarray*}
Again we apply the induction hypothesis to cancel the last  term in each expression above.  Apply the Laplacian again and obtain
\begin{eqnarray*}
(A.8.4)& \ & \!\!\!\!\!\!\!\Delta_o \ \left\{   4 \cdot 3 \ y_{|_1}^{\,2 } \cdot [ \,\cdot \cdot \cdot \ {\mbox{degree}} \ \,=\, \   {\cal o} \,+\, 2 \ \cdot \cdot \cdot \,] \ \right\} \\[0.1in]
 &\ & \!\!\!\!\!\!\!\!\!\!\!\!\!\!\!\!\!\!\!\!\!\!\!\!\!\!\!\!=\ 4 \cdot 3 \cdot 2 \cdot  [ \,\cdot \cdot \cdot \ {\mbox{degree}} \ \,=\,  {\cal o} \,+\, 2 \ \cdot \cdot \cdot \,] \ \,+\, \ 4 \cdot 3 \cdot y_{|_1}^2 \cdot \left\{ \Delta_o \ [ \,\cdot \cdot \cdot \ {\mbox{degree}} \ \,=\, \ {\cal o}\ +\ 2 \ \cdot \cdot \cdot \,] \, \right\}\ ;\\[0.25in]
(A.8.5)& \ &  \Delta_o \ \left\{   3 \cdot 2 \ [\, y_{|_1}^{\,2 } \,+\, y_{|_n}^2\,] \cdot [ \,\cdot \cdot \cdot \ {\mbox{degree}} \ \,=\,  {\cal o} \,+\, 2 \ \cdot \cdot \cdot \,] \,\right\}\\[0.1in]
 &  \,=\, &     3 \cdot 2 \cdot 4 \cdot  [ \,\cdot \cdot \cdot \ {\mbox{degree}} \ \,=\,  \ {\cal o}\ \,+\, \ 2 \ \cdot \cdot \,\cdot \,] \\[0.1in]
   & \ & \ \ \ \ \ \ \ \ \,+\, \  3 \cdot 2 \cdot [\, y_{|_1}^2 \ \,+\, \  y_{|_n}^2\,] \cdot \left\{ \Delta_o \ [ \,\cdot \cdot \cdot \ {\mbox{degree}} \ \,=\, \  {\cal o} \ +\  2 \ \cdot \cdot \cdot \,] \, \right\}\,.
\end{eqnarray*}

\newpage

\begin{eqnarray*}
{\mbox{As}} \  \ & \ & 4 \cdot 3 \cdot  y_{|_1}^2\,  \cdot \left\{ \Delta_o \ [ \,\cdot \cdot \cdot \ {\mbox{degree}} \ \,=\, \ {\cal o} \,+\, 2 \ \cdot \cdot \cdot \,] \ \right\} \\[0.1in]
& \,=\, & 2\cdot 3 \cdot  y_{|_1}^2\,  \cdot \left\{ \Delta_o \ [ \,\cdot \cdot \cdot \ {\mbox{degree}} \ \,=\, \ {\cal o} \ +\ 2 \ \cdot \cdot \cdot \,] \ \right\} \ \,+\, \\[0.1in]
& \ & \ \ \ \ \ \ \,+\, \  2\cdot 3 \cdot  y_{|_1}^2\,  \cdot \left\{ \Delta_o \ [ \,\cdot \cdot \cdot \ {\mbox{degree}} \ \,=\,  \ {\cal o} \,+\, 2 \ \cdot \cdot \cdot \,] \ \right\}\ ,\\[0.25in]
& \ &   3 \cdot 2 \cdot [\, y_{|_1}^2 \,+\, y_{|_n}^2\,] \cdot \left\{ \Delta_o \ [ \,\cdot \cdot \cdot \ {\mbox{degree}} \ \,=\,  \  {\cal o} \,+\, 2 \ \cdot \cdot \cdot \,] \ \right\}\\[0.1in]
 & \,=\, &   3 \cdot 2 \cdot  y_{|_1}^2   \cdot \left\{ \Delta_o \ [ \,\cdot \cdot \cdot \ {\mbox{degree}} \ \,=\,  \ {\cal o} \ \,+\, \ 2 \ \cdot \cdot \cdot \,] \ \right\} \\[0.1in]
   & \ & \ \ \ \ \,+\, \ 3 \cdot 2 \cdot  y_{|_n}^2\,  \cdot \left\{ \Delta_o \ [ \,\cdot \cdot \cdot \ {\mbox{degree}} \ \,=\, \  {\cal o} \,+\, 2 \ \cdot \cdot \cdot \,] \ \right\}\ ,\ \ \ \ \ \ \ \ \ \ \ \ \ \ \ \ \ \ \ \ \ \ \ \ \ \ \ \\[0.1in]
{\mbox{and}} \,& \ & \!\!\!\!\!\!\Delta_o^{\left({{2 \,+\, ( \,{\cal o} \,+\, 2)}\over 2} \right) } \ \left[\  y_{|_1}^2   \cdot \left\{ \Delta_o \ [ \,\cdot \cdot \cdot \ {\mbox{degree}} \ \,=\,  \ {\cal o} \,+\, 2 \ \cdot \cdot \cdot \,] \ \right\} \right] \\[0.1in]
 &\ &  \!\!\!\!\!\!\!\!\!\!\!\!\!\!\!\!\!\!\!\!\!\!\!\!\!=\ \Delta_o^{\left({{2 \,+\, ( \,{\cal o} \,+\, 2)}\over 2} \right) } \ \left[\  y_{|_n}^2   \cdot \left\{ \Delta_o \ [ \,\cdot \cdot \cdot \ {\mbox{degree}} \ =\   {\cal o} \ \,+\, \ 2 \ \cdot \cdot \cdot \,] \ \right\} \right] \   (\leftarrow \   {\mbox{equal \   to \   a \   number}})\,,
\end{eqnarray*}
we apply the remaining order of Laplacian on (A.8.4) and (A.8.5), yielding the same numbers. Hence we verify (A.8.3), and so (A.8.2). This  completes the induction step.

\vspace*{0.3in}

{\large{\bf \S\,A.9.}{\bf\ \  Verification of\,} (A.6.44).}\\[0.2in]
Let ${\cal A}$\,,\, $\,{\cal B}$\, and $\,{\cal C}$\,  be numbers such that
$$
{\cal A} \,>\, 0\,, \ \ \  \ {\cal A} \,+ \,{\cal B} \,+\, {\cal C} \,>\, 0\,, \ \ \ \ {\mbox{and}} \ \  \ \ \beta \ \,=\, \ {{2n}\over {n \,-\, 2}}\ .\leqno (A.9.1)
$$
 We show that for any  number $\,\varepsilon\, > \,0$\, small enough,  there exists a positive number ${\bar C}_{\beta}$\, so that

 \vspace*{-0.25in}

$$
 |\, ({\cal A} \,+\, {\cal B} \,+\, {\cal C})^\beta \,-  \,{\cal A}^{\beta}\,|
  \ \le \ \varepsilon\,   {\cal A}^\beta \ +\  {{{\bar C}_{\beta}}\over{\varepsilon^{{2n}\over {n \,-\, 2}} }} \cdot  \left( |\, {\cal B} |^{\,\beta}\, \,+\, \,|\, {\cal C} |^{\,\beta} \,\right)\,.\leqno (A.9.2)
 $$
In this presentation, the proof of this statement relies on the following ($\beta$ being relaxed).

 \vspace*{0.2in}

{\bf Lemma A.9.3.} \ \ {\it Let $\,\beta \ge 1\,$ be given\,.\, There is a positive number $\,c_o\,$   such that for any number  $\,\varepsilon \,\in\, (0, \ c_o)\,$, we have}
$$
|\,(1 \,+\, t)^\beta \,- 1| \ \le \ \varepsilon \ \,+\, \ {{C_\beta}\over {\varepsilon^\beta}} \cdot |\,t|^\beta \ \ \ \ \ \ \  \ {\it{for}} \ \ \
 t \in [\,-1, \ \infty)\,. \leqno (A.9.4)
 $$
{\it   Here $\,c_o\,$ and $\,C_\beta$\, do not depend on $\,\varepsilon\,$ or $\,t\,$\,.\,}\\[0.2in]
{\bf Proof.} \ \ Via Taylor expansion, there is a positive number $\,c_o\,$ ($\,< \,1$)\,  such that
$$
|\, (1 \,+\, t)^\beta  \,-\, 1|\ \le \  ( \beta \,+\, 1)\, |\,t| \ \ \ \ \mfor \ \ |\,t|\,\le\, {{c_o}\over {\beta \,+\, 1}}\,.\leqno (A.9.5)
$$
For $\,\varepsilon \in (0, \ c_o)\,,\,$ consider the case
\begin{eqnarray*}
 |\,t|\ \le\  {{\varepsilon}\over {\beta \,+\, 1}} &  \Longrightarrow & {{|\,(1 \,+\, t)^\beta \,-\, 1|   }\over { \varepsilon    }} \ \le \  {{  ( \beta \,+\, 1) |\,t|  }\over { \varepsilon  }} \ \le \ 1  \ \ \ \ \  \ \ \ \ \  \ \ \ \ [\,{\mbox{using \ \ (A.9.4)}}]\\[0.1in]
&  \Longrightarrow &   |\,(1 \,+\, t)^\beta \,-\, 1|     \ \le \    \varepsilon    \ \ \Longrightarrow \ \  (A.9.4) \ \ {\mbox{holds}} \mfor |\,t|\ \le\  {{\varepsilon}\over {\beta \,+\, 1}} \ .
\end{eqnarray*}

\vspace*{-0.25in}

$$
1 \ \ge \ |\,t| \ \ge\  {{\varepsilon}\over {\beta \,+\, 1}}\,,\, \leqno{\mbox{When}}
$$
 we find $\,C_\varepsilon\,$ to be large enough so that
$$
{{ 2^\beta+1    }\over { C_\varepsilon \cdot \left[\,{{\varepsilon}\over {\beta \,+\, 1}} \right]^\beta    }} \ \le \  1\,, \ \ \  \ \ \ {\mbox{that \ \ is}}\,, \ \ C_\varepsilon \ \,=\, \ {1\over {\varepsilon^\beta}}\cdot \left(2^\beta \,+\, 1\right)\, (\,\beta \,+\, 1)^\beta\ . \leqno (A.9.6)
$$
It follows that
$$
|\,(1 \,+\, t)^\beta \,-\, 1| \ \le \ \left(2^\beta \,+\, 1\right)\  \le \  C_\varepsilon |\,t|^\beta \mfor 1 \ \ge \ |\,t| \ \ge\  {{\varepsilon}\over {\beta \,+\, 1}}\,.
$$

When $t > 1$\,,\,  from (A.9.6)
\begin{eqnarray*}
C_\varepsilon\ > \ 2^\beta & \Longrightarrow &  {{|\,(1 \,+\, t)^\beta \,-\, 1|   }\over { C_\varepsilon  \,t^\beta  }} \ \le \ {{|\,(2\,t)^\beta  |   }\over { C_\varepsilon  \,t^\beta  }}\ \le \ 1\\[0.1in]
& \Longrightarrow &  |\,(1 \,+\, t)^\beta \,-\, 1|  \ \le \ C_\varepsilon  \,t^\beta  \ \ \ \ \ \ \ \ \mfor t > 1\,.
\end{eqnarray*}
Combining the three cases, we have (A.9.4) with the choice of $\,C_\beta\ \,=\,  (2^\beta \,+\, 1) (\,\beta \,+\, 1)^\beta\,$ as specified in (A.9.6)\,.\, \qedwh
{\it Proof of\ } (A.9.2)\,.\, Using (A.8.4), we obtain
\begin{eqnarray*}
& \ & |\, ({\cal A}  \,+\, {\cal B} \,+\, {\cal C} )^\beta \,-\, {\cal A} ^\beta|\ \,=\, \
 {\cal A} ^\beta   \cdot \bigg\vert \left(1 \,+\, \left[\,{{{\cal B}+ {\cal C} }\over {{\cal A} }} \right]\,\right)^\beta  \ \,-\, \ 1 \ \bigg\vert   \\[0.1in]
 &  \ & \!\!\!\!\!\!\!\!\!\!\!\!\le \ {\cal A} ^\beta  \left( \varepsilon \ \,+\, \ {{C_\beta}\over {\varepsilon^\beta}} \cdot   \bigg\vert\,{{{\cal B}+ {\cal C} }\over { {\cal A} }}\bigg\vert ^\beta   \,\right)\ \ \    \ \ \ \   \left( {\cal A} > 0 \ \ \& \ \ {\cal A} +{\cal B} \,+\, {\cal C} > 0 \   \Longrightarrow \   {{{\cal B} \,+\, {\cal C}}\over {{\cal A} }} \,\ge\, \,-\, 1  \right) \\[0.1in]
 &  \ & \!\!\!\!\!\!\!\!\!\!\!\!\le \ {\cal A} ^\beta  \left[\ \varepsilon \ \,+\, \ {{C_\beta}\over {\varepsilon^\beta}} \cdot   \left( \bigg\vert\,{{{\cal B} }\over { {\cal A} }}\bigg\vert \ \,+\, \  \bigg\vert\,{{ {\cal C} }\over { {\cal A} }}\bigg\vert \, \right) ^\beta  \,\right] \, \le \, {\cal A}^\beta\!\cdot  \! \left\{ \varepsilon \, \,+\, \, {{C_\beta}\over {\varepsilon^\beta}} \cdot   \left[  \, 2^{\,\beta \,-\, 1} \left( \ \bigg\vert\,{{{\cal B} }\over { {\cal A} }}\bigg\vert^{\,\beta}   \,+\,  \bigg\vert\,{{ {\cal C} }\over { {\cal A} }}\bigg\vert^{\,\beta} \right)  \! \right]     \right\} \\[0.1in]
 &  \ & \!\!\!\!\!\!\!\!\!\!\!\!\le \  \varepsilon\,  {\cal A} ^\beta \,+\, 2^{\,\beta \,-\, 1} \cdot {{C_\beta}\over {\varepsilon^\beta}} \cdot  \left[\,|\,{\cal B} |^{\,\beta} \,+\, |\,{\cal C} |^{\,\beta}\, \right]\ \ \   [ \,{\mbox{using}} \ (\,|\,{\cal B} | \,+\, |\,{\cal C} |\,)^\beta \ \le \ 2^{\,\beta -1} \, (\,|\,{\cal B} |^{\,\beta} \,+\,  |\,{\cal C} |^{\,\beta} ) ].
\end{eqnarray*}
We can take  $\,{\bar C}_\beta \,=\,  2^{\,\beta \,-\, 1} \cdot C_\beta$\, to obtain (A.9.2)\,,\,  and take $\,\beta\, =\,  {{2n}\over {n \,-\, 2}}\ $ to obtain 7.23.\qed

\vspace*{0.15in}

{\large {\bf \S\,A.10\,.  \  Linear approximation to}\, $\left( A_1^{{n \,+\, 2}\over {n \,-\, 2}} \,-\, {\cal V}_{\,i}^{{n \,+\, 2}\over {n \,-\, 2}}  \right)\,$ {\bf in case of }} \\[0.075in]
\hspace*{0.83in} {\large {\bf simple blow\,-\,up\,.}}\\[0.2in]
The Taylor expansion
\begin{eqnarray*}
  (1 \,+\, t)^{\,p} & \,=\, & 1 \ \,+\, \  p\, t \ \,+\, \    O (t^2) \cdot (1 \,+\, \tau)^{ p \,-\, 2}     \ \ \ \ \ \ \ \  \ \ \  (\,{\mbox{here}} \ \ 0 \ \le\ \tau\ \le\  t)
\end{eqnarray*}
 tells us that

 \vspace*{-0.3in}

 \begin{eqnarray*}
  \ \ \ \ \ \  \ \ \ \ \ \  a^{\,p}   & \,=\, & [\,b \ +\ (a \,-\, b)]^{\,p} \ =\  b^{\,p}\, \left[\, 1 \ +\  {{(a \,-\, b)}\over b} \right]^{\,p} \\[0.1in]
 & \,=\, &   b^{\,p} \ +\   p\cdot (a \,-\, b)\, b^{\,p-1} \ +\   O \,(1) \, (1 \,+\, \tau)^{p-2} \,(a \,-\, b)^2 \cdot b^{ \,p -2} \ \ \ \ \ \ \ \ \ \ \ \ \ \ \ \ \ \
\end{eqnarray*}
for numbers with $\,a \,>\, b\, >\, 0\,$ and $\,p \,>\, 1\,$.\,
It follows that
$$
  a^{{n+2}\over {n \,-\, 2}} \,-\,  b^{{n+2}\over {n \,-\, 2}} \,=\,   \left[\, {{n+2}\over {n \,-\, 2}}\right] (a \,-\, b)\, b^{4\over {n \,-\, 2}} \ +\ O \,(1) \,(1 \,+\, \tau)^{ {4\over {n \,-\, 2}} \,-\, 1} \cdot ( a \,-\, b )^2    \cdot b^{\, {{n \,+\, 2}\over {n \,-\, 2}} \,-\,2 }\ . \leqno (A.10.1)
$$
Here

\vspace*{-0.3in}

$$
  0 \ \le \ \,\tau \  \le \ {{ a \,-\, b }\over b}\ .
$$
For $\,A_1\, ({\cal Y}) \ge \,{\cal V}_1\, ({\cal Y})\,,\,$ (A.10.1) implies that
 \begin{eqnarray*}
(A.10.2)   \!\!\!\!\!\!\!\!& \ & [\,A_1\, ({\cal Y})]^{{n \,+\, 2}\over {n \,-\, 2}} \ -\ [\,{\cal V}_i\, ({\cal Y})]^{{n \,+\, 2}\over {n \,-\, 2}}   \,=\,  \left(\, {{n+2}\over {n \,-\, 2}} \right)   [\, A_1 \, ({\cal Y}) \,-\, {\cal V}_{\,i} \, \, ({\cal Y})] \cdot [\,A_1\, ({\cal Y})]^{4\over {n \,-\, 2}} \\[0.15in]
 & \ & \ \ \ \ \ \ \ \ \ \ \ \ \ \ \ \  \ \,+\, \ \left[\,O\, (1) \,(1 \,+\, \tau)^{ {4\over {n \,-\, 2}} \,-\, 1}\, \right]\,[\, A_1 \, ({\cal Y}) \,-\, {\cal V}_{\,i} \, \, ({\cal Y})]^{\,2} \cdot [\,A_1\, ({\cal Y})]^{{4\over {n \,-\, 2}} -1} \,. \ \ \ \ \ \  \ \
\end{eqnarray*}
$$
n \ge 6 \ \   \Longrightarrow \ \  1 \,-\, {4\over {n \,-\, 2}}\ \ge \ 0 \ \ \Longrightarrow \ \ (1 \,+\, \tau)^{ {4\over {n \,-\, 2}} \,-\, 1} \,=\, {1\over {(1 \,+\, \tau)^{ 1 \,-\, {4\over {n \,-\, 2}}  }}} \le 1\,, \leqno {\mbox{Since}}
$$

\vspace*{-0.075in}

and
$$
{{A_1\, ({\cal Y}) \,-\, \,{\cal V}_i \, ({\cal Y})}\over { {\cal V}_i\, ({\cal Y})}} \le 1 \,+\, {{A_1\, ({\cal Y})  }\over { {\cal V}_i \, ({\cal Y})}} \le 1 \ \,+\, \ C\,, \ \  3 \ \le \ n \ \le \  5  \ \   \Longrightarrow \ \    (1 \,+\, \tau)^{ {4\over {n \,-\, 2}} \,-\, 1}   \le C_1\,.
$$
Recalling (2.26) in Proposition 2.24\,,\, and also    \,\S\,2\,f\,,\, we have
$$
C^{-1} \cdot A_1 \,({\cal Y})\ \le \ {\cal V}_{\,i}\,({\cal Y}) \ \le \ C\, A_1\,({\cal Y}) \ \ \ \ \ \ \mfor \ \ \  \,|\,{\cal Y}| \ \le \ \rho_o \,\lambda_i^{-1}\,. \leqno (A.10.3)
$$

It follows that
\begin{eqnarray*}
(A.10.4)\!\!\!\!\!\! & \ &  \,A_1\, ({\cal Y}) \ \ge \ {\cal V}_i\, ({\cal Y}) \ \ \Longrightarrow \ \
[\,A_1\, ({\cal Y})]^{{n \,+\, 2}\over {n \,-\, 2}} \ -\ [\,{\cal V}_i\, ({\cal Y})]^{{n \,+\, 2}\over {n \,-\, 2}}\\[0.1in]
   & \ & \hspace*{2in} \,=\,  \left(\, {{n+2}\over {n \,-\, 2}} \right)  [\,A_1\, ({\cal Y})]^{4\over {n \,-\, 2}} \cdot  [\, A_1 \, ({\cal Y}) \,-\, {\cal V}_{\,i} \, \, ({\cal Y})] \\[0.1in]
 & \ &  \ \ \ \ \ \ \ \ \ \ \ \ \ \  \ \ \ \ \ \ \ \ \ \ \ \ \ \ \ \ \ \ \ \ \ \ \ \ \ \ \ \ \ \, \,+\, \ O \,(1) \,[\, A_1 \, ({\cal Y}) \,-\, {\cal V}_{\,i} \, \, ({\cal Y})]^{\,2}\cdot  [\,A_1\, ({\cal Y})]^{\,{{4\over {n \,-\, 2}} -1}} \ . \ \ \ \
\end{eqnarray*}

For the opposite case
$$
 {\cal V}_i\, ({\cal Y}) \  >  \ A_1\, ({\cal Y})
$$
we obtain a similar expression as in (A.10.4), with the last term
$$
[\,A_1\, ({\cal Y})]^{\,{{4\over {n \,-\, 2}}\, -\,1}} \ \ \ \   {\mbox{changed \ \ to}} \ \ \ [\,{\cal V}_i\, ({\cal Y})]^{\,{{4\over {n \,-\, 2}} \,-\,1}}\ .   \leqno (A.10.5)
$$
We obtain
$$
{{1}\over { C_2 \cdot A_1 \, ({\cal Y})  }} \ \le \  {1\over { {\cal V}_i\, ({\cal Y})}} \ \le \  {{C_2}\over { A_1 \, ({\cal Y})  }} \ \ \mfor \ \  |\, {\cal Y}|\le \rho_o \,\lambda_i^{-1}\ .
$$
Thus the two terms in (A.10.5) have the same order. We come to the conclusion that
\begin{eqnarray*}
(A.10.6)\!\!\!\!\!\!& \ &
[\,A_1\, ({\cal Y})]^{{n \,+\, 2}\over {n \,-\, 2}} \ -\ [\,{\cal V}_1\, ({\cal Y})]^{{n \,+\, 2}\over {n \,-\, 2}}   \,=\,  \left(\, {{n+2}\over {n \,-\, 2}} \right)  [\,A_1\, ({\cal Y})]^{4\over {n \,-\, 2}} \cdot  [\, A_1 \, ({\cal Y}) \,-\, {\cal V}_{\,i} \, \, ({\cal Y})] \\[0.1in]
 & \ &  \ \ \ \ \ \ \ \ \ \ \ \ \ \  \ \ \ \ \ \ \ \ \ \ \ \ \ \ \ \ \ \ \ \ \ \ \ \ \ \   \,+\, \ O \,(1) \,[\, A_1 \, ({\cal Y}) \,-\, {\cal V}_{\,i} \, \, ({\cal Y})]^{\,2}\cdot  [\,A_1\, ({\cal Y})]^{\,{{4\over {n \,-\, 2}} -1}} \ , \ \ \ \
\end{eqnarray*}
which holds for $|\,{\cal Y}|\ \le \ \rho_o\, \lambda_i^{-1}\,.\,$ (A.10.6) is
independent on which value is bigger, $\,A_1\, ({\cal Y})\,$ or $\,{\cal V}_i\, ({\cal Y})\,$.

\end{document}